%% file: curves-and-surfaces.tex
\documentclass[twoside]{book}
\usepackage{lectures}

\newcommand{\arxiv}[2]{#1} %for arxiv

\newcommand{\spell}[2]{#1} %for text

\arxiv{\geometry{top=0.9in, bottom=0.9in,inner=0.9in, outer=0.7in, paperwidth=6in, paperheight=9in}}{\geometry{top=1.025in, bottom=1.025in,inner=0.9in, outer=0.825in, paperwidth=6.125in, paperheight=9.25in}}
\makeindex

\spell{}{\makeatletter\usepackage{environ}\RenewEnviron{figure}[1][]{}\RenewEnviron{figure*}[1][]{}\RenewEnviron{wrapfigure}[1][]{}\RenewEnviron{Figure}[1][]{}\pagestyle{empty}\let\ps@plain\ps@empty\makeatother\usepackage[none]{hyphenat}\def\emph{\textit}\renewcommand{\footnote}[1]{\textup{\ (#1)}}\renewcommand{\frac}[2]{(#1)/(#2)}\renewcommand{\tfrac}[2]{(#1/#2)}\everydisplay{\textstyle}}

\begin{document}

\addtocontents{toc}{\cftpagenumbersoff{part}}

\input{title.tex}

\input{intro.tex}

\input{part-curves.tex}

\input{curves.tex}

\input{length-new.tex}

\input{curve-curvature-new.tex}

\input{poly.tex}
\input{torsion-new.tex}

\input{plane-curves.tex}

\input{supporting-curves.tex}

\input{part-surfaces.tex}

\input{surfaces-def.tex}

\input{first-order.tex}

\input{surfaces-curvature.tex}

\input{saddle.tex}

\input{part-geodesics.tex}

\input{shortest-path.tex}
\input{geodesics-new.tex}
\input{parallel-new.tex}

\input{gauss-bonnet.tex}
\input{polar-chart.tex}
\input{comparison.tex}

\arxiv{\newgeometry{top=0.83in, bottom=0.83in,inner=0.52in, outer=0.52in}}{\newgeometry{top=0.955in, bottom=0.955in,inner=0.52in, outer=0.645in}}

{
\appendix
\input{prelim.tex}
}

\backmatter

{

\footnotesize

\input{sol-new.tex}

}

\arxiv{\newgeometry{top=0.9in, bottom=0.9in,inner=0.9in, outer=0.7in}}{\newgeometry{top=1.025in, bottom=1.025in,inner=0.9in, outer=0.825in}}

\input{afterword.tex}

\arxiv{\newgeometry{top=0.83in, bottom=0.83in,inner=0.52in, outer=0.52in}}{\newgeometry{top=0.955in, bottom=0.955in,inner=0.52in, outer=0.645in}}

{
\newpage
\phantomsection
\footnotesize\sloppy
\input{curves-and-surfaces.ind}

}

{

\sloppy

\def\emph{\textit}
\printbibliography[heading=bibintoc]

\fussy

}
\end{document}

%% file: title.tex
\hypersetup
{
pdftitle={What is differential geometry?: curves and surfaces},
pdfauthor={Anton Petrunin and Sergio Zamora Barrera}
}

\title{\tt What is differential geometry?:\\
curves and surfaces}

\author{\tt Anton Petrunin and \\ \tt Sergio Zamora Barrera}

\date{}

\maketitle

\thispagestyle{empty}

%% file: intro.tex
Artwork by Ana Cristina Chávez Cáliz.
\null
\vfill
\noindent{\includegraphics[scale=0.5]{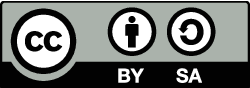}
\vspace*{1mm}
\\
\hbox{\parbox{1\textwidth}
{This work is licensed under CC BY-SA 4.0.
To view a copy of this license, visit
\texttt{https://creativecommons.org/licenses/by-sa/4.0}}}}

\thispagestyle{empty}
\newpage
\clearpage
\phantomsection
\pdfbookmark[1]{\contentsname}{bm:toc}
\tableofcontents

\vfill

\begin{figure}[h!]
\centering
\begin{tikzpicture}[->,>=stealth',shorten >=1pt,auto,scale=.25,
  thick,main node/.style={circle,draw,font=\sffamily\bfseries,minimum size=8mm}]

  \node[main node] (1) at (0,0) {\ref{chap:curves-def}};
  \node[main node] (2) at (3,-5) {\ref{chap:length}};
  \node[main node] (3) at (6,-10) {\ref{chap:curve-curvature}};
  \node[main node] (4) at (9,-5) {\ref{chap:poly}};
  \node[main node] (5) at (9,-15) {\ref{chap:torsion}};
  \node[main node] (6) at (12,-10) {\ref{chap:signed-curvature}};
  \node[main node] (7) at (15,-15) {\ref{chap:supporting-curves}};
  \node[main node] (8) at (18,0) {\ref{chap:surfaces-def}};
  \node[main node] (9) at (15,-5) {\ref{chap:first-order}};
  \node[main node] (10) at (18,-10) {\ref{chap:surface-curvature}};
  \node[main node] (11) at (21,-15) {\ref{chap:Curves in a surface}};  
  \node[main node] (12) at (18,-20) {\ref{chap:surface-support}};
  \node[main node] (13) at (21,-5) {\ref{chap:shortest}};
  \node[main node] (14) at (24,-10) {\ref{chap:geodesics}};
  \node[main node] (15) at (27,-15) {\ref{chap:parallel-transport}};
  \node[main node] (16) at (33,-15) {\ref{chap:gauss-bonnet}};
  \node[main node] (17) at (30,-10) {\ref{chap:semigeodesic}};
  \node[main node] (18) at (36,-10) {\ref{chap:comparison}};
  
  \path[every node/.style={font=\sffamily\small}]
   (1) edge node{}(2)
   (2) edge node{}(3)
   (3) edge node{}(5)
   (3) edge node{}(4)
   (3) edge node{}(6)
   (4) edge[dashed] node{}(6)
   (6) edge node{}(7)
   (6) edge node{}(10)
   (7) edge node{}(12)
   (8) edge node{}(9)
   (8) edge node{}(13)
   (9) edge node{}(10)
   (10) edge node{}(11)
   (11) edge node{}(12)
   (10) edge node{}(14)
   (13) edge node{}(14)
   (14) edge node{}(15)   
   (14) edge[bend left=15] node{}(17)
(17) edge[dashed, bend left=15] node{}(14)
   (15) edge node{}(16)
   (16) edge[dashed] node{}(18)
   (17) edge[dashed] node{}(15)
   (17) edge[dashed] node{}(16)
   (17) edge node{}(18);
\end{tikzpicture}
\end{figure}

\vfill

\newpage
%\phantomsection
\chapter*{Preface}
\addcontentsline{toc}{chapter}{Preface}
\thispagestyle{myheadings}
\markboth{PREFACE}{PREFACE}

These notes are designed for those who either plan to work in differential geometry,
or at least want to have a good reason \textit{not} to do~it.
It should be more than sufficient for a semester-long course. 

Differential geometry draws on several branches of mathematics, including
real analysis,
measure theory,
the calculus of variations,
differential equations,
topology,
and elementary and convex geometry.
Physical intuition also helps in understanding many of its aspects.
This subject is wide even at the beginning. 
For that reason, it is fun and painful both to teach and to study.

In this book, we discuss smooth curves and surfaces --- the main gate to differential geometry.
This subject provides a collection of examples and ideas critical for further study.
It is wise to become a master in this subject before making further steps --- there is no need to rush.

We give a general overview of the subject, keeping it
problem-centered,
elementary, 
visual, 
and virtually rigorous; we allow gaps that belong to other branches of mathematics, most of these subjects are discussed briefly in the appendix.

We focus on the techniques that are absolutely essential for further study.
For  that reason we omit a number of topics that are traditionally included in the introductory texts;
for example, we almost do not touch 
minimal surfaces
and the Peterson--Codazzi formulas.
%(In an ideal world, we could also remove the chapter about torsion, but we decided to yield to tradition this once.) 

At the same time, we get to applications that are not in the scope of typical introductory texts.
 
The first example is the theorem of Vladimir Ionin and German Pestov about the Moon in a puddle (\ref{thm:moon-orginal}).
This theorem might be the simplest meaningful example of the so-called {}\emph{local-to-global} theorems which lie at the heart of differential geometry;
for that reason, it is a good answer to the main question of this book --- ``What is differential geometry?''.

Other examples include the theorem of Sergei Bernstein on saddle graphs (\ref{thm:bernshtein}) and the theorem of Stefan Cohn-Vossen
on a two-sided infinite geodesic (\ref{thm:cohn-vossen}).

These notes are based on the lectures given at the MASS program (Mathematics Advanced Study Semesters at Pennsylvania State University) in the fall of 2018.
Many of the topics were presented by Yurii Burago in his lectures teaching the first author at Leningrad University.
We extensively used textbooks
by Wilhelm Blaschke and Kurt Leichtweiss~\cite{blaschke-leichtweiss},
by Aleksei Chernavskii~\cite{chernavsky},
and by Victor Toponogov~\cite{toponogov-book}, and also the lecture notes of Sergei Ivanov \cite{ivanov};
many advanced exercises are taken from \cite{petrunin2020}.
The last chapter is based on the introductory material in the book by Stephanie Alexander, Vitali Kapovitch, and the first author~\cite{alexander-kapovitch-petrunin2027}.
We want to thank
Stephanie Alexander,
Yurii Burago, 
Berk Ceylan,
Nina Lebedeva,
Alexander Lytchak,
Benjamin McKay,
and all the students in our class
for their help.

The present work is partially supported by NSF grant DMS-2005279
and by the Simons Foundation under grant \#584781.

\begin{flushright}
Anton Petrunin and\\
Sergio Zamora Barrera.
\end{flushright}

%% file: part-curves.tex
\arxiv{\cleardoublepage
\phantomsection
\AddToShipoutPictureBG*{\includegraphics[width=\paperwidth]{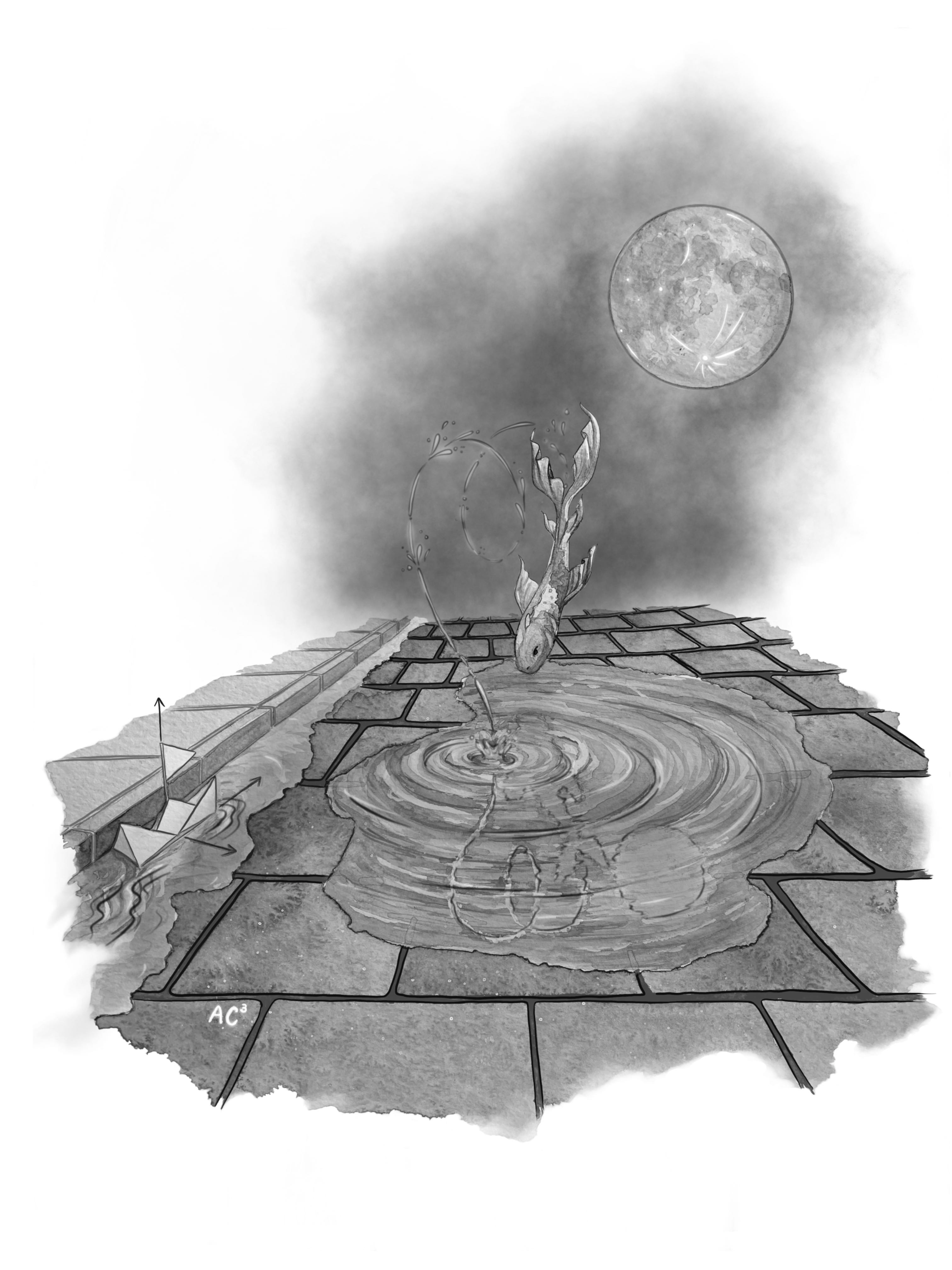}}
\cleardoublepage
\thispagestyle{empty}
\stepcounter{part}
\begin{center}
{\Huge\textbf{Part \Roman{part}:\qquad Curves}}\
\end{center}
\addcontentsline{toc}{part}{\Roman{part} Curves}
\clearpage}
{\backgroundsetup{
scale=.95,
opacity=1,
angle=0,
contents={
\includegraphics[width=\paperwidth%,height=\paperheight
]{pics/Curves}
}%
}
\cleardoublepage
\phantomsection
\thispagestyle{empty}
\BgThispage
\stepcounter{part}
\addcontentsline{toc}{part}{\Roman{part} Curves}
\begin{center}
{\Huge\textbf{Part \Roman{part}:\qquad Curves}}\
\end{center}
}

%% file: curves.tex
\chapter{Definitions}
\label{chap:curves-def}

In most cases, by looking at tire tracks of a bicycle,
it is possible to tell in which direction it traveled and to tell apart the trails made by the front and rear wheels;
plus one can find the distance between the wheels.

Try to figure out how to do this.
This question might motivate the reader to reinvent a considerable part of the differential geometry of curves.
To learn more about this problem check the reference given after Exercise~\ref{ex:bike}.

\section{Before we start}

The notion of curve comes with many variations.
Some of them are described by nouns (path, arc, and so on) and  
others by adjectives (closed, open, proper, simple, smooth, and so on).
The following picture gives an idea about four of these variations.

\vskip-0mm
\begin{figure}[h!]
\begin{minipage}{.48\textwidth}
\centering
\includegraphics{mppics/pic-110}
\end{minipage}\hfill
\begin{minipage}{.48\textwidth}
\centering
\includegraphics{mppics/pic-115}
\end{minipage}
\bigskip
\begin{minipage}{.48\textwidth}
\centering
\includegraphics{mppics/pic-120}
\end{minipage}\hfill
\begin{minipage}{.48\textwidth}
\centering
\includegraphics{mppics/pic-125}
\end{minipage}
\end{figure}
\vskip-0mm

This chapter covers all these definitions.
One may skip this chapter and use it as a quick reference while reading the book.

\section{Simple curves}

We assume that the reader is familiar with metric spaces (see Appendix \ref{sec:metric-spcaes});
their main examples will be the Euclidean plane $\mathbb{R}^2$ and the Euclidean space $\mathbb{R}^3$.

Recall that a \index{real interval}\emph{real interval} is a connected subset of the real line.
A bijective continuous map $f\:X\to Y$ between subsets of some metric spaces is called a {}\emph{homeomorphism} if its inverse $f^{-1}\:Y\z\to X$ is continuous.  

\begin{thm}{Definition} 
A connected subset $\gamma$ in a metric space is called a \index{simple curve}\emph{simple curve} if it is \index{locally homeomorphic}\emph{locally} homeomorphic to a real interval;
that is, any point $p\in\gamma$ has a neighborhood in $\gamma$ that is homeomorphic to a real interval.
\end{thm}

It turns out that any simple curve $\gamma$ can be \index{parametrization}\emph{parametrized} by a real interval or circle.
That is, there is a homeomorphism $\GG\to\gamma$ 
where $\GG$ is a real interval (open, closed, or semi-open) or the circle
\[\mathbb{S}^1=\set{(x,y)\in\mathbb{R}^2}{x^2+y^2=1}.\] 
A complete proof of the latter statement is given by David Gale \cite{gale}.
The proof is not hard, but it would take us away from the main subject.
%A~finicky reader may add this property in the definition of curves.

The parametrization $\GG\to\gamma$ describes the curve completely.
We will denote a curve and its parametrization by the same letter;
for example, we may say a curve $\gamma$ is given with a parametrization $\gamma\:(a,b]\z\to \mathbb{R}^2$.
Note, however, that any simple curve admits many parametrizations. 

\begin{thm}{Exercise}\label{ex:9}

\begin{subthm}{ex:9:compact}
Show that the image of any continuous injective map $\gamma\:[0,1]\to\mathbb{R}^2$ is a simple arc.
\end{subthm}

\begin{subthm}{ex:9:9}
Find a continuous injective map $\gamma\:(0,1)\to\mathbb{R}^2$ such that its image is \textit{not} a simple curve.
\end{subthm}

\end{thm}

\section{Parametrized curves}\label{sec:Parametrized curves}

Let $\GG$ be a circle or a real interval (open, closed, or semi-open), and let $\spc{X}$ be a metric space.
A \index{parametrized curve}\emph{parametrized curve} is defined as a continuous map $\gamma\:\GG\to\spc{X}$. 
For a parametrized curve, we do \textit{not} assume that the map is injective; in other words, a parametrized curve might have {}\emph{self-intersections}.

From the previous section, we know that a simple curve can be parametrized.
The term \index{curve}\emph{curve} can be used if we do not want to specify whether it is parametrized or simple.

If $\GG$ is an open interval or a circle, we say that $\gamma$ is a {}\emph{curve without endpoints};
otherwise, it is called a {}\emph{curve with endpoints}.
In the case when $\GG$ is a circle, we say that the curve is \index{closed!curve}\emph{closed}. 
If the domain $\GG$ is a closed interval $[a,b]$, then the curve is called an \index{arc}\emph{arc}.
Further, if it is the unit interval $[0,1]$, then it is also called a \index{path}\emph{path}.

\begin{wrapfigure}{o}{15 mm}
\vskip-0mm
\centering
\includegraphics{mppics/pic-130}
\end{wrapfigure}

%%loops are not necesarily parametrized by the interval [0,1].

If an arc $\gamma \: [a,b] \to \spc{X}$ satisfies $p=\gamma (a)=\gamma (b)$, then $\gamma$ is called a \index{loop}\emph{loop};
in this case, the point $p$ is called the \index{base of loop}\emph{base} of the loop.

Suppose $\GG_1$ and $\GG_2$ are either real intervals or circles.
A continuous onto map $\tau\:\GG_1\to\GG_2$ is called \index{monotone map}\emph{monotone} if, for any $t\in \GG_2$, the set $\tau^{-1}\{t\}$ is connected.
If $\GG_1$ and $\GG_2$ are intervals, then, by the intermediate value theorem, a monotone map is either nondecreasing or nonincreasing;
that is, our definition agrees with the standard one when $\GG_1$ and $\GG_2$ are intervals.

\begin{thm}{Exercise}\label{ex:mono}
Construct a monotone (and, in particular, onto) map $(0,1)\to [0,1]$.
\end{thm}

Suppose $\gamma_1\:\GG_1\to \spc{X}$ and $\gamma_2\:\GG_2\to \spc{X}$ are two parametrized curves such that 
$\gamma_1=\gamma_2\circ\tau$ for a monotone map $\tau\:\GG_1\to\GG_2$.
Then we say that $\gamma_2$ is a \index{reparametrization}\emph{reparametrization}%
\footnote{In general, $\gamma_1$ is \textit{not} a reparametrization of $\gamma_2$.
In other words, according to our definition, the described relation \textit{being a reparametrization} is not symmetric;
in particular, it is not an equivalence relation.
Usually, it is fixed by extending it to the minimal equivalence relation that includes ours \cite[2.5.1]{burago-burago-ivanov}.
But we will stick to our version.}
of $\gamma_1$ by $\tau$.

\begin{thm}{Advanced exercise}\label{aex:simple-curve}
Let $X$ be a subset of the plane.
Suppose two distinct points $p,q\in X$ can be connected by a path in~$X$.
Show that there is a simple arc in~$X$ connecting $p$ to~$q$.
\end{thm}

Any loop and any simple closed curve can be described by a {}\emph{periodic} parametrized curve $\gamma\: \mathbb{R}\to \spc{X}$;
that is, a curve such that $\gamma(t+\ell)=\gamma(t)$ for a fixed period $\ell>0$ and all~$t$.
For example, the unit circle in the plane can be described by the $2{\cdot}\pi$-periodic parametrization $\gamma(t)=(\cos t,\sin t)$.

In the reverse direction, any curve with periodic parametrization yields a loop that can be also considered as a closed curve. 

\section{Smooth curves}\label{sec:Smooth curves}

Curves in the Euclidean space or plane are called \index{space curve}\emph{space curves} or, respectively, \index{plane curve}\emph{plane curves}.

A parametrized space curve can be described by its coordinate functions 
\[\gamma(t)=(x(t),y(t),z(t)).\]
Plane curves can be considered as a special case of space curves with $z(t)\equiv 0$.

Recall that a real-to-real function is called \index{smooth!function}\emph{smooth} if its derivatives of all orders are defined everywhere in the domain of definition.  
If each of the coordinate functions $x(t), y(t)$, and $z(t)$ are smooth, then the parametrization is called 
\index{smooth!parametrization}\emph{smooth}.

If the \index{velocity vector}\emph{velocity vector} 
\[\gamma'(t)=(x'(t),y'(t),z'(t))\] 
does not vanish at any point, then the parametrization $\gamma$ is called \index{regular!parametrization}\emph{regular}.

A parametrized curve is called {}\emph{smooth} if its parametrization is smooth and regular.
A simple space curve is called \index{smooth!curve}\emph{smooth} if it admits a regular smooth parametrization;
for closed curves the parametrization is assumed to be periodic.
These curves are the main objects in the first part of the book.
Pedantically, one could call them {}\emph{regular smooth curves}. 

Notice that the closed curve described by a smooth loop might fail to be smooth at its base; an example is shown on the picture.

\begin{wrapfigure}{o}{17 mm}
\vskip-4mm
\centering
\includegraphics{mppics/pic-51}
\bigskip
\includegraphics{mppics/pic-140}
\vskip-8mm
\end{wrapfigure}

The following exercise shows that curves with smooth parametrizations might have corners and might fail to be smooth curves.

\begin{thm}{Exercise}\label{ex:L-shape}
From Appendix~\ref{sec:analysis}, we know that the following function is smooth:
\[f(t)=
\begin{cases}
0&\text{if}\ t\le 0,
\\
\frac{t}{e^{1\!/\!t}}&\text{if}\ t> 0.
\end{cases}
\]

Show that $\alpha(t)=(f(t),f(-t))$ gives a smooth parametrization of the curve shown on the picture;
it is a simple curve formed by the union of two half-axes in the plane.

Show that any smooth parametrization of this curve has a vanishing velocity vector at the origin.
Conclude that this curve is not smooth;
that is, it does not admit a regular smooth parametrization.
\end{thm}

\begin{thm}{Exercise}\label{ex:cycloid}
Describe the set of real numbers $\ell$
such that the parametrization $\gamma_\ell (t)= (t+\ell \cdot \sin t,\ell \cdot \cos t)$, $t\in\mathbb{R}$ is

\begin{minipage}{.30\textwidth}
\begin{subthm}{ex:cycloid:smooth}
smooth; 
\end{subthm}
\end{minipage}
\hfill
\begin{minipage}{.30\textwidth}
\begin{subthm}{ex:cycloid:regular}
regular;
\end{subthm}
\end{minipage}
\hfill
\begin{minipage}{.30\textwidth}
\begin{subthm}{ex:cycloid:simple}
simple.
\end{subthm}
\end{minipage}

\end{thm}

\begin{thm}{Exercise}\label{ex:nonregular}
Find a parametrization of the cubic parabola $y=x^3$ in the plane that is smooth, but \textit{not} regular.
\end{thm}

\section{Implicitly defined curves}\label{sec:implicit-curves}

Suppose $f\:\mathbb{R}^2\to \mathbb{R}$ is a smooth function; 
that is, all its partial derivatives are defined everywhere.
Let $\gamma\subset \mathbb{R}^2$ be its level set described by  the equation $f(x,y)=0$.

Assume $0$ is a \index{regular!value}\emph{regular value} of~$f$; that is, the gradient $\nabla_p f$ does not vanish at any point $p\in \gamma$.
In other words, if $f(p)=0$, then   
$f_x(p)\ne 0$ or $f_y(p)\ne 0$.%
\footnote{Here $f_x$ is a shortcut notation for the partial derivative
$\tfrac{\partial f}{\partial x}$.\index{10f@$f_x$ (partial derivative)}}
If $\gamma$ is connected, then by the implicit function theorem (\ref{thm:imlicit}), it is a smooth simple curve. 

The described condition is sufficient but \textit{not necessary}.
For example, zero is \textit{not} a regular value of the function $f(x,y)=y^2$, but the equation $f(x,y)=0$ describes a smooth curve --- the $x$-axis.

Similarly, assume $(f,h)$ is a pair of smooth functions defined in $\mathbb{R}^3$.
Again, the implicit function theorem (\ref{thm:imlicit}) implies that the system of equations
\[\begin{cases}
   f(x,y,z)=0,
   \\
   h(x,y,z)=0
  \end{cases}
\]
defines a smooth space curve if the set $\gamma$ of solutions is connected, and $0$ is a regular value of the map $F\:\mathbb{R}^3\to\mathbb{R}^2$ defined by
\[F\:(x,y,z)\mapsto (f(x,y,z),h(x,y,z)).\]
It means that the gradients $\nabla f$ and $\nabla h$ are linearly independent at any point $p\in \gamma$.
In other words, the Jacobian matrix
\[
\Jac_pF=\begin{pmatrix}
f_x&f_y&f_z\\
h_x&h_y&h_z
\end{pmatrix}
\]
for the map $F\:\mathbb{R}^3\to\mathbb{R}^2$ has rank $2$ at any point $p \in \gamma$.

If a curve $\gamma$ is described in such a way,
then we say that it is \index{implicitly defined curve}\emph{implicitly defined}.

The implicit function theorem guarantees the existence of regular smooth parametrizations for any implicitly defined curve.
However, when it comes to calculations, it is usually easier to work directly with implicit representations. 

\begin{thm}{Exercise}\label{ex:y^2=x^3}
Consider the set in the plane described by the equation
$y^2=x^3$.
Is it a simple curve?
Is it a smooth curve?
\end{thm}

\begin{thm}{Exercise}\label{ex:viviani}
Describe the set of real numbers $\ell$
such that the system of equations
\[\begin{cases}
x^2+y^2+z^2&=1
\\
x^2+\ell\cdot x+y^2&=0
\end{cases}\]
describes a smooth curve.
\end{thm}

\section{Proper, closed, and open curves}\label{sec:proper-curves}

A parametrized curve $\gamma$ in a metric space $\spc{X}$ is called \index{proper!curve}\emph{proper} if, for any compact set $K \subset \spc{X}$, the inverse image $\gamma^{-1}(K)$ is compact.

For example, the curve $\gamma(t)=(e^t,0,0)$ defined on the real line is not proper.
Indeed, the half-line $(-\infty,0]$ is not compact, but it is the inverse image of the closed unit ball around the origin.

\begin{thm}{Exercise}\label{ex:open-curve}
Suppose $\gamma\:\mathbb{R}\to\mathbb{R}^3$ is a proper curve.
Show that  $|\gamma(t)|\z\to\infty$ as $t\to\pm\infty$.
\end{thm}

Recall that closed bounded intervals are compact (\ref{thm:Heine--Borel}) and closed subsets of a compact set are compact;
see Appendix~\ref{sec:topology}.
Since circles and closed intervals are compact,
it follows that closed curves and arcs are proper.

A simple curve is called proper if it admits a proper parametrization.

\begin{thm}{Exercise}\label{ex:proper-closed}
Show that a simple space curve is proper if and only if it is a closed set.

\end{thm}

A proper simple curve is called \index{open!curve}\emph{open} if it is not closed and has no endpoints.
So any simple proper curve without endpoints is either closed or open.
The terms \textit{open curve} and \textit{closed curve} have nothing to do with open and closed sets.

\begin{thm}{Exercise}\label{ex:proper-curve}
Use Jordan's theorem (\ref{thm:jordan}) to show that any simple open plane curve divides the plane in two connected components.  
\end{thm}

%% file: length-new.tex
\chapter{Length}
\label{chap:length}

\section{Definitions}

A sequence 
\[a=t_0 < t_1 < \cdots < t_k=b.\]
is called a \index{partition}\emph{partition} of the interval $[a,b]$.

\begin{thm}{Definition}\label{def:length}
Let $\gamma\:[a,b]\to \spc{X}$ be a curve in a metric space.
The \index{length of curve}\emph{length} of $\gamma$ is defined as
\begin{align*}
\length \gamma
&= 
\sup
\set{\dist{\gamma(t_0)}{\gamma(t_1)}{\spc{X}}
+\dots+
\dist{\gamma(t_{k-1})}{\gamma(t_k)}{\spc{X}}}{},
\end{align*}
where the least upper bound is taken over all partitions $t_0,\dots,t_k$ of $[a,b]$.

The length of a closed curve is defined as the length of the corresponding loop.
If a curve is parametrized by an open or semi-open interval, then its length is defined as the least upper bound of the lengths of all its restrictions to closed intervals.
 
\end{thm}

A curve is called \index{rectifiable curve}\emph{rectifiable} if its length is finite.

\begin{thm}{Exercise}\label{ex:integral-length-0}
Suppose a curve $\gamma_1\:[a_1,b_1] \to\mathbb{R}^3$ is a reparametrization of $\gamma_2\:[a_2,b_2] \to\mathbb{R}^3$.
Show that
\[\length \gamma_1 = \length \gamma_2.\]

\end{thm}

\begin{wrapfigure}[4]{r}{33 mm}
\vskip-4mm
\centering
\includegraphics{mppics/pic-224}
\end{wrapfigure}

Suppose $\gamma\:[a,b]\to \mathbb{R}^3$ is a parametrized space curve.
For a partition $a=t_0 < t_1 < \z\cdots < t_k=b$, set $p_i=\gamma(t_i)$.
Then the polygonal line $p_0\dots p_k$ is said to be \index{inscribed polygonal line}\emph{inscribed} in~$\gamma$.
If $\gamma$ is closed, then $p_0=p_k$, so the inscribed polygonal line is also closed.

Note that the length of a space curve can be defined as the least upper bound of the lengths of its inscribed polygonal lines.

\begin{thm}{Exercise}\label{ex:length-chain}
Let $\gamma\:[0,1]\to\mathbb{R}^3$ be a path.
Suppose that $\beta_n$ is a sequence of polygonal lines inscribed in $\gamma$ with vertices $\gamma(\tfrac in)$ for $i\z\in\{0,\dots,n\}$.
Show that 
\[\length\beta_n\to\length\gamma
\quad\text{as}\quad
n\to \infty.
\]
\end{thm}

\begin{thm}{Exercise}\label{ex:length-image}
Let $\gamma\:[0,1]\to\mathbb{R}^3$ be a simple path.
Suppose that a path $\beta\:[0,1]\to\mathbb{R}^3$ has the same image as $\gamma$;
that is, $\beta([0,1])=\gamma([0,1])$.
Show that 
\[\length \beta\ge \length \gamma.\]
Try to prove the same assuming only that $\beta([0,1])\supset\gamma([0,1])$.
\end{thm}

\begin{thm}{Exercise}\label{ex:integral-length}
Assume $\gamma\:[a,b]\to\mathbb{R}^3$ is a smooth curve.
Show that
\vskip1mm
\begin{minipage}{.48\textwidth}
\begin{subthm}{ex:integral-length>}
$\length \gamma
\le
\int_a^b|\gamma'(t)|\cdot dt$;
\end{subthm}
\end{minipage}
\hfill
\begin{minipage}{.48\textwidth}
\begin{subthm}{ex:integral-length<}
$\length \gamma
\ge
\int_a^b|\gamma'(t)|\cdot dt$.
\end{subthm}
\end{minipage}

\vskip1mm
Conclude that 
\[\length \gamma
=
\int_a^b|\gamma'(t)|\cdot dt.\eqlbl{eq:length}\]
\end{thm}

\begin{thm}{Advanced exercises}\label{adex:integral-length}

\begin{subthm}{adex:integral-length:a}Show that the formula \ref{eq:length} holds for any Lipschitz curve $\gamma\:[a,b]\z\to\mathbb{R}^3$.
\end{subthm}

\begin{subthm}{adex:integral-length:b}
Construct a non-constant curve $\gamma\:[a,b]\to\mathbb{R}^3$ such that $\gamma'(t)=0$ almost everywhere.
(In this case, the formula \ref{eq:length} does not hold for~$\gamma$ despite that both sides are defined.)
\end{subthm}

\end{thm}

\section{Nonrectifiable curves}

Let us describe the so-called \index{Koch snowflake}\emph{Koch snowflake} ---
a classical example of a nonrectifiable curve.

\begin{figure}[ht!]
\centering
\includegraphics{mppics/pic-225}
\end{figure}

Start with an equilateral triangle.
For each side, divide it into three segments of equal length, and then add an equilateral triangle with the middle segment as its base.
Repeat this construction recursively with the obtained polygons.
The Koch snowflake is the boundary of the union of all the polygons.
Two iterations and the resulting Koch snowflake are on the picture.

\begin{thm}{Exercise}\label{ex:nonrectifiable-curve}

\begin{subthm}{ex:nonrectifiable-curve:a} Show that the Koch snowflake is a simple closed curve; in particular, it can be parametrized by a circle.
\end{subthm}

\begin{subthm}{ex:nonrectifiable-curve:b} Show that the Koch snowflake is not rectifiable. 
\end{subthm}
\end{thm}
  
\section{Semicontinuity of length}

The lower limit of a sequence of real numbers $x_n$ will be denoted by
\[\liminf_{n\to\infty} x_n.\] 
It is defined as the lowest partial limit; that is, the lowest possible limit of a subsequence of $x_n$.
The lower limit is defined for any sequence of real numbers, and it lies in the extended real line $[-\infty,\infty]$.

\begin{thm}{Theorem}
Assume a sequence
of curves $\gamma_n\:[a,b]\to \spc{X}$ in a metric space $\spc{X}$ converges pointwise 
to a curve $\gamma_\infty\:[a,b]\to \spc{X}$;
that is, for any fixed $t \in [a,b]$, we have $\gamma_n(t)\z\to\gamma_\infty(t)$ as $n\to\infty$. 
Then 
$$\liminf_{n\to\infty} \length\gamma_n \ge \length\gamma_\infty.\eqlbl{eq:semicont-length}$$
\end{thm}

\begin{thm}{Corollary}\label{thm:length-semicont}
Length is lower semicontinuous with respect to the pointwise convergence of curves. 
\end{thm}

\parbf{Proof.}
Fix a partition $a=t_0<t_1<\dots<t_k=b$.
Set 
\begin{align*}\Sigma_n
&\df
\dist{\gamma_n(t_0)}{\gamma_n(t_1)}{}
+\dots+
\dist{\gamma_n(t_{k-1})}{\gamma_n(t_k)}{},
\\
\Sigma_\infty
&\df
\dist{\gamma_\infty(t_0)}{\gamma_\infty(t_1)}{}
+\dots+
\dist{\gamma_\infty(t_{k-1})}{\gamma_\infty(t_k)}{}.
\end{align*}

For each $i$ we have 
\[\dist{\gamma_n(t_{i-1})}{\gamma_n(t_i)}{}
\to
\dist{\gamma_\infty(t_{i-1})}{\gamma_\infty(t_i)}{},\]
and therefore
$\Sigma_n\to \Sigma_\infty$
as $n\to\infty$.
Note that 
$\Sigma_n\le\length\gamma_n$
for each~$n$.
Hence,
$$\liminf_{n\to\infty} \length\gamma_n \ge \Sigma_\infty.$$

Since the partition was arbitrary, the definition of length implies inequality \ref{eq:semicont-length}.
\qeds

\begin{wrapfigure}{o}{20 mm}
\vskip3mm
\centering
\includegraphics{mppics/pic-6}
\end{wrapfigure}

The inequality \ref{eq:semicont-length} might be strict.
For example, the diagonal $\gamma_\infty$ of the unit square 
can be approximated by stairs-like polygonal lines $\gamma_n$
with sides parallel to the sides of the square ($\gamma_6$ and $\gamma_\infty$ are in the picture).
In this case,
\[\length\gamma_\infty=\sqrt{2}
\quad
\text{and}\quad
\length\gamma_n=2
\quad
\text{for any}\quad
n.\]

\section{Arc-length parametrization}

We say that a curve $\gamma$ has an \index{arc-length parametrization}\emph{arc-length parametrization} (also called \index{natural parametrization}\emph{natural parametrization})
if 
\[t_2-t_1=\length \gamma|_{[t_1,t_2]}\]
for any two parameter values $t_1<t_2$;
that is, the arc of $\gamma$ from $t_1$ to $t_2$ has length $t_2-t_1$.

\begin{thm}{Exercise}\label{ex:cont-length}
Let  $\gamma\:[a,b]\to \spc{X}$ be a rectifiable curve in a metric space.
Given $t\in [a,b]$, denote by $s(t)$ the length of the arc $\gamma|_{[a,t]}$.
Show that the function $t\mapsto s(t)$ is continuous.

Conclude that $\gamma$ admits an arc-length parametrization.
\end{thm}

By Exercise~\ref{ex:integral-length}, a smooth curve $\gamma(t)=(x(t),y(t),z(t))$ is an arc-length parametrization if and only if it has unit velocity vector at all times;
that is, 
\[|\gamma'(t)|=\sqrt{x'(t)^2+y'(t)^2+z'(t)^2}=1\]
for all $t$; by that reason smooth curves equipped with an arc-length parametrization are also called \index{unit-speed curve}\emph{unit-speed} curves.
Observe that smooth unit-speed parametrizations are automatically regular (see Section~\ref{sec:Smooth curves}).

\begin{thm}{Proposition}\label{prop:arc-length-smooth}
If $t\mapsto \gamma(t)$ is a smooth curve, 
then its arc-length parametrization is also smooth and regular.
Moreover, the arc-length parameter $s$ of $\gamma$ can be written as an integral
\[s(t)=\int_{t_0}^t |\gamma'(\tau)|\cdot d\tau.\eqlbl{s(t)}\]
\end{thm}

Most of the time we will use $s$ for an arc-length parameter of a curve.

\parbf{Proof.}
Since $\gamma$ is smooth, $|\gamma'(t)|>0$ for any~$t$.
Therefore, the function $t\mapsto|\gamma'(t)|$ is smooth.

By the fundamental theorem of calculus, $s'(t)=|\gamma'(t)|$.
Therefore, $t\mapsto s(t)$ is a smooth increasing function with positive derivative.

By the inverse function theorem (\ref{thm:inverse}), the inverse function $s^{-1}(t)$ is also smooth
and $|(\gamma\circ s^{-1})'|\equiv1$.
Therefore, $\gamma\circ s^{-1}$ is a unit-speed reparametrization  of $\gamma$ by~$s$.
By construction, $\gamma\circ s^{-1}$ is smooth, and since $|(\gamma\circ s^{-1})'|\equiv1$, it is regular.
\qeds

\begin{thm}{Exercise}\label{ex:arc-length-helix}
Reparametrize the \index{helix}\emph{helix} 
\[\gamma_{a,b}(t)=(a\cdot\cos t,a\cdot \sin t, b\cdot t)\]
by its arc-length.
\end{thm}

We will be interested in the properties of curves that are invariant under reparametrizations.
Therefore, we can always assume that any given smooth curve comes with an arc-length parametrization.
A nice property of arc-length parametrizations is that they are almost canonical --- these parametrizations differ only by a sign and an additive constant.
For that reason, they make it easier to define parametrization-independent quantities.
This observation will be used in the definitions of curvature and torsion.

On the other hand, it is usually impossible to find an explicit arc-length parametrization.
Therefore, when it comes to calculations, it is often more convenient to use the original parametrization.

\section{Convex curves}

A simple plane curve is called \index{convex!curve}\emph{convex} if it bounds a convex region.
Since the boundary of any region is closed, any convex curve is either closed or open (see Section~\ref{sec:proper-curves}).

\begin{thm}{Proposition}\label{prop:convex-curve}
Assume that a closed convex curve $\alpha$ lies inside the domain bounded by a simple closed plane curve $\beta$.
Then
\[\length\alpha\le \length\beta.\]
\end{thm}

To prove Proposition \ref{prop:convex-curve} it is sufficient to show that the perimeter of any polygon inscribed in $\alpha$ is less or equal than the length of $\beta$.
Since any polygon inscribed in $\alpha$ is convex, it is sufficient to prove the following lemma.

\begin{thm}{Lemma}\label{lem:perimeter}
Assume that a convex polygon $P$ lies in a figure $F$ bounded by a simple closed curve.
Then 
\[\perim P\le \perim F,\]
where $\perim F$ denotes the perimeter of~$F$.
\end{thm}

\parbf{Proof.} 
A \index{chord}\emph{chord} in $F$ is defined to be a line segment in $F$ with endpoints in its boundary.
Suppose $F'$ is a figure obtained from $F$ by cutting it along a chord and removing one side.
By the triangle inequality, we have
\[\perim F'\le \perim F.\]

\begin{wrapfigure}{o}{24 mm}
\vskip-10mm
\centering
\includegraphics{mppics/pic-7}
\vskip3mm
\end{wrapfigure}

Observe that there is a decreasing sequence of figures 
\[F=F_0\supset F_1\supset\dots\supset F_n=P\]
such that $F_{i+1}$ is obtained from $F_{i}$ by cutting along a chord.
Therefore, 
\begin{align*}
\perim F=\perim F_0&\ge\perim F_1\ge\dots\ge\perim F_n=\perim P.
\end{align*}
\qedsf

\parit{Comment.}
Two other proofs of \ref{lem:perimeter} can be obtained by applying Crofton's formulas (see \ref{ex:convex-croftons}) and the nearest-point projection (see Lemma~\ref{lem:nearest-point-projection}).  

\begin{thm}{Corollary}\label{cor:convex=>rectifiable}
Any convex closed plane curve is rectifiable.  
\end{thm}

\parbf{Proof.}
Any closed curve is bounded.
Indeed, the curve can be described as an image of a loop $\alpha\:[0,1]\to\mathbb{R}^2$, $\alpha(t)=(x(t),y(t))$.
The coordinate functions $t\mapsto x(t)$ and $t\mapsto y(t)$ are continuous functions defined on $[0,1]$.
This implies that both coordinate functions are bounded by some constant~$C$.
Therefore, $\alpha$ lies in the square defined by the inequalities $|x|\le C$ and $|y|\le C$.

By Proposition~\ref{prop:convex-curve}, the length of the curve cannot exceed the perimeter of the square; hence the result.
\qeds

Recall that the \index{convex!hull}\emph{convex hull} of a set $X$ is the smallest convex set that contains $X$;
equivalently, the convex hull of $X$ is the intersection of all convex sets containing~$X$.

\begin{thm}{Exercise}\label{ex:convex-hull}
Let $\alpha$ be a simple closed plane curve.
Denote by $K$ the convex hull of $\alpha$; let $\beta$ be the boundary curve of~$K$.
Show that 
\[\length \alpha\ge \length \beta.\]

Try to show that the statement holds for arbitrary closed plane curves, assuming only that $K$ has a nonempty interior.
\end{thm}

\section{Crofton's formulas}
\label{sec:crofton}
\index{Crofton's formula}

For a function $f\: \mathbb{S}^1 \to \mathbb{R}$, we will denote its average value as $\overline{f(\vec u)}$;
that is,
\[\overline{f(\vec u)}=\frac1{2\cdot \pi}\cdot\int_{\vec u \in\mathbb{S}^1} f(\vec u).\]
For a vector $\vec w$ and a unit vector $\vec u$, we will denote by $\vec w_{\vec u}$ the orthogonal projection of $\vec w$ to the line in the direction of  $\vec u$;
that is,
\[\vec w_{\vec u}=\langle\vec u,\vec w\rangle\cdot\vec u.\] 

\begin{thm}{Theorem}
For any plane curve $\gamma$ we have
\[
\length\gamma
=\tfrac\pi2\cdot \overline{\length\gamma_{\vec u}}, \eqlbl{crofton-formula}
\]
where $\gamma_{\vec u}$ is the curve defined by $\gamma_{\vec u}(t) \df (\gamma (t))_{\vec u}$.
\end{thm}

\parbf{Proof.}
The magnitude of any vector ${\vec w}$ is proportional to the average magnitude of its projections; that is,
\[|{\vec w}|=k\cdot \overline{|{\vec w}_{\vec u}|}\]
for some $k \in \mathbb{R}$.
(The exact value of $k$ can be found by integration\footnote{It is the average value of $|\cos x|$.}, but we will find it differently.)
Assume $\gamma\:[a,b]\to\mathbb{R}^2$ is a smooth curve.
Note that 
\[\gamma_{\vec u}'(t)=(\gamma'(t))_{\vec u}
\quad\text{and}\quad
|\gamma'_{\vec u}(t)|=|\langle\vec u,\gamma'(t)\rangle|\]
for any $t \in [a,b]$. Then, according to Exercise~\ref{ex:integral-length},
\begin{align*}
\length\gamma
&=\int_a^b|\gamma'(t)|\cdot dt=
\\
&=\int_a^b  k\cdot \overline{|\gamma_{\vec u}'(t)|}\cdot dt=
\\
&=k\cdot \overline{\length\gamma_{\vec u}}.
\end{align*}

Since $k$ is a universal constant, we can compute it by taking $\gamma$ to be the unit circle.
In this case,
\[\length \gamma=2\cdot\pi.\]
For any unit plane vector ${\vec u}$, the curve $\gamma_{\vec u}$ runs back and forth along an interval of length~$2$.
Hence $\length\gamma_{\vec u}=4$ for any $\vec u$, and 
\[\overline{\length\gamma_{\vec u}} =4.\]
It follows that $2\cdot \pi =k\cdot 4$.
Therefore, \ref{crofton-formula} holds for smooth curves.

Applying the same argument together with \ref{adex:integral-length}, we get that \ref{crofton-formula} holds for arbitrary Lipschitz curves.
Further, since an arc length parametrization of any rectifiable curve is Lipschitz, \ref{ex:cont-length} implies \ref{crofton-formula} for arbitrary rectifiable curves.

It remains to consider the nonrectifiable case;
we have to show that 
\[\length\gamma=\infty
\quad\Longrightarrow\quad
\overline{\length\gamma_{\vec u}}=\infty.
\]

By the definition of length, 
\[\length\gamma_{\vec u}+\length\gamma_{\vec v}\ge \length\gamma\]
for any plane curve $\gamma$ and any pair $(\vec u , \vec v )$ of orthonormal vectors in $\mathbb{R}^2$.
Therefore, if $\gamma$ has infinite length, then the average length of $\gamma_{\vec u}$ is infinite as well.
\qeds

\begin{thm}{Exercise}\label{ex:convex-croftons}
Suppose a simple closed plane curve $\gamma$ bounds a figure~$F$.
Let $s$ be the average length of the projections of $F$ to lines.
Show that $\length\gamma\ge \pi \cdot s$.
Moreover, equality holds if and only if $\gamma$ is convex.

Use this statement to give another solution to Exercise~\ref{ex:convex-hull}.
\end{thm}

The following exercise gives analogous formulas in the Euclidean space.

As before, we denote by $\vec w_{\vec u}$ the orthogonal projection of $\vec w$ to the line passing thru the origin with direction $\vec u$.
Further, let us denote by $\vec w_{\vec u}^\bot$ the projection of $\vec w$ to the plane orthogonal to $\vec u$;
that is,
\[\vec w_\vec u^\bot=\vec w - \vec w_{\vec u}.\]

We will use the notation 
$\overline{f(\vec u)}$ for the average value
of a function $f$ defined on $\mathbb{S}^2$.

\begin{thm}{Advanced exercise}\label{adex:more-croftons}
Show that the length of a space curve is proportional to 
\begin{subthm}{}
the average length of its projections to all lines; that is,
\[\length\gamma=k_1\cdot\overline{\length\gamma_{\vec u}}\]
for some $k_1 \in \mathbb{R}$.
\end{subthm}
\begin{subthm}{}the average length of its projections to all planes; that is,
\[\length\gamma=k_2\cdot\overline{\length\gamma_{\vec u}^\bot}\]
for some $k_2 \in \mathbb{R}$.
\end{subthm}
Find the values $k_1$ and $k_2$.
\end{thm}

\section{Length metric}\label{sec:Length metric}

Let $\spc{X}$ be a metric space.
Given two points $x,y$ in $\spc{X}$, denote by $\ell(x,y)$ the greatest lower bound of lengths of all paths connecting $x$ to $y$; if there is no such path, then $\ell(x,y)=\infty$.

It is straightforward to see that the function $\ell$ satisfies all the axioms of a metric except it might take infinite values.
Therefore, if any two points in $\spc{X}$ can be connected by a rectifiable curve, then $\ell$ defines a new metric on $\spc{X}$;
in this case, $\ell$ is called the \index{induced length-metric}\emph{induced length-metric}.

Evidently, $\ell(x,y)\ge \dist{x}{y}{}$ for any pair of points $x,y\in \spc{X}$.
If the equality holds for all pairs, then the metric $\dist{{*}}{{*}}{}$ is said to be a \index{length-metric}\emph{length-metric}, and the corresponding metric space is called a \index{}\emph{length-metric space}.

\begin{thm}{Exercise}\label{ex:induced-is-length}
Let $\spc{X}$ be a metric space with a well-defined induced length-metric $(x,y)\mapsto \ell(x,y)$.
Show that $\ell$ is a length-metric.
\end{thm}

Most of the time we consider length-metric spaces.
In particular, the Euclidean space is a length-metric space.
A subspace $A$ of a length-metric space $\spc{X}$ is not necessarily a length-metric space;
the induced length distance between points $x$ and $y$ in the subspace $A$ will be denoted as $\dist{x}{y}A$;
that is, $\dist{x}{y}A$ is the greatest lower bound of the lengths of paths in $A$ from $x$ to~$y$.\index{10aaa@$\lvert x-y\rvert_A$ (length distance)}

\begin{thm}{Exercise}\label{ex:intrinsic-convex}
Let $A\subset \mathbb{R}^3$ be a closed subset.
Show that $A$ is convex if and only if
\[\dist{x}{y}A=\dist{x}{y}{\mathbb{R}^3}\]
for any $x,y\in A$
\end{thm}

\section{Spherical curves}

Let us denote by $\mathbb{S}^2$ the unit sphere in the space; that is,
\[\mathbb{S}^2=\set{(x,y,z)\in\mathbb{R}^3}{x^2+y^2+z^2=1}.\]
A space curve $\gamma$ is called \index{spherical!curve}\emph{spherical} if it runs in $\mathbb{S}^2$;
that is, $|\gamma(t)|=1$ for any~$t$.

Recall that $\measuredangle(u,v)$ denotes the angle between vectors $u$ and~$v$.

\begin{thm}{Observation}\label{obs:S2-length}
For any $u,v\in \mathbb{S}^2$, we have
\[\dist{u}{v}{\mathbb{S}^2}=\measuredangle(u,v).\]

\end{thm}

\parbf{Proof.}
Let $\gamma$ be the short arc of a great circle%
\footnote{A great circle is the intersection of the sphere with a plane passing thru its center.}
from $u$ to $v$ in $\mathbb{S}^2$.
Note that $\length\gamma=\measuredangle(u,v)$.
Therefore,
\[\dist{u}{v}{\mathbb{S}^2}\le\measuredangle(u,v).\]

It remains to prove the opposite inequality.
In other words, we need to show that given a polygonal line $\beta=p_0\dots p_n$ inscribed in $\gamma$ there is a polygonal line
$\beta_1=q_0\dots q_n$ inscribed in any given spherical path $\gamma_1$ connecting $u$ to $v$ such that 
\[\length\beta_1\ge \length \beta.\eqlbl{eq:length beta=<length beta}\]

Define $q_i$ as the first point on $\gamma_1$ such that $|u-p_i|=|u-q_i|$, but set $q_n=v$.
Clearly, $\beta_1$ is inscribed in $\gamma_1$, and, according to the triangle inequality for angles (\ref{thm:spherical-triangle-inq}), we have that 
\begin{align*}
 \measuredangle(q_{i-1},q_i) &\ge  \measuredangle (u, q_i) - \measuredangle ( u , q_{i-1})  =
\\
&= \measuredangle (u,p_i) - \measuredangle (u,p_{i-1}) =
\\
& = \measuredangle(p_{i-1},p_i).
\end{align*}
By the angle monotonicity (\ref{lem:angle-monotonicity}),
\[|q_{i-1}-q_i|\ge|p_{i-1}-p_i|\]
and \ref{eq:length beta=<length beta} follows.
\qeds

\begin{thm}{Hemisphere lemma}\label{lem:hemisphere}
Any closed spherical curve of length less than $2\cdot \pi$ lies in an open hemisphere. 
\end{thm}

This lemma will play a key role in the proof of Fenchel's theorem (\ref{thm:fenchel}).
The following proof is due to Stephanie Alexander;
it is not as simple as one may think.
Try to prove the lemma before reading further.

\parbf{Proof.}
Let $\gamma$ be a closed curve in $\mathbb{S}^2$ of length $2\cdot\ell$.
Suppose $\ell<\pi$.

\begin{wrapfigure}[8]{o}{35mm}
\vskip-0mm
\centering
\includegraphics{mppics/pic-52}
\end{wrapfigure}

Let us subdivide $\gamma$ into two arcs $\gamma_1$ and $\gamma_2$ of length $\ell$;
denote their common endpoints by $p$ and~$q$. 
By \ref{obs:S2-length},
\begin{align*}
\measuredangle(p,q)&\le \length \gamma_1=
\\
&= \ell<
\\
&<\pi.
\end{align*}

Denote by $z$ the midpoint between $p$ and $q$ in $\mathbb{S}^2$;
that is, $z$ is the midpoint of the short arc of a great circle from $p$ to $q$ in $\mathbb{S}^2$. 
We claim that $\gamma$ lies in the open hemisphere with the pole at~$z$.
If not, $\gamma$ intersects the equator at some point~$r$.
Without loss of generality, we may assume that $r$ lies on~$\gamma_1$. 

Rotate the arc $\gamma_1$ by the angle $\pi$ around the line thru $z$ and the center of the sphere.
The obtained arc $\gamma_1^{*}$ together with $\gamma_1$ forms a closed curve of length $2\cdot \ell$ passing thru $r$ and its antipodal point $r^{*}$.
Applying \ref{obs:S2-length} again, we get
\[\tfrac12\cdot\length \gamma=\ell\ge \measuredangle(r,r^{*})=\pi\] 
--- a contradiction.
\qeds

\begin{thm}{Exercise}\label{ex:antipodal}
Describe a simple closed spherical curve that does not pass thru a pair of antipodal points and does not lie in any hemisphere.
\end{thm}

\begin{thm}{Exercise}\label{ex:bisection-of-S2}
Suppose a simple closed spherical curve $\gamma$ divides $\mathbb{S}^2$ in two regions of equal area.
Show that 
\[\length\gamma\ge2\cdot\pi.\]
\end{thm}

\begin{thm}{Exercise}\label{ex:flaw}
Find a flaw in the solution of the following problem.
Come up with a correct argument.
\end{thm}

\parbf{Problem.}
Suppose a closed plane curve $\gamma$ has length at most 4.
Show that $\gamma$ lies in a unit disc.

\parbf{Wrong solution.}
It is sufficient to show that the \index{diameter}\emph{diameter} of $\gamma$ is at most 2;
that is, 
\[|p-q|\le 2\eqlbl{eq:|pq|=<2}\]
for any two points $p$ and $q$ on~$\gamma$.

The length of $\gamma$ cannot be smaller than the closed inscribed polygonal line which goes from $p$ to $q$ and back to~$p$.
Therefore, 
\[2\cdot |p-q|\le\length \gamma\le 4;\]
whence \ref{eq:|pq|=<2} follows.
\qedsf

\begin{thm}{Advanced exercises} \label{adex:crofton}
Given unit vectors ${\vec u},{\vec w}\in\mathbb{S}^2$, denote by ${\vec w}^*_{\vec u}$ the nearest point to ${\vec w}$ on the equator with the pole at ${\vec u}$;
in other words, if ${\vec w}^\perp_{\vec u}$ is the projection of ${\vec w}$ to the plane perpendicular to ${\vec u}$, then ${\vec w}^*_{\vec u}$ is the unit vector in the direction of ${\vec w}^\perp_{\vec u}$.
The vector ${\vec w}^*_{\vec u}$ is defined if ${\vec w}\ne\pm {\vec u}$.

\begin{subthm}{adex:crofton:crofton}
Show that for any rectifiable 
%one requires rectifiable for the right hand side to be defined for almost all U 
spherical curve $\gamma$ we have
\[\length\gamma=\overline{\length\gamma^*_{\vec u}},\]
where $\overline{\length\gamma^*_{\vec u}}$ denotes the average length of $\gamma^*_{\vec u}$ with ${\vec u}$ varying in~$\mathbb{S}^2$.
(This is a spherical analog of Crofton's formula.)
\end{subthm}

\begin{subthm}{adex:crofton:hemisphere}
Use \ref{SHORT.adex:crofton:crofton} to give another proof of the hemisphere lemma (\ref{lem:hemisphere}). 
\end{subthm}
 
\end{thm}

Spherical Crofton's formula can be rewritten the following way:
\[\length\gamma=\overline n\cdot \pi,\]
where $\overline n$ denotes the average number of intersection points of $\gamma$ with equators.
The equivalence can be proved using Levi's monotone convergence theorem.

%% file: curve-curvature-new.tex
\chapter{Curvature}
\label{chap:curve-curvature}

The term curvature is used for anything that measures how much a geometric object deviates from being \textit{straight}, whatever that may mean.
For curves, the curvature is a pointwise quantitative measure of how much the curve differs from a straight line.

\section{Acceleration of a unit-speed curve}

Recall that any smooth curve can be reparametrized by its arc-length (\ref{prop:arc-length-smooth}).
The obtained parametrized curve, say $\gamma$, remains smooth, and it has unit speed;
that is, $|\gamma'|=1$.
From mechanics, you might know that acceleration and velocity are perpendicular if the speed is constant.
Let us restate it.

\begin{thm}{Proposition}\label{prop:a'-pertp-a''}
Assume $\gamma$ is a smooth unit-speed space curve.
Then $\gamma'(s)\perp \gamma''(s)$ for any~$s$.
\end{thm}

The \index{scalar product}\emph{scalar product} (also known as \textit{dot product}) of two vectors $\vec v$ and $\vec w$ will be denoted by $\langle \vec v,\vec w\rangle$.
Recall that the derivative of a scalar product satisfies the product rule;
that is, if $\vec v=\vec v(t)$ and $\vec w=\vec w(t)$ are smooth vector-valued functions of a real parameter $t$, then
\[\langle \vec v,\vec w\rangle'=\langle \vec v',\vec w\rangle+\langle \vec v,\vec w'\rangle.\]

\parbf{Proof.}
The identity $|\gamma'|=1$ can be rewritten as $\langle\gamma',\gamma'\rangle=1$.
Differentiating both sides, we get
$2\cdot\langle\gamma'',\gamma'\rangle=\langle\gamma',\gamma'\rangle'=0$;
whence $\gamma''\perp\gamma'$.
\qeds

\section{Curvature}\label{sec:curvature}

For a unit-speed smooth space curve $\gamma$, the magnitude of its acceleration $|\gamma''(s)|$ is called its \index{10k@$\kur$ (curvature)}\index{curvature}\emph{curvature} at  time~$s$.
If $\gamma$ is simple, then we can say that $|\gamma''(s)|$ is the curvature at the point $p=\gamma(s)$ without ambiguity.
The curvature is usually denoted by $\kur(s)$ or $\kur(s)_\gamma$, and, in the case of simple curves, it might be also denoted by $\kur(p)$ or $\kur(p)_\gamma$.

\begin{thm}{Exercise}\label{ex:zero-curvature-curve}
Show that a smooth simple space curve has zero curvature at each point if and only if it is a segment of a straight line.
\end{thm}

\begin{thm}{Exercise}\label{ex:scaled-curvature}
Let $\gamma$ be a smooth simple space curve, and let $\gamma_{\lambda}$ be a scaled copy of $\gamma$ with factor $\lambda >0$;
that is, $\gamma_{\lambda}(t)=\lambda \cdot\gamma(t)$ for any~$t$.
Show that 
\[\kur(\lambda \cdot p )_{\gamma_{\lambda}}
= \frac{\kur(p)_{\gamma}}\lambda\]
for any $p \in \gamma$.
\end{thm}

\begin{thm}{Exercise}\label{ex:curvature-of-spherical-curve}
Show that any smooth spherical curve has curvature at least 1.
\end{thm}

\section{Tangent indicatrix}\label{sec:Tangent indicatrix}

Let $\gamma$ be a smooth space curve.
The curve 
\[\tan(t)=\tfrac{\gamma'(t)}{|\gamma'(t)|};
\eqlbl{eq:tantrix}\] 
it is called a \index{tangent!indicatrix}\emph{tangent indicatrix} of~$\gamma$.
Note that $|\tan(t)|=1$ for any $t$;
that is, $\tan$ is a spherical curve.

If $s\mapsto \gamma(s)$ is a unit-speed parametrization, then $\tan(s)=\gamma'(s)$.
In this case, we have the following expression for the curvature:
\[\kur(s)\z=|\tan'(s)|\z=|\gamma''(s)|.\]

For a general parametrization $t\mapsto \gamma(t)$,
we have instead
\[ \kur(t)=\frac{|\tan'(t)|}{|\gamma'(t)|}.\eqlbl{eq:curvature}\]
Indeed, for an arc-length reparametrization by $s(t)$, we have $s'(t)=|\gamma'(t)|$.
Therefore,
\begin{align*}
\kur&=\left|\frac{d\tan}{ ds}\right|=
\left|\frac{d\tan}{ dt}\right|/\left|\frac{ds}{ dt}\right|=
\frac{|\tan'|}{|\gamma'|}.
\end{align*}

Recall that smooth curves are required to have a smooth and regular parametrization; see \ref{sec:Smooth curves}.
If $\gamma$ is smooth, then $t\mapsto \tan(t)$ is a smooth parametrization,
which is regular if $|\tan'|=\kur\cdot|\gamma'|$ does not vanish.
It follows that \textit{if the curvature of a smooth curve does not vanish, then its tangent indicatrix is a smooth curve}.

\begin{thm}{Exercise}\label{ex:curvature-formulas}
Use the formulas \ref{eq:tantrix} and \ref{eq:curvature} to show that 
for any smooth space curve $\gamma$ we have the following expressions for its curvature:
\setlength{\columnseprule}{0.4pt}
\begin{multicols}{2}

\begin{subthm}{ex:curvature-formulas:a} 
\[\kur=\frac{|\vec w|}{|\gamma'|^2},\]
where $\vec w$ is the projection of $\gamma''(t)$ to the plane normal to $\gamma'(t)$;
\end{subthm}

\begin{subthm}{ex:curvature-formulas:b}
\[\kur=\frac{|\gamma''\times \gamma'|}{|\gamma'|^{3}},\]
where $\times$ denotes the {}\emph{vector product} (also known as {}\emph{cross product}).
\end{subthm}
\end{multicols}
\end{thm}

\begin{thm}{Exercise}\label{ex:curvature-graph}
Apply the formulas in the previous exercise to show that if $f$ is a smooth real function,
then its graph $y=f(x)$  has curvature
\[\kur=\frac{|f''(x)|}{(1+f'(x)^2)^{\frac32}}\]
at the point $(x,f(x))$.
\end{thm}

\begin{thm}{Advanced exercise}\label{ex:approximation-const-curvature}
Show that any smooth curve $\gamma\:\mathbb{I}\z\to\mathbb{R}^3$ with curvature at most 1 can be approximated by smooth curves with constant curvature~1.

In other words, construct a sequence of smooth curves $\gamma_n\:\mathbb{I}\to\mathbb{R}^3$ with constant curvature $1$ such that $\gamma_n(t)\to \gamma(t)$ for any $t$ as $n\to\infty$.
\end{thm}

\section{Tangents}

{

\begin{wrapfigure}{r}{32 mm}
\vskip-10mm
\centering
\includegraphics{mppics/pic-3400}
\vskip0mm
\end{wrapfigure}

Let $\gamma$ be a smooth space curve, and let $\tan$ be its tangent indicatrix.
The line thru $\gamma(t)$ in the direction of $\tan(t)$ is called the \index{tangent!line}\emph{tangent line} to $\gamma$  at~$t$.
Any vector proportional to $\tan(t)$ is called the \index{tangent!vector}\emph{tangent} to $\gamma$ at~$t$.

The tangent line could be also defined as a unique line that has \index{order of contact}\emph{first-order contact} with $\gamma$ at $t$;
that is, $\rho(\ell)=o(\ell)$, where $\rho(\ell)$ denotes the distance from $\gamma(t+\ell)$ to the line.

}

\begin{thm}{Advanced exercise}\label{ex:no-parallel-tangents}
Construct a smooth closed space curve without parallel tangent lines.
\end{thm}

We say that smooth curves $\gamma_1$ and $\gamma_2$ are \index{tangent!curves}\emph{tangent} at $s_1$ and $s_2$
if $\gamma_1(s_1)=\gamma_2(s_2)$ and the tangent line of $\gamma_1$ at $s_1$ coincides with the tangent line of $\gamma_2$ at $s_2$;
if both curves are simple we can also say that they are tangent at the point $p=\gamma_1(s_1)\z=\gamma_2(s_2)$ without ambiguity.

\section{Total curvature}\label{sec:Total curvature}

Let $\gamma\:\mathbb{I}\to\mathbb{R}^3$ be a smooth unit-speed curve.
The integral 
\[\tc\gamma\df\int_{\mathbb{I}}\kur(s)\cdot ds\]
is called the \index{total!curvature}\emph{total curvature of}\label{page:total curvature of:smooth-def}
$\gamma$, and it measures the total change of direction along $\gamma$.

Rewriting the above integral using a change of variables produces a formula for a general parametrization $t\mapsto \gamma(t)$:
\[\tc\gamma\df\int_{\mathbb{I}}\kur(t)\cdot|\gamma'(t)| \cdot dt.
\eqlbl{eq:tocurv}\]

\begin{thm}{Exercise}\label{ex:helix-curvature}
Find the curvature of the helix 
\[\gamma_{a,b}(t)=(a\cdot \cos t,a\cdot \sin t,b\cdot t),\]
its tangent indicatrix, and the total curvature of its arc $\gamma_{a,b}|_{[0,2\cdot\pi]}$.
\end{thm}

Note that for a unit-speed smooth curve,
the speed of its tangent indicatrix equals its curvature.
Therefore, we get the following.

\begin{thm}{Observation}\label{obs:tantrix}
The total curvature of a smooth curve is the length of its tangent indicatrix.
\end{thm}

\begin{thm}{Fenchel's theorem}
\label{thm:fenchel}
\index{Fenchel's theorem}
The total curvature of any closed smooth space curve is at least $2\cdot\pi$.
\end{thm}

\parbf{Proof.}
Fix a closed smooth space curve~$\gamma$.
We can assume that $\gamma$ is described by a unit-speed loop $\gamma\:[a,b]\to \mathbb{R}^3$;
in this case, $\gamma(a)=\gamma(b)$ and $\gamma'(a)=\gamma'(b)$.

Consider its tangent indicatrix $\tan=\gamma'$.
Recall that $|\tan(s)|=1$ for any $s$; that is, $\tan$ is a closed spherical curve.

Let us show that $\tan$ cannot lie in a hemisphere.
Arguing by contradiction, we can assume that it lies in the hemisphere defined by the inequality $z>0$ in $(x,y,z)$-coordinates.
In other words, if $\gamma(t)=(x(t), y(t), z(t))$, then $z'(t)>0$ for any~$t$.
Therefore,
\[z(b)-z(a)=\int_a^b z'(s)\cdot ds>0.\]
In particular, $\gamma(a)\ne \gamma(b)$ --- a contradiction.

Applying the observation (\ref{obs:tantrix}) and the hemisphere lemma (\ref{lem:hemisphere}), we get  
\[\tc\gamma=\length \tan\ge2\cdot\pi.\]
\qedsf

\begin{thm}{Exercise}\label{ex:length>=2pi}
Show that a closed space curve $\gamma$ with curvature at most~$1$ cannot be shorter than the unit circle;
that is, 
\[\length\gamma\ge 2\cdot \pi.\]

\end{thm}

\begin{thm}{Advanced exercise}\label{ex:gamma/|gamma|}
Suppose $\gamma$ is a smooth space curve that does not pass thru the origin.
Consider the spherical curve $\sigma$ defined by $\sigma(t)\z\df\frac{\gamma(t)}{|\gamma(t)|}$.
Show that 
\[\length \sigma< \tc\gamma+\pi.\]
Moreover, if $\gamma$ is closed, then
\[\length \sigma\le \tc\gamma.\]
\end{thm}

The last inequality gives an alternative proof of Fenchel's theorem.
Indeed, without loss of generality, we can assume that the origin lies on a chord of~$\gamma$.
In this case, the closed spherical curve $\sigma$ goes from a point to its antipode and comes back; 
it takes length at least $\pi$ each way, 
whence 
\[\length\sigma\ge 2\cdot\pi.\]

Recall that the curvature of a spherical curve is at least $1$
(see \ref{ex:curvature-of-spherical-curve}).
In particular, the length of a spherical curve cannot exceed its total curvature.
The following theorem shows that the same inequality holds for \textit{closed} curves in a unit ball.

\begin{thm}{DNA theorem}\label{thm:DNA}
Let $\gamma$ be a smooth closed curve that lies in a unit ball.
Then 
\[\tc\gamma\ge \length\gamma.\]

\end{thm}

Several proofs of this theorem are collected by Serge Tabachnikov~\cite{tabachnikov}.
The 2-dimensional case of this theorem was proved by Istv\'{a}n F\'{a}ry \cite{fary1950}.
It was generalized by Don Chakerian \cite{chakerian1962} to higher dimensions.
The following exercise guides you thru another proof by him \cite{chakerian1964}.
Yet another proof is given in Section~\ref{sec:DNA-poly}.

\begin{thm}{Exercise}\label{ex:DNA}
Let $\gamma\:[0,\ell]\to\mathbb{R}^3$ be a smooth unit-speed closed curve that lies in the unit ball; that is, $|\gamma|\le 1$.

\begin{subthm}{ex:DNA:c''c>=k}
Show that 
\[\langle\gamma''(s),\gamma(s)\rangle\ge-\kur(s)\]
for any~$s$.
\end{subthm}

\begin{subthm}{ex:DNA:int>=length-tc}
Use part \ref{SHORT.ex:DNA:c''c>=k} to show that 
\[\int_0^\ell\langle\gamma(s),\gamma'(s)\rangle'\cdot ds\ge
\ell-\tc\gamma.\]

\end{subthm}

\begin{subthm}{ex:DNA:end}
Suppose $\gamma(0)=\gamma(\ell)$ and $\gamma'(0)=\gamma'(\ell)$.
Show that 
\[\int_0^\ell\langle\gamma(s),\gamma'(s)\rangle'\cdot ds=0.\]
Use this equality together with  part \ref{SHORT.ex:DNA:int>=length-tc} to prove \ref{thm:DNA}.
\end{subthm}
\end{thm}

\section{Convex curves}

In this section, we show that the tangent indicatrix of a convex curve rotates monotonically. 
The following exercise provides a key observation for the proof.

\begin{thm}{Exercise}\label{ex:tangent-support}
Let $\gamma$ be a smooth convex plane curve;
denote by $F$ the convex set bounded by $\gamma$.
Show that a line $\ell$ is tangent to $\gamma$ at point $p$ if and only if it \index{supporting!plane}\emph{supports} $F$ at $p$;
that is, $\ell\ni p$ and $F$ lies in a closed half-plane bounded by $\ell$.
\end{thm}

Recall that a map is monotone if the inverse image of any point in the target space is a connected set (in particular, nonempty).

\begin{thm}{Proposition}\label{prop:convex-monotone}
Let $\gamma$ be a smooth convex plane curve.

\begin{subthm}{prop:convex-monotone:closed}
Suppose $\gamma$ is closed and it is parametrized as $\gamma\:\mathbb{S}^1\to \mathbb{R}^2$.
Then its tangent indicatrix $\tan\:\mathbb{S}^1\to\mathbb{S}^1$ is a monotone map.
\end{subthm}

\begin{subthm}{prop:convex-monotone:open}
Suppose $\gamma$ is open, it is parametrized as $\gamma\:\mathbb{R}\to \mathbb{R}^2$.
Then its tangent indicatrix $\tan\:\mathbb{R}\to \mathbb{S}^1$ defines a monotone map to an interval in a closed semicircle.
\end{subthm}

\end{thm}

The following corollary says that \textit{for convex curves we have equality in Fenchel's theorem} (\ref{thm:fenchel}).
Later, in \ref{prop:fenchel=}, we will show that equality holds \textit{only} for convex curves.

\begin{thm}{Corollary}\label{cor:fenchel=convex}
Let $\gamma$ be a convex plane curve.

\begin{subthm}{}
If $\gamma$ is closed, then $\tc\gamma=2\cdot\pi$.
\end{subthm}

\begin{subthm}{}
If $\gamma$ is open, then $\tc\gamma\le\pi$.
\end{subthm}

\end{thm}

\parbf{Proof.}
Follows from \ref{prop:convex-monotone:closed}, \ref{obs:tantrix} and \ref{ex:integral-length-0}.
\qeds

\begin{wrapfigure}{r}{32 mm}
\vskip-0mm
\centering
\includegraphics{mppics/pic-3500}
\vskip0mm
\end{wrapfigure}

\parbf{Proof of \ref{prop:convex-monotone};} \ref{SHORT.prop:convex-monotone:closed}.
Since $\gamma$ is closed, it bounds a compact convex set $F$.
We may assume that $F$ lies on the left from $\gamma$.

Choose a unit vector $\vec u$.
Let $(x,y)$ be the coordinates of the plane with $x$-axis in the direction of $\vec u$.

From Exercise~\ref{ex:tangent-support}, it follows that  $\vec u=\tan(s)$ if and only if $p=\gamma(s)$ is a minimum point of the $y$-coordinate function on~$F$.
Indeed suppose $p$ is a minimum point.
Let $\ell$ be the horizontal line thru $p$.
Then $\ell$ supports $F$ at $p$.
By the exercise, $\ell$ is tangent to $\gamma$ at $p$.
Since $F$ lies on the left from $\gamma$, we get that $\tan(s)=\vec u$.
And the other way around, if $\tan(s)=\vec u$ then, the tangent line at $p=\gamma(s)$ is horizontal,
and, by the exercise, it supports $F$.
Since $F$ lies on the left from $\gamma$, we get that the $y$-coordinate on $F$ has a minimum at $p$.

Since $F$ is compact, there is a minimum point $p=\gamma(s)$ for the $y$-coordinate.
The point $p$ might be unique, in this case, $\tan^{-1}\{\vec u\}=\{s\}$,
or $F$ might have a line segment of minimal points in $\gamma$, in this case, $\tan^{-1}\{\vec u\}$ is an arc of $\mathbb{S}^1$.
It follows that $\tan\:\mathbb{S}^1\to \mathbb{S}^1$ is monotone.

\parit{\ref{SHORT.prop:convex-monotone:open}} 
The same argument in \ref{SHORT.prop:convex-monotone:closed} shows that $\tan$ is a monotone map to its image.
Evidently, the image is a connected set in $\mathbb{S}^1$.
It remains to show that the image lies in a semicircle;
in other words, 
\[\measuredangle(\vec w,\tan(s))\ge\tfrac\pi2
\quad\text{for some}\quad \vec w
\quad\text{and any}\quad
s.
\eqlbl{eq:<(w,tan).pi/2}
\]

\begin{figure}[ht!]
\centering
\includegraphics{mppics/pic-3502}
\end{figure}

Since $\gamma$ is open, it bounds an unbounded convex closed region $F$.
As before we assume that $F$ lies on the left from $\gamma$.

There is a half-line, say $h$, that lies in~$F$.
Indeed, we can assume that the origin $o$ lies in~$F$.
Consider a sequence of points $q_n\in F$ such that $|q_n|\z\to \infty$ as $n\to \infty$.
Denote by $\vec v_n$ the unit vector in the direction of $q_n$; that is $\vec v_n=\tfrac{q_n}{|q_n|}$.

Since the unit circle is compact, we can pass to a subsequence of $q_n$ such that $\vec v_n$ converges to a unit vector, say $\vec v$.
Let us draw a half-line $h$ from $o$ in the direction of $\vec{v}$.
Any point on $h$ can be approximated by points from the segments $[o, q_n]$ as $n \to \infty$.
Since the set $F$ is closed, the half-line $h$ lies in $F$

Let $\vec w$ be the counterclockwise rotation of $\vec v$ by an angle $\tfrac\pi 2$, let
$\ell$ be the tangent line to $\gamma$ at the point $p=\gamma(s)$,
and let $H$ be the closed left half-plane bounded by $\ell$;
that is, $H$ lies to the left of $\ell$ in the direction of $\tan(s)$.
The same reasoning as in \ref{SHORT.prop:convex-monotone:closed} shows that $F$, and therefore $h$,
lie in $H$.
In particular, $\vec v$ points from $p$ into $H$, which is equivalent to~\ref{eq:<(w,tan).pi/2}.
\qeds

\section{Bow lemma}

The following lemma proved by Erhard Schmidt \cite{schmidt}; it generalizes a result by Axel Schur \cite{shur}.

\begin{wrapfigure}{r}{39 mm}
\vskip-6mm
\centering
\includegraphics{mppics/pic-3510}
\vskip-6mm
\end{wrapfigure}

This lemma is a differential-geometric analog of the so-called {}\emph{arm lemma} of Augustin-Louis Cauchy.
The arm lemma says that 
\textit{if $p_0\dots p_n$ is a convex plane polygon and $q_0\dots q_n$ is a polygonal line in the space such that 
\begin{align*}
|p_i-p_{i-1}|&=|q_i-q_{i-1}|,
\\
\measuredangle\hinge{p_i}{p_{i+1}}{p_{i-1}}&\le \measuredangle\hinge{q_i}{q_{i+1}}{q_{i-1}}
\end{align*}
for each $i$, then $|p_0-p_n|\le |q_0-q_n|$.}
(Intuitively, if you extend all the joints in your arm, then the distance from your shoulder to the tip of your middle finger increases.) 

\begin{thm}{Lemma}\label{lem:bow}\index{bow lemma}
Let $\gamma_1\:[a,b]\to\mathbb{R}^2$ and $\gamma_2\:[a,b] \to\mathbb{R}^3$ be two smooth unit-speed curves.
Suppose $\kur(s)_{\gamma_1}\ge\kur(s)_{\gamma_2}$ for any $s$ 
and the curve
$\gamma_1$ is an arc of a convex curve; that is, it runs in the boundary of a convex plane figure.
Then the distance between the endpoints of $\gamma_1$ cannot exceed the  distance between the endpoints of $\gamma_2$; that is,
\[|\gamma_1(b)-\gamma_1(a)|\le |\gamma_2(b)-\gamma_2(a)|.\]

\end{thm}

The following exercise states that the condition that $\gamma_1$ is a convex arc is necessary.
It is instructive to do this exercise before reading the proof of the lemma.

{\sloppy 

\begin{thm}{Exercise}\label{ex:anti-bow}
Construct two simple smooth unit-speed plane curves $\gamma_1,\gamma_2\:[a,b]\to\mathbb{R}^2$ such that $\kur(s)_{\gamma_1}>\kur(s)_{\gamma_2}>0$ for any $s$ and
\[|\gamma_1(b)-\gamma_1(a)|> |\gamma_2(b)-\gamma_2(a)|.\]
\end{thm}

}

\parbf{Proof of \ref{lem:bow}.}
We can assume that $\gamma_1(a)\ne \gamma_1(b)$;
otherwise the statement is evident.
By convexity, the curve $\gamma_1$ lies on one side of the line $\ell$ thru $\gamma_1(a)$ and $\gamma_1(b)$;
we may assume $\ell$ is horizontal and $\gamma_1$ lies below $\ell$.

Let $\gamma_1(s_0)$ be the lowest point on $\gamma_1$;
that is, $\gamma_1(s_0)$ has minimal $y$-coordinate.

Denote by $\tan_1$ and $\tan_2$ the tangent indicatrixes of $\gamma_1$ and $\gamma_2$, respectively.
Consider two unit vectors 
\[\vec u_1=\tan_1(s_0)=\gamma_1'(s_0)
\quad\text{and}\quad
\vec u_2=\tan_2(s_0)=\gamma_2'(s_0).\]
Note that $\gamma_1(b)$ lies in the direction of $\vec u_1$ from $\gamma_1(a)$.

Let us show that
\[\measuredangle(\gamma'_1(s),\vec u_1)\ge \measuredangle(\gamma'_2(s),\vec u_2)
\eqlbl{<gamma',u}
\]
for any $s$.
We will prove it for $s\le s_0$; the case $s\ge s_0$ is analogous.

Note that
\[\measuredangle(\gamma'_1(s),\vec u_1)=\measuredangle(\tan_1(s),\vec u_1)=\length (\tan_1|_{[s,s_0]}).\eqlbl{<=length}\]
for any $s\le s_0$.
Indeed, by \ref{ex:tangent-support}, the $y$-coordinate of $\gamma_1$ is nonincreasing in the interval $[a,s_0]$.
Therefore, the arc $\tan_1|_{[a,s_0]}$ lies in one of the unit semicircles with endpoints $\vec u_1$ and $-\vec u_1$.
It remains to apply \ref{obs:tantrix}, \ref{prop:convex-monotone}, and \ref{ex:integral-length-0}.

By \ref{obs:S2-length}, we also have 
\[\measuredangle(\gamma'_2(s),\vec u_2)=\measuredangle(\tan_2(s),\vec u_2)\le \length (\tan_2|_{[s,s_0]}).\eqlbl{<=<length}\]

{

\begin{wrapfigure}{r}{39 mm}
\vskip0mm
\centering
\includegraphics{mppics/pic-57}
\vskip0mm
\end{wrapfigure}

Further,
\begin{align*}
\length (\tan_1|_{[s,s_0]})
&=\int_s^{s_0}|\tan_1'(t)|\cdot d t=
\\
&=\int_s^{s_0}\kur_1(t)\cdot d t\ge
\\
&\ge
\int_s^{s_0}\kur_2(t)\cdot d t=
\\
&=\int_s^{s_0}|\tan_2'(t)|\cdot d t= 
\\
&=\length (\tan_2|_{[s,s_0]}).
\end{align*}
This inequality, together with \ref{<=length} and \ref{<=<length}, implies \ref{<gamma',u}.
}

Since $1=|\gamma_1'(s)|=|\gamma_2'(s)|=|\vec u_1|=|\vec u_2|$,
we have 
\[\langle\gamma'_1(s),\vec u_1\rangle=\cos \measuredangle(\gamma'_1(s),\vec u_1)
\quad\text{and}\quad
\langle\gamma'_2(s),\vec u_2\rangle=\cos \measuredangle(\gamma'_2(s),\vec u_2).
\]
The cosine is decreasing in the interval $[0,\pi]$; therefore, \ref{<gamma',u} implies 
\[\langle\gamma'_1(s),\vec u_1\rangle\le \langle\gamma'_2(s),\vec u_2\rangle\eqlbl{<gamma',u>}\]
for any~$s$.

Further, since $\gamma_1(b)$ lies in the direction of $\vec u_1$ from $\gamma_1(a)$, we have that
\[|\gamma_1(b)-\gamma_1(a)|=\langle \vec u_1,\gamma_1(b)-\gamma_1(a)\rangle.\]
Since $\vec u_2$ is a unit vector, we have that
\[|\gamma_2(b)-\gamma_2(a)|\ge\langle \vec u_2,\gamma_2(b)-\gamma_2(a)\rangle.\]

Integrating \ref{<gamma',u>}, we get 
\begin{align*}
|\gamma_1(b)-\gamma_1(a)|&=\langle \vec u_1,\gamma_1(b)-\gamma_1(a)\rangle=
\\
&=
\int_a^b\langle \vec u_1,\gamma'_1(s)\rangle\cdot ds \le 
\int_a^b\langle \vec u_2,\gamma'_2(s)\rangle\cdot ds 
=
\\
&=\langle \vec u_2,\gamma_2(b)-\gamma_2(a)\rangle
\le |\gamma_2(b)-\gamma_2(a)|.
\end{align*}
\qedsf

\begin{thm}{Advanced exercise}\label{ex:bow'}
Let $\gamma_1$ and $\gamma_2$ be as in the bow lemma (\ref{lem:bow});
denote by $\tan_1$ and $\tan_2$ their tangent indicatrixes.

{

\begin{wrapfigure}{r}{50 mm}
\vskip-0mm
\centering
\includegraphics{mppics/pic-251}
\vskip-0mm
\end{wrapfigure}

Set
%\[\vec w_i=\gamma_i(b)-\gamma_i(a),\quad \alpha_i=\measuredangle(\tan_i(a),\vec w_i)\quad\text{и}\quad\beta_i&=\measuredangle(\tan_i(b),\vec w_i).\]
\begin{align*}
\vec w_i&=\gamma_i(b)-\gamma_i(a),
\\
\alpha_i&=\measuredangle(\tan_i(a),\vec w_i),
\\
\beta_i&=\measuredangle(\tan_i(b),\vec w_i).
\end{align*}
%\begin{align*}
%\vec w_1&=\gamma_1(b)-\gamma_1(a),
%&
%\vec w_2&=\gamma_2(b)-\gamma_2(a),
%\\
%\alpha_1&=\measuredangle(\tan_1(a),\vec w_1),
%&
%\alpha_2&=\measuredangle(\tan_2(a),\vec w_2),
%\\
%\beta_1&=\measuredangle(\tan_1(b),\vec w_1),
%&
%\beta_2&=\measuredangle(\tan_2(b),\vec w_2).
%\end{align*}

}

\begin{subthm}{ex:bow'+}
Suppose that $\beta_1\le\tfrac\pi2$.
Show that $\alpha_1\ge \alpha_2$.
\end{subthm}

\begin{subthm}{ex:bow'-} Construct an example, showing that the inequality $\alpha_1\ge \alpha_2$ does not hold in general.
\end{subthm}

\end{thm}

\begin{thm}{Exercise}\label{ex:length-dist}
Let $\gamma\:[a,b]\to \mathbb{R}^3$ be a smooth curve and $0\z<\theta\z\le\tfrac\pi2$.
Assume 
\[\tc\gamma\le 2\cdot\theta.\]

\begin{subthm}{ex:length-dist:>} Show that
\[|\gamma(b)-\gamma(a)|> \cos\theta\cdot\length\gamma.\]
\end{subthm}

\begin{subthm}{ex:length-dist:self-intersection:>pi}
Show that if a smooth curve $\gamma\:[a,b]\to\mathbb{R}^3$ has a self-intersection, then $\tc\gamma>\pi$.
Draw a smooth plane curve $\gamma$ with a self-intersection, such that $\tc\gamma<2\cdot\pi$.
\end{subthm}

\begin{subthm}{ex:length-dist:=} Show that the inequality in \ref{SHORT.ex:length-dist:>} is optimal; that is, given 
$\theta$ there is a smooth curve $\gamma$ such that $\tc\gamma\le 
2\cdot\theta$, and $\frac{|\gamma(b)-\gamma(a)|}{\length\gamma}$ is arbitrarily 
close to $\cos\theta$.
\end{subthm}

\end{thm}

\begin{thm}{Exercise}\label{ex:schwartz}
Let $p$ and $q$ be points on a unit circle dividing it in two arcs with lengths $\ell_1<\ell_2$.
Suppose a space curve $\gamma$ connects $p$ to $q$ and has curvature at most $1$.
Show that either
\[\length \gamma\le \ell_1
\quad\text{or}\quad
\length \gamma\ge \ell_2.
\]
\end{thm}

The following exercise generalizes \ref{ex:length>=2pi}.

\begin{thm}{Exercise}\label{ex:loop}
Suppose $\gamma\:[a,b]\to \mathbb{R}^3$ is a smooth loop with curvature at most~1.
Show that 
\[\length\gamma\ge2\cdot\pi.\]

\end{thm}

\begin{thm}{Exercise}\label{ex:bow-upper}
Let $\kur$ be a smooth nonnegative function defined on $[0,\ell]$.
Show that there is a smooth unit-speed curve $\gamma\:[0,\ell]\to\mathbb{R}^3$ with curvature $\kur(s)$ for any $s$ such that the distance $|\gamma(\ell)-\gamma(0)|$ is arbitrarily close to~$\ell$.
\end{thm}

\begin{wrapfigure}{r}{39 mm}
\vskip-0mm
\centering
\includegraphics{mppics/pic-283}
\vskip0mm
\end{wrapfigure}

\begin{thm}{Advanced exercise}\label{ex:gromov-twist}
Let $\gamma$ be a closed smooth space curve with curvature at most $2$.
Suppose $|\gamma(t)|\le 1$ for any $t$.
Show that if $\gamma(t)\ne 0$, then 
\[|\gamma(t)| \le \sin (\alpha(t))\]
where $\alpha(t)$ denotes the angle between  $\gamma(t)$ and $\gamma'(t)$.
\end{thm}

%% file: poly.tex
\chapter{Polygonal lines}
\label{chap:poly}

This chapter reinterprets the curvature of curves via angles of inscribed polygonal lines;
it should help to build a firm geometric intuition about curvature.

\section{Piecewise smooth curves}

\begin{wrapfigure}{o}{25 mm}
\vskip-0mm
\centering
\includegraphics{mppics/pic-54}
\end{wrapfigure}

Assume $\alpha\:[a,b]\to \mathbb{R}^3$ and $\beta\:[b,c]\z\to \mathbb{R}^3$ are two curves such that $\alpha(b)\z=\beta(b)$.
These two curves can be combined into one $\gamma\:[a,c]\z\to \mathbb{R}^3$ by the rule 
\[\gamma(t)=
\begin{cases}
\alpha(t)&\text{if}\quad t\le b,
\\
\beta(t)&\text{if}\quad t> b.
\end{cases}
\]
The obtained curve $\gamma$ is called the 
\index{concatenation}\emph{concatenation} of $\alpha$ and $\beta$.
(The condition $\alpha(b)=\beta(b)$ ensures that the map $t\mapsto\gamma(t)$ is continuous.)

The same definition of concatenation can be applied if $\alpha$ and/or $\beta$ are defined on semiopen intervals 
$(a,b]$ and/or $[b,c)$.

The assumption that the time intervals of $\alpha$ and $\beta$ fit together is not essential --- one can concatenate if the endpoint of $\alpha$ coincides with the starting point of~$\beta$.
To do this, the time intervals of the curves should be shifted so that they fit together. 

If in addition $\beta(c)=\alpha(a)$, then we can do cyclic concatenation of these curves;
this way we obtain a closed curve.

If $\alpha'(b)$ and $\beta'(b)$ are defined, then the angle $\theta\z=\measuredangle(\alpha'(b),\beta'(b))$ is called the \index{external angle}\emph{external angle} of $\gamma$ at time~$b$.
If $\theta=\pi$, then we say that $\gamma$ has a \index{cusp}\emph{cusp} at  time~$b$.

A space curve $\gamma$ is called \index{piecewise smooth curve}\emph{piecewise smooth} if it can be presented as an iterated concatenation of a finite number of smooth curves; if $\gamma$ is closed, then the  concatenation is assumed to be cyclic.

If $\gamma$ is a concatenation of smooth arcs $\gamma_1,\dots,\gamma_n$, then the total curvature of $\gamma$ is defined as a sum of the total curvatures of $\gamma_i$ and the external angles;
that is, 
\[\tc\gamma=\tc{\gamma_1}+\dots+\tc{\gamma_n}+\theta_1+\dots+\theta_{n-1}\]
where $\theta_i$ is the external angle at the joint between $\gamma_i$ and $\gamma_{i+1}$.

If $\gamma$ is closed, then the total curvature of $\gamma$ is defined by
\[\tc\gamma=\tc{\gamma_1}+\dots+\tc{\gamma_n}+\theta_1+\dots+\theta_{n},\]
where $\theta_n$ is the external angle at the joint between $\gamma_n$ and $\gamma_1$.

{

\begin{wrapfigure}{r}{23 mm}
\vskip-3mm
\centering
\includegraphics{mppics/pic-354}
\end{wrapfigure}

In particular, for a smooth loop $\gamma\:[a,b] \z\to \mathbb{R}^3$, the total curvature of the corresponding closed curve $\hat\gamma$ is
\[\tc{\hat\gamma}\df\tc\gamma + \theta,\]
where $\theta=\measuredangle(\gamma'(a),\gamma'(b))$.

}

\section{Generalized Fenchel's theorem}

\begin{thm}{Theorem}\label{thm:gen-fenchel}
Let $\gamma$ be a closed piecewise smooth space curve.
Then 
\[\tc\gamma\ge2\cdot\pi.\]

\end{thm}

\parbf{Proof.}
Suppose $\gamma$ is a cyclic concatenation of smooth arcs $\gamma_1,\dots,\gamma_n$.
Denote by $\theta_1,\dots,\theta_n$ its external angles.
We need to show that \index{10phi@$\tc{\gamma}$ (total curvature)}
\[\tc{\gamma_1}+\dots+\tc{\gamma_n}+\theta_1+\dots+\theta_n\ge2\cdot\pi.\eqlbl{eq:gen-fenchel}\]

Consider the tangent indicatrix $\tan_i$ for each arc $\gamma_i$;
these are spherical arcs.

The same argument as in the proof of Fenchel's theorem (\ref{thm:fenchel}) shows that the curves $\tan_1,\dots,\tan_n$ cannot lie in an open hemisphere.

The spherical distance from the endpoint of $\tan_i$ to the starting point of $\tan_{i+1}$ is equal to the external angle $\theta_i$ (we enumerate the arcs modulo $n$, so $\gamma_{n+1}=\gamma_1$).
Let us connect the endpoint of $\tan_i$ to the starting point of $\tan_{i+1}$ by a short arc of a great circle in the sphere.
This way we get a closed spherical curve that is $\theta_1+\dots+\theta_n$ longer than the total length of $\tan_1,\dots,\tan_n$.

Applying the hemisphere lemma (\ref{lem:hemisphere}) to the obtained closed curve, we get that
\[\length\tan_1+\dots+\length\tan_n+\theta_1+\dots+\theta_n\ge 2\cdot\pi.\]
By \ref{obs:tantrix}, the statement follows.
\qedsf

\begin{thm}{Chord lemma}\label{lem:chord}
Let $\gamma\:[a,b]\z\to\mathbb{R}^3$
be a smooth arc, and
$\ell$ be its chord.
Assume $\gamma$ meets $\ell$ at angles $\alpha$ and $\beta$ at $\gamma (a)$ and $\gamma (b)$, respectively;
that is,
\[\alpha=\measuredangle(\vec w,\gamma'(a))\quad\text{and}\quad \beta=\measuredangle(\vec w,\gamma'(b)),\]
where $\vec w=\gamma(b)-\gamma(a)$.
Then 
\[\tc\gamma\ge \alpha+\beta.\eqlbl{tc>a+b}\] 

\end{thm}

\parbf{Proof.}
Let us parametrize the chord $\ell$ from $\gamma(b)$ to $\gamma(a)$ and consider the cyclic concatenation $\hat\gamma$ of $\gamma$ and $\ell$.
The closed curve $\hat\gamma$ has two external angles $\pi-\alpha$ and $\pi-\beta$.

\begin{wrapfigure}{r}{45 mm}
\vskip-5mm
\centering
\includegraphics{mppics/pic-53}
\vskip0mm
\end{wrapfigure}

Since the curvature of $\ell$ vanishes, we get 
\[\tc{\hat\gamma}=\tc\gamma+(\pi-\alpha)+(\pi-\beta).\]
According to the generalized Fenchel's theorem (\ref{thm:gen-fenchel}),
$\tc{\hat\gamma}\ge 2\cdot\pi$;
hence \ref{tc>a+b} follows.
\qeds

\begin{thm}{Exercise}\label{ex:chord-lemma-optimal}
Show that the estimate in the chord lemma is optimal.

More precisely, given two distinct points $p, q$ and two unit vectors $\vec u,\vec v$ in $\mathbb{R}^3$,
construct a smooth curve $\gamma$ that starts at $p$ in the direction $\vec u$ and ends at $q$ in the direction $\vec v$ such that 
$\tc\gamma$ is arbitrarily close to $\measuredangle(\vec w,\vec u)+\measuredangle(\vec w,\vec v)$, where vector $\vec w$ points from $p$ to $q$.

\end{thm}

\section{Polygonal lines} 

Polygonal lines are a particular case of piecewise smooth curves;
each arc in its concatenation is a line segment.
Since the curvature of a line segment vanishes, the total curvature of a polygonal line is the sum of its external angles.

\begin{thm}{Exercise}\label{ex:monotonic-tc}
Let $a$, $b$, $c$, $d$, and $x$ be distinct points in $\mathbb{R}^3$.
Show that the total curvature of the polygonal line $abcd$ cannot exceed the total curvature of $abxcd$; that is, 
\[\tc {abcd} \le \tc {abxcd}.\]

Use this statement to show that any closed polygonal line has curvature at least $2\cdot\pi$.
\end{thm}

\begin{thm}{Proposition}\label{prop:inscribed-total-curvature}
Assume a polygonal line $p_1\dots p_n$ is inscribed in a smooth curve~$\gamma$.
Then 
\[\tc\gamma\ge \tc{p_1\dots p_n}.\]
Moreover, if $\gamma$ is closed, we allow the inscribed polygonal line $p_1\dots p_n$ to be closed.

\end{thm}

{

\begin{wrapfigure}[7]{o}{40 mm}
\vskip-8mm
\centering
\includegraphics{mppics/pic-55}
\vskip0mm
\end{wrapfigure}

\parbf{Proof.}
Assume that $\gamma$ is closed.
Set 
\begin{align*}
p_i&=\gamma(t_i),
&
\alpha_i&=\measuredangle(\vec w_i,\vec v_i),
\\
\vec w_i&=p_{i+1}-p_i,
& 
\beta_i&=\measuredangle(\vec w_{i-1},\vec v_i),
\\
\vec v_i&=\gamma'(t_i),
&
\theta_i&=\measuredangle(\vec w_{i-1},\vec w_i).
\end{align*}
Let us use indexes modulo $n$;
so $p_{n+1}\z=p_1$.

}

Since the curvature of line segments vanishes, 
the total curvature of the polygonal line $p_1\dots p_n$ is the sum of external angles $\theta_i$.

By triangle inequality for angles \ref{thm:spherical-triangle-inq}, we get that
\[\theta_i\le \alpha_i+\beta_i.\]
By the chord lemma (\ref{lem:chord}), the total curvature of the arc of $\gamma$ from $p_i$ to $p_{i+1}$ is at least $\alpha_i+\beta_{i+1}$.
Therefore, if $\gamma$ is a closed curve, we have
\begin{align*}
\tc{p_1\dots p_n}&=\theta_1+\dots+\theta_n\le
\\
&\le\beta_1+\alpha_1+\dots+\beta_n+\alpha_n = 
\\
&=(\alpha_1+\beta_2)+\dots+(\alpha_n+\beta_1) \le 
\\
&\le \tc\gamma.
\end{align*}

If $\gamma$ is not closed, the calculations are analogous:
\begin{align*}
\tc{p_1\dots p_n}&=\theta_2+\dots+\theta_{n-1}\le
\\
&\le\beta_2+\alpha_2+\dots+\beta_{n-1}+\alpha_{n-1} \le
\\
&\le (\alpha_1+\beta_2)+\dots+(\alpha_{n-1}+\beta_n) \le
\\
&\le \tc\gamma.
\end{align*}
\qedsf

\begin{thm}{Exercise}\label{ex:sef-intersection}\label{ex:sef-intersection:>pi}  
Use the above results to give another solution to Exercise \ref{ex:length-dist:self-intersection:>pi}. 
\end{thm}

\begin{wrapfigure}{r}{30 mm}
\vskip-0mm
\centering
\includegraphics{mppics/pic-20}
\vskip0mm
\end{wrapfigure}

\begin{thm}{Exercise}\label{ex:quadrisecant}
Suppose a closed curve $\gamma$ crosses a line at four points $a$, $b$, $c$, and~$d$.
Assume these points appear on the line in the order $a$, $b$, $c$, $d$
and they appear on the curve $\gamma$ in the order $a$, $c$, $b$,~$d$.
Show that 
\[\tc\gamma\ge 4\cdot\pi.\]

\end{thm}

\section[\texorpdfstring{What if $\Phi(\gamma)=2\cdot \pi$?}{What if Φ(γ)=2·π?}]{What if $\bm{\Phi(\gamma)=2\cdot \pi}$?}

\begin{thm}{Proposition}\label{prop:fenchel=}
The equality case in Fenchel's theorem holds only for convex plane curves;
that is, the total curvature of a smooth space curve $\gamma$ equals $2\cdot\pi$ if and only if $\gamma$ is a convex plane curve.
\end{thm}

\parbf{Proof of \ref{prop:fenchel=}.}
The if part is proved in Corollary~\ref{cor:fenchel=convex};
it remains to prove the only-if part.

Consider an inscribed quadrangle $abcd$ in~$\gamma$.
By the definition of total curvature, we have that
\begin{align*}
\tc{abcd}&=(\pi-
\measuredangle\hinge adb)+(\pi-
\measuredangle\hinge bac)+(\pi-
\measuredangle\hinge cbd)+(\pi-
\measuredangle\hinge dca)=
\\
&=4\cdot\pi -(
\measuredangle\hinge adb
+
\measuredangle\hinge bac
+
\measuredangle\hinge cbd
+
\measuredangle\hinge dca)
\end{align*}

By the triangle inequality for angles (\ref{thm:spherical-triangle-inq}),
\[
\measuredangle\hinge bac
\le
\measuredangle\hinge bad
+ 
\measuredangle\hinge bdc
\quad\text{and}\quad
\measuredangle\hinge dca\le
\measuredangle\hinge dcb
+ 
\measuredangle\hinge dba.
\eqlbl{eq:spheric-triangle}
\]

\begin{wrapfigure}{r}{30 mm}
\vskip-5mm
\centering
\includegraphics{mppics/pic-56}
\vskip0mm
\end{wrapfigure}

The sum of angles in any triangle is $\pi$, so combining these inequalities, we get that 
\begin{align*}
\tc{abcd}\ge 4\cdot \pi 
&- (\measuredangle\hinge adb+\measuredangle\hinge bad+ 
\measuredangle\hinge dba)-
\\
&-(\measuredangle\hinge cbd+\measuredangle\hinge dcb 
+\measuredangle\hinge  bdc)=
\\
=2\cdot\pi.&
\end{align*}

By \ref{prop:inscribed-total-curvature},
\[\tc{abcd}\le \tc\gamma\le 2\cdot\pi.\]
Therefore, we have equalities in \ref{eq:spheric-triangle}.
It means that point $d$ lies in the angle $abc$ 
and point $b$ lies in the angle $cda$.
The latter implies that $abcd$ is a convex plane quadrangle.

It follows that any quadrangle inscribed in $\gamma$ is a convex plane quadrangle.
Therefore, all points of $\gamma$ lie in one plane defined by three points on~$\gamma$.
Further, since any quadrangle inscribed in $\gamma$ is convex,
we get that $\gamma$ is convex as well. 
\qeds

\section{Generalized DNA theorem}\label{sec:DNA-poly}

\begin{thm}{Theorem}\label{thm:DNA-poly}
Let $p_1\dots p_n$ be a closed polygonal line in a unit ball.
Then 
\[\tc{p_1\dots p_n}>\length(p_1\dots p_n).\]
\end{thm}

Exercise \ref{ex:total-curvature=} and this theorem imply the original DNA theorem (\ref{thm:DNA}).
Therefore, \ref{thm:DNA-poly} generalizes \ref{thm:DNA}.

\parbf{Proof.}
We assume that $p_n=p_0$, $p_{n+1}=p_1$, and so on.
Denote by $\theta_i$ the external angle at $p_i$.

\begin{figure}[ht!]
\vskip-0mm
\centering
\includegraphics{mppics/pic-16}
\vskip0mm
\end{figure}

Denote by $o$ the center of the ball.
Consider a sequence of triangles 
\[\triangle q_0q_1s_0\cong \triangle p_0p_1o,\ \ \triangle q_1q_2s_1\cong \triangle p_1p_2o,\ \dots\]
such that the points $q_0,q_1,\dots$ lie on one line in that order and all the $s_i$'s lie on one side from this line.

Note that $s_0s_nq_nq_0$ is a parallelogram; hence
\[|s_n-s_0|=|q_n-q_0|=\length (p_1\dots p_n).\]
Therefore 
\[|s_0-s_1|+\dots+|s_{n-1}-s_n|\ge \length (p_1\dots p_n).\]

Furthermore,
\[|q_i-s_{i-1}|=|q_i-s_i|=|p_i-o|\le 1\]
and
\[\theta_i\ge\measuredangle \hinge{q_i}{s_{i-1}}{s_i}\]
for each $i$.
It follows that
\[\theta_i>|s_{i-1}-s_i|\]
for each $i$.
Therefore,
\begin{align*}
\tc {p_1\dots p_n}
&=\theta_1+\dots+\theta_n>
\\
&> |s_{0}-s_1|+\dots +|s_{n-1}-s_n|\ge 
\\
&\ge\length (p_1\dots p_n).
\end{align*}
\qedsf

Let us mention the following  generalization of the DNA inequality;
it was proved by Jeffrey Lagarias and Thomas Richardson \cite{lagarias-richardso}; another proof is given by Alexander Nazarov and Fedor Petrov \cite{nazarov-petrov}.
Both proofs are annoyingly tricky.

\begin{thm}{Theorem}
Let $\alpha$ be a closed curve that lies in a convex plane figure bounded by a curve $\gamma$.
Then the average curvature of $\alpha$ is not less than the average curvature of $\gamma$.

\end{thm}

\section{Generalized curvature}

The following exercise states that the inequality in \ref{prop:inscribed-total-curvature} is optimal.

\begin{thm}{Exercise}\label{ex:total-curvature=}
Show that for any smooth space curve $\gamma$ we have 
\[\tc\gamma=\sup\{\tc\beta\},\]
where the least upper bound is taken over all polygonal lines~$\beta$ inscribed in $\gamma$
(if $\gamma$ is closed, then we assume so is $\beta$).
\end{thm}

This equality can be used to define the total curvature of an arbitrary curve~$\gamma$.
Namely, it can be defined as {}\textit{the least upper bound on the total curvatures of nondegenerate polygonal lines inscribed in~$\gamma$.}

This theory was reinvented and developed several times; see \cite[III §~1]{pogorelov}, \cite{aleksandrov-reshetnyak}, and \cite{sullivan-curves}.
It is possible to generalize most of the statements in this chapter to the curves of finite total curvature in this generalized sense.

\begin{thm}{Exercise}\label{ex:tc-rectifiable}
Suppose that a curve $\gamma\:[0,1]\to\mathbb{R}^3$ has bounded total curvature in the generalized sense;
that is, there is an upper bound on the total curvatures of polygonal lines inscribed in~$\gamma$.

Show that $\gamma$ is rectifiable.
Construct an example showing that the converse does not hold. 
\end{thm}

%% file: torsion-new.tex
\chapter{Torsion}
\label{chap:torsion}

This chapter provides mostly practice in computations.
Except for the definitions in Section~\ref{sec:frenet-frame},
it will not be used in the sequel.

Torsion is a close relative of curvature.
The curvature measures how much a curve fails to be a straight line, and the torsion measures how much a space curve fails to be a plane curve (see \ref{ex:lancret}).

\section{Frenet frame}\label{sec:frenet-frame}

Let $\gamma$ be a smooth space curve.
Without loss of generality, we may assume that $\gamma$ has an arc-length parametrization,
so the velocity vector $\tan(s)\z=\gamma'(s)$ is unit.

Assume its curvature does not vanish at time $s$;
in other words, $\gamma''(s)\z\ne 0$.
Then we can define the so-called \index{normal!vector}\emph{normal vector} at $s$ as
\[\norm(s)=\frac{\gamma''(s)}{|\gamma''(s)|}.\]
Note that \index{10tnb@$\tan$, $\norm$, $\bi$ (Frenet frame)}
\[\tan'(s)=\gamma''(s)=\kur(s)\cdot\norm(s).\]

According to \ref{prop:a'-pertp-a''}, $\norm(s)\perp \tan(s)$.
Therefore, the vector product 
\[\bi(s)=\tan(s)\times \norm(s)\]
is a unit vector.
Moreover, the triple $\tan(s),\norm(s),\bi(s)$ is an oriented orthonormal basis in $\mathbb{R}^3$;
it is called the \index{Frenet frame}\emph{Frenet frame} of $\gamma$ at~$s$.
In particular, we have that
\[\begin{aligned}
\langle\tan,\tan\rangle&=1,
&
\langle\norm,\norm\rangle&=1,
&
\langle\bi,\bi\rangle&=1,
\\
\langle\tan,\norm\rangle&=0,
&
\langle\norm,\bi\rangle&=0,
&
\langle\bi,\tan\rangle&=0.
\end{aligned}
\eqlbl{eq:orthornomal}
\]
$\tan(s)$, $\norm(s)$, and $\bi(s)$ are called \index{tangent}\emph{tangent}, \index{normal}\emph{normal}, and \index{binormal}\emph{binormal} vectors of the Frenet frame, respectively.
The frame $\tan(s),\norm(s),\bi(s)$ is defined only if $\kur(s)\z\ne 0$.

The plane $\Pi_s$ thru $\gamma(s)$ spanned by vectors $\tan(s)$ and $\norm(s)$ is called the \index{osculating!plane}\emph{osculating plane} at $s$;
equivalently it can be defined as a plane thru $\gamma(s)$ that is perpendicular to the binormal vector $\bi(s)$.
This is the unique plane that has a \index{order of contact}\emph{second-order contact} with $\gamma$ at $s$;
that is, $\rho(\ell)=o(\ell^2)$, where $\rho(\ell)$ denotes the distance from $\gamma(s+\ell)$ to $\Pi_s$.

\section{Torsion}

Let $\gamma$ be a smooth unit-speed space curve,
and let $\tan,\norm,\bi$ be its Frenet frame.
The value \index{10tau@$\tor$ (torsion)}
\[\tor(s)=\langle \norm'(s),\bi(s)\rangle\]
is called the \index{torsion}\emph{torsion} of $\gamma$ at~$s$.

The torsion $\tor(s_0)$ is defined if $\kur(s_0)\z\ne0$.
Indeed, since the function $s\mapsto \kur(s)$ is continuous, 
$\kur(s_0)\z\ne 0$ implies that $\kur(s)\z\ne 0$ for all $s$ near $s_0$.
Therefore, the Frenet frame is also defined in an open interval containing $s_0$.
Clearly, $\tan(s)$, $\norm(s)$, and $\bi(s)$ depend smoothly on $s$ in their domains of definition.
Therefore, $\norm'(s_0)$ is defined, and so is the torsion.

%\begin{thm}{Exercise}\label{ex:torsion-plane}
%Let $\gamma$ be a smooth unit-speed space curve and assume that $\gamma''$ never vanishes.
%Show that $\gamma$ is contained in a plane if and only if $\tor \equiv 0$.
%\end{thm} !!!It is a copy of \ref{ex:torsion=0}

\begin{thm}{Exercise}\label{ex:helix-torsion}
Given real numbers $a$ and $b$, calculate the curvature and the torsion of the helix
$\gamma_{a,b}(t)=(a\cdot \cos t,a\cdot\sin t, b\cdot t)$.

Conclude that for any $\kur>0$ and $\tor$ there is a helix with constant curvature $\kur$ and torsion $\tor$.
\end{thm}

\section{Frenet formulas}

Assume the Frenet frame $\tan(s),\norm(s),\bi(s)$ of a curve $\gamma$ is defined at~$s$.
Recall that 
\[\tan'=\kur\cdot \norm.
\eqlbl{eq:frenet-tau}\]
Let us write the remaining derivatives $\norm'$ and $\bi'$ in the frame $\tan,\norm,\bi$.

Let us show that
\[\norm'=-\kur\cdot\tan+\tor\cdot\bi.\eqlbl{eq:frenet-nu}\]

Since the frame $\tan,\norm,\bi$ is orthonormal, the above formula is equivalent to the following three identities:
\[\begin{aligned}
\langle \norm',\tan\rangle&=-\kur,
&
\langle \norm',\norm\rangle&=0,
&
\langle \norm',\bi\rangle&=\tor,
\end{aligned}\eqlbl{eq:<N',?>}\]
The last identity follows from the definition of torsion.
The second one is a consequence of the identity $\langle \norm,\norm\rangle\z=1$ in \ref{eq:orthornomal}. 
By differentiating the identity $\langle\tan,\norm\rangle\z=0$ in \ref{eq:orthornomal}
we get 
\[\langle\tan',\norm\rangle+\langle\tan,\norm'\rangle=0.\]
Applying \ref{eq:frenet-tau}, we get the first equation in \ref{eq:<N',?>}.

Differentiating the third identity in \ref{eq:orthornomal}, we get that $\bi'\perp\bi$.
Taking further derivatives of the other identities with $\bi$ in \ref{eq:orthornomal}, we get that 
\begin{align*}
\langle\bi',\tan\rangle&=-\langle\bi,\tan'\rangle=-\kur\cdot\langle\bi,\norm\rangle=0,
\\
\langle\bi,\norm'\rangle&=-\langle\bi',\norm\rangle=\tor.
\end{align*}
Since the frame $\tan,\norm,\bi$ is orthonormal, it follows that
\[\bi'=-\tor\cdot\norm.\eqlbl{eq:frenet-beta}\]

The equations \ref{eq:frenet-tau}, \ref{eq:frenet-nu}, and \ref{eq:frenet-beta} are called \index{Frenet formulas}\emph{Frenet formulas}.
All three can be written as one matrix identity:
\[
\begin{pmatrix}
\tan'
\\
\norm'
\\
\bi'
\end{pmatrix}
=
\begin{pmatrix}
0&\kur&0
\\
-\kur&0&\tor
\\
0&-\tor&0
\end{pmatrix}
\cdot
\begin{pmatrix}
\tan
\\
\norm
\\
\bi
\end{pmatrix}.
\]

Since $\bi$ is the normal vector to the osculating plane, equation \ref{eq:frenet-beta} shows that the torsion measures how fast the osculating plane rotates when one travels along~$\gamma$.

\begin{thm}{Exercise}\label{ex:beta-from-tau+nu}
Deduce the formula \ref{eq:frenet-beta} from  \ref{eq:frenet-tau} and \ref{eq:frenet-nu} by differentiating the identity
$\bi=\tan\times \norm$.
\end{thm}

\begin{thm}{Exercise}\label{ex:torsion=0}
Let $\gamma$ be a smooth space curve with nonvanishing curvature.
Show that $\gamma$ lies in a plane if and only if its torsion vanishes.
\end{thm}

\begin{thm}{Exercise}\label{ex:+B}
Let $\gamma_0\:[a,b]\to \mathbb{R}^3$ be a smooth space curve with Frenet frame $\tan,\norm,\bi$.
Consider the curve $\gamma_1(t)=\gamma_0(t)+\bi(t)$.
Show that
\[\length\gamma_1\ge\length\gamma_0.\]
\end{thm}

\begin{thm}{Exercise}\label{ex:frenet}
Suppose $\gamma$ is a smooth space curve.
Let $\tan,\norm,\bi$ be its Frenet frame and $\tor$ be its torsion.
Show that 
\[\bi=\frac{\gamma'\times\gamma''}{|\gamma'\times\gamma''|}
\quad\text{and}\quad
\tor=\frac{\langle\gamma'\times\gamma'',\gamma'''\rangle}{|\gamma'\times\gamma''|^2}.
\]

\end{thm}

\begin{thm}{Exercise}\label{ex:moment-curve}
Find curvature $\kur(t)$ and torsion $\tor(t)$ of the \index{moment curve}\emph{moment curve} $\gamma\:t\z\mapsto (t,t^2,t^3)$ at $\gamma(t)$.
\end{thm}

The following exercise is closely related to the bow lemma (\ref{lem:bow}).

\begin{thm}{Exercise}\label{ex:bow-converse}
Let $\gamma_1,\gamma_2\:[a,b]\to\mathbb{R}^3$ be two smooth unit-speed curves.
Assume that 
\[\dist{\gamma_1(t_1)}{\gamma_1(t_2)}{}\ge \dist{\gamma_2(t_1)}{\gamma_2(t_2)}{}\]
for any $t_1$ and $t_2$.
Show that $\kur(t_0)_{\gamma_1}\le \kur(t_0)_{\gamma_2}$ for any $t_0$.
\end{thm}

\begin{thm}{Advanced exercise}\label{ex:torsion-indicatrix}
Let $\gamma$ be a closed smooth space curve with positive torsion.
Show that the tangent indicatrix of $\gamma$ has a self-intersection.
\end{thm}

\section{Curves of constant slope}

We say that a smooth space curve $\gamma$ has \index{constant slope}\emph{constant slope} if its velocity vector makes a constant angle with a fixed direction.
The following theorem was proved by Michel Ange Lancret~\cite{lancret}.

\begin{thm}{Theorem}\label{thm:const-slope}
Let $\gamma$ be a smooth curve;
denote by $\kur$ and $\tor$ its curvature and torsion.
Suppose $\kur(s)>0$ for all~$s$.
Then $\gamma$ has a constant slope if and only if the ratio $\tfrac\tor\kur$ is constant.
\end{thm}

\begin{thm}{Proof-guided exercise}\label{ex:lancret}
Let $\gamma$ be a smooth space curve with nonvanishing curvature, $\tan,\norm,\bi$ 
its Frenet frame, and $\kur$, $\tor$ its curvature and torsion.

\begin{subthm}{ex:lancret:a}
Assume $\langle \vec w,\tan\rangle$ is constant for a fixed nonzero vector $\vec w$.
Show that $\langle \vec w, \norm\rangle =0$.
Conclude that $\langle \vec w,\bi\rangle$ is constant.
Show that \[\tor\cdot\langle \vec w,\bi\rangle -\kur\cdot\langle \vec w,\tan\rangle =0.\]
Conclude that $\tfrac\tor\kur$ is constant.
\end{subthm}

\begin{subthm}{ex:lancret:b} Assume $\tfrac\tor\kur$ is constant, show that the vector $\vec w=\tfrac\tor\kur\cdot \tan+\bi$ is constant.
Conclude that $\gamma$ has a constant slope.
\end{subthm}

\end{thm}

Let $\gamma$ be a smooth unit-speed curve and $s_0$ a fixed real number. 
Then the curve 
\[\alpha(s)=\gamma(s)+(s_0-s)\cdot \gamma'(s)\]
is called the \index{evolvent}\emph{evolvent} of~$\gamma$.
Note that if $\ell(s)$ denotes the tangent line to $\gamma$ at $s$,
then $\alpha(s)\in \ell(s)$ and $\alpha'(s)\perp \ell$ for all~$s$.

\begin{thm}{Exercise}\label{ex:evolvent-constant-slope}
Show that the evolvent of a constant slope curve is a plane curve.
\end{thm}

\section{Spherical curves}

\begin{thm}{Theorem}
Suppose $\gamma$ is a smooth space curve with nonvanishing torsion $\tor$ and (therefore) curvature $\kur$.
Then $\gamma$ lies on a unit sphere if and only if 
the following identity holds true:
\[\left|\frac{\kur'}{\tor}\right|=\kur\cdot\sqrt{\kur^2-1}.\]
\end{thm}

\begin{thm}{Proof-guided exercise}\label{ex:spherical-frenet}
Suppose $\gamma$ is a smooth unit-speed space curve.
Denote by $\tan,\norm,\bi$ its Frenet frame and by $\kur$, $\tor$ its curvature and torsion.

\smallskip

Assume $\gamma$ is spherical; that is, $|\gamma(s)|=1$ for any~$s$.
Show that

\begin{subthm}{ex:spherical-frenet:tau} $\langle\tan,\gamma\rangle=0$; conclude that $\langle\norm,\gamma\rangle^2+\langle\bi,\gamma\rangle^2=1$.
\end{subthm}

\begin{subthm}{ex:spherical-frenet:nu} $\langle\norm,\gamma\rangle=-\tfrac1\kur$;
\end{subthm}

\begin{subthm}{ex:spherical-frenet:beta} $\langle\bi,\gamma\rangle'=\tfrac\tor\kur$.
\end{subthm}

\begin{subthm}{ex:spherical-frenet:beta+}
Use \ref{SHORT.ex:spherical-frenet:beta} to show that if $\gamma$ is closed, then $\tor(s)=0$ for some~$s$.
\end{subthm}

\begin{subthm}{ex:spherical-frenet:kur-tor} Assume the torsion of $\gamma$ does not vanish.
Use \ref{SHORT.ex:spherical-frenet:tau}--\ref{SHORT.ex:spherical-frenet:beta} to show that 
\[\left|\frac{\kur'}{\tor}\right|=\kur\cdot\sqrt{\kur^2-1}.\]
\end{subthm}
Now assume $\gamma$ is a space curve that satisfies the identity in \ref{SHORT.ex:spherical-frenet:kur-tor}.
\begin{subthm}{ex:spherical-frenet:f} Show that $\gamma+\tfrac1\kur\cdot \norm+\tfrac{\kur'}{\kur^2\cdot\tor}\cdot\bi$ is a constant point; conclude that $\gamma$ lies in the unit sphere centered at this point.
\end{subthm}

\end{thm}

For a unit-speed curve $\gamma$ with nonzero curvature and torsion at~$s$,
the sphere $\Sigma_s$ that passes thru $\gamma(s)$ and has the center at
\[p(s)=\gamma(s)+\tfrac1{\kur(s)}\cdot \norm(s)+\tfrac{\kur'(s)}{\kur^2(s)\cdot\tor(s)}\cdot\bi(s)\]
 is called the \index{osculating!sphere}\emph{osculating sphere} of $\gamma$ at~$s$.
This is the unique sphere that has \index{order of contact}\emph{third-order contact} with $\gamma$ at~$s$;
that is, $\rho(\ell)=o(\ell^3)$, where $\rho(\ell)$ denotes the distance from $\gamma(s+\ell)$ to $\Sigma_s$.
 
\section{Fundamental theorem of space curves}

\begin{thm}{Theorem}\label{thm:fund-curves}
Let $s\mapsto \kur(s)$ and $s\mapsto \tor(s)$ be two smooth real-valued functions defined on a real interval $\mathbb{I}$.
Suppose $\kur(s)>0$ for all~$s$.
Then there is a smooth unit-speed curve $\gamma\:\mathbb{I}\to\mathbb{R}^3$ with curvature $\kur(s)$ and torsion $\tor(s)$ for every~$s\in \mathbb{I}$.
Moreover, $\gamma$ is uniquely defined up to an orientation-preserving isometry of the space.
\end{thm}

The proof is an application of the theorem on the existence and uniqueness of solutions of ordinary differential equations (\ref{thm:ODE}).

\parbf{Proof.}
Fix a parameter value $s_0$, a point $\gamma(s_0)$, and an oriented orthonormal frame $\tan(s_0)$, $\norm(s_0)$, $\bi(s_0)$.
Consider the following system of differential equations
\[
\begin{cases}
\gamma'=\tan,
\\
\tan'=\kur\cdot\norm,
\\
\norm'=-\kur\cdot\tan+\tor\cdot\bi,
\\
\bi'=-\tor\cdot\norm
\end{cases}
\eqlbl{eq:gamma'tan'norm'bi'}
\]
with the initial condition formed by $\gamma(s_0)$, $\tan(s_0)$, $\norm(s_0)$, $\bi(s_0)$.
(The system of equations has four vector equations, so it can be rewritten as a system of $12$ scalar equations.)

By \ref{thm:ODE}, this system has a unique solution which is defined in a maximal subinterval $\mathbb{J}\subset \mathbb{I}$ containing $s_0$.
Let us show that actually $\mathbb{J}=\mathbb{I}$.

First note that 
\[\begin{aligned}
\langle\tan,\tan\rangle&=1,
&
\langle\norm,\norm\rangle&=1,
&
\langle\bi,\bi\rangle&=1,
\\
\langle\tan,\norm\rangle&=0,&
\langle\tan,\norm\rangle&=0,&
\langle\bi,\tan\rangle&=0
\end{aligned}
\eqlbl{eq:111000}
\]
at any parameter value $s$.

Indeed, by \ref{eq:gamma'tan'norm'bi'}, we have the following system of scalar equations:
\[
\begin{cases}
\langle\tan,\tan\rangle'
&=
2\cdot\langle\tan,\tan'\rangle
=
2\cdot\kur\cdot \langle\tan,\norm\rangle,
\\
\langle\norm,\norm\rangle'
&=
2\cdot\langle\norm,\norm'\rangle
=
-
2\cdot\kur\cdot\langle\norm,\tan\rangle
+
2\cdot\tor\cdot\langle\norm,\bi\rangle,
\\
\langle\bi,\bi\rangle'
&=
2\cdot\langle\bi,\bi'\rangle
=
-2\cdot\tor\langle\bi,\norm\rangle,
\\
\langle\tan,\norm\rangle'
&=
\langle\tan',\norm\rangle
+
\langle\tan,\norm'\rangle
=
\kur\cdot\langle\norm,\norm\rangle
-
\kur\cdot\langle\tan,\tan\rangle
+
\tor\cdot\langle\tan,\bi\rangle,
\\
\langle\norm,\bi\rangle'
&=
\langle\norm',\bi\rangle+\langle\norm,\bi'\rangle
=\kur\cdot\langle\tan,\bi\rangle+\tor\cdot\langle\bi,\bi\rangle-\tor\cdot\langle\norm,\norm\rangle,
\\
\langle\bi,\tan\rangle'
&=
\langle\bi',\tan\rangle+\langle\bi,\tan'\rangle
=
-\tor\cdot \langle\norm,\tan\rangle
+\kur\cdot\langle\bi,\norm\rangle.
\end{cases}
\eqlbl{eq:<gamma'tan'norm'bi'>}
\]
The constants in \ref{eq:111000} give a solution of this system.
Moreover, since we choose $\tan(s_0)$, $\norm(s_0)$, $\bi(s_0)$ to be an oriented orthonormal frame,
\ref{eq:111000} solves our initial value problem for \ref{eq:<gamma'tan'norm'bi'>}.

Assume $\mathbb{J} \varsubsetneq \mathbb{I}$.
Then an end of $\mathbb{J}$, say $b$, lies in the interior of $\mathbb{I}$.
The Theorem~\ref{thm:ODE} is applicable for $\Omega=\mathbb{R}^{12}\times\mathbb{I}$;
therefore, at least one of the values $\gamma(s)$, $\tan(s)$, $\norm(s)$, $\bi(s)$
escapes to infinity as $s\to b$.
But this is impossible --- the vectors $\tan(s)$, $\norm(s)$, $\bi(s)$ remain unit, and $|\gamma'(s)|\z=|\tan(s)|=1$;
so $\gamma$ travels only a finite distance as $s\to b$ --- a contradiction.
Hence, $\mathbb{J}= \mathbb{I}$, and the first statement follows.

Now assume there are two curves $\gamma_1$ and $\gamma_2$ with the given curvature and torsion functions.
Applying an isometry of the space, we can assume that $\gamma_1(s_0)=\gamma_2(s_0)$, and the Frenet frames of the curves coincide at $s_0$.
Then $\gamma_1=\gamma_2$ by the uniqueness of  solutions to the system (\ref{thm:ODE}).
Hence, the last statement follows.
\qeds

\begin{thm}{Exercise}\label{ex:cur+tor=helix}
Assume a curve $\gamma\:\mathbb{R}\to\mathbb{R}^3$ has constant speed, curvature, and torsion.
Show that $\gamma$ is a helix, possibly degenerating to a circle;
that is, in a suitable coordinate system we have
$\gamma(t)=(a\cdot \cos t,a\cdot\sin t, b\cdot t)$
for some constants $a$ and~$b$.
\end{thm}

\begin{thm}{Advanced exercise}\label{ex:const-dist}
Let $\gamma$ be a smooth space curve such that the distance $|\gamma(t)-\gamma(t+\ell)|$ depends only on $\ell$.
Show that $\gamma$ is a helix, possibly degenerate to a line or a circle.
\end{thm}

%% file: plane-curves.tex
\chapter{Signed curvature}

On the plane, it makes sense to think about left (right) turns as positive (respectively negative).
This leads to the so-called \textit{signed curvature} of plane curves.
If you drive a car along a plane curve, then the signed curvature describes the position of the steering wheel at a given %point. 
time.

\label{chap:signed-curvature}

\section{Definitions}\label{sec:def(skur)}

Suppose $\gamma$ is a smooth unit-speed plane curve,
so $\tan(s)=\gamma'(s)$ is its unit tangent vector for any~$s$.

Let us rotate $\tan(s)$ by the angle $\tfrac\pi 2$ counterclockwise; 
denote the obtained vector by $\norm(s)$.
The pair $\tan(s),\norm(s)$ is an oriented orthonormal frame in the plane which is analogous to the Frenet frame
defined in Section~\ref{sec:frenet-frame};
we will keep the name \index{Frenet frame}\emph{Frenet frame} for it.

Recall that $\gamma''(s)\perp \gamma'(s)$ (see \ref{prop:a'-pertp-a''}).
Therefore, \index{10k@$\skur$ (signed curvature)}
\[\tan'(s)=\skur(s)\cdot \norm(s).\eqlbl{eq:tau'}\]
for a real number $\skur(s)$;
the value $\skur(s)$ is called \index{curvature}\index{signed curvature}\emph{signed curvature} of $\gamma$ at~$s$.
We may use the notation $\skur(s)_\gamma$ if we need to specify the curve~$\gamma$.

Note that 
\[\kur(s)=|\skur(s)|;\]
that is, up to sign, the signed curvature $\skur(s)$ equals the curvature $\kur(s)$  of $\gamma$ at $s$ defined in Section~\ref{sec:curvature};
the sign prescribes the direction of turn --- if $\gamma$ turns left at time $s$, then $\skur (s)>0$.
If we want to emphasize that we are working with the \textit{non-signed} curvature of a curve, 
we call it \index{curvature}\index{absolute curvature}\emph{absolute curvature}.

If we reverse the parametrization of the curve or change the orientation of the plane, then
the signed curvature changes its sign.

Since $\tan(s),\norm(s)$ is an orthonormal frame, we have 
\begin{align*}
\langle\tan,\tan\rangle&=1,
&
\langle\norm,\norm\rangle&=1, 
&
\langle\tan,\norm\rangle&=0.
\end{align*}
Differentiating these identities we get 
\begin{align*}
\langle\tan',\tan\rangle&=0,
&
\langle\norm',\norm\rangle&=0,
&
\langle\tan',\norm\rangle+\langle\tan,\norm'\rangle&=0.
\end{align*}
By \ref{eq:tau'}, $\langle\tan',\norm\rangle=\skur$. 
Therefore $\langle\tan,\norm'\rangle=-\skur$.
Hence, we get 
\[\norm'(s)=-\skur(s)\cdot \tan(s).\eqlbl{eq:nu'}\]
The equations \ref{eq:tau'} and \ref{eq:nu'} are the Frenet formulas for plane curves. 
They can be written in a matrix form as:
\[
\begin{pmatrix}
\tan'
\\
\norm'
\end{pmatrix}
=
\begin{pmatrix}
0&\skur
\\
-\skur&0
\end{pmatrix}
\cdot
\begin{pmatrix}
\tan
\\
\norm
\end{pmatrix}.
\]

\begin{thm}{Exercise}\label{ex:bike}
Let $\gamma_0\:[a,b]\to\mathbb{R}^2$ be a smooth curve and $\tan$ its tangent indicatrix.
Consider another curve $\gamma_1\:[a,b]\to\mathbb{R}^2$ defined by $\gamma_1(t)\z\df\gamma_0(t)+\tan(t)$.
Show that

\begin{minipage}{.47\textwidth}
\begin{subthm}{ex:bike:length}$\length\gamma_0\le \length\gamma_1$;
\end{subthm}
\end{minipage}
\hfill
\begin{minipage}{.47\textwidth}
\begin{subthm}{ex:bike:tc}$\tc{\gamma_0}\le \length\gamma_1$.
\end{subthm}
\end{minipage}

\end{thm}

The curves $\gamma_0$ and $\gamma_1$ in the exercise above describe the tracks of an idealized bicycle with  distance 1 from the rear to the front wheel.
Thus, by the exercise, the front wheel must have a longer track.
For more on the geometry of bicycle tracks, see the survey of Robert Foote, Mark Levi, and Serge Tabachnikov \cite{foote-levi-tabachnikov} and the references therein.

\section{Fundamental theorem of plane curves}

\begin{thm}{Theorem}\label{thm:fund-curves-2D}
Let $s\mapsto \skur(s)$ be a smooth real-valued function defined on a real interval $\mathbb{I}$.
Then there is a smooth unit-speed curve $\gamma\:\mathbb{I}\to\mathbb{R}^2$ with signed curvature $\skur(s)$.
Moreover, $\gamma$ is uniquely defined up to an orientation-preserving isometry of the plane.
\end{thm}

This theorem is a partial case of its $3$-dimensional analog (\ref{thm:fund-curves}), but we present a direct proof.

\parbf{Proof.} 
Fix $s_0\in\mathbb{I}$.
Consider the function
\[\theta(s)
=
\int_{s_0}^s\skur(t)\cdot dt.\]
By the fundamental theorem of calculus, we have $\theta'(s)\z=\skur(s)$ for all~$s$.

Set 
$\tan(s)=(\,\cos(\theta(s)),\,\sin(\theta(s))\,)$,
and let $\norm(s)$ be its counterclockwise rotation by angle $\tfrac\pi2$; so 
$\norm(s)=(\,-\sin(\theta(s)),\,\cos(\theta(s))\,)$.
Consider the curve 
\[\gamma(s)=\int_{s_0}^s\tan(s)\cdot ds.\]
Since $|\gamma'|=|\tan|=1$, the curve $\gamma$ is unit-speed, and its Frenet frame is~$\tan,\norm$. 

Note that
\begin{align*}
\gamma''(s)&=\tan'(s)=
\\
&=(\,\cos(\theta(s))',\,\sin(\theta(s))'\,)=
\\
&=\theta'(s)\cdot (\,-\sin(\theta(s)),\,\cos(\theta(s))\,)=
\\
&=\skur(s)\cdot \norm(s).
\end{align*}
So, $\skur(s)$ is the signed curvature of $\gamma$ at~$s$. 

This proves the existence;
it remains to prove the uniqueness.

Assume $\gamma_1$ and $\gamma_2$ are two curves that satisfy the assumptions of the theorem.
Applying an isometry, we can assume that $\gamma_1(s_0)\z=\gamma_2(s_0)$, and both curves have the same Frenet frame at $s_0$.
Let us denote by $\tan_1,\norm_1$ and $\tan_2,\norm_2$ the Frenet frames of $\gamma_1$ and $\gamma_2$, respectively.
Both triples $\gamma_i,\tan_i,\norm_i$ satisfy the following system of ordinary differential equations 
\[
\begin{cases}
\gamma_i'=\tan_i,
\\
\tan_i'=\skur\cdot\norm_i,
\\
\norm_i'=-\skur\cdot\tan_i.
\end{cases}
\]
Moreover, they have the same initial values at $s_0$.
By the uniqueness of solutions of ordinary differential equations (\ref{thm:ODE}), we have $\gamma_1=\gamma_2$.
\qeds

Suppose $\gamma\:\mathbb{I}\to\mathbb{R}^2$ is a unit-speed curve.
A continuous function $\theta\:\mathbb{I}\to\mathbb{R}$ is called \index{continuous!argument}\emph{continuous argument} of $\gamma$ if 
\[\gamma'(s)=(\,\cos (\theta(s)),\,\sin(\theta(s))\,)\]
for any $s$.
The proof of the theorem implies the following.

\begin{thm}{Corollary}\label{cor:2D-angle}
For any smooth unit-speed curve $\gamma\:\mathbb{I}\to\mathbb{R}^2$ there is a continuous argument $\theta\:\mathbb{I}\to\mathbb{R}$.
Moreover 
\[\theta'(s)=\skur(s),\]
where $\skur$ denotes the signed curvature of~$\gamma$.
\end{thm}

\section{Total signed curvature}\label{sec:Total signed curvature}

Let $\gamma\:\mathbb{I}\to\mathbb{R}^2$ be a smooth unit-speed plane curve.
The \index{curvature}\index{total!signed curvature}\emph{total signed curvature} of $\gamma$, denoted by $\tgc\gamma$, is defined as the integral \index{10psi@$\tgc\gamma$ (total signed curvature)}
\[\tgc\gamma
=
\int_\mathbb{I} \skur(s)\cdot ds,\eqlbl{eq:tsc-k}\]
where $\skur$ denotes the signed curvature of~$\gamma$.

If $\mathbb{I}=[a,b]$, then 
\[\tgc\gamma=\theta(b)-\theta(a),\eqlbl{eq:tsc-theta}\]
where $\theta$ is a continuous argument of $\gamma$ (see \ref{cor:2D-angle}).

If $\gamma$ is a piecewise smooth plane curve, then we define its total signed curvature as the sum of the total signed curvatures of its arcs plus the sum of the \textit{signed} external angles at its joints;
they are positive where $\gamma$ turns left, negative where $\gamma$ turns right, and 0 where $\gamma$ goes straight.
It is undefined if $\gamma$ turns exactly backwards;
that is, if it has a cusp.

In other words, if $\gamma$ is a concatenation of smooth arcs $\gamma_1,\dots,\gamma_n$, then 
\[\tgc\gamma=\tgc{\gamma_1}+\dots+\tgc{\gamma_n}+\theta_1+\dots+\theta_{n-1},\]
where $\theta_i$ is the signed external angle at the joint between $\gamma_i$ and $\gamma_{i+1}$.
If $\gamma$ is closed, then the concatenation is cyclic, and
\[\tgc\gamma=\tgc{\gamma_1}+\dots+\tgc{\gamma_n}+\theta_1+\dots+\theta_{n},\]
where $\theta_n$ is the signed external angle at the joint between $\gamma_n$ and $\gamma_1$.

Since $\left|\int \skur(s)\cdot ds\right|\le \int|\skur(s)|\cdot ds$, we have
\[|\tgc\gamma|\le \tc\gamma\eqlbl{eq:tsc-tc}\] 
for any smooth plane curve $\gamma$;
that is, the total signed curvature $\tgc{}$ cannot exceed the total curvature $\tc{}$ in absolute value.
The equality holds if and only if the signed curvature does not change the sign.

\begin{thm}{Exercise}\label{ex:trochoids}
A trochoid is a curve traced out by a point fixed to a wheel as it rolls along a straight line.
\begin{figure}[!ht]
\centering
\begin{lpic}[t(-0mm),b(0mm),r(0mm),l(0mm)]{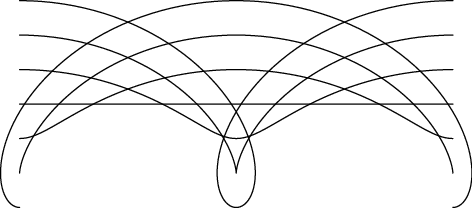}
\lbl[l]{4,0;{\footnotesize $-\tfrac32$}}
\lbl[l]{4,6;{\footnotesize $-1$}}
\lbl[tl]{10,15;{\footnotesize $-\tfrac12$}}
\lbl[t]{22,17;{\footnotesize $0$}}
\lbl[r]{3,23.6;{\footnotesize $\tfrac12$}}
\lbl[r]{3,29.4;{\footnotesize $1$}}
\lbl[r]{3,35.2;{\footnotesize $\tfrac32$}}

\end{lpic}
\end{figure}
A family of \index{trochoid}\emph{trochoids} $\gamma_a\:[0,2\cdot\pi]\to \mathbb{R}^2$ (see the picture) can be parametrized as
\[\gamma_a(t)=(t+a\cdot \sin t, a\cdot \cos t).\]
\begin{enumerate}[(a)]
\item Given $a\in \mathbb{R}$, find $\tgc{\gamma_a}$ if it is defined.
\item Given $a\in \mathbb{R}$, find $\tc{\gamma_a}$.
\end{enumerate}
\end{thm}

\begin{thm}{Proposition}\label{prop:total-signed-curvature}
Any simple closed smooth plane curve $\gamma$ has total signed curvature  $\pm2\cdot\pi$; it is $+2\cdot\pi$
if the region bounded by $\gamma$ lies on the left from it and  $-2\cdot\pi$ otherwise.

Moreover, the same statement holds for any simple closed piecewise smooth plane curve $\gamma$ if its total signed curvature is defined.
\end{thm}

This proposition is called sometimes {}\emph{Umlaufsatz}; it is a differential-geometric analog of the theorem about the sum of the internal angles of a polygon (\ref{thm:sum=(n-2)pi}) which we use in the proof.
A more conceptual proof was given by Heinz Hopf \cite{hopf1935}, \cite[p. 42]{hopf1989}.

\parbf{Proof.}
Without loss of generality, we may assume that $\gamma$ is oriented in such a way that the region bounded by $\gamma$ lies on its left.
We can also assume that $\gamma$ is unit-speed.

Consider a closed polygonal line $p_1\dots p_n$ inscribed in~$\gamma$.
We can assume that the arcs between the vertices are sufficiently small 
so that the polygonal line is simple, and each arc $\gamma_i$ from $p_i$ to $p_{i+1}$ has small total absolute curvature, say  $\tc{\gamma_i}<\pi$ for each~$i$.

Assume $p_i=\gamma(t_i)$.
As usual, we use indexes modulo $n$; in particular, $p_{n+1}\z=p_1$.
Set 
\begin{align*}
\vec w_i&=p_{i+1}-p_i,& \vec v_i&=\gamma'(t_i),
\\
\alpha_i&=\measuredangle(\vec v_i,\vec w_i),&\beta_i&=\measuredangle(\vec w_{i-1},\vec v_i),
\end{align*}
where $\alpha_i,\beta_i\in(-\pi,\pi)$ are signed angles --- $\alpha_i$ is positive if $\vec w_i$ points to the left from~$\vec v_i$.

{

\begin{wrapfigure}[13]{o}{43 mm}
\vskip-4mm
\centering
\includegraphics{mppics/pic-59}
\vskip0mm
\end{wrapfigure}

By \ref{eq:tsc-theta}, the value
\[\tgc{\gamma_i}-\alpha_i-\beta_{i+1}\eqlbl{eq:Psi-alpha-beta}\]
is a multiple of $2\cdot\pi$.
Since $\tc{\gamma_i}\z<\pi$, the chord lemma (\ref{lem:chord}) implies that $|\alpha_i|\z+|\beta_{i+1}|\z<\pi$.
By \ref{eq:tsc-tc}, we have that $|\tgc{\gamma_i}|\z\le\tc{\gamma_i}$;
therefore the value in \ref{eq:Psi-alpha-beta} vanishes.
In other words, for each $i$ we have
\[\tgc{\gamma_i}=\alpha_i+\beta_{i+1}.\]

}

Note that 
\[\delta_i=\pi-\alpha_i-\beta_i\eqlbl{eq:delta=pi-alpha-beta}\] 
is the internal angle of the polygonal line at $p_i$;
$\delta_i\in (0,2\cdot\pi)$ for each~$i$.
Recall that the sum of the internal angles of an $n$-gon is $(n-2)\cdot \pi$ (see \ref{thm:sum=(n-2)pi}); that is,
\[\delta_1+\dots+\delta_n=(n-2)\cdot \pi.\]
Therefore, 
\[
\begin{aligned}
\tgc\gamma&=\tgc{\gamma_1}+\dots+\tgc{\gamma_n}=
\\
&=(\alpha_1+\beta_2)+\dots+(\alpha_n+\beta_1)=
\\
&=(\beta_1+\alpha_1)+\dots+(\beta_n+\alpha_n)=
\\
&=(\pi-\delta_1)+\dots+(\pi-\delta_n)=
\\
&=n\cdot\pi-(n-2)\cdot \pi=
\\
&=2\cdot\pi.
\end{aligned}\eqlbl{eq:delta=pi-alpha-beta-sum}\]

The case of piecewise smooth curves is done the same way;
we need to subdivide the arcs in the cyclic concatenation further to meet the requirement above, and, instead of equation \ref{eq:delta=pi-alpha-beta}, we have 
\[\delta_i=\pi-\alpha_i-\beta_i-\theta_i,\]
where $\theta_i$ is the signed external angle of $\gamma$ at $p_i$; it vanishes if the curve $\gamma$ is smooth at $p_i$.
Therefore, instead of equation \ref{eq:delta=pi-alpha-beta-sum}, we have
\begin{align*}
\tgc\gamma&=\tgc{\gamma_1}+\dots+\tgc{\gamma_n}+\theta_1+\dots+\theta_n=
\\
&=(\alpha_1+\beta_2)+\dots+(\alpha_n+\beta_1)+\theta_1+\dots+\theta_n=
\\
&=(\beta_1+\alpha_1+\theta_1)+\dots+(\beta_n+\alpha_n+\theta_n)=
\\
&=(\pi-\delta_1)+\dots+(\pi-\delta_n)=
\\
&=n\cdot\pi-(n-2)\cdot \pi=
\\
&=2\cdot\pi.
\end{align*}
\qedsf

\begin{thm}{Exercise}\label{ex:zero-tsc}
Draw a smooth closed plane curve $\gamma$ such that 

\begin{subthm}{ex:zero-tsc:0}
$\tgc\gamma=0$;
\end{subthm}
 
\begin{subthm}{ex:zero-tsc:5}
$\tgc\gamma=\tc\gamma=10\cdot\pi$;
\end{subthm}

\begin{subthm}{ex:zero-tsc:2-4}
$\tgc\gamma=2\cdot \pi$ and $\tc\gamma=4\cdot\pi$.
\end{subthm}

\end{thm}

\begin{thm}{Exercise}\label{ex:length'}
Let $\gamma\:[a,b]\to\mathbb{R}^2$ be a smooth plane curve with Frenet frame $\tan,\norm$.
Given a real parameter $\ell$, consider
the curve $\gamma_\ell(t)\z=\gamma(t)\z+\ell\cdot\norm(t)$; it is called a \index{parallel!curve}\emph{parallel curve} of $\gamma$ at signed distance~$\ell$.

\begin{subthm}{ex:length':reg}
Show that the parametrization $\gamma_\ell$ is regular if $\ell\cdot \skur(t)_\gamma\ne 1$ for all~$t$.
\end{subthm}
 
\begin{subthm}{ex:length':formula}
Set $L(\ell)=\length\gamma_\ell$.
Show that 
\[L(\ell)=L(0)-\ell\cdot\tgc\gamma\eqlbl{eq:length(parallel-curve)}\]
for all $\ell$ sufficiently close to $0$. 
\end{subthm}

\begin{subthm}{ex:length':antiformula}
Describe an example showing that formula \ref{eq:length(parallel-curve)} may not hold for some values~$\ell$. 
\end{subthm}

\end{thm}

\section{Osculating circle}

\begin{thm}{Proposition}\label{prop:circline}
Given a point $p\in\mathbb{R}^2$,
a unit vector $\tan$ 
and a real number $\skur$, there is a unique smooth unit-speed curve $\sigma\:\mathbb{R}\to\mathbb{R}^2$ 
that starts at $p$ in the direction of $\tan$ and signed curvature $\skur$.

Moreover, if $\skur=0$, then it is a line $\sigma(s)=p+s\cdot \tan$;
if $\skur\ne 0$, then $\sigma$ runs around a circle of radius $\tfrac1{|\skur|}$ with center at $p+\tfrac1\skur\cdot \norm$, where $\tan,\norm$ is an oriented orthonormal frame.
\end{thm}

\parbf{Proof.}
Choose a coordinate system such that $p$ is its origin and $\tan$ points in the direction of the $x$-axis.
Therefore, $\norm$ points in the direction of the $y$-axis.
Then
\begin{align*}\theta(s)&=\int_{0}^s\skur\cdot dt=\skur\cdot s
\end{align*}
is a continuous argument of $\sigma$, see \ref{cor:2D-angle}.
Therefore,
\[\sigma'(s)=(\,\cos(\skur\cdot s),\,\sin(\skur\cdot s)\,).\]
It remains to integrate the last identity.
If $\skur=0$, we get $\sigma(s)=(s,0)$
which describes the line $\sigma(s)=p+s\cdot \tan$.

If $\skur\ne 0$, we get
\[\sigma(s)=(\,\tfrac1\skur\cdot\sin(\skur\cdot s),\,\tfrac1\skur\cdot(1-\cos(\skur\cdot s))\,);\]
it is the circle of radius $r=\tfrac1{|\skur|}$ centered at $(0,\tfrac1\skur)=p+\tfrac1\skur\cdot\norm$.
\qeds

\begin{thm}{Definition}
Let $\gamma$ be a smooth unit-speed plane curve;
denote by $\skur(s)$ the signed curvature of $\gamma$ at~$s$.

The unit-speed curve $\sigma_s$ of constant signed curvature $\skur(s)$ that starts at $\gamma(s)$ in the direction $\gamma'(s)$ is called the \index{osculating!circle}\emph{osculating circle} of $\gamma$ at~$s$.

The center and radius of the osculating circle at a given point are called the \index{center of curvature}\emph{center of curvature} and \index{radius of curvature}\emph{radius of curvature} of the curve at that point.
\end{thm}

{

\begin{wrapfigure}{o}{31 mm}
\vskip-0mm
\centering
\includegraphics{mppics/pic-21}
\vskip0mm
\end{wrapfigure}

The \textit{osculating circle} might be a circle or a line.
In the latter case, the center of curvature is undefined and the radius of curvature is infinite.

The osculating circle $\sigma_s$ can be also defined as the unique circle (or line) that has \index{order of contact}\emph{second-order contact} with $\gamma$ at $s$;
that is, $\rho(\ell)\z=o(\ell^2)$, where $\rho(\ell)$ denotes the distance from $\gamma(s+\ell)$ to $\sigma_s$.

The following exercise is recommended to the reader familiar with the notion of \index{inversion}\emph{inversion}.

}

\begin{thm}{Advanced exercise}\label{ex:inverse}
Suppose $\gamma$ is a smooth plane curve that does not pass thru the origin.
Let $\hat \gamma$ be the inversion of $\gamma$ with respect to a circle centered at the origin.
Show that the osculating circle of $\hat\gamma$ at $s$ is the inversion of the osculating circle of $\gamma$ at~$s$.
\end{thm}

\section{Spiral lemma}
\label{spiral}
\index{spiral lemma}

The following lemma was proved by Peter Tait \cite{tait}
and later rediscovered by Adolf Kneser \cite{kneser}.

\begin{thm}{Lemma}\label{lem:spiral}
Assume $\gamma$ is a smooth plane curve with strictly decreasing positive signed curvature. Then the osculating circles of $\gamma$ are nested; that is, if $\sigma_s$ denotes the osculating circle of $\gamma$ at $s$,
then $\sigma_{s_0}$ lies in the open disc bounded by $\sigma_{s_1}$ for any $s_0<s_1$. \index{10sigma@$\sigma_s$ (osculating circle)}
\end{thm}

{

\begin{wrapfigure}{r}{31 mm}
\vskip-6mm
\includegraphics{mppics/pic-61}
\end{wrapfigure}

The picture shows a curve $\gamma$ as in the theorem and its osculating circles.
They form a peculiar foliation of the annulus; it meets the following property:
if a smooth function is constant on each osculating circle, then it must be constant in the annulus \cite[Lecture 10]{fuchs-tabachnikov}.
Also, the curve $\gamma$ is tangent to a circle of the foliation at each of its points.
However, it does not run along any of those circles.

}

\parbf{Proof.}
Let $\tan(s),\norm(s)$ be the Frenet frame,
$\omega(s)$, $r(s)$
the center and radius of curvature of~$\gamma$.
By \ref{prop:circline},  we have
\[\omega(s)=\gamma(s)+r(s)\cdot \norm(s).\]

Since $\skur>0$, we have that $r(s)\cdot\skur(s)=1$.
Therefore, applying Frenet formula \ref{eq:nu'}, we get that
\begin{align*}
\omega'(s)&=\gamma'(s)+r'(s)\cdot \norm(s)+r(s)\cdot \norm'(s)=
\\
&=\tan(s)+r'(s)\cdot \norm(s)-r(s)\cdot \skur(s)\cdot \tan(s)=
\\
&=r'(s)\cdot \norm(s).
\end{align*}
Since $\skur(s)$ is decreasing, $r(s)$ is increasing;
therefore $r'\ge 0$.
It follows that $|\omega'(s)|= r'(s)$ and $\omega'(s)$ points in the direction of $\norm(s)$.

Since $\norm'(s)=-\skur(s)\cdot\tan(s)$, the direction of $\omega'(s)$ cannot be constant on a nontrivial interval;
in other words, the curve $s\mapsto \omega(s)$ contains no line segments.

\begin{wrapfigure}[5]{o}{35 mm}
\vskip-0mm
\centering
\includegraphics{mppics/pic-84}
\end{wrapfigure}

In particular, $|\omega(s_1)-\omega(s_0)|\z<\length(\omega|_{[s_0,s_1]})$ for any $s_0<s_1$.
Therefore, 
\begin{align*}
|\omega(s_1)-\omega(s_0)|&<\length(\omega|_{[s_0,s_1]})=
\\
&=\int_{s_0}^{s_1}|\omega'(s)|\cdot ds=
\\
&=\int_{s_0}^{s_1}r'(s)\cdot ds=
\\
&=r(s_1)-r(s_0).
\end{align*}
In other words, the distance between the centers of $\sigma_{s_1}$ and $\sigma_{s_0}$
is strictly less than the difference between their radii, hence the result.
\qeds

The curve $s\mapsto \omega(s)$ is called the \index{evolute}\emph{evolute} of $\gamma$; 
it traces the centers of curvature of the curve. 
The evolute of $\gamma$ can be written as 
\[\omega(t)=\gamma(t)+\tfrac1{\skur(t)}\cdot \norm(t);\] 
in the proof, we showed that $(\tfrac1{\skur})'\cdot\norm$ is its velocity vector.

\begin{thm}{Exercise}\label{ex:evolute}
Let $\omega$ be the evolute of a smooth plane curve~$\gamma$.
Suppose $\gamma$ has positive signed curvature $\skur$ and $\skur' \neq 0$ at all points.
%%%%%%% Vertices haven't been defined yet.
Find the Frenet frame and curvature of \(\omega\) in terms of $\skur$ and Frenet frame $(\tan,\norm)$ of $\gamma$.
\end{thm}

The following theorem states formally that 
\textit{if you drive on the plane and turn the steering wheel to the left all the time,
then you will not be able to come back to the place you started.}

\begin{thm}{Theorem}\label{thm:spiral}
Assume $\gamma$ is a smooth plane curve with positive strictly monotone signed curvature. 
Then $\gamma$ is simple.
\end{thm}

The same statement holds true without assuming positivity of curvature; the proof requires only minor modifications.

\parbf{Proof.}
Note that $\gamma(s)$ lies on the osculating circle $\sigma_s$ of $\gamma$ at~$s$.
If $s_1\ne s_0$, then by \ref{lem:spiral}, $\sigma_{s_0}$ does not intersect $\sigma_{s_1}$.
Therefore, $\gamma(s_1)\ne \gamma(s_0)$,
hence the result.\qeds

\begin{thm}{Advanced exercise}\label{ex:3D-spiral}
Show that a 3-dimensional analog of the theorem does not hold.
That is, there are self-intersecting smooth space curves with strictly monotone curvature.
\end{thm}

\begin{thm}{Exercise}\label{ex:double-tangent}
Assume $\gamma$ is a smooth plane curve with positive strictly monotone signed curvature.

\begin{subthm}{ex:double-tangent:a}Show that no line can be tangent to $\gamma$ at two distinct points.
\end{subthm}

\begin{subthm}{ex:double-tangent:b}Show that no circle can be tangent to $\gamma$ at three distinct points. 
\end{subthm}

\end{thm}

{

\begin{wrapfigure}{o}{25 mm}
\vskip-6mm
\centering
\includegraphics{mppics/pic-25}
\vskip0mm
\end{wrapfigure}

Part \ref{SHORT.ex:double-tangent:a} does not hold if we allow the curvature to be negative; an example is shown on the picture.

}

\begin{thm}{Advanced exercise}\label{ex:spherical-spiral}
Show that a smooth spherical curve with nonvanishing torsion has no self-intersections.
\end{thm}

%% file: supporting-curves.tex
\chapter{Supporting curves}
\label{chap:supporting-curves}

When two plane curves touch each other without crossing, it is possible to control the signed curvature of one of them in terms of the signed curvature of the other.
This will be proved and used to study the global behavior of plane curves.

\section{Cooriented tangent curves}

Suppose $\gamma_1$ and $\gamma_2$ are smooth plane curves.
Recall that the curves $\gamma_1$ and $\gamma_2$ are tangent at the  time parameters $t_1$ and $t_2$
if $\gamma_1(t_1)=\gamma_2(t_2)$
and they share the tangent line at these time parameters.
In this case, the point $p=\gamma_1(t_1)=\gamma_2(t_2)$ is called a \index{point!of tangency}\emph{point of tangency} of the curves.

In this case, the velocity vectors $\gamma_1'(t_1)$ and $\gamma_2'(t_2)$ are parallel.
\begin{figure}[!ht]
\vskip-0mm
\centering
\includegraphics{mppics/pic-85}
\vskip-0mm
\end{figure}
If $\gamma_1'(t_1)$ and $\gamma_2'(t_2)$ have the same direction, we say that the curves are \index{cooriented and counteroriented!curves}\emph{cooriented},
if their directions are opposite, the curves are called {}\emph{counteroriented}.

If we reverse the parametrization of one of the curves, then cooriented curves become counteroriented and vice versa; so we can always assume the curves are cooriented at any given point of tangency.

\section{Supporting curves}

\begin{wrapfigure}[8]{o}{43 mm}
\vskip-4mm
\centering
\includegraphics{mppics/pic-86}
\vskip0mm
\end{wrapfigure}

Let $\gamma_1$ and $\gamma_2$ be two smooth plane curves that share a point 
\[p=\gamma_1(t_1)=\gamma_2(t_2);\] 
we assume that $p$ is not an endpoint for $\gamma_1$ nor $\gamma_2$.
Suppose that there is $\epsilon>0$ such that the arc $\gamma_2|_{[t_2-\epsilon, t_2+\epsilon]}$ lies in a closed plane region $R$ with the arc $\gamma_1|_{[t_1-\epsilon, t_1+\epsilon]}$ in its boundary,
then we say that $\gamma_1$ \index{supporting!curve}\emph{locally supports} $\gamma_2$ at the time parameters $t_1$ and $t_2$.
Furthermore, if the curves on the picture are oriented according to the arrows, then $\gamma_1$ supports $\gamma_2$ from the right at $p$ (as well as $\gamma_2$ supports $\gamma_1$ from the left at~$p$).

Suppose $\gamma_1$ is a simple curve that cuts the plane into two closed regions, one lies on the left and the other on the right from~$\gamma_1$.
We say that $\gamma_1$ \index{supporting!curve}\emph{globally supports} $\gamma_2$ at point $p=\gamma_2(t_2)$ 
if $\gamma_2$ runs in one of these closed regions, and 
$p$ lies on~$\gamma_1$.

Further, suppose $\gamma_2$ is a simple closed plane curve.
By Jordan's theorem (\ref{thm:jordan}), $\gamma_2$ cuts from the plane two closed regions, one is bounded and the other is unbounded.
We say that a point $p$ lies {}\emph{inside} (respectively, {}\emph{outside}) $\gamma_2$ if $p$ lies in the bounded region (respectively, unbounded) region.
In this case, we say that $\gamma_1$ supports $\gamma_2$ \index{supporting!curve}\emph{from the inside} (\index{supporting!curve}\emph{from the outside}) if $\gamma_1$ supports $\gamma_2$ and lies inside $\gamma_2$ (respectively outside $\gamma_2$). 

If $\gamma_1$ and $\gamma_2$ share a point $p=\gamma_1(t_1)=\gamma_2(t_2)$ and are not tangent at $t_1$ and $t_2$, then at time $t_2$ the curve $\gamma_2$ crosses $\gamma_1$  moving from one of its sides to the other.
It follows that $\gamma_1$ cannot locally support $\gamma_2$ at the time parameters $t_1$ and $t_2$.
Whence we get the following.

\begin{thm}{Definition-Observation}
Let $\gamma_1$ and $\gamma_2$ be two smooth plane curves.
Suppose $\gamma_1$ locally supports $\gamma_2$ at time parameters $t_1$ and $t_2$.
Then $\gamma_1$ is tangent to $\gamma_2$ at $t_1$ and $t_2$.

In such a case, if the curves are cooriented, and the region $R$ in the definition of supporting curves lies on the right (left) from the arc of $\gamma_1$, then we say that 
$\gamma_1$ supports $\gamma_2$ from the left (respectively right).
\end{thm}

We say that a smooth plane curve $\gamma$ has a \index{vertex of  curve}\emph{vertex} at $s$
if the signed curvature function is critical at $s$;
that is, if $\skur'(s)_\gamma=0$.
If in addition to that, $\gamma$ is simple, we could say that the point $p=\gamma(s)$ is a vertex of~$\gamma$.

\begin{thm}{Exercise}\label{ex:vertex-support}
Assume the osculating circle $\sigma_s$ of a smooth plane curve $\gamma$ at $s$ locally supports $\gamma$ at $p=\gamma(s)$.
Show that $p$ is a vertex of~$\gamma$.
\end{thm}

\section{Supporting test}

The following proposition resembles the second derivative test. 

\begin{thm}{Proposition}\label{prop:supporting-circline}
Let $\gamma_1$ and $\gamma_2$ be two smooth plane curves.

Suppose $\gamma_1$ locally supports $\gamma_2$ from the left (right) at the time parameters $t_1$ and $t_2$.
Then 
\[\skur_1(t_1)\ge \skur_2(t_2)\quad(\text{respectively}\quad \skur_1(t_1)\le \skur_2(t_2)),\]
where $\skur_1$ and $\skur_2$ denote the signed curvature of $\gamma_1$ and $\gamma_2$, respectively.

A partial converse also holds.
Namely, if $\gamma_1$ and $\gamma_2$ are tangent and cooriented at the time parameters $t_1$ and $t_2$,
and 
\[\skur_1(t_1)> \skur_2(t_2)\quad(\text{respectively}\quad \skur_1(t_1)< \skur_2(t_2)),\]
then $\gamma_1$ locally supports $\gamma_2$ from the left (right) at the time parameters $t_1$ and $t_2$. 
\end{thm}

\parbf{Proof.} Without loss of generality, we can assume that $t_1=t_2=0$, the shared point $\gamma_1(0)=\gamma_2(0)$ is the origin, and the velocity vectors $\gamma'_1(0)$, $\gamma'_2(0)$ point in the direction of $x$-axis.
Then small arcs $\gamma_1|_{[-\epsilon,+\epsilon]}$ and $\gamma_2|_{[-\epsilon,+\epsilon]}$ can be described as a graph 
$y=f_1(x)$ and $y=f_2(x)$ for smooth functions $f_1$ and $f_2$ such that $f_i(0)=0$ and $f_i'(0)=0$.
Note that $f_1''(0)=\skur_1(0)$, and $f_2''(0)=\skur_2(0)$ (see \ref{ex:curvature-graph})

Clearly, $\gamma_1$ supports $\gamma_2$ from the left (right) if 
\[f_1(x)\ge f_2(x)\quad(\text{respectively}\quad f_1(x)\le f_2(x))\]
for all sufficiently small values of~$x$.
Applying the second derivative test to the function $f_1-f_2$, we get the result.

The partial converse can be proved along the same lines.
\qeds

\begin{thm}{Advanced exercise}\label{ex:support}
Suppose that two smooth unit-speed simple plane curves $\gamma_0$ and $\gamma_1$ are tangent and cooriented at the point $p\z=\gamma_0(0)\z=\gamma_1(0)$.
Assume $\skur_0(s)\le\skur_1(s)$ for any~$s$.
Show that $\gamma_0$ locally supports $\gamma_1$ from the right at~$p$.

Give an example of two simple curves $\gamma_0$ and $\gamma_1$ satisfying the above condition such that $\gamma_0$ is closed, but does not support $\gamma_1$ at $p$ globally.
\end{thm}

According to \ref{thm:DNA}, for any closed smooth curve that runs in a unit disc, the average of its absolute curvature is at least~1; in particular, there is a point where the absolute curvature is at least~1.
The following exercise says that this statement also holds for loops.

\begin{thm}{Exercise}\label{ex:in-circle}
Assume a smooth plane loop $\gamma$ runs in a unit disc.
Show that there is a point on $\gamma$ with absolute curvature at least~1.
\end{thm}

\begin{thm}{Exercise}\label{ex:between-parallels-1}
Assume a closed smooth plane curve $\gamma$ runs between two parallel lines at distance 2 from each other.
Show that there is a point on $\gamma$ with absolute curvature at least~1.

Try to prove the same for a smooth plane loop.
\end{thm}

\begin{thm}{Exercise}\label{ex:in-triangle}
Assume a closed smooth plane curve $\gamma$ runs inside a triangle $\triangle$ with inradius~$1$; that is, a unit circle is tangent to all three sides of $\triangle$.
Show that $\gamma$ has curvature at least~$1$ at some point.
\end{thm}

The three exercises above are baby cases of \ref{ex:moon-rad}, but try to find a direct solution.

{

\begin{wrapfigure}{r}{32 mm}
\vskip-4mm
\centering
\includegraphics{mppics/pic-70}
\vskip0mm
\end{wrapfigure}

\begin{thm}{Exercise}\label{ex:lens}
Let $F$ be a plane figure bounded by two circle arcs $\sigma_1$ and $\sigma_2$ of signed curvature 1 that run from $x$ to~$y$.
Suppose $\sigma_1$ is shorter than~$\sigma_2$.
Assume a smooth arc $\gamma$ runs in $F$ and has both endpoints on $\sigma_1$.
Show that the absolute curvature of $\gamma$ is at least 1 at some parameter value.

\end{thm}

}

\section{Convex curves}

Recall that a plane curve is \index{convex!curve}\emph{convex} if it bounds a convex region.

\begin{thm}{Proposition}\label{prop:convex}
Suppose a simple closed smooth plane curve $\gamma$ bounds a figure~$F$.
Then $F$ is convex if and only if the signed curvature of $\gamma$ does not change the sign.
\end{thm}

\begin{thm}{Lens lemma}\label{lem:lens}
Let $\gamma$ be a smooth simple plane curve that runs from $x$ to~$y$.
Assume $\gamma$ runs strictly on the right side (left side) of the oriented line $xy$; only its endpoints $x$ and $y$ lie on the line.
Then $\gamma$ has a point with positive (respectively negative) signed curvature.
\end{thm}

{

\begin{wrapfigure}{o}{35 mm}
\vskip-4mm
\centering
\includegraphics{mppics/pic-22}
\vskip0mm
\end{wrapfigure}

The lemma fails for curves with self-intersections.
For example, the curve $\gamma$ on the picture always turns right, 
so it has negative curvature everywhere, but it lies on the right side of the line $xy$.

}

\begin{wrapfigure}[6]{i}{50 mm}
\vskip-3mm
\centering
\includegraphics{mppics/pic-24}
\end{wrapfigure}

\parbf{Proof.}
Choose points $p$ and $q$ on the line $xy$
so that the points $p, x, y, q$ appear in that order.
We can assume that $p$ and $q$ lie sufficiently far from $x$ and $y$, so the half-disc with diameter $pq$ contains~$\gamma$.

Consider the smallest disc segment with chord $[p,q]$ that contains~$\gamma$.
Its arc $\sigma$ supports $\gamma$ at some point $w=\gamma(t_0)$.

Let us parametrize $\sigma$ from $p$ to~$q$.
Note that $\gamma$ and $\sigma$ are tangent and cooriented at~$w$.
If not, then the arc of $\gamma$ from $w$ to $y$ would be trapped in the curvilinear triangle $xwp$ bounded by the line segment $[p,x]$ and the arcs of $\sigma$,~$\gamma$.
But this is impossible since $y$ does not belong to this triangle.

{

\begin{wrapfigure}{o}{50 mm}
\vskip-4mm
\centering
\includegraphics{mppics/pic-23}
\bigskip
\includegraphics{mppics/pic-230}
\end{wrapfigure}

It follows that $\sigma$ supports $\gamma$ at $t_0$ from the right.
By \ref{prop:supporting-circline}, 
\[\skur(w)_\gamma\ge \skur_\sigma >0.\]
\qedsf

\parit{Remark.}
Instead of taking the minimal disc segment, one can take a point $w$ on $\gamma$ that maximizes the distance to the line $xy$.
The same argument shows that the curvature at $w$ is nonnegative, which is slightly weaker than the required positive curvature.

}

\parbf{Proof of \ref{prop:convex};} \textit{only-if part.}
If $F$ is convex, then every tangent line of $\gamma$ supports~$\gamma$.
If a point moves along $\gamma$, the figure $F$ has to stay on one side from its tangent line;
that is, we can assume that each tangent line supports $\gamma$ on one side, say on the right.
Since a line has vanishing curvature, the supporting test (\ref{prop:supporting-circline}) implies that $\skur\ge 0$ at each point.

\begin{wrapfigure}{r}{35 mm}
\vskip-3mm
\centering
\includegraphics{mppics/pic-68}
\vskip0mm
\end{wrapfigure}

\parit{If part.}
Denote by $K$ the convex hull of~$F$.
If $F$ is not convex, then $F$ is a proper subset of~$K$.
Therefore, the boundary $\partial K$ contains a line segment that is not a part of $\partial F$.
In other words, there is a line that supports $\gamma$ at two points, say $x$ and $y$.
These points divide $\gamma$ in two arcs $\gamma_1$ and $\gamma_2$, both distinct from the line segment $[x,y]$.

One of the arcs $\gamma_1$ or $\gamma_2$ is parametrized from $x$ to $y$ and the other from $y$ to~$x$.
Passing to a smaller arc if necessary we can ensure that only its endpoints lie on the line. 
Applying the lens lemma, we get that the arcs $\gamma_1$ and $\gamma_2$ contain points with signed curvatures of opposite signs.
\qeds

\begin{thm}{Exercise}\label{ex:convex small}
Suppose $\gamma$ is a smooth simple closed plane curve of diameter larger than~$2$.
Show that $\gamma$ has a point with absolute curvature less than~$1$.
\end{thm}

\begin{wrapfigure}{r}{45 mm}
\vskip-8mm
\centering
\includegraphics{mppics/pic-713}
\vskip0mm
\end{wrapfigure}

\begin{thm}{Exercise}\label{ex:convex-lens}
Let $\gamma$ be a simple smooth plane arc with endpoints $p$ and~$q$.
Suppose that $\gamma$ has nonnegative signed curvature and $|p-q|$ is the maximal distance between points on~$\gamma$.
Show that the arc $\gamma$ and its chord $[p,q]$ bound a convex figure of the plane.
\end{thm}

\begin{thm}{Exercise}\label{ex:diameter-of-simple-curve}
Show that any simple smooth plane curve $\gamma$ with curvature at least 1 has diameter at most 2.

Try to prove that $\gamma$ lies in a unit disc.
\end{thm}

\section{Moon in a puddle}

The following theorem is a slight generalization of the theorem proved by Vladimir Ionin and German Pestov \cite{ionin-pestov}.\index{Ionin--Pestov theorem}
For convex curves, this result was known earlier \cite[§ 24]{blaschke}.

\begin{wrapfigure}{r}{18 mm}
\vskip-8mm
\centering
\includegraphics{mppics/pic-67}
\vskip-2mm
\end{wrapfigure}

\begin{thm}{Theorem}\label{thm:moon-orginal}
Assume $\gamma$ is a simple smooth plane loop with absolute curvature bounded by~$1$.
Then it surrounds a unit disc.
\end{thm}

This theorem gives a simple but nontrivial example of the so-called \index{local-to-global}\emph{local-to-global} theorems --- based on some local data (in this case, the curvature of a curve) we conclude a global property (in this case, the existence of a unit disc surrounded by the curve).

{

\begin{wrapfigure}{r}{33 mm}
\vskip-4mm
\centering
\includegraphics{mppics/pic-62}
\vskip0mm
\end{wrapfigure}

A straightforward approach would be to start with a disc in the region bounded by the curve and blow it up to maximize its radius.
However, as one may see from the spinner-like example on the picture it does not always lead to a solution --- a closed plane curve of curvature at most $1$ may surround a disc of radius smaller than $1$ that cannot be enlarged continuously.

}

\begin{thm}{Key lemma}\label{thm:moon}
Assume $\gamma$ is a simple smooth plane loop.
Then at one point of $\gamma$ (distinct from its base), its osculating circle globally supports $\gamma$ from the inside.
\end{thm}

First, let us show that the theorem follows from the lemma.

\parbf{Proof of \ref{thm:moon-orginal} modulo \ref{thm:moon}.}
Since $\gamma$ has absolute curvature at most 1, each osculating circle has radius at least~1.
According to the key lemma, one of the osculating circles $\sigma$ globally supports $\gamma$ from the inside.
In particular, $\sigma$ lies inside $\gamma$, whence the result.
\qeds

\parbf{Proof of \ref{thm:moon}.}
Denote by $F$ the closed region surrounded by $\gamma$.
We can assume that $F$ lies on the left from~$\gamma$.
Arguing by contradiction,
assume that the osculating circle at each point $p\in \gamma$ does not lie in~$F$.

\begin{figure}[!ht]
\vskip-0mm
\centering
\includegraphics{mppics/pic-32}
\vskip-2mm
\end{figure}

Given a point $p\in\gamma$ let us consider the maximal circle that lies completely in $F$ and tangent to $\gamma$ at~$p$.
This circle, say $\sigma$, will be called the {}\emph{incircle} of $F$ at~$p$;
its curvature $\skur_\sigma$ has to be larger than $\skur(p)_\gamma$.
Indeed, since $\sigma$ supports $\gamma$ from the left, by \ref{prop:supporting-circline} we have $\skur_\sigma\ge \skur(p)_\gamma$; in the case of equality, $\sigma$ is the osculating circle at~$p$.
The latter is impossible by our assumption.

It follows that $\sigma$ has to touch $\gamma$ at another point.
Otherwise, we can increase $\sigma$ slightly while keeping it inside~$F$.

Indeed, since $\skur_\sigma> \skur(p)_\gamma$, 
by \ref{prop:supporting-circline} we can choose a neighborhood $U$ of $p$ such that after a slight increase of $\sigma$, the intersection $U\cap \sigma$ is still in~$F$.
On the other hand, if $\sigma$ does not touch $\gamma$ at another point, then after some (maybe smaller) increase of $\sigma$ the complement $\sigma\setminus U$ is still in~$F$.
That is, a slightly increased $\sigma$ is still in $F$ --- a contradiction.

Choose a point $p_1$ on $\gamma$ that is distinct from its base point. 
Let $\sigma_1$ be the incircle at $p_1$.
Denote by $\gamma_1$ an arc of $\gamma$ from $p_1$ to a first point $q_1$ on $\sigma_1$.
Denote by $\hat\sigma_1$ and $\check\sigma_1$ two arcs of $\sigma_1$ from $p_1$ to $q_1$ such that the cyclic concatenation of $\hat\sigma_1$ and $\gamma_1$ surrounds~$\check\sigma_1$. 

\begin{wrapfigure}{r}{32 mm}
\vskip-8mm
\centering
\includegraphics{mppics/pic-64}
\caption*{Two ovals pretend to be circles.}
\vskip-2mm
\end{wrapfigure}

Let $p_2$ be the midpoint of $\gamma_1$.
Denote by $\sigma_2$ the incircle at $p_2$.

The circle $\sigma_2$ cannot intersect $\hat\sigma_1$.
Otherwise, if $\sigma_2$ intersects $\hat\sigma_1$ at some point $s$, then $\sigma_2$ has to have two more common points with $\check\sigma_1$, say $x$ and $y$ --- one for each arc of $\sigma_2$ from $p_2$ to~$s$.
Therefore, $\sigma_1\z=\sigma_2$ since these two circles have three common points: $s$, $x$, and~$y$. 
On the other hand, by construction, $p_2\in \sigma_2$ and $p_2\notin \sigma_1$ --- a contradiction.

Recall that $\sigma_2$ has to touch $\gamma$ at another point.
From above, it follows that $\sigma_2$ can only touch $\gamma_1$. 
Therefore we can choose an arc $\gamma_2\subset \gamma_1$ that runs from $p_2$ to a first point $q_2$ on $\sigma_2$.
Since $p_2$ is the midpoint of $\gamma_1$, we have that
\[\length \gamma_2< \tfrac12\cdot\length\gamma_1.\eqlbl{eq:length<length/2}\]

Repeating this construction recursively,
we get an infinite sequence of arcs $\gamma_1\supset \gamma_2\supset\dots$;
by \ref{eq:length<length/2}, we also get that 
\[\length\gamma_n\to0\quad\text{as}\quad n\to\infty.\] 
Therefore, the intersection $\gamma_1\cap\gamma_2\cap\dots$
contains a single point; denote it by $p_\infty$.

Let $\sigma_\infty$ be the incircle at $p_\infty$; it has to touch $\gamma$ at another point, say $q_\infty$.
The same argument as above shows that $q_\infty\in\gamma_n$ for any~$n$.
It follows that $q_\infty =p_\infty$ --- a contradiction.
\qeds

\begin{thm}{Exercise}\label{ex:moon-rad}
Assume a closed smooth curve $\gamma$ lies in a figure $F$ bounded by a simple closed plane curve.
Suppose $R$ is the maximal radius of discs that lie in~$F$.
Show that the absolute curvature of $\gamma$ is at least $\tfrac1R$ at some parameter value.
\end{thm}

\section{Four-vertex theorem}
\index{Four-vertex theorem}
{

\begin{wrapfigure}{r}{20 mm}
\vskip-8mm
\centering
\includegraphics{mppics/pic-26}
\vskip0mm
\end{wrapfigure}

Recall that a vertex of a smooth curve is defined as a critical point of its signed curvature;
in particular, any local minimum (or maximum) of the signed curvature is a vertex.
For example, every point of a circle is a vertex.

\begin{thm}{Theorem}\label{thm:4-vert}
Any smooth simple closed plane curve has at least four vertices.
\end{thm}

}

Evidently, any closed smooth curve has at least two vertices --- where the minimum and the maximum of the curvature are attained.
On the picture the vertices are marked;
the first curve has one self-intersection and exactly two vertices;
the second curve has exactly four vertices and no self-intersections.

The four-vertex theorem was first proved by Syamadas Mukhopadhyaya \cite{mukhopadhyaya} for convex curves.
One of our favorite proofs was given by Robert Osserman \cite{osserman}.
Our proof of the following stronger statement is based on the key lemma in the previous section.
For more on the subject, see \cite{petrunin-zamora:moon} and the references therein.

{

\begin{wrapfigure}[8]{r}{33 mm}
\vskip-8mm
\centering
\includegraphics{mppics/pic-63}
\vskip0mm
\end{wrapfigure}

\begin{thm}{Theorem}\label{thm:4-vert-supporting}
Any smooth simple closed plane curve is globally supported by its osculating circle at least at 4 distinct points;
two from the inside and two from the outside.
\end{thm}

\parbf{Proof of \ref{thm:4-vert} modulo \ref{thm:4-vert-supporting}.}
It is sufficient to show that if an osculating circle $\sigma$ at a point $p$ supports $\gamma$ locally, then $p$ is a vertex.

}

If not, then a small arc around $p$ has monotone curvature.
Applying the spiral lemma (\ref{lem:spiral}) we get that the osculating circles at this arc are nested.
In particular, the curve $\gamma$ crosses $\sigma$ at~$p$. 
Therefore $\sigma$ does not support $\gamma$ locally at~$p$.
\qeds

\parbf{Proof of \ref{thm:4-vert-supporting}.}
According to the key lemma (\ref{thm:moon}), there is a point $p\in\gamma$ such that its osculating circle supports $\gamma$ from the inside.
The curve $\gamma$ can be considered as a loop with the base at~$p$.
Therefore, the key lemma implies the existence of another point $q\in\gamma$ with the same property.

It gives a pair osculating circles that support $\gamma$ from the inside;
it remains to find another pair.

In order to get the osculating circles supporting $\gamma$ from the outside, one can repeat the proof of the key lemma taking instead of incircle the circle (or line) of maximal signed curvature that supports the curve from the outside, assuming that $\gamma$ is oriented so that the region on the left from it is bounded.\qeds

\parbf{Alternative end of proof.}
If one applies to $\gamma$ an inversion with respect to a circle whose center lies inside~$\gamma$, then the obtained curve $\gamma_1$ also has  two osculating circles that support $\gamma_1$ from the inside.
According to \ref{ex:inverse}, these osculating circles are inverses of the osculating circles of~$\gamma$.
Note that the region lying inside $\gamma$ is mapped to the region outside $\gamma_1$, and the other way around.
Therefore, these two circles (or lines) correspond to the osculating circles supporting $\gamma$ from the outside.\qeds

{

\begin{wrapfigure}{r}{20 mm}
\vskip-6mm
\centering
\includegraphics{mppics/pic-725}
\vskip0mm
\end{wrapfigure}

\begin{thm}{Exercise}\label{ex:2-squares}
Assume that a smooth simple closed plane curve \(\gamma\) is contained in a square with side $2$ and surrounds a square with diagonal~$2$.
Show that \(\gamma\) contains a point with curvature~$1$.
\end{thm}

}

\begin{thm}{Exercise}\label{ex:moon-area}
Suppose a smooth simple plane loop bounds a region of area $a$.
Show that it has curvature $\sqrt{\pi/a}$ at some point.
\end{thm}

\begin{wrapfigure}[5]{r}{25 mm}
\vskip-7mm
\centering
\includegraphics{mppics/pic-65}
\vskip0mm
\end{wrapfigure}

\begin{thm}{Advanced exercise}\label{ex:curve-crosses-circle}
Suppose $\gamma$ is a simple closed smooth plane curve, and $\sigma$ is a circle.
Assume $\gamma$ crosses $\sigma$ at the points $p_1,\dots,p_{2{\cdot} n}$, and these points appear in the same cyclic order on $\gamma$ and on $\sigma$.
Show that $\gamma$ has at least $2\cdot n$ vertices.

Construct an example of a simple closed smooth plane curve $\gamma$ with only $4$ vertices that crosses a given circle at arbitrarily many points.
\end{thm}

\begin{thm}{Advanced exercise}\label{ex:berk}
Let $\gamma$ be a simple smooth plane curve with absolute curvature bounded by~$1$.
Show that $\gamma$ surrounds two disjoint open unit discs if and only if its diameter is at least $4$;
that is, $\dist{p}{q}{}\ge 4$ for some points $p,q\in\gamma$. 
\end{thm}

%\begin{thm}{Advanced exercise}\label{ex:order} Show that the points $a,b,c,d$ guaranteed by \ref{thm:4-vert-supporting} can be chosen so that they appear in the same order on the curve and the osculating circles at $a$ and $c$ support the curve from the inside and $b$ and $d$ support the curve from the outside. \end{thm}

The following exercise is a version of the four-vertex theorem for space curves without parallel tangents; they exist by \ref{ex:no-parallel-tangents}.

\begin{thm}{Advanced exercise}\label{ex:4x0-torsion}
Let $\gamma$ be a closed smooth space curve that has no pair of points with parallel tangent lines.
Assume the curvature of $\gamma$ does not vanish at any point.
Show that $\gamma$ has at least four points with zero torsion.
\end{thm}

%% file: part-surfaces.tex
\arxiv{\cleardoublepage
\phantomsection
\AddToShipoutPictureBG*{\includegraphics[width=\paperwidth]{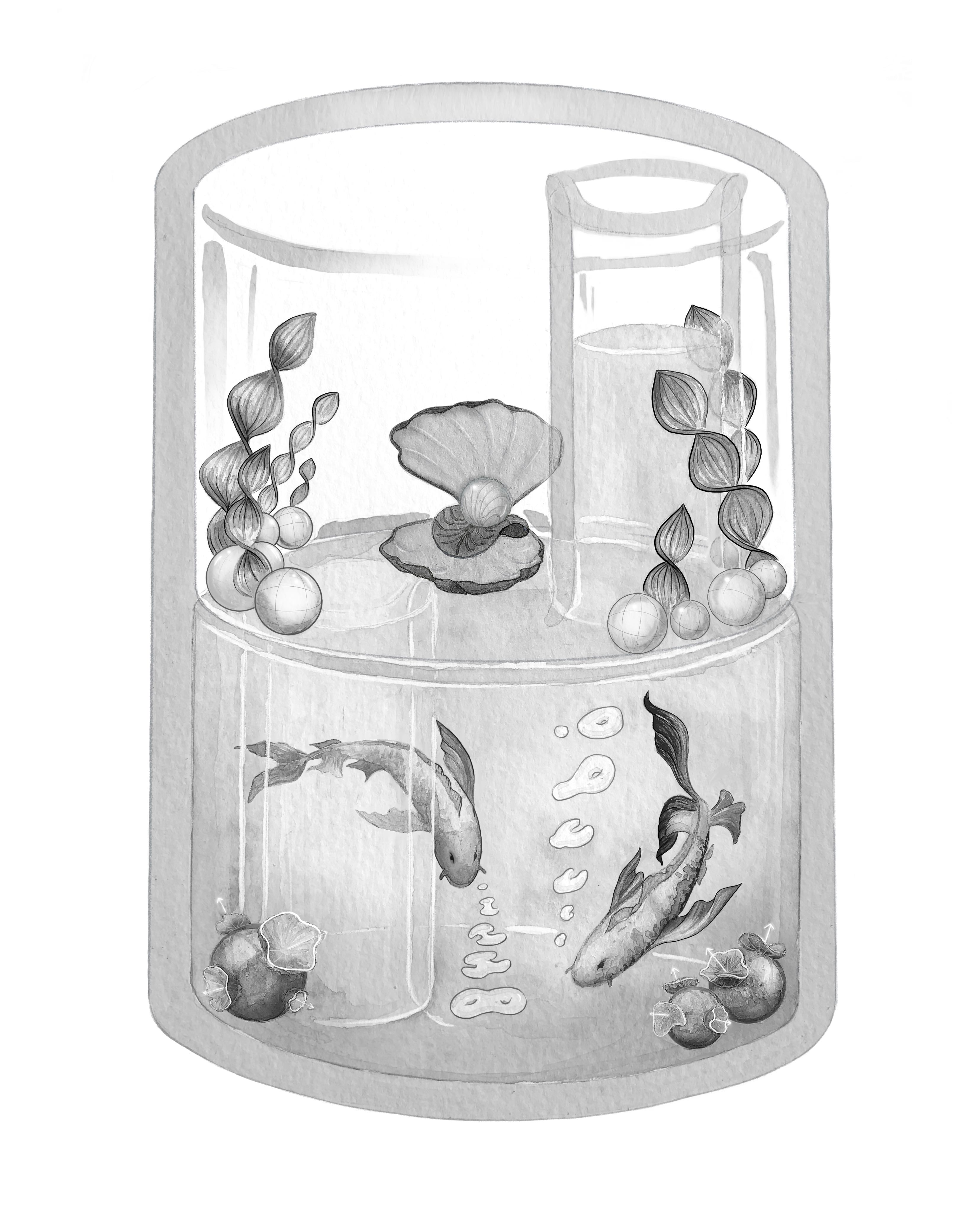}}
\cleardoublepage
\thispagestyle{empty}
\stepcounter{part}
\begin{center}
{\Huge\textbf{Part \Roman{part}:\qquad Surfaces}}\
\end{center}
\addcontentsline{toc}{part}{\Roman{part} Surfaces}
\clearpage
}
{\backgroundsetup{scale=.95,opacity=1,angle=0,vshift=10mm,contents={%
\includegraphics[width=\paperwidth]{pics/Surfaces}
}%
}

\cleardoublepage
\phantomsection
\stepcounter{part}
\addcontentsline{toc}{part}{\Roman{part} Surfaces}
\thispagestyle{empty}
\begin{center}
{\Huge\textbf{Part \Roman{part}:\qquad Surfaces}}\
\end{center}
\BgThispage
}

%% file: surfaces-def.tex
\chapter{Definitions}
\label{chap:surfaces-def}

\section{Topological surfaces}

We will be most interested in smooth surfaces defined in the following section.
The following general definition will be used only a few times.

A connected subset $\Sigma$ in the Euclidean space $\mathbb{R}^3$
is called a \index{surface}\index{topological!surface}\emph{topological surface} (more precisely an {}\emph{embedded topological surface without boundary}) 
if any point $p\in \Sigma$ admits a neighborhood $W$ in $\Sigma$
that can be parametrized by an open subset in the Euclidean plane; 
that is, there is a homeomorphism $V\to W$ from an open set $V\subset \mathbb{R}^2$; see Appendix~\ref{sec:topology}.

\section{Smooth surfaces}\label{sec:def-smooth-surface}

Recall that a function $f$ of two variables $x$ and $y$ is called \index{smooth!function}\emph{smooth} if all its partial derivatives $\frac{\partial^{m+n}}{\partial x^m\partial y^n}f$ are defined and are continuous in the domain of definition of~$f$. 

A connected set $\Sigma \subset \mathbb{R}^3$ is called a \index{surface}\index{smooth!surface}\emph{smooth surface}%
\footnote{We use it as a shortcut for the more precise term {}\emph{smooth regular embedded surface}.} if it can be described locally as a graph of a smooth function in an appropriate coordinate system.

More precisely, for any point $p\in \Sigma$ one can choose a coordinate system $(x,y,z)$ and a neighborhood $U\ni p$ such that
the intersection $W=U\cap \Sigma$ is a graph $z=f(x,y)$ of a smooth function $f$ defined in an open domain $V$ of the $(x,y)$-plane.

Note that $(x,y)\mapsto (x,y,f(x,y))$ defines a homeomorphism $V\to W$.
Therefore, a smooth surface is an example of a topological surface.

\parbf{Examples.}
The simplest example of a smooth surface is the $(x,y)$-plane 
$\Pi=\set{(x,y,z)\in\mathbb{R}^3}{z=0}$.
The plane $\Pi$ is a surface since
it can be described as the graph of the function $f(x,y)=0$.

All other planes are smooth surfaces as well since one can choose a coordinate system so that it becomes the $(x,y)$-plane.
We may also present a plane as a graph of a linear function 
$f(x,y)=a\cdot x+b\cdot y+c$ for some constants $a$, $b$, and $c$
(assuming the plane is not perpendicular to the $(x,y)$-plane, in which case a different coordinate system is needed).

A more interesting example is the unit sphere 
\[\mathbb{S}^2=\set{(x,y,z)\in\mathbb{R}^3}{x^2+y^2+z^2=1}.\]
This set is not a graph,
but it is locally a graph ---
it can be covered by the following 6 graphs:
\begin{align*}
z&=f_\pm(x,y)=\pm \sqrt{1-x^2-y^2},
\\
y&=g_\pm(x,z)=\pm \sqrt{1-x^2-z^2},
\\
x&=h_\pm(y,z)=\pm \sqrt{1-y^2-z^2},
\end{align*}
where each function $f_+$, $f_-$, $g_+$, $g_-$, $h_+$, and $h_-$ is defined in an open unit disc.
Any point $p\in\mathbb{S}^2$ lies in one of these graphs therefore $\mathbb{S}^2$ is a smooth surface.

\section{Surfaces with boundary}
A connected subset in a surface that is bounded by one or more piecewise
smooth curves is called a \index{surface!with boundary}\emph{surface with boundary}; such curves form the \index{boundary!line}\emph{boundary line} of the surface.

When we say {}\emph{surface} we usually mean a {}\emph{smooth surface without boundary}.
If needed, one may use the term {}\emph{surface with possibly nonempty boundary}.

\section{Proper, closed, and open surfaces}
If the surface $\Sigma$ is formed by a closed set in $\mathbb{R}^3$, then it is called \index{proper!surface}\emph{proper}.
For example, for any smooth function $f$ defined on the whole plane, its graph $z=f(x,y)$ is a proper surface.
The sphere $\mathbb{S}^2$ gives another example of a proper surface.

On the other hand, the open disc 
\[\set{(x,y,z)\in\mathbb{R}^3}{x^2+y^2<1,\  z=0}\]
is not a proper surface; this set is neither open nor closed in $\mathbb{R}^3$.

A compact surface without boundary is called \index{closed!surface}\emph{closed}
(this term is closely related to the closed curve, but has nothing to do with closed set).

A proper noncompact surface without boundary is called \index{open!surface}\emph{open} (again, the term open curve is relevant, but open set is not).

For example, the paraboloid $z=x^2+y^2$
is an open surface; the 
sphere $\mathbb{S}^2$ is a closed surface.
Note that \textit{any proper surface without boundary is either closed or open}.

The following claim is a three-dimensional analog of the plane separation theorem (\ref{ex:proper-curve}).
It follows from the so-called {}\emph{Alexander's duality} \cite{hatcher}.
We omit its proof; it would take us far away from the main subject.

\begin{thm}{Claim}\label{clm:proper-divides}
The complement of any proper topological surface without boundary (or, equivalently any open or closed topological surface) has exactly two connected components. 
\end{thm}

\section{Implicitly defined surfaces}

\begin{thm}{Proposition}\label{prop:implicit-surface}
Let $f\:\mathbb{R}^3\to \mathbb{R}$ be a smooth function.
Suppose $0$ is a regular value of $f$;
that is, $\nabla_p f\ne 0$ at any point $p$ such that $f(p)=0$.
Then any connected component $\Sigma$ of the level set  $f(x,y,z)=0$ is a smooth surface.
\end{thm}

\parbf{Proof.}
Fix $p\in\Sigma$.
Since $\nabla_p f\ne 0$ we have 
$f_x(p)\ne 0$,
$f_y(p)\ne 0$, or
$f_z(p)\ne 0$.
We may assume that $f_z(p)\ne 0$;
otherwise, permute the coordinates $x,y,z$.

The implicit function theorem (\ref{thm:imlicit}) implies that a neighborhood of $p$ in $\Sigma$ is the graph $z=h(x,y)$ of a smooth function $h$ defined on an open domain in $\mathbb{R}^2$.
It remains to apply the definition of smooth surface (Section~\ref{sec:def-smooth-surface}).
\qeds

\begin{thm}{Exercise}\label{ex:hyperboloids}
For which constants $\ell$ is the level set $x^2+y^2-z^2=\ell$
a~smooth surface?
\end{thm}

\section{Local parametrizations}
\index{parametrization}

Let $U$ be an open domain in $\mathbb{R}^2$, and $s\:U\to \mathbb{R}^3$ be a smooth map.
We say that $s$ is \index{regular!parametrization}\emph{regular} if its Jacobian matrix has maximal rank;
in this case, it means that the derivatives $s_u$ and $s_v$ are linearly independent at any $(u,v)\in U$;
equivalently $s_u\times s_v\ne 0$, where $\times$ denotes the vector product.

\begin{thm}{Proposition}\label{prop:graph-chart}
If $s\:U\to \mathbb{R}^3$ is a smooth regular embedding of an open connected set $U\subset \mathbb{R}^2$, then its image $\Sigma=s(U)$ is a smooth surface.
\end{thm}

\parbf{Proof of \ref{prop:graph-chart}.}
Let $s(u,v)=(x(u,v),y(u,v),z(u,v))$.
Since $s$ is regular, its Jacobian matrix
\[\Jac s=
\renewcommand\arraystretch{1.3}
\begin{pmatrix}
x_u&x_v\\
y_u&y_v\\
z_u&z_v
\end{pmatrix}
\]
has rank two at any point $(u,v)\in U$.

Choose a point $p\in \Sigma$; by shifting the $(x,y,z)$ and $(u,v)$ coordinate systems, we may assume that $p=(0,0,0)=s(0,0)$.
Permuting the coordinates $x,y,z$ if necessary, we may assume that 
the matrix $\left(\begin{smallmatrix}
x_u&x_v\\
y_u&y_v
\end{smallmatrix}\right)$
is invertible at the origin.
Note that this is the Jacobian matrix of the map $(u,v)\mapsto (x(u,v),y(u,v))$.

The inverse function theorem (\ref{thm:inverse}) implies that there is a smooth regular map
$w\:(x,y)\mapsto (u,v)$ defined on an open set $W\ni 0$ in the $(x,y)$-plane
such that $w(0,0)=(0,0)$ and  $s\circ w(x,y)=(x,y,f(x,y))$ where $f=z\circ w$.
That is, the subset $s\circ w(W)\subset \Sigma$ is 
the graph of $f$.

Again, by the inverse function theorem, $w(W)$ is an open subset of $U$.
Since $s$ is an embedding, our graph is open in~$\Sigma$;
that is, there is an open set $V\subset \mathbb{R}^3$ such that $s\circ w(W)=V\cap \Sigma$ is a graph of a smooth function.

Since $p$ is arbitrary, we get that $\Sigma$ is a smooth surface.
\qeds

\begin{thm}{Exercise}\label{ex:9-surf}
Construct a smooth regular injective map $s\:\mathbb{R}^2\to\mathbb{R}^3$ such that its image is \textit{not} a surface.
\end{thm}

If we have $s$ and $\Sigma$ as in the proposition, then we say that $s$ is a \index{smooth!parametrization}\emph{smooth parametrization} of the surface~$\Sigma$. 

Not all smooth surfaces can be described by such a parametrization;
for example, the sphere $\mathbb{S}^2$ cannot.
However, \textit{any smooth surface $\Sigma$ admits a local parametrization at any point $p\in\Sigma$}; that is,  $p$ admits an open neighborhood $W\subset \Sigma$ with a smooth regular parametrization~$s$.
In this case, any point in $W$ can be described by two parameters, usually denoted by $u$ and $v$, 
which are called \index{local coordinates}\emph{local coordinates} at~$p$.
The map $s$ is called a \index{chart}\emph{chart} of~$\Sigma$.

If $W$ is a graph $z=h(x,y)$ of a smooth function $h$, then the map
\[s\:(u,v)\mapsto (u,v,h(u,v))\] is a chart.
Indeed, $s$ has an inverse $(u,v,h(u,v))\mapsto (u,v)$ which is continuous;
that is, $s$ is an embedding.
Further,
$s_u=(1,0,h_u)$, and $s_v=(0,1,h_v)$. 
Whence the partial derivatives $s_u$ and $s_v$ are linearly independent;
that is, $s$ is a regular map.

\begin{thm}{Corollary}\label{cor:reg-parmeterization}
A connected set $\Sigma\subset \mathbb{R}^3$ is a smooth surface if and only if a neighborhood of any point in $\Sigma$ can be covered by a chart.
\end{thm}

A function $g\: \Sigma \to \mathbb{R}$ defined on a smooth surface $\Sigma$ is said to be \index{smooth!function}\emph{smooth} if for any chart $s \: U\to \Sigma$,
the composition $g\circ s$ is smooth; that is, all partial derivatives $\frac{\partial^{m+n}}{\partial u^m\partial v^n}(g\circ s)$ are defined and continuous in the domain of definition.

\begin{thm}{Exercise}\label{ex:smooth-fun(surf)}
Let $\Sigma\subset \mathbb{R}^3$ be a smooth surface.
Show that a function $g\:\Sigma\to\mathbb{R}$ is smooth if and only if for any point $p\in \Sigma$ there is a smooth function $h\:N\to\mathbb{R}$ defined in a neighborhood $N\subset \mathbb{R}^3$ of $p$ such that the equality $g(q)=h(q)$ holds for $q\in \Sigma\cap N$.

Construct a smooth surface $\Sigma$ and a smooth function $g\:\Sigma\to\mathbb{R}$ that cannot be extended to a smooth function $h\:\mathbb{R}^3\to\mathbb{R}$.
\end{thm}

\begin{thm}{Exercise}\label{ex:inversion-chart}
Consider the following map 
\[s(u,v)=(\tfrac{2\cdot u}{1+u^2+v^2},\tfrac{2\cdot v}{1+u^2+v^2},\tfrac{2}{1+u^2+v^2}).\]
Show that $s$ is a chart of the unit sphere centered at $(0,0,1)$; describe the image of~$s$.
\end{thm}

\begin{wrapfigure}{o}{31 mm}
\vskip-3mm
\centering
\includegraphics{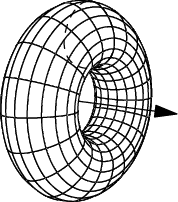}
\vskip0mm
\end{wrapfigure}

Let $\gamma(t)=(x(t),y(t))$ be a plane curve.
Recall that the \index{surface!of revolution}\emph{surface of revolution} of the curve $\gamma$ around the $x$-axis can be described as the 
image of the map 
\[(t, s)\mapsto (x(t), y(t)\cdot\cos s,y(t)\cdot\sin s).\]
The parameters $t$ and $s$ are called \index{latitude}\emph{latitude} and \index{longitude}\emph{longitude} respectively;
for fixed $t$ or $s$ the obtained curves are called \index{parallel}\emph{parallels} or
\index{meridian}\emph{meridians} respectively. 
The parallels are formed by circles in the plane perpendicular to the axis of rotation.
The curve $\gamma$ is called the \index{generatrix}\emph{generatrix} of the surface.

\begin{thm}{Exercise}\label{ex:revolution}
Assume $\gamma$ is a simple closed smooth plane curve that does not intersect the $x$-axis.
Show that the surface of revolution around the $x$-axis with generatrix $\gamma$ is a smooth surface.
\end{thm}

\section{Global parametrizations}\label{sec:global-parametrizations}
\index{parametrization}

A surface can be described by an embedding from a known surface.

For example, consider the ellipsoid
\[\Theta=\set{(x,y,z)\in\mathbb{R}^3}{\tfrac{x^2}{a^2}+\tfrac{y^2}{b^2}+\tfrac{z^2}{c^2}=1}\]
for positive constants $a$, $b$, and~$c$.
By \ref{prop:implicit-surface}, $\Theta$ is a smooth surface.
Indeed, let $h(x,y,z)=\tfrac{x^2}{a^2}+\tfrac{y^2}{b^2}+\tfrac{z^2}{c^2}$,
then
\[\nabla h(x,y,z)=(\tfrac{2}{a^2}\cdot x,\tfrac{2}{b^2}\cdot y,\tfrac{2}{c^2}\cdot z).\]
Therefore, $\nabla h\ne0$ if $h=1$; that is, $1$ is a regular value of~$h$.
It remains to observe that $\Theta$ is connected.

The surface $\Theta$ can be defined as the image of the map $ \mathbb{S}^2\to\mathbb{R}^3$, defined as the restriction of the following map to the unit sphere $\mathbb{S}^2$:
\[(x,y,z)\z\mapsto (a\cdot x, b\cdot y,c\cdot z).\]

A map $f\:\Sigma \to \mathbb{R}^3$ (or $f\:\Sigma \to \mathbb{R}^2$) is smooth if each of its coordinate functions is smooth.
Further, a smooth map $f \: \Sigma \to \mathbb{R}^3$ is called a 
\emph{smooth parametrized surface} if it is an embedding, and for any chart $s \:U\to \Sigma$,
the composition $f\circ s$ is regular;
that is,  the two vectors 
$\frac{\partial}{\partial u}(f\circ s)$ and $\frac{\partial}{\partial v}(f\circ s)$ are linearly independent.
In this case, the image $\Sigma^{*}=f(\Sigma)$ is a smooth surface.
The latter follows since for any chart $s\:U\to \Sigma$ the composition $f\circ s\:U\to \Sigma^{*}$ is a chart of $\Sigma^{*}$. 

The map $f$ is called a \index{diffeomorphism}\emph{diffeomorphism} from $\Sigma$ to $\Sigma^{*}$;
the surfaces $\Sigma$ and $\Sigma^{*}$ are said to be {}\emph{diffeomorphic} if there is a diffeomorphism $f\:\Sigma\to\Sigma^{*}$.
The following exercise implies that \textit{being diffeomorphic} is an equivalence relation for surfaces.

\begin{thm}{Exercise}\label{ex:inv-diffeomorphism}
Show that the inverse of a diffeomorphism is a diffeomorphism.
\end{thm}

\begin{thm}{Advanced exercise}\label{ex:star-shaped-disc}
Show the following.

\begin{subthm}{ex:plane-n}
Complements of $n$-point sets in the plane are diffeomorphic to each other.
\end{subthm}

\begin{subthm}{ex:star-shaped-disc:smooth}
Open convex subsets of the plane bounded by smooth closed curves are diffeomorphic to each other.
\end{subthm}

\begin{subthm}{ex:star-shaped-disc:nonsmooth}
Any pair of open convex subsets of the plane are diffeomorphic to each other.
\end{subthm}

\begin{subthm}{ex:star-shaped-disc:star-shaped}
Open star-shaped subsets of the plane are diffeomorphic to each other.
\end{subthm}
\end{thm}

%% file: first-order.tex
\chapter{First-order structure}
\label{chap:first-order}
\section{Tangent plane}

\begin{thm}{Definition}\label{def:tangent-vector}
Let $\Sigma$ be a smooth surface.
A vector $\vec w$ is \index{tangent!vector}\emph{tangent} to $\Sigma$ at $p$ if and only if there is a curve $\gamma$ that runs in $\Sigma$ and has $\vec w$ as a velocity vector at $p$;
that is, $p=\gamma(t)$ and $\vec w=\gamma'(t)$ for some~$t$.
\end{thm}

\begin{thm}{Proposition and definition}\label{def:tangent-plane}
Let $\Sigma$ be a smooth surface and $p\z\in \Sigma$.
Then the set of tangent vectors of $\Sigma$ at $p$ forms a plane;
this plane is called the \index{tangent!plane}\emph{tangent plane} of $\Sigma$ at~$p$.

Moreover, if $s\:U\to \Sigma$ is a local chart and $p=s(u_0,v_0)$, then 
the tangent plane of $\Sigma$ at $p$ is spanned by vectors $s_u(u_0,v_0)$ and $s_v(u_0,v_0)$.
\end{thm}

The tangent plane to $\Sigma$ at $p$ is usually denoted by $\T_p$ or $\T_p\Sigma$.
This plane might be considered as a linear subspace of $\mathbb{R}^3$ or as a parallel plane passing thru $p$;
the latter is sometimes called the \index{tangent!plane}\emph{affine tangent plane}.
The affine tangent plane can be interpreted as the best approximation at~$p$ of the surface $\Sigma$ by a plane.
More precisely, 
it has \index{order of contact}\emph{first-order contact} with $\Sigma$ at $p$;
that is, $\rho(q)\z=o(|p-q|)$, where $q\in \Sigma$, and $\rho(q)$ denotes the distance from $q$ to $\T_p$.

An assignment of a tangent vector $\vec v_p$ to each point $p$ of the surface $\Sigma$ is called a \emph{tangent vector field} if $\vec v_p$ depends smoothly on $p$.
More precisely, for any chart $s$ we have
\[\vec v_{s(u,v)}=a(u,v)\cdot s_u+b(u,v)\cdot s_v\]
for smooth functions $a$ and $b$ defined in the domain of $s$.

\parbf{Proof.}
Fix a chart $s$ at~$p$.
Assume $\gamma$ is a smooth curve that starts at~$p$.
We can assume that $\gamma$ is covered by the chart;
in particular, there are smooth functions $t\mapsto u(t)$ and $t\mapsto v(t)$ such that 
\[\gamma(t)=s(u(t),v(t)).\]
Applying the chain rule, we get
\[\gamma'=s_u\cdot u'+ s_v\cdot v';\]
that is, $\gamma'$ is a linear combination of $s_u$ and $s_v$.

The smooth functions $t\mapsto u(t)$ and $t\mapsto v(t)$ can be chosen arbitrarily.
Therefore any linear combination of $s_u$ and $s_v$ is a tangent vector at~$p$. 
\qeds

\begin{thm}{Exercise}\label{ex:tangent-normal}
Let $f:\mathbb{R}^3\to\mathbb{R}$ be a smooth function with $0$ as a regular value, and $\Sigma$ be a surface described as a connected component of the level set $f(x,y,z)=0$.
Show that the tangent plane $\T_p\Sigma$ is perpendicular to the gradient $\nabla_pf$ at any point $p\in\Sigma$.
\end{thm}

\begin{thm}{Exercise}\label{ex:vertical-tangent}
Let $\Sigma$ be a smooth surface and $p\in\Sigma$.
Choose $(x,y,z)$-coordinates.
Show that a neighborhood of $p$ in $\Sigma$ is a graph $z=f(x,y)$ of a smooth function $f$ defined on an open subset in the $(x,y)$-plane if and only if the tangent plane $\T_p$ is not {}\emph{vertical}; that is, if $\T_p$ is not perpendicular to the $(x,y)$-plane.
\end{thm}

\begin{thm}{Exercise}\label{ex:tangent-single-point}
Show that if a smooth surface $\Sigma$ meets a plane $\Pi$ at a single point $p$, then $\Pi$ is tangent to $\Sigma$ at~$p$.
\end{thm}

\section{Directional derivative}\label{sec:dirder}

In this section, we extend the definition of directional derivative to smooth functions defined on smooth surfaces.
First, let us recall the standard definition for the plane.

Given a point $p\in \mathbb{R}^2$, a vector $\vec w\in \mathbb{R}^2$ and a function $f\:\mathbb{R}^2\to\mathbb{R}$, 
consider the function
$h(t)=f(p+t\cdot\vec w)$.
Then the directional derivative of $f$ at $p$ along $\vec w$ is defined as \index{10d@$D_{\vec{w}}f$ (directional derivative)}
\[D_{\vec w}f(p)\df h'(0).\]

Recall that a function $f\: \Sigma \to \mathbb{R}$ is said to be smooth if for any chart $s\: U \to \Sigma$, the composition $f \circ s$ is smooth.

\begin{thm}{Proposition-Definition}\label{def:directional-derivative}
Let $f$ be a smooth function defined on a smooth surface $\Sigma$.
Suppose $\gamma$ is a smooth curve in $\Sigma$ that starts at $p$ with velocity vector $\vec{w}\in \T_p$;
that is, $\gamma(0)=p$, and $\gamma'(0)=\vec{w}$.
Then the derivative $(f\circ\gamma)'(0)$
depends only on $f$, $p$, and $\vec{w}$;
it is called the \index{directional derivative}\emph{directional derivative} of $f$ along $\vec{w}$ at $p$
and is denoted by
\[D_{\vec{w}}f,\quad D_{\vec{w}}f(p), \quad\text{or}\quad D_{\vec{w}}f(p)_\Sigma\] 
--- we may omit $p$ and $\Sigma$ if it is clear from the context.

Moreover, if $(u,v)\mapsto s(u,v)$ is a local chart and $\vec{w}=a\cdot s_u +b\cdot s_v$ at $p$, then 
\[D_{\vec{w}}f=a\cdot (f\circ s)_u+b\cdot (f\circ s)_v.\]

\end{thm}

Our definition agrees with the standard definition of the directional derivative if $\Sigma$ is a plane.
Indeed, in this case, $\gamma(t)=p+\vec w\cdot t$ is a curve in $\Sigma$ that starts at $p$ with velocity vector $\vec{w}$.
For a general surface, the point $p+\vec w\cdot t$ might not lie on the surface, and the value $f(p+\vec w\cdot t)$ might be undefined.
In this case, the standard definition does not work.

\parbf{Proof.}
Without loss of generality, we may assume that $p=s(0,0)$ and the curve $\gamma$ is covered by the chart $s$;
if not we can chop~$\gamma$.
In this case, 
\[\gamma(t)=s(u(t),v(t))\]
for smooth functions $u,v$ defined in a neighborhood of $0$ such that 
$u(0)\z=v(0)\z=0$.

Applying the chain rule, we get that
\begin{align*}
\gamma'(0)&=u'(0)\cdot s_u+v'(0)\cdot s_v
\end{align*}
at $(0,0)$.
Since $\vec{w}=\gamma'(0)$ and the vectors $s_u$, $s_v$ are linearly independent, we get that $a=u'(0)$ and $b=v'(0)$.

Applying the chain rule again, we get that
\[
(f\circ\gamma)'(0)=a\cdot (f\circ s)_u+b\cdot (f\circ s)_v.
\eqlbl{eq:f-gamma}
\]
at $(0,0)$.

Notice that the left-hand side in \ref{eq:f-gamma} does not depend on the choice of the chart $s$, and the right-hand side depends only on $p$, $\vec w$, $f$, and~$s$. 
It follows that $(f\circ\gamma)'(0)$ depends only on $p$, $\vec w$, and~$f$.

The last statement follows from \ref{eq:f-gamma}.
\qeds

\begin{thm}{Advanced exercise}\label{ex:lin-ind-chart}
Let $\vec x$ and $\vec y$ be vector fields on a smooth surface $\Sigma$.
Suppose that $\vec x_p$ and $\vec y_p$ are linearly independent at some point $p\in \Sigma$.
Construct two functions $u$ and $v$ in a neighborhood of $p$ such that 
\begin{align*}
D_{\vec x} u&>0,
&
D_{\vec y} u&=0,
&
D_{\vec x} v&=0,
&
D_{\vec y} v&>0.
\end{align*}

Conclude that there is a chart of $\Sigma$ at $p$ such that $\vec x$ and $\vec y$ are tangent to its coordinate lines.
\end{thm}

\section{Tangent vectors as functionals}

In this section, we introduce a more conceptual way to define tangent vectors.
We will not use this approach in the sequel, but it is better to know about it.

A tangent vector $\vec w\in \T_p$ to a smooth surface $\Sigma$ 
defines a linear functional%
\footnote{The term {}\emph{functional} is used for functions that take a function as an argument and return a number.} $D_{\vec w}$ that swallows a smooth function $\phi$ on $\Sigma$ and spits its directional derivative $D_{\vec w}\phi$.
The functional $D_{\vec w}$ obeys the product rule:
\[D_{\vec w}(\phi\cdot\psi)=(D_{\vec w}\phi)\cdot \psi(p)+\phi(p)\cdot(D_{\vec w}\psi).
\eqlbl{eq:tangent-functional}\]

It is not hard to show that the tangent vector $\vec w$ is completely determined by the functional $D_{\vec w}$.
Moreover, tangent vectors at $p$ can be \textit{defined} as linear functionals on the space of smooth functions
that satisfy the product rule \ref{eq:tangent-functional}.

This definition grabs the key algebraic property of tangent vectors.
It might be a less intuitive way to think about tangent vectors, but it is often convenient to use in the proofs. 
For example, \ref{def:directional-derivative} becomes a tautology.

\section{Differential of map}\label{sec:differential}

Recall that a smooth map from a domain on the plane to space is defined in \ref{sec:Multivariable calculus}.
Let \( \Sigma_0 \) be a smooth surface in \( \mathbb{R}^3 \).
A map \(s\:\Sigma_0 \to \mathbb{R}^3 \) is called \index{smooth!map}\emph{smooth} if, for any chart \( w\:U \to \Sigma_0 \), the composition \( s \circ w \) defines a smooth map \( U \to \mathbb{R}^3 \).

Any smooth map $s$ from a surface $\Sigma_0$ to $\mathbb{R}^3$ can be described by its coordinate functions 
$ s(p)=(x(p),y(p),z(p))$. 
To take a directional derivative of the map we should take the  directional derivative of each of its coordinate functions.
\[D_{\vec{w}} s\df(D_{\vec{w}}x,D_{\vec{w}}y,D_{\vec{w}}z).\]

Assume $s$ maps one smooth surface $\Sigma_0$ to another $\Sigma_1$.
Let $p_0\in \Sigma_0$ and $p_1=s(p_0)$.
Choose a curve $\gamma_0$ in $\Sigma_0$ such that $\gamma_0(0)=p_0$ and $\gamma_0'(0)=\vec w$.
Observe that $\gamma_1= s\circ \gamma_0$ is a smooth curve in $\Sigma_1$. 
Curve $\gamma_1$ lies in $\Sigma_1$ and $\gamma_1(0)\z= p_1$;
hence its velocity vector $\gamma_1'(0)$ is in $\T_{ p_1}\Sigma_1$.
By the definition of the directional derivative, we have $D_{\vec w} s=\gamma_1'(0)$.
Therefore, $D_{\vec w} s\in \T_{p_1}\Sigma_1$ for any $\vec w\in \T_{p_0}$.

By \ref{def:directional-derivative},
$\vec w \mapsto D_{\vec w} s$ defines a linear map $\T_{p_0}\Sigma_0\z\to \T_{ p_1}\Sigma_1$;
that is,
\[D_{c\cdot \vec w} s=c\cdot D_{\vec w} s
\quad\text{and}\quad D_{\vec v+ \vec w} s=D_{\vec v} s+ D_{\vec w} s\]
for any $c\in\mathbb{R}$ and $\vec v, \vec w\in\T_{p_0}$.
The map $d_{p_0} s\:\T_{p_0}\Sigma_0\z\to \T_{ p_1}\Sigma_1$ defined by
\[d_{p_0} s\:\vec w \mapsto D_{\vec w} s\]
is called the \index{differential}\emph{differential} of $s$ at~$p_0$.

The differential $d_{p_0} s$ can be described by a $2{\times}2$-matrix $M$ in orthonormal bases of $\T_{p_0}\Sigma_0$ and $\T_{p_1}\Sigma_1$.
Set $\jac_{p_0} s=|\det M|$; this value  
does not depend on the choice of orthonormal bases in $\T_{p_0}\Sigma_0$ and $\T_{p_1}\Sigma_1$.%
\label{page:|L|}\index{10d@$d_p f$ (differential)}%
\footnote{The value $\jac_{p_0} s$ has the following geometric meaning:
if $R_0$ is a region in $\T_{p_0}$, and $R_1=(d_{p_0} s)(R_0)$ is its image, then
\[\area R_1=\jac_{p_0} s \cdot \area R_0.\]
This identity will play a key role in the definition of surface area.}

Let $r\:\Sigma_1\to\Sigma_2$ be another smooth map between smooth surfaces $\Sigma_1$ and $\Sigma_2$,
then 
\[d_{p_0}( r\circ s)=d_{p_1} r \circ d_{p_0} s,\]
and therefore,
\[\jac_{p_0}( r\circ s)
=
\jac_{p_1} r\cdot\jac_{p_0} s .\eqlbl{eq:jac-composition}\]

The constructions above can be applied to a chart $s$;
in this case, the surface $\Sigma_0$ is an open domain in $\mathbb{R}^2$.
Then, the value $\jac_{p_0} s$ can be found using the following formulas:
\begin{align*}
\jac s
&=|s_v\times s_u|=
\\
&=\sqrt{\langle s_u, s_u\rangle\cdot\langle s_v, s_v\rangle -\langle s_u, s_v\rangle^2}=
\\
&=\sqrt{\det[(\Jac s)^\top\cdot \Jac s ]},
\end{align*}
where $\Jac s$ denotes the Jacobian matrix of $s$; see Appendix~\ref{sec:Multivariable calculus}.

\section{Surface integral and area}

Let $\Sigma$ be a smooth surface, and $h\:\Sigma\to\mathbb{R}$ be a smooth function.
Let us define the integral of $h$ along a Borel set $R\subset \Sigma$;
most of the time we will apply this definition to surfaces with boundary.

Recall that $\jac_ps$ is defined in the previous section.
Assume there is a chart $(u,v)\mapsto s(u,v)$ of $\Sigma$ defined on an open set $U\subset\mathbb{R}^2$ such that $R\subset s(U)$.
In this case, set
\[\iint_R h\df \iint_{s^{-1}(R)} h\circ s(u,v)\cdot \jac_{(u,v)}s  \cdot du\cdot dv.\eqlbl{eq:area-def}\]

By the substitution rule (\ref{thm:mult-substitution}), the right-hand side in \ref{eq:area-def} does not depend on the choice of~$s$.
That is, if $s_1\:U_1\to \Sigma$ is another chart such that $s_1(U_1)\supset R$, then 
\[\iint_{s^{-1}(R)} h\circ s(u,v)\cdot \jac_{(u,v)}s  \cdot du\cdot dv=\iint_{s_1^{-1}(R)} h\circ s_1(u,v)\cdot \jac_{(u,v)}s_1  \cdot du\cdot dv.\]
In other words, the defining identity \ref{eq:area-def} makes sense.

A general region $R$ can be subdivided into regions $R_1,R_2\dots$ such that each $R_i$ lies in the image of a chart.
After that one could define the integral along $R$ as the sum
\[\iint_Rh
\df
\iint_{R_1}h+\iint_{R_2}h+\dots\]
In case $R$ is compact or $h \geq 0$, it is straightforward to check that the value $\iint_Rh$ does not depend on the choice of such subdivision. 

The area of a region $R$ in a smooth surface $\Sigma$ is defined as the surface integral 
\[\area R=\iint_R1.\]

The following proposition is the substitution rule for surface integral.

\begin{thm}{Area formula}\label{prop:surface-integral}
Suppose $s\:\Sigma_0\to \Sigma_1$ is a diffeomorphism between smooth surfaces.
Consider a region $R\subset \Sigma_0$ and a smooth function $f\:\Sigma_1\to\mathbb{R}$.
Assume that either $R$ is compact or $f$ is non-negative.
Then 
\[\iint_R (f\circ s)\cdot \jac s
=
\int_{s(R)}f.
\]
In particular, if $f\equiv 1$, we have
\[\iint_R \jac s=\area (s(R)).\]
\end{thm}

\parbf{Proof.}
Follows from \ref{eq:jac-composition} and the definition of surface integral.
\qeds

Let $\Sigma_1$ and $\Sigma_2$ be two smooth surfaces.
A map $f\:\Sigma_1\to \Sigma_2$ is called \index{length-nonincreasing}\emph{length-nonincreasing} if, for any curve $\gamma$ in $\Sigma_1$ we have $\length\gamma\z\ge \length (f\circ\gamma)$.
The following theorem provides a more natural definition of area.
Despite its intuitive statement, the proof is well beyond the scope of our book;
it is based on a generalization of the area formula that works for Lipschitz maps \cite[3.2.3]{federer}.

\begin{thm}{Theorem}\label{thm:area-axioms}
The area functional satisfies the following properties:

\begin{subthm}{thm:area-axioms:aditivity}
Sigma-additivity:
Let $R_1,R_2,\dots$ be a sequence of disjoint Borel sets in a smooth surface.
Then 
\[\area (R_1\cup R_2\cup \dots)=\area R_1+\area R_2+{}\dots\]
\end{subthm}

\begin{subthm}{thm:area-axioms:monotonicity}
Monotonicity:
Let $f\:\Sigma_1\to \Sigma_2$ be length-nonincreasing map between two smooth surfaces.
Suppose that $R_1\subset \Sigma_1$ and $R_2\subset \Sigma_2$ are Borel sets such that $f(R_1)\supset R_2$.
Then 
\[\area R_1\ge \area R_2.\]
\end{subthm}

\begin{subthm}{thm:area-axioms:unit}
Unit square has unit area. 
\end{subthm}

Moreover, the area functional is uniquely defined by these properties on all Borel sets.
\end{thm}

\parit{Remark.}
The notion of the area of a surface is closely related to the length of a curve.
However, to define length we use a different idea --- it was defined as the least upper bound on the lengths of inscribed polygonal lines.
It turns out that an analogous definition does not work even for very simple surfaces.
The latter is shown by a classical example --- \textit{Schwarz's boot}.
This example and different approaches to the notion of the area are discussed in a popular article by Vladimir Dubrovsky~\cite{dubrovsky}.

\section{Normal vector and orientation}

A surface $\Sigma$ is called \index{oriented surface}\emph{oriented} if it is equipped with a unit normal vector field $\Norm$ or \index{10nu@$\Norm$ (normal field)}
\index{normal!field}\emph{normal field};
that is, a continuous map $p\mapsto \Norm(p)$ such that $\Norm(p)\perp\T_p$ and $|\Norm(p)|=1$ for any~$p$.
The choice of the field $\Norm$ is called the {}\emph{orientation} of~$\Sigma$.
A surface $\Sigma$ is called {}\emph{orientable} if it can be oriented.
Each orientable surface admits two orientations: $\Norm$ and $-\Norm$.

Let $\Sigma$ be a smooth oriented surface with unit normal field $\Norm$.
The map $\Norm\:\Sigma\to \mathbb{S}^2$ defined by $p\mapsto \Norm(p)$ is called the \index{spherical!map}\emph{spherical map} or \index{Gauss map}\emph{Gauss map}.

\begin{wrapfigure}{r}{42 mm}
\vskip-7mm
\centering
\includegraphics{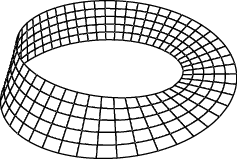}
\vskip-1mm
\end{wrapfigure}

For surfaces, the spherical map plays essentially the same role as the tangent indicatrix for curves.

The Möbius strip on the picture gives an example of a nonorientable surface --- there is no choice of a normal vector field that is continuous along the middle of the strip (it changes the sign if you try to go around).

Each surface is locally orientable.
In fact, each chart $s(u,v)$ admits an orientation 
\[\Norm=
\frac{s_u\times s_v}
{\left|s_u\times s_v\right|}.\]
Indeed, the vectors $s_u$ and $s_v$ are tangent vectors at $p$; 
since they are linearly independent, their vector product does not vanish, and it is perpendicular to the tangent plane.
Evidently, $(u,v)\mapsto \Norm(u,v)$ is a continuous map.
Therefore, $\Norm$ is a unit normal field. 

\begin{thm}{Exercise}\label{ex:const-normal}
Suppose a smooth surface $\Sigma$ has constant unit normal vector $\nu_0$.
Show that $\Sigma$ lies in a plane perpendicular to $\nu_0$.
\end{thm}

\begin{thm}{Exercise}\label{ex:implicit-orientable}
Let $0$ be a regular value of a smooth function $h:\mathbb{R}^3\z\to\mathbb{R}$.
Show that a connected component of the level set $h(x,y,z)=0$ is an oriented surface.
\end{thm}

Recall that any proper surface $\Sigma$ without boundary in the Euclidean space divides it into two connected components (\ref{clm:proper-divides}).
Therefore, we can choose the unit normal field on $\Sigma$ that points into one of the components of the complement; so, we obtain the following.

\begin{thm}{Observation}
Any smooth proper surface without boundary in the Euclidean space is orientable.
\end{thm}

In particular, the Möbius strip does not lie on a proper smooth surface without boundary.

\section{Sections}

\begin{thm}{Lemma}\label{lem:reg-section}
Let $\Sigma$ be a smooth surface.
Suppose $f\:\Sigma\z\to\mathbb{R}$ is a smooth function.
Then for any constant $r_0$ there is an arbitrarily close value $r$ such that 
each connected component of the level set $L_r=\set{x\in\Sigma}{f(x)=r}$ is a smooth curve.
\end{thm}

\parbf{Proof.}
The surface $\Sigma$ can be covered by a countable set of charts $s_i\:U_i\z\to \Sigma$.
The composition $f\circ s_i$ is a smooth function for any~$i$.
By Sard's lemma (\ref{lem:sard}), $r$ is a regular value of each $f\circ s_i$ for almost all $r$.

Fix such a value $r$ sufficiently close to $r_0$, and consider its level set $L_r\subset \Sigma$.
Any point on $L_r$ lies in the image of one of the charts.
From above, it admits a neighborhood which is a smooth curve;
hence the result.
\qeds

\begin{thm}{Advanced exercise}\label{ex:plane-section}
Let $\Pi$ be the $(x,y)$-plane and $A \subset \Pi$ be any closed subset.
Construct an open smooth surface $\Sigma$ such that $\Sigma \cap \Pi = A$.
\end{thm}

The exercise above says that sections of smooth surfaces might look ugly.
The following corollary makes it possible to perturb the plane so that the section becomes nice.

\begin{thm}{Corollary}
Let $\Sigma$ be a smooth surface.
Then for any plane $\Pi$ there is a parallel plane $\Pi^{*}$ that lies arbitrarily close to $\Pi$ and such that the intersection $\Sigma\cap\Pi^{*}$ is a union of disjoint smooth curves.
\end{thm}

%% file: surfaces-curvature.tex
\chapter{Curvatures}
\label{chap:surface-curvature}

In the previous chapter, we learned that the tangent plane to any smooth surface has first-order contact with it.
This means that up to first order, all surfaces locally look the same at a given point.
This is not the case for second-order approximations, and this is the subject of this chapter.

\section{Tangent-normal coordinates} \label{sec:lmn}
\index{tangent-normal coordinates}

Fix a point $p$ in a smooth oriented surface~$\Sigma$;
let $\Norm$ be its normal field.
Consider a coordinate system $(x,y,z)$ with origin at $p$ such that the $(x,y)$-plane coincides with $\T_p$, and the $z$-axis points in the direction of the normal vector $\Norm(p)$. Such coordinates are called \emph{tangent-normal coordinates}. 
By \ref{ex:vertical-tangent}, we can present $\Sigma$ locally at $p$ as a graph $z=f(x,y)$ of a smooth function. 
Note that 
\begin{align*}
f(0,0)&=0,
&
f_x(0,0)&=0,
&
f_y(0,0)&=0.
\end{align*}
The first equality holds since $p=(0,0,0)$ lies on the graph, and the other two mean that the tangent plane at $p$ is horizontal.

\begin{wrapfigure}[7]{o}{42 mm}
\vskip-4mm
\centering
\includegraphics{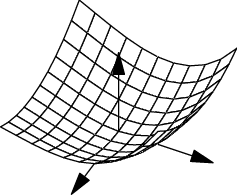}
\vskip-3mm
\end{wrapfigure}

Set 
\begin{align*}
\ell&=f_{xx}(0,0),
\\
m&=f_{xy}(0,0)=f_{yx}(0,0),
\\
n&=f_{yy}(0,0).
\end{align*}
The \textit{Taylor series} 
for $f$ at $(0,0)$ up to the second-order term can be then written as
\[f(x,y)=\tfrac12(\ell\cdot x^2+2\cdot m\cdot x\cdot y+n\cdot y^2)+o(x^2+y^2).\]
The values $\ell$, $m$, and $n$ are completely determined by this equation.\index{10lmn@$\ell$, $m$, $n$ (Hessian-matrix components)}
The graph of Taylor polynomial
\[z=\tfrac12\cdot(\ell\cdot x^2+2\cdot m\cdot x\cdot y+n\cdot y^2)\]
is called \index{osculating!paraboloid}\emph{osculating paraboloid}.
It has \index{order of contact}\emph{second-order contact} with $\Sigma$ at~$p$.

Note that 
\[\ell\cdot x^2+2\cdot m\cdot x\cdot y+n\cdot y^2=\langle M_p\cdot (\begin{smallmatrix}
x\\y
\end{smallmatrix}), (\begin{smallmatrix}
x\\y
\end{smallmatrix})\rangle,\]
where 
\[M_p=\begin{pmatrix}
 \ell
 &m
 \\
 m
 &n
 \end{pmatrix}
 =\begin{pmatrix}
 f_{xx}(0,0)
 &f_{xy}(0,0)
 \\
 f_{yx}(0,0)
 &f_{yy}(0,0)
 \end{pmatrix}.
\eqlbl{eq:hessian}
\]
is the so-called \index{Hessian matrix}\emph{Hessian matrix} of $f$ at $(0,0)$.\index{10m@$M_p$ (Hessian matrix)}

\section{Principal curvatures}\label{sec:Principal curvatures}

The tangent-normal coordinates are almost canonical; they depend on the point $p$ and are unique up to a rotation of the $(x,y)$-plane.
Rotating the $(x,y)$-plane results in rewriting 
the Hessian matrix $M_p$ in the new basis.

Since the matrix $M_p=\left(\begin{smallmatrix}
 \ell
 &m
 \\
 m
 &n
 \end{smallmatrix}\right)$ is symmetric, by the spectral theorem (\ref{thm:spectral}) it is diagonalizable by orthogonal matrices.
That is, by rotating the $(x,y)$-plane we can assume that $m=0$.
In this case,
\[M_p=\begin{pmatrix}
 k_1
 &0
 \\
 0
 &k_2
 \end{pmatrix},
\]
the diagonal components $k_1$ and $k_2$ of $M_p$ are called the \index{principal curvatures and directions}\emph{principal curvatures} of $\Sigma$ at $p$;\index{10k@$k_1$, $k_2$ (principal curvatures)}
they might also be denoted as $k_1(p)$ and $k_2(p)$, or $k_1(p)_\Sigma$ and $k_2(p)_\Sigma$;
if we need to emphasize that we compute them at the point $p$ for the surface~$\Sigma$.
We will assume 
\[k_1\le k_2\]
unless indicated otherwise.

If $z=f(x,y)$ is a local graph representation of $\Sigma$ in these coordinates, then 
\[f(x,y)=\tfrac12\cdot(k_1\cdot x^2+k_2\cdot y^2)+o(x^2+y^2).\]

The principal curvatures can be also defined as the eigenvalues of the Hessian matrix $M_p$.
The eigendirections of $M_p$ are called the {}\emph{principal directions} of $\Sigma$ at~$p$.
If $k_1(p)\ne k_2(p)$, then $p$ has exactly two principal directions, which are perpendicular to each other;
if $k_1(p)\z= k_2(p)$, then all tangent directions at $p$ are principal.

If we reverse the orientation of $\Sigma$, then the principal curvatures switch their signs and indexes: $k_1$ converts into $-k_2$ and $k_2$ into $-k_1$.

A smooth curve on a surface $\Sigma$ that always runs in the principal directions is called a \index{line of curvature}\emph{line of curvature} of~$\Sigma$.
If $k_1(p)\ne k_2(p)$, then there is a chart $(u,v)\mapsto s(u,v)$ at $p$ with coordinate lines formed by lines of curvature (it follows from \ref{ex:lin-ind-chart}).
Since principal directions are orthogonal we have $s_u\perp s_v$.

\begin{thm}{Exercise}\label{ex:line-of-curvature}
Assume a smooth surface $\Sigma$ is mirror-symmetric with respect to a plane $\Pi$.
Suppose $\Sigma$ and $\Pi$ intersect along a smooth curve~$\gamma$.
Show that $\gamma$ is a line of curvature of~$\Sigma$.
\end{thm}

\section{More curvatures}\label{sec:More curvatures}

Let $p$ be a point on an oriented smooth surface~$\Sigma$.
Recall that $k_1(p)$ and $k_2(p)$ denote the principal curvatures at~$p$.

The product 
\[K(p)=k_1(p)\cdot k_2(p)\]
is called the \index{10k@$K$ (Gauss curvature)}\index{Gauss curvature}\emph{Gauss curvature} at~$p$.
We may denote it by $K$, $K(p)$, or $K(p)_\Sigma$ if we need to specify the point and the surface.

Since the determinant is equal to the product of the eigenvalues, we get
\[K=\ell\cdot n-m^2,\]
where 
$M_p\z=
(\begin{smallmatrix}
\ell&m
\\
m&n
\end{smallmatrix}
)
$ is the Hessian matrix.

Reversing the orientation of $\Sigma$ does not change the Gauss curvature.
In particular, the Gauss curvature is well-defined for nonoriented surfaces.

\begin{thm}{Exercise}\label{ex:gauss+orientable}
Show that any surface with positive Gauss curvature is orientable. 
\end{thm}

The sum \index{10h@$H$ (mean curvature)}
\[H(p)=k_1(p)+ k_2(p)\] 
is called the \index{curvature}\index{mean curvature}\emph{mean curvature}%
\footnote{Some authors define it as $\tfrac12\cdot(k_1(p)+ k_2(p))$ --- the mean value of the principal curvatures. It suits the name better, but it is not as convenient when it comes to computations.}
at~$p$.
We may also denote it by $H(p)_\Sigma$.
The mean curvature can also be interpreted as the trace of the Hessian matrix $M_p\z=
(\begin{smallmatrix}
\ell&m
\\
m&n
\end{smallmatrix}
)$;
that is,
\[H=\ell+n.\] 
Reversing the orientation of $\Sigma$ changes the sign of the mean curvature.

A surface with vanishing mean curvature is called \index{minimal surface}\emph{minimal}.

\begin{thm}{Exercise}\label{ex:re-scale-surface-curvature}
Let $\Sigma$ be an oriented surface, and let $\Sigma_{\lambda}$ be a scaled copy of $\Sigma$ with factor $\lambda > 0$; that is, $\Sigma_{\lambda}$ consists of the points $\lambda \cdot x$ with $x \in \Sigma$. Show that
\[K(\lambda\cdot p)_{\Sigma_{\lambda}}
= \tfrac{1}{\lambda^2}\cdot K(p)_{\Sigma}
\quad\text{and}\quad
H(\lambda \cdot p)_{\Sigma_{\lambda}} = \tfrac1\lambda\cdot H(p)_{\Sigma}\]
for any $p\in \Sigma$.  
\end{thm}

\section{Shape operator}\label{sec:shape}

In the following definitions, we use the notion of directional derivative and differential defined in \ref{sec:dirder} and \ref{sec:differential}.

Let $p$ be a point on a smooth surface $\Sigma$ with orientation defined by the unit normal field $\Norm$.
Given $\vec w\in \T_p$,
the \index{shape operator}\emph{shape operator}%
\footnote{
The following bilinear forms on a tangent plane  
\begin{align*}
\mathrm{I}(\vec v,\vec w)&=\langle\vec v,\vec w\rangle,
&
\mathrm{II}(\vec v,\vec w)&=\langle\Shape\vec v,\vec w\rangle,
&
\mathrm{III}(\vec v,\vec w)&=\langle\Shape\vec v,\Shape\vec w\rangle
\end{align*}
are called the \index{fundamental form}\emph{first, second, and third fundamental forms}, respectively.
Historically, these forms were introduced before the shape operator, but we will not touch them.
} of $\vec w$ is defined by
\[\Shape_p\vec w=-D_{\vec w}\Norm.\]
Equivalently, the shape operator can be defined by
\[\Shape=-d\Norm,\eqlbl{eq:shape=-L}\] 
where $d\Norm$ denotes the differential of the spherical map $\Norm\:\Sigma\to\mathbb{S}^2$; that is, $d_p\Norm(\vec v)=(D_{\vec v}\Norm)(p)$.

Recall that $d_p\Norm$ is a linear map $\T_p\Sigma\to \T_{\Norm(p)}\mathbb{S}^2$.
Note that $\T_p\Sigma$ coincides with $\T_{\Norm(p)}\mathbb{S}^2$ --- both of them are normal subspaces to $\Norm(p)$.
Therefore, $\Shape_p$ is indeed a linear operator $\T_p\Sigma\to \T_p\Sigma$ (the latter also will follow from \ref{thm:shape-chart}).

The value of the shape operator on a tangent vector $\vec w\in \T_p\Sigma$ can be denoted by $\Shape\vec w$
if it is clear from the context which base point $p$ and which surface $\Sigma$ we are working with;
otherwise, we may use the notations 
\[\Shape_p(\vec w)\quad\text{or}\quad \Shape_p(\vec w)_\Sigma.\]

\begin{thm}{Theorem}\label{thm:shape-chart}
Suppose $(u,v)\mapsto s(u,v)$ is a smooth map to a smooth surface $\Sigma$ with unit normal field $\Norm$.
Then 
\begin{align*}
\langle \Shape(s_u), s_u\rangle 
&=\langle s_{uu},\Norm\rangle,
&
\langle \Shape(s_v), s_u\rangle 
&=\langle s_{uv},\Norm\rangle,
\\
\langle \Shape(s_u), s_v\rangle 
&=\langle s_{uv},\Norm\rangle,
&
\langle \Shape(s_v), s_v\rangle 
&=\langle s_{vv},\Norm\rangle,
\\
\langle \Shape(s_u), \Norm\rangle 
&=0,
&
\langle \Shape(s_v), \Norm\rangle 
&=0
\end{align*}
for any $(u,v)$.

\end{thm}

\parbf{Proof.}
We will use the shortcut $\Norm=\Norm(u,v)$ for $\Norm(s(u,v))$,
so 
\[
\begin{aligned}
\Shape(s_u)&=-D_{s_u}\Norm=-\Norm_u,
&
\Shape(s_v)&=-D_{s_v}\Norm=-\Norm_v.
\end{aligned}
\eqlbl{eq:shape=norm_u}
\]

Note that $\Norm$ is a unit vector orthogonal to $s_u$ and $s_v$;
therefore
\begin{align*}
\langle \Norm,s_u\rangle&\equiv0,
&
\langle \Norm,s_v\rangle&\equiv0,
&
\langle \Norm,\Norm\rangle&\equiv1.
\end{align*}
Taking partial derivatives of these two identities we get
\begin{align*}
\langle \Norm_u,s_u\rangle+\langle \Norm,s_{uu}\rangle&=0,
&
\langle \Norm_v,s_u\rangle+\langle \Norm,s_{uv}\rangle&=0,
\\
\langle \Norm_u,s_v\rangle+\langle \Norm,s_{uv}\rangle&=0,
&
\langle \Norm_v,s_v\rangle+\langle \Norm,s_{vv}\rangle&=0,
\\
2\cdot\langle \Norm_u,\Norm\rangle&=0,
&
2\cdot\langle \Norm_v,\Norm\rangle&=0.
\end{align*}
It remains to plug in the expressions from \ref{eq:shape=norm_u}.
\qeds

\begin{thm}{Exercise}\label{ex:self-adjoint}
Show that the shape operator is \index{self-adjoint operator}\emph{self-adjoint}; that is,
\[\langle \Shape\vec u,\vec v\rangle=\langle \vec u,\Shape\vec v\rangle\]
for any $\vec u,\vec v\in\T_p$.
\end{thm}

Recall that the components $\ell$, $m$, and $n$ of the Hessian matrix are defined in Section~\ref{sec:lmn}.

\begin{thm}{Corollary}\label{cor:Shape(ij)}
Let $z=f(x,y)$ be a local representation of a smooth surface $\Sigma$ in the tangent-normal coordinates at~$p$.
Suppose that its Hessian matrix at $p$ is $(\begin{smallmatrix}
\ell&m\\ m&n
\end{smallmatrix})$.
Then 
\begin{align*}
\Shape\vec i&=\ell\cdot \vec i+m\cdot \vec j,
&
\Shape\vec j&=m\cdot \vec i+n\cdot\vec j,
\end{align*}
where $\vec i,\vec j,\vec k$\index{10i@$\vec i$, $\vec j$, $\vec k$ (standard basis)} is the standard basis in $\mathbb{R}^3$.
That is, the multiplication by the Hessian matrix describes its shape operator.
\end{thm}

This corollary illustrates the close relationship between the curvatures of a surface and its shape operator; the principal curvatures of $\Sigma$ at $p$ are the eigenvalues of $\Shape_p$, the principal directions are the eigendirections of $\Shape_p$, the Gaussian curvature is the determinant of $\Shape_p$, and the mean curvature is the trace of $\Shape_p$.

Since the Hessian matrix is symmetric, the corollary also implies that $\Shape$ is self-adjoint which gives a way to solve \ref{ex:self-adjoint}.

\parbf{Proof.}
The map $s\:(u,v)\mapsto (u,v,f(u,v))$ defines a chart of $\Sigma$ at~$p$.
Furthermore, 
\begin{align*}
s_u(0,0)&=\vec i,
&
s_v(0,0)&=\vec j,
&
\Norm(0,0)&=\vec k,
\\
s_{uu}(0,0)&=\ell\cdot \vec k,
&
s_{uv}(0,0)&=m\cdot \vec k,
&
s_{vv}(0,0)&=n\cdot \vec k.
\end{align*}
It remains to apply \ref{thm:shape-chart}.
\qeds

\begin{thm}{Corollary}\label{cor:intK}
Let $\Sigma$ be a smooth surface whose Gaussian curvature never vanishes.
Assume that its spherical map $\Norm\:\Sigma\to\mathbb{S}^2$ is injective.
Then
\[\iint_\Sigma |K| = \area(\Norm(\Sigma)).\]
\end{thm}

\parbf{Proof.}
Choose an orthonormal basis of $\T_p$ consisting of principal directions,
so the shape operator can be expressed by the matrix
$(\begin{smallmatrix}
 k_1
 &0
 \\
 0
 &k_2
 \end{smallmatrix})$.

Since $\Shape_p=-d_p\Norm$, it follows from \ref{cor:Shape(ij)} that
\[
\jac_p\Norm
=
\left|\det\bigl(\begin{smallmatrix}
k_1 & 0\\
0   & k_2
\end{smallmatrix}\bigr)\right|
= |K(p)| \ne 0
\]
for every point $p\in\Sigma$.
In particular, the map $\Norm$ is regular and, since it is injective,
it defines a diffeomorphism $\Sigma\to\Norm(\Sigma)$.
It remains to apply the area formula (\ref{prop:surface-integral}).
\qeds

\begin{thm}{Exercise}\label{ex:normal-curvature=const}
Let $\Sigma$ be a smooth surface with orientation defined by a unit normal field $\Norm$.
Suppose the principal curvatures of $\Sigma$ are 1 at all points.

\begin{subthm}{ex:normal-curvature=const:a} Show that $\Shape_p(\vec w)=\vec w$ for any $p\in\Sigma$ and $\vec w\in \T_p\Sigma$.
\end{subthm}

\begin{subthm}{ex:normal-curvature=const:b}
Show that the point $q=p+\Norm(p)$ does not depend on $p\in\Sigma$.
Conclude that $\Sigma$ is a subset of the unit sphere centered at~$q$.
\end{subthm}

\end{thm}

\begin{thm}{Advanced exercise}\label{ex:normal-curvature=0}
Let $\Sigma$ be a smooth surface with orientation defined by a unit normal field $\Norm$,
and let $Z_0\subset \Sigma$ be a connected set with vanishing shape operator.
Show that $Z_0$ lies in a plane.
\end{thm}

We define the {}\emph{angle} between two oriented surfaces at a point of their intersection $p$ as the angle between their normal vectors at~$p$.

\begin{thm}{Exercise}\label{ex:shape-curvature-line}
Assume two smooth oriented surfaces $\Sigma_1$ and $\Sigma_2$ intersect at constant angle along a smooth curve~$\gamma$.
Show that if $\gamma$ is a curvature line in $\Sigma_1$, then it is also a curvature line in $\Sigma_2$.

Conclude that if a smooth surface $\Sigma$ intersects a plane or sphere along a smooth curve $\gamma$ at a constant angle,
then $\gamma$ is a curvature line of~$\Sigma$.
\end{thm}

\begin{thm}{Exercise}\label{ex:equidistant}
Let $\Sigma$ be a closed smooth surface with orientation defined by a unit normal field $\Norm$.

\begin{subthm}{ex:equidistant:smooth}
Show that if $t$ is sufficiently close to zero, then the set 
\[\Sigma_t=\set{p+t\cdot \Norm(p)}{p\in\Sigma}\] 
is a smooth surface.
\end{subthm}

\begin{subthm}{ex:equidistant:area}
Show that for all $t$ sufficiently close to zero we have
\[\area\Sigma_t=\area\Sigma-t\cdot \iint_\Sigma H+t^2\cdot \iint_\Sigma K,\]
where $H$ and $K$ denote mean and Gauss curvature of~$\Sigma$.
\end{subthm}

\end{thm}

\begin{thm}{Advanced exercise}\label{ex:flat-plane}
Let $\Sigma$ be a smooth oriented surface parametrized by a coordinate rectangle as $(u,v)\mapsto s(u,v)$.
Set $\vec u=\tfrac{s_u}{|s_u|}$, $\vec v=\tfrac{s_v}{|s_v|}$;
denote by $\Norm(u,v)$ the normal unit vector at $s(u,v)$.

Suppose $\vec u$ and $\vec v$ are principal directions at each point;
let $0$ and $k\z=k(u,v)$ be their principal curvatures.
Assume $k$ does not vanish and $|s_u|=1$ on some $u$-coordinate line.

\begin{subthm}{ex:flat-plane:orthonormal}
Show that $\Norm(u,v)$, $\vec u(u,v)$, and $\vec v(u,v)$ form an orthonormal frame for any $(u,v)$.
\end{subthm}

\begin{subthm}{ex:flat-plane:depend}
Show that $\Norm(u,v)$, $\vec u(u,v)$, and $\vec v(u,v)$ depend only on $v$.
Conclude that $u$-coordinate lines are line segments.
\end{subthm}

\begin{subthm}{ex:flat-plane:depend-u}
Show that $s_{uu}$ is proportional to $s_u$ at all points.
Use it to show that $|s_u|=1$ at all points.
\end{subthm}

\begin{subthm}{ex:flat-plane:linear}
Show that for any fixed $v_0$ the value $\tfrac1{k(u,v_0)}$ depends linearly on $u$.
\end{subthm}

\end{thm}

This exercise is the key part in the proof of the following theorem:
\textit{Any open surface with vanishing Gauss curvature is \index{cylindrical surface}\emph{cylindrical}}\,;
that is, the surface is swept out by a family of parallel lines.
Exercises \ref{ex:lin-ind-chart} and \ref{ex:line-cylinder} illustrate other parts of this proof.
The mentioned theorem was originally proved by Aleksey Pogorelov \cite[II §3 Thm 2]{pogorelov1956} in more general settings; it was rediscovered couple of times after that \cite{hartman-nirenberg,massey1962}.

\chapter{Curves in a surface}\label{chap:Curves in a surface}

\section{Darboux frame}\label{sec:Darboux}

\begin{wrapfigure}{r}{42 mm}
\vskip-10mm
\centering
\begin{lpic}[t(-0mm),b(0mm),r(0mm),l(0mm)]{asy/paraboloid+curve(1)}
\lbl[ul]{34,14;$\tan$}
\lbl[b]{20,43;$\Norm$}
\lbl[bl]{38,35;$\mu$}
\end{lpic}
\vskip-0mm
\end{wrapfigure}

Suppose $\gamma$ is a smooth curve in a smooth oriented surface~$\Sigma$.
As usual, we denote by $\Norm$ the unit normal field on~$\Sigma$.
Without loss of generality, we may assume that $\gamma$ is unit-speed;
in this case, $\tan(s)=\gamma'(s)$ is its tangent indicatrix.
Further we will use the shortcut notation $\Norm(s)\z=\Norm(\gamma(s))$.

The unit vectors $\tan(s)$ and $\Norm(s)$ are orthogonal.
Therefore there is a unique unit vector $\mu(s)$\index{10tmn@$\tan$, $\mu$, $\Norm$ (Darboux frame)} such that 
$\tan(s),\mu(s),\Norm(s)$ is an oriented orthonormal basis;
it is called the \index{Darboux frame}\emph{Darboux frame} of $\gamma$ in~$\Sigma$.

Since $\T_{\gamma (s)}\z\perp\Norm(s)$, the vector $\mu(s)$ is tangent to $\Sigma$ at $\gamma(s)$.
In fact, $\mu(s)$ is a counterclockwise rotation of $\tan(s)$ by the angle $\tfrac\pi2$ in the tangent plane $\T_{\gamma(s)}$.
This vector can also be defined as the vector product $\mu(s)\z\df\Norm(s)\times \tan(s)$.

Since $\gamma$ is unit-speed, we have that $\gamma''\perp \gamma'$ (see \ref{prop:a'-pertp-a''}).
Therefore, the acceleration of $\gamma$ can be written as a linear combination of $\mu$ and $\Norm$;
that is, \index{10k@$k_g$ (geodesic curvature)}\index{10k@$k_n$ (normal curvature)}
\[\gamma''(s)=k_g(s)\cdot \mu(s)+k_n(s)\cdot\Norm(s).\]
The values $k_g(s)$ and $k_n(s)$ are called \index{curvature}\index{geodesic!curvature}\emph{geodesic} and \index{normal!curvature}\emph{normal curvatures} of $\gamma$ at $s$, respectively.
Since the frame $\tan(s),\mu(s),\Norm(s)$ is orthonormal, these curvatures can also be written as the following scalar products:
\begin{align*}
k_g(s)&=\langle \gamma''(s),\mu(s)\rangle= 
&
k_n(s)&=\langle \gamma''(s),\Norm(s)\rangle=
\\
&=\langle \tan'(s),\mu(s)\rangle,
&
&=\langle \tan'(s),\Norm(s)\rangle.
\end{align*}

Since $0=\langle \tan(s),\Norm(s)\rangle$ we have 
that 
\begin{align*}
0&=\langle \tan(s),\Norm(s)\rangle'=
\\
&=\langle \tan'(s),\Norm(s)\rangle+\langle \tan(s),\Norm'(s)\rangle=
\\
&=k_n(s)+\langle \tan(s),D_{\tan(s)}\Norm\rangle.
\end{align*}
Applying the definition of the shape operator,
we get the following.

\begin{thm}{Proposition}\label{prop:normal-shape}
Assume $\gamma$ is a smooth unit-speed curve in a smooth surface~$\Sigma$.
Let $p=\gamma(s_0)$ and $\vec v=\gamma'(s_0)$.
Then 
\[k_n(s_0)=\langle \Shape_p(\vec v),\vec v\rangle,\]
where $k_n$ denotes the normal curvature of $\gamma$ at $s_0$, and $\Shape_p$ is the shape operator at~$p$.
\end{thm}

According to the proposition, the normal curvature of a smooth curve in $\Sigma$ is completely determined by the velocity vector $\vec v$.
For that reason, the normal curvature is also denoted by $k_{\vec v}$,\index{10k@$k_{\vec v}$ (normal curvature)}
and
\[k_{\vec v}=\langle \Shape_p(\vec v),\vec v\rangle\]
for any unit vector $\vec v$ in $\T_p$.

\section{Euler's formula}

Let $p$ be a point on a smooth surface~$\Sigma$.
Assume we choose tangent-normal coordinates at $p$ so that the Hessian matrix is diagonalized, then we have
\[M_p=\begin{pmatrix}
 k_1(p)
 &0
 \\
 0
 &k_2(p)
 \end{pmatrix}.
\]
Consider a vector ${\vec w}=a\cdot\vec i+b\cdot\vec j$ in the $(x,y)$-plane.
Then by \ref{cor:Shape(ij)}, we have
\[
\langle \Shape\vec w,\vec w\rangle
=a^2\cdot k_1(p) +b^2\cdot k_2(p). 
\]
If ${\vec w}$ is unit, then $a^2+b^2=1$ which implies the following.

\begin{thm}{Observation}\label{obs:k1-k2}
For any point $p$ on an oriented smooth surface $\Sigma$,
the principal curvatures $k_1(p)$ and $k_2(p)$ are respectively the minimum and maximum of the normal curvatures at~$p$.
Moreover, if $\theta$ is the angle between a unit vector ${\vec w}\in\T_p$ and the first principal direction at $p$, then 
\[k_{\vec w}(p)=k_1(p)\cdot(\cos\theta)^2+k_2(p)\cdot(\sin\theta)^2.\]

\end{thm}

The last identity is called \index{Euler's formula}\emph{Euler's formula}.

\begin{thm}{Exercise}\label{ex:mean-curvature}
Let $p$ be a point on a smooth surface.
Show that the sum of the normal curvatures for any pair of orthogonal tangent directions, at $p$ is $H(p)$ --- the mean curvature at~$p$. 
\end{thm}

\begin{thm}{Exercise}\label{ex:average}
Let $p$ be a point on a smooth surface.
Show that $\tfrac38\cdot H(p)^2-\tfrac12\cdot K(p)$ is the average value squared normal curvatures at $p$;
that is, the average value of $k_{\vec w}^2$ for all unit vectors ${\vec w}\in\T_p$.
\end{thm}

\begin{thm}{Meusnier's theorem}
\label{thm:meusnier}
\index{Meusnier's theorem}
Let $\gamma$ be a smooth curve on a smooth oriented surface~$\Sigma$.
Let $p=\gamma(t_0)$, ${\vec v}\z=\gamma'(t_0)$, and $\alpha\z=\measuredangle(\Norm(p),\norm(t_0))$;
that is, $\alpha$ is the angle between the normal vector to $\Sigma$ at $p$ and the normal vector in the Frenet frame of $\gamma$ at~$t_0$.
Then 
\[\kur(t_0)\cdot\cos\alpha=k_{n}(t_0);\]
here $\kur(t_0)$ and $k_n(t_0)$ are the curvature and the normal curvature of $\gamma$ at $t_0$, respectively. 
\end{thm}

\parbf{Proof.} Since $\gamma''=\tan'=\kur\cdot \norm$, we get that
\begin{align*}
k_{n}(t_0)&=\langle\gamma''(t_0),\Norm(p)\rangle=
\\
&=\kur(t_0)\cdot\langle\norm(t_0),\Norm(p)\rangle=
\\
&=\kur(t_0)\cdot\cos\alpha.
\end{align*}
\qedsf

\begin{thm}{Exercise}\label{ex:meusnier}
Let $\Sigma$ be a smooth surface, $p\in\Sigma$, and ${\vec v}\in \T_p\Sigma$ a unit vector.
Assume that the normal curvature $k_{\vec v}(p)$ does not vanish.

Show that the osculating circles at $p$ of smooth curves in $\Sigma$ that run in the direction ${\vec v}$ sweep out a sphere $S$ with center at $p+\tfrac1{k_{\vec v}}\cdot\Norm$ and radius $r=\tfrac1{|k_{\vec v}|}$.
\end{thm}

\begin{thm}{Exercise}\label{ex:principal-revolution}
Let $\gamma(s)=(x(s),y(s))$ be a smooth unit-speed simple plane curve in the upper half-plane,
and $\Sigma$ be the surface of revolution around the $x$-axis with generatrix~$\gamma$.
Assume that $x'\ne 0$.

\begin{subthm}{ex:principal-revolution:a}
Show that the parallels and meridians are lines of curvature on~$\Sigma$.
\end{subthm}

\begin{subthm}{ex:principal-revolution:formula}
Show that 
\[\frac{|x'(s)|}{y(s)}
\quad
\text{and}
\quad
\frac{-y''(s)}{|x'(s)|}
\]
are the principal curvatures of $\Sigma$ at $(x(s),y(s),0)$ in the direction of the corresponding parallel and meridian, respectively.
\end{subthm}

\begin{subthm}{ex:principal-revolution:pseudosphere}
Show that $\Sigma$ has Gauss curvature $-1$ at all points if and only if $y$ satisfies the differential equation $y''=y$. 
\end{subthm}

\end{thm}

{

\begin{wrapfigure}{r}{31 mm}
\vskip-0mm
\includegraphics{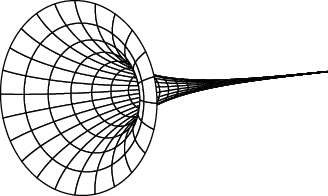}
\vskip-3mm
\end{wrapfigure}

The case \(y = e^{-s}\) is shown.
In this case, the curve \( \gamma \) is called a \index{tractrix}\emph{tractrix}, and the surface is called a \index{pseudosphere}\emph{pseudosphere};
by \ref{SHORT.ex:principal-revolution:pseudosphere}, its Gauss curvature is $-1$ at all points.

\begin{thm}{Exercise}\label{ex:catenoid-is-minimal}
Show that the \index{catenoid}\emph{catenoid} defined implicitly by the equation
\[(\cosh z)^2=x^2+y^2\]
is a minimal surface.
\end{thm}

}

\begin{thm}{Exercise}\label{ex:helicoid-is-minimal}
Show that the \index{helicoid}\emph{helicoid} defined by the following parametrization
\[s(u,v)=(u\cdot \sin v,u\cdot \cos v,v)\]
is a minimal surface.
\end{thm}

\begin{wrapfigure}{r}{51 mm}
\vskip-6mm
\centering
\includegraphics{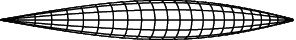}
\vskip0mm
\end{wrapfigure}

\begin{thm}{Exercise}\label{ex:rev(sin)}
Let $\Sigma$ be the surface of revolution around the $x$-axis
with generatrix $y=a\cdot \sin x$ for a constant $a>0$ and $x\in (0,\pi)$.
Show that the Gauss curvature of $\Sigma$ does not exceed 1.
\end{thm}

\begin{thm}{Exercise}\label{ex:rev(lin)}
Let $f\:(a,b)\to\mathbb{R}$ be a smooth positive function and $\Sigma$ be a surface of revolution of the graph $y=f(x)$ around $x$-axis.
Suppose $\Sigma$ has vanishing Gauss curvature.
Show that $f$ is a linear function; that is $f(x)=c\cdot x+d$ for some constants $c$ and $d$.
\end{thm}

\section{Lagunov's fishbowl}
\index{Lagunov's fishbowl}

The following question is a 3-dimensional analog of the moon in a puddle problem (\ref{thm:moon-orginal}).

\begin{thm}{Question}\label{quest:lagunov}
Assume a set $V\subset \mathbb{R}^3$ is bounded by a closed surface $\Sigma$ with 
principal curvatures bounded in absolute value by $1$.
Is it true that $V$ contains a ball of radius~$1$?
\end{thm}

\begin{thm}{Exercise}\label{ex:moon-revolution}
Show that Question \ref{quest:lagunov} has an affirmative answer if $V$ is a body of revolution.
\end{thm}

Later (see \ref{ex:convex-lagunov} and \ref{ex:rad=2})
we will show that the answer is ``yes'' for convex sets and surfaces in an open ball of radius $2$.
Now we are going to show by an example that the answer is ``no'' in the general case;
it was constructed by Vladimir Lagunov \cite{lagunov-1961}.

\parbf{Construction.}
Let us start with a body of revolution $V_1$ whose cross-section is shown on the next picture.
The boundary curve of the cross-section consists of $6$ long horizontal line segments included in $3$ simple closed smooth curves.
(To make the curves smooth, one has to use cutoffs and mollifiers from Appendix~\ref{sec:analysis}.)
The boundary of $V_1$ has $3$ components, each of which is a smooth sphere.

\begin{figure}[t!]
\centering
\includegraphics{mppics/pic-33}
\vskip0mm
\end{figure}

{

\begin{wrapfigure}{r}{21 mm}
\vskip-4mm
\centering
\includegraphics{mppics/pic-910}
\vskip0mm
\end{wrapfigure}

We assume that the curves have curvature at most~1.
Moreover, except for the almost horizontal parts, the curves have absolute curvature close to 1 all the time.
The only thick part in $V_1$ is the place where all three boundary components come close together;
the remaining part of $V_1$ is assumed to be very thin.
It could be arranged that the radius $r$ of the maximal ball in $V_1$ is just a bit above 
$r_2=\tfrac2{\sqrt{3}}-1$.
(The small black disc on the picture has radius $r_2$,
assuming that the three big circles are unit.)
In particular, we may assume that $r<\tfrac16$.

}

Exercise \ref{ex:principal-revolution} gives formulas for the principal curvatures of the boundary of $V_1$;
which imply that both principal curvatures are at most 1 in absolute value. 

It remains to modify $V_1$ to make its boundary connected without increasing the bounds on its principal curvatures and without allowing larger balls inside.

Each sphere in the boundary contains two flat discs;
they come into three pairs closely lying to each other. 
Let us drill thru two of such pairs and reconnect the holes by a body of revolution whose 
axis is shifted but stays parallel to the axis of $V_1$.
Denote the obtained body by $V_2$; its cross-section is shown on the picture. 

Let us repeat the operation for the other two pairs of discs.
Denote the obtained body by $V_3$.

Note that the boundary of $V_3$ is connected.
Assuming that the holes are large, its boundary can be made so that its principal curvatures are still at most $1$; the latter can be proved in the same way as for~$V_1$.
\qeds

The bound $r_2=\tfrac2{\sqrt{3}}-1$ is optimal;
this and more was proved by Vladimir Lagunov and Abram Fet \cite{lagunov-1960, lagunov-1961, lagunov-fet-1963, lagunov-fet-1965}.

The surface of $V_3$ in Lagunov's fishbowl has genus $2$;
that is, it can be parametrized by a sphere with two handles.

\begin{wrapfigure}{o}{55 mm}
\centering
\vskip-0mm
\includegraphics{mppics/pic-920}
\vskip0mm
\end{wrapfigure}

Indeed, the boundary of $V_1$ consists of three smooth spheres.

When we drill, we make four holes --- two spheres get one hole in each, and one sphere gets two holes.
We reconnect two spheres with a tube and obtain one sphere.
By connecting the two holes of the other sphere with a tube we get a torus;
it is on the left side in the picture of $V_2$.
That is, the boundary of $V_2$ is formed by one sphere and one torus.

To construct $V_3$ from $V_2$, we make a torus from the remaining sphere and connect it to the other torus by a tube.
This way we get a sphere with two handles.

\begin{thm}{Exercise}\label{ex:lagunov-genus4}
Modify Lagunov's construction to make the boundary surface a sphere with $4$ handles.

Try to make the boundary surface a sphere with $3$ handles.
\end{thm}

%Recall that Lagunov's fishbowl contains a ball of radius $r_2=\tfrac2{\sqrt{3}}-1$.
%It turns out that this radius is optimal.
%Moreover, \textit{suppose a connected body $V\subset \mathbb{R}^3$ is bounded by a finite number of closed smooth surfaces with principal curvatures bounded by~$1$.
%Assume $V$ does not contain a ball of radius $r_2$.
%Then the boundary of $V$ has two diffeomorphic connected components.}
%Such $V$ could be the region squeezed between two close surfaces; for example, a region between two concentric spheres.

%For bodies bounded by a sphere, there is a better bound $r_3=\sqrt{\tfrac32}-1$;
%it is the radius of the smallest sphere tangent to four unit spheres that are tangent to each other.
%Notice that $r_3>r_2$.
%\textit{If a body $V\subset \mathbb{R}^3$ is bounded by a smooth sphere with principal curvatures bounded  by~$1$.
%Then $V$ contains a ball of radius~$r_3$.
%Moreover, this bound is sharp; that is, for every $\epsilon>0$, there are examples of such $V$ not containing a ball of radius $r_3+\epsilon$.}

%% file: saddle.tex
\chapter{Supporting surfaces}
\label{chap:surface-support}

\section{Definitions}

Two smooth surfaces $\Sigma_1$ and $\Sigma_2$ are \index{tangent!surfaces}\emph{tangent} at a point $p \z\in \Sigma_1 \cap \Sigma_2$ if they share the tangent plane at $p$.
In this case their normals at $p$ either coincide or are opposite.
If the surfaces are oriented, we say they are \index{cooriented and counteroriented!surfaces}\emph{cooriented} if they have a common unit normal vector at $p$.
If their normal vectors at $p$ point in opposite directions, we say the surfaces are \index{}\emph{counteroriented}.

A surface $\Sigma_1$ \index{local support}\emph{locally supports} a surface $\Sigma_2$ at a point $p$ if
$p$ is a common point of the two surfaces and there is a neighborhood $U$ of $p$ such that
$\Sigma_2\cap U$ lies entirely on one side of $\Sigma_1$ in $U$.
According to the next exercise, such surfaces must be tangent at $p$; after choosing their orientations,
we may assume that they are cooriented at $p$.

\begin{thm}{Exercise}\label{ex:supp>tan}
Let a smooth surface $\Sigma_1$ locally support a smooth surface $\Sigma_2$ at a point~$p$.
Prove that $\Sigma_1$ is tangent to $\Sigma_2$ at $p$.
\end{thm}

Now, suppose $\Sigma_1$ and $\Sigma_2$ are cooriented at $p$.
Let us describe $\Sigma_1$ and $\Sigma_2$ locally at $p$ in tangent-normal coordinates as the graphs of functions $f_1$ and $f_2$ respectively.
$\Sigma_2$ locally supports $\Sigma_1$ at $p$ if and only if either
\[   f_1(x,y)\ge f_2(x,y)
\quad\text{(\,or}\quad
f_1(x,y)\le f_2(x,y)\,)\]
holds for all $(x,y)$ sufficiently close to the origin;
in this case, we say $\Sigma_2$ locally supports $\Sigma_1$ \index{supporting!surface}\emph{from the outside} (respectively, \emph{from the inside}).

\begin{thm}{Proposition}\label{prop:surf-support}
Let $\Sigma_1$ and $\Sigma_2$ be smooth oriented surfaces.
Assume $\Sigma_1$ locally supports $\Sigma_2$ from the inside at the point $p$ (equivalently $\Sigma_2$ locally supports $\Sigma_1$ from the outside).
Then 
\[k_1(p)_{\Sigma_1}\ge k_1(p)_{\Sigma_2}\quad\text{and}\quad k_2(p)_{\Sigma_1}\z\ge k_2(p)_{\Sigma_2}.\]
\end{thm}

\begin{thm}{Exercise}\label{ex:surf-support}
Construct two surfaces $\Sigma_1$ and $\Sigma_2$ that have a  common point $p$ with a common unit normal vector that
do not support each other, but
$k_1(p)_{\Sigma_1}\z> k_1(p)_{\Sigma_2}$ and $k_2(p)_{\Sigma_1}\z> k_2(p)_{\Sigma_2}$.
\end{thm}

\parbf{Proof of \ref{prop:surf-support}.}
We can assume that $\Sigma_1$ and $\Sigma_2$ are graphs $z=f_1(x,y)$  and $z=f_2(x,y)$ in common tangent-normal coordinates at $p$, so we have $f_1\ge f_2$.

Choose a common tangent unit vector ${\vec w}\in \T_p\Sigma_1=\T_p\Sigma_2$.
Consider the plane $\Pi$ passing thru $p$ and spanned by the normal vector $\Norm_p$ and ${\vec w}$.
Let $\gamma_1$ and $\gamma_2$ be the curves of intersection of $\Sigma_1$ and $\Sigma_2$ with $\Pi$.

Let us orient $\Pi$ so that the common normal vector $\Norm_p$ for both surfaces at $p$ points to the left from ${\vec w}$.
Further, let us parametrize both curves so that they are running in the direction of ${\vec w}$ at $p$ and are therefore cooriented.
In this case, the curve $\gamma_1$ supports the curve $\gamma_2$ from the left.

\begin{wrapfigure}{o}{40 mm}
\vskip-4mm
\centering
\includegraphics{mppics/pic-80}
\vskip-0mm
\end{wrapfigure}

By \ref{prop:supporting-circline}, we have the following inequality for the normal curvatures of $\Sigma_1$ and $\Sigma_2$ at $p$ in the direction of ${\vec w}$:
\[k_{\vec w}(p)_{\Sigma_1}\ge k_{\vec w}(p)_{\Sigma_2}.\eqlbl{kw>=kw}\]

\setlength{\columnseprule}{0.4pt}
\begin{multicols}{2}
According to \ref{obs:k1-k2},
\[k_1(p)_{\Sigma_i}=\min_{\vec w}\set{k_{\vec w}(p)_{\Sigma_i}}{},\]
where $\vec w\in\T_p$ and $|\vec w|=1$.
Fix ${\vec w}$ such that 
\[k_1(p)_{\Sigma_1}\z=k_{\vec w}(p)_{\Sigma_1}.\]
By \ref{kw>=kw}, 
\begin{align*}
k_1(p)_{\Sigma_1}&=k_{\vec w}(p)_{\Sigma_1}\ge
\\
&\ge k_{\vec w}(p)_{\Sigma_2}\ge
\\
&\ge\min_{\vec v}\set{k_{\vec v}(p)_{\Sigma_2}}{}=
\\
&=k_1(p)_{\Sigma_2}.
\end{align*}

\columnbreak

Similarly, by \ref{obs:k1-k2},
\[k_2(p)_{\Sigma_i}=\max_{\vec w}\set{k_{\vec w}(p)_{\Sigma_i}}{},\]
where $\vec w\in\T_p$ and $|\vec w|=1$.
Fix ${\vec w}$ such that 
\[k_2(p)_{\Sigma_2}=k_{\vec w}(p)_{\Sigma_2}.\]
By \ref{kw>=kw},
\begin{align*}
k_2(p)_{\Sigma_2}&=k_{\vec w}(p)_{\Sigma_2}\le
\\
&\le k_{\vec w}(p)_{\Sigma_1}\le
\\
&\le\max_{\vec v}\set{k_{\vec v}(p)_{\Sigma_1}}{}=
\\
&=k_2(p)_{\Sigma_1}.
\end{align*}
\end{multicols}

Both times we assume  that  $\vec v\in\T_p$ and $|\vec v|=1$.\qeds

\begin{thm}{Corollary}\label{cor:surf-support}
Let $\Sigma_1$ and $\Sigma_2$ be oriented smooth surfaces.
Assume $\Sigma_1$ locally supports $\Sigma_2$ from the inside at the point~$p$.
Then:

\begin{subthm}{cor:surf-support:mean}$H(p)_{\Sigma_1}\ge H(p)_{\Sigma_2}$;
\end{subthm}

\begin{subthm}{cor:surf-support:gauss} If $k_1(p)_{\Sigma_2}\ge 0$, then $K(p)_{\Sigma_1}\ge K(p)_{\Sigma_2}$.
\end{subthm}
 
\end{thm}

\parbf{Proof;} \ref{SHORT.cor:surf-support:mean}.
The statement follows from  \ref{prop:surf-support} since $H=k_1+k_2$.
%\[H(p)_{\Sigma_i}=k_1(p)_{\Sigma_i}+k_2(p)_{\Sigma_i}.\]

\parit{\ref{SHORT.cor:surf-support:gauss}.} Since $k_2(p)_{\Sigma_i}\ge k_1(p)_{\Sigma_i}$, and $k_1(p)_{\Sigma_2}\ge 0$, we get that all the principal curvatures 
$k_1(p)_{\Sigma_1}$, 
$k_1(p)_{\Sigma_2}$, 
$k_2(p)_{\Sigma_1}$, and 
$k_2(p)_{\Sigma_2}$ are nonnegative.
By \ref{prop:surf-support}, it implies that
\begin{align*}
K(p)_{\Sigma_1}&=k_1(p)_{\Sigma_1}\cdot k_2(p)_{\Sigma_1}\ge 
\\
&\ge k_1(p)_{\Sigma_2}\cdot k_2(p)_{\Sigma_2}=K(p)_{\Sigma_2}.
\end{align*}
\qedsf

\begin{thm}{Exercise}\label{ex:positive-gauss-0}
Show that any closed surface in a unit ball has a point with Gauss curvature at least~1.
Conclude that any closed surface has a point with strictly positive Gauss curvature.
\end{thm}

\begin{thm}{Exercise}\label{ex:positive-gauss}
Show that any closed surface that lies at distance at most 1 from a straight line has a point with Gauss curvature at least~1.
\end{thm}

\section{Convex surfaces}

A surface that bounds a convex region is called \index{convex!surface}\emph{convex}.

\begin{thm}{Exercise}\label{ex:convex-surf}
Show that the Gauss curvature of any convex smooth surface is nonnegative at each point.
\end{thm}

\begin{thm}{Advanced exercise}\label{ex:convex-lagunov}
Assume $R$ is a convex body in $\mathbb{R}^3$ bounded by a surface with principal curvatures at most~1.
Show that $R$ contains a unit ball.
\end{thm}

Recall that a region $R$ in the Euclidean space is called  {}\emph{strictly convex} if, for any two points $x,y\in R$, any point $z$ between $x$ and $y$ lies in the interior of~$R$.

For example, any open convex set is strictly convex.
The cube (as well as any convex polyhedron) gives an example of a non-strictly convex set.
Evidently, \textit{a closed convex region is strictly convex if and only if its boundary does not contain a line segment}.

\begin{thm}{Lemma}\label{lem:gauss+=>convexity}
Let $z=f(x,y)$ be the local description of a smooth surface $\Sigma$ in tangent-normal coordinates at a point $p\in\Sigma$.
Assume both principal curvatures of $\Sigma$ are positive at~$p$.
Then the function $f$ is strictly convex in a neighborhood of the origin and has a local minimum at the origin.

In particular, the tangent plane $\T_p$ locally supports $\Sigma$ from the outside at~$p$.
\end{thm}

\parbf{Proof.}
Since both principal curvatures are positive, by \ref{cor:Shape(ij)}, we have 
\[D^2_{\vec w}f(0,0)=\langle \Shape_p({\vec w}),{\vec w}\rangle\ge k_1(p)>0\] 
for any unit tangent vector ${\vec w}\in\T_p\Sigma$ (in our coordinates, ${\vec w}$ is horizontal).

By continuity of the function $(x,y,{\vec w})\mapsto D^2_{\vec w}f(x,y)$,
we have that $D^2_{\vec w}f(x,y)>0$ if $\vec w\ne 0$ and $(x,y)$ lies in a sufficiently small neighborhood of the origin.
This property implies that $f$ is a strictly convex function in a neighborhood of the origin in the $(x,y)$-plane (see \ref{thm:Jensen}).

Finally, since $\nabla f(0,0)=0$, and $f$ is strictly convex in a neighborhood of the origin, $f$ has a strict local minimum at the origin.
\qeds

\begin{thm}{Exercise}\label{ex:section-of-convex}
Let $\Sigma$ be a smooth surface (without boundary) with positive Gauss curvature.
Show that any connected component of the intersection of $\Sigma$ with a plane is either a single point or a smooth plane curve whose signed curvature has constant sign.
\end{thm}

\section{Convexity test}

\begin{thm}{Theorem}\label{thm:convex-embedded}
Suppose $\Sigma$ is an open or closed smooth surface with positive Gauss curvature.
Then $\Sigma$ bounds a strictly convex region.
\end{thm}

In the proof we will have to use that the surface is connected (by the definition);
otherwise, a pair of spheres would give a counterexample.
It is true that \textit{any closed smooth surface with nonnegative Gauss curvature bounds a convex region}, but the proof requires more work \cite{hadamard,gomes,sacksteder}.

The same holds for surfaces with possible self-intersections.
This was proved by James Stoker \cite{stoker}, but he attributed the result to Jacques Hadamard who proved a closely relevant statement \cite[§ 23]{hadamard}.

From the next proof, one can also extract the following statement: \textit{any closed connected locally convex region in the Euclidean space is convex}.

\parbf{Proof.}
Since the Gauss curvature is positive, we can assume that the principal curvatures are positive at any point.
Denote by $R$ the region bounded by $\Sigma$ that lies on the side of $\Norm$;
that is, $\Norm$ points inside $R$ at any point of~$\Sigma$.
(The region $R$ exists by \ref{clm:proper-divides}.)

Let us show that $R$ is {}\emph{locally strictly convex};
that is, for any point $p\in \Sigma$, the intersection of $R$ with a small ball centered at $p$ is strictly convex.

Indeed, suppose $z=f(x,y)$ is a local description of $\Sigma$ in the tangent-normal coordinates at~$p$.
By \ref{lem:gauss+=>convexity}, $f$ is strictly convex in a neighborhood of the origin.
In particular, the intersection of a small ball centered at $p$ with the epigraph $z\ge f(x,y)$ is strictly convex.

Since $\Sigma$ is connected, so is $R$;
moreover, any two points in the interior of $R$ can be connected by a polygonal line in the interior of~$R$.

Assume the interior of $R$ is not convex;
that is, there are points $x,y\in R$, and a point $z$ between $x$ and $y$ that does not lie in the interior of~$R$.
Consider a polygonal  line $\beta$ from $x$ to $y$ in the interior of~$R$.
Let $y_0$ be the first point on $\beta$ such that the chord $[x,y_0]$ touches $\Sigma$, say at~$z_0$.

\begin{wrapfigure}{o}{43 mm}
\vskip-0mm
\centering
\includegraphics{mppics/pic-37}
\vskip-0mm
\end{wrapfigure}

Since $R$ is locally strictly convex, $R\z\cap B(z_0,\epsilon)$ is strictly convex for all sufficiently small $\epsilon>0$.
On the other hand, $z_0$ lies between two points in the intersection $[x,y_0]\cap B(z_0,\epsilon)$.
Since $[x,y_0]\subset R$, we arrived to a contradiction.

Therefore, the interior of $R$ is a convex set.
The region $R$ is the closure of its interior, therefore $R$ is convex as well.

Since $R$ is locally strictly convex, its boundary $\Sigma$ contains no line segments.
Therefore, $R$ is strictly convex.
\qeds

\begin{thm}{Exercise}\label{ex:surrounds-disc}
Assume a closed surface $\Sigma$ surrounds a unit circle.
Show that there is a point  $p \in \Sigma$ with $K(p)\le 1$.
\end{thm} 

\begin{thm}{Exercise}\label{ex:small-gauss}
Let $\Sigma$ be a closed convex smooth surface of diameter at least $\pi$;
that is, there is a pair of points $p,q\in\Sigma$ such that $|p-q|\ge \pi$.
Show that $\Sigma$ has a point with Gauss curvature at most~1.
\end{thm}

\section{More convexity}

\begin{thm}{Theorem}\label{thm:convex-closed}
Suppose $\Sigma$ is a closed smooth convex surface.
Then it is a smooth sphere; that is, $\Sigma$ admits a smooth regular parametrization by $\mathbb{S}^2$.\end{thm}

{

\begin{wrapfigure}[8]{r}{33 mm}
\vskip-6mm
\centering
\includegraphics{mppics/pic-78}
\end{wrapfigure}

\begin{thm}{Proof-guided exercise}\label{ex:convex-proper-sphere}
Assume a convex compact region $R$ contains the origin in its interior and is bounded by a smooth surface~$\Sigma$.

\begin{subthm}{ex:convex-proper-sphere:single}
Show that any half-line that starts at the origin intersects $\Sigma$ at a single point;
that is, there is a positive function $\rho\:\mathbb{S}^2\z\to\mathbb{R}$ such that $\Sigma$ consists of the points $q\z=\rho(\xi)\cdot \xi$ for $\xi\in \mathbb{S}^2$.
\end{subthm}

\begin{subthm}{ex:convex-proper-sphere:smooth}
Show that $\rho\:\mathbb{S}^2\to\mathbb{R}$ is a smooth function.
Conclude that $\xi\z\mapsto \rho(\xi)\cdot \xi$ is a smooth regular parametrization $\mathbb{S}^2\z\to \Sigma$.
\end{subthm}

\end{thm}

}

\begin{thm}{Theorem}\label{thm:convex-open}
Suppose $\Sigma$ is an open smooth surface that bounds a strictly convex closed region.
Then there is a coordinate system in which $\Sigma$ is the graph $z\z=f(x,y)$ of a convex function $f$ defined on a convex open region $\Omega$ of the $(x,y)$-plane.
Moreover, $f(x_n,y_n)\to\infty$ as $(x_n,y_n)\z\to(x_\infty,y_\infty)\z\in \partial\Omega$.

\end{thm}

{

\begin{wrapfigure}{r}{33 mm}
\vskip-0mm
\centering
\includegraphics{mppics/pic-1181}
\vskip2mm
\end{wrapfigure}

\begin{thm}{Proof-guided exercise}\label{ex:convex-proper-plane}
Assume a strictly convex closed noncompact region $R$ contains the origin in its interior and is bounded by a smooth surface~$\Sigma$.

\begin{subthm}{ex:convex-proper-plane:a}
Show that $R$ contains a half-line $\ell$.
\end{subthm}

\begin{subthm}{ex:convex-proper-plane:b}
Show that any line $m$ parallel to $\ell$ intersects $\Sigma$ at most at one point.
\end{subthm}

\begin{subthm}{ex:convex-proper-plane:c}
Consider an $(x,y,z)$-coordinate system such that the $z$-axis points in the direction of $\ell$.
Show that the projection of $\Sigma$ to the $(x,y)$ plane is an open convex set; denote it by $\Omega$.
\end{subthm}

\begin{subthm}{ex:convex-proper-plane:d}
Conclude that $\Sigma$ is a graph $z=f(x,y)$ of a convex function $f$ defined on $\Omega$.
\end{subthm}

\begin{subthm}{ex:convex-proper-plane:e}
Prove the last statement in \ref{thm:convex-open}.
\end{subthm}

\end{thm}

}

\begin{thm}{Exercise}\label{ex:open+convex=plane}
Show that any open surface $\Sigma$ with positive Gauss curvature is a topological plane;
that is, there is an embedding $\mathbb{R}^2\hookrightarrow\mathbb{R}^3$ with image~$\Sigma$.

Try to show that $\Sigma$ is a smooth plane;
that is, the embedding $f$ can be made smooth and regular.
\end{thm}

\begin{thm}{Exercise}\label{ex:circular-cone}
Show that any open smooth surface $\Sigma$ with positive Gauss curvature
lies inside an infinite circular cone.
In other words, there is an $(x,y,z)$-coordinate system in which $\Sigma$ lies in the region $z \z\ge m \cdot\sqrt{x^2 + y^2}$ for a positive constant $m$.
\end{thm} 

\begin{thm}{Exercise}\label{ex:intK}
Assume $\Sigma$ is a smooth convex surface with positive Gauss curvature
and $\Norm\:\Sigma\to \mathbb{S}^2$ is its spherical map.

\begin{multicols}{2}

\begin{subthm}{ex:intK:4pi}
Show that if $\Sigma$ is closed, then $\Norm\:\Sigma\to \mathbb{S}^2$ is a bijection.
Conclude that 
\[\iint_\Sigma K=4\cdot\pi.\]
\end{subthm}

\columnbreak

\begin{subthm}{ex:intK:2pi}
Show that if $\Sigma$ is open, then $\Norm\:\Sigma\to \mathbb{S}^2$
maps $\Sigma$ bijectively into a subset of a hemisphere.
Conclude that 
\[\iint_\Sigma K\le 2\cdot\pi.\]
\end{subthm}

\end{multicols}

\end{thm}

\chapter{Saddle surfaces}

\section{Definitions}\label{sec:saddle}

A surface is called \index{saddle surface}\emph{saddle} if its Gauss curvature at each point is nonpositive;
in other words, the principal curvatures at each point have opposite signs or at least one of them is zero.

\begin{wrapfigure}{r}{43 mm}
\vskip-0mm
\centering
\includegraphics{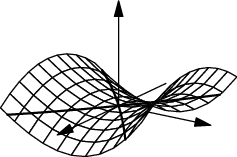}
\vskip0mm
\end{wrapfigure}

If the Gauss curvature is negative at each point,
then the surface is said to be {}\emph{strictly saddle};
equivalently, this means that the principal curvatures have opposite signs at each point.
In this case, the tangent plane cannot support the surface even locally --- moving along the surface in the principal directions at a given point, one goes above and below the tangent plane.

\begin{thm}{Exercise}\label{ex:convex-revolution}
Let $f\:\mathbb{R}\to\mathbb{R}$ be a smooth positive function.
Show that the surface of revolution of the graph $y\z=f(x)$ around the $x$-axis
 is saddle if and only if $f$ is convex; that is, if $f''(x)\ge0$ for any~$x$.
\end{thm}

A surface $\Sigma$ is called \index{ruled surface}\emph{ruled} if, for every point $p\in \Sigma$, there is a line $\ell_p\subset \Sigma$ passing thru $p$.

\begin{thm}{Exercise}\label{ex:ruled=>saddle}
Show that any ruled surface is saddle.
\end{thm}

\begin{thm}{Exercise}\label{ex:saddle-convex}
Let $\Sigma$ be an open strictly saddle surface and $f\:\mathbb{R}^3\z\to\mathbb{R}$ be a smooth convex function.
Show that the restriction of $f$ to $\Sigma$ does not have a point of local maximum.
\end{thm}

A tangent direction on a smooth surface with vanishing normal curvature is called \index{asymptotic direction and line}\emph{asymptotic}.
A smooth curve that always runs in an asymptotic direction is called an
{}\emph{asymptotic line}.\label{page:asymptotic line}

Recall that a set $R$ in the plane is called \index{star-shaped}\emph{star-shaped} if there is a point $p\in R$ such that for any $x\in R$ the line segment $[p,x]$ belongs to~$R$.

\begin{thm}{Advanced exercise}\label{ex:panov}
Let $\gamma$ be a closed smooth asymptotic line
in the graph $z\z=f(x,y)$ of a smooth function~$f$.
Assume the graph is strictly saddle in a neighborhood of~$\gamma$.
Show that the region in the $(x,y)$-plane bounded by the projection $\bar \gamma$ of $\gamma$ cannot be star-shaped.
\end{thm}

\begin{thm}{Advanced exercise}\label{ex:crosss}
Let $\Sigma$ be a smooth surface and $p\in \Sigma$.
Assume $K(p)<0$.
Show that there is a neighborhood $\Omega$ of $p$ in $\Sigma$
such that the intersection of $\Omega$ with the tangent plane $\T_p$ is a union of two smooth curves  \index{transversality}\emph{intersecting transversally} at~$p$;
that is, their velocity vectors at $p$ are linearly independent.
\end{thm}

\section{Hats}

A \textit{closed surface cannot be saddle}.
Indeed, suppose $\Sigma$ is a closed surface.
Consider the smallest sphere that surrounds~$\Sigma$.
The sphere supports $\Sigma$ at some point, and at this point, the principal curvatures must have the same sign.
The following more general statement is proved using the same idea.

{

\begin{wrapfigure}[6]{r}{45 mm}
\vskip-6mm
\centering
\includegraphics{mppics/pic-73}
\vskip0mm
\end{wrapfigure}

\begin{thm}{Lemma}\label{lem:convex-saddle}
Assume $\Sigma$ is a compact saddle surface, and its boundary line lies in a closed convex  region~$R$.
Then the entire surface $\Sigma$ lies in~$R$.
\end{thm}

\parit{Remark.}
In the case of strictly saddle surface, the lemma can be deduced from \ref{ex:saddle-convex}.

\parbf{Proof.}
Arguing by contradiction,
assume there is a point $p\in \Sigma$ that does not lie in~$R$.
Let $\Pi$ be a plane that separates $p$ from $R$; it exists by \ref{lem:separation}.
Denote by $\Sigma'$ the part of $\Sigma$ that lies with $p$ on the same side of~$\Pi$.

}

Since $\Sigma$ is compact, it is surrounded by a sphere;
let $\sigma$ be the circle of intersection of this sphere and $\Pi$.
Consider the smallest spherical dome $\Sigma_0$ with boundary $\sigma$ that surrounds~$\Sigma'$.

$\Sigma_0$ supports $\Sigma$ at some point~$q$.
Without loss of generality, we may assume that $\Sigma_0$ and $\Sigma$ are cooriented at $q$, and $\Sigma_0$ has positive principal curvatures.
In this case, $\Sigma_0$ supports $\Sigma$ from the outside.
By \ref{cor:surf-support}, we have $K(q)_\Sigma\z\ge K(q)_{\Sigma_0}>0$, a contradiction.
\qeds

\begin{thm}{Exercise}\label{ex:proper-saddle}
Construct a proper saddle surface that does not lie in the convex hull of its boundary line.
(By \ref{lem:convex-saddle} it cannot be compact.)
\end{thm}

\begin{thm}{Exercise}\label{ex:length-of-bry}
Let $\Delta$ be a compact smooth saddle surface with boundary and $p\in \Delta$.
Assume the boundary curve of $\Delta$ lies in the unit sphere centered at~$p$.
Show that if $\Delta$ is a topological disc, then $\length(\partial\Delta)\ge 2\cdot\pi$.

Show that the statement does not hold without assuming that $\Delta$ is a topological disc.
\end{thm}

\parit{Remark.}
In fact, $\area \Delta\ge \pi$;
that is, the unit plane disc has the minimal possible area.
This follows from the exercise by applying the so-called \index{coarea formula}\emph{coarea formula}.

\begin{thm}{Exercise}\label{ex:circular-cone-saddle}
Show that an open saddle surface
cannot lie inside an infinite circular cone.
\end{thm}

A topological disc $\Delta$ in a surface $\Sigma$ is called a \index{hat}\emph{hat} of $\Sigma$
if its boundary line $\partial\Delta$ lies in a plane $\Pi$ and $\Delta \setminus \partial \Delta$ lies strictly on one side of $\Pi$.

\begin{thm}{Proposition}\label{prop:hat}
A smooth surface $\Sigma$ is saddle if and only if it has no hats.
\end{thm}

A saddle surface can contain a closed plane curve.
For example, the hyperboloid $x^2+y^2-z^2=1$ contains the unit circle in the $(x,y)$-plane.
However, according to the proposition (as well as the lemma), a plane curve cannot bound a topological disc (as well as any compact set) in a saddle surface.

\parbf{Proof.}
Since a plane is convex, the only-if part follows from \ref{lem:convex-saddle};
it remains to prove the if part.

Assume $\Sigma$ is not saddle; that is, it has a point $p$ with strictly positive Gauss curvature;
or equivalently, the principal curvatures $k_1(p)$ and $k_2(p)$ have the same sign.
From now on we assume $0 < k_1 (p) \leq k_2 (p) $, the other case being analogous.

Let $z=f(x,y)$ be a graph representation of $\Sigma$ in tangent-normal coordinates at~$p$.
Consider the set $F_\epsilon$ in the $(x,y)$-plane defined by the inequality $f(x,y)\le \epsilon$.
By \ref{lem:gauss+=>convexity}, $f$ is convex in a small neighborhood of $(0,0)$.
Therefore, $F_\epsilon$ is convex, for sufficiently small $\epsilon>0$.
In particular, $F_\epsilon$ is a topological disc.

The map $(x,y)\mapsto (x,y,f(x,y))$ defines a homeomorphism from $F_\epsilon$
to
\[\Delta_\epsilon=\set{(x,y,f(x,y))\in \mathbb{R}^3}{f(x,y)\le \epsilon};\]
so $\Delta_\epsilon$ is a topological disc for any small $\epsilon>0$.
The boundary line of $\Delta_\epsilon$ lies on the plane $z=\epsilon$, and the entire disc lies below it;
that is, $\Delta_\epsilon$ is a hat of~$\Sigma$.
\qeds

The following exercise shows that $\Delta_\epsilon$ is in fact a smooth disc.
This can be used to prove a slightly stronger version of \ref{prop:hat};
namely in the definition of hats one can assume that the disc is smooth.

\begin{thm}{Exercise}\label{ex:disc-hat}
Let $f\:\mathbb{R}^2\to\mathbb{R}$ be a smooth strictly convex function with a minimum at the origin.
Show that the set $F_\epsilon$ in the graph $z=f(x,y)$ defined by the inequality $f(x,y)\le \epsilon$ is a smooth disc for any $\epsilon>0$;
that is, there is a diffeomorphism from $F_\epsilon$ to the unit disc $\Delta\z=\set{(x,y)\in\mathbb{R}^2}{x^2+y^2\le 1}$.
\end{thm}

\begin{thm}{Exercise}\label{ex:saddle-linear}
Let $L\:\mathbb{R}^3\to\mathbb{R}^3$ be an \index{affine transformation}\emph{affine transformation}; that is, a bijection $\mathbb{R}^3\to\mathbb{R}^3$ that sends any plane to a plane.
Show that for any saddle surface $\Sigma$ the image $L(\Sigma)$ is also a saddle surface.
\end{thm}

%???Monkey saddle

\section{Saddle graphs}

The following theorem was proved by Sergei Bernstein \cite{bernstein}.
\index{Bernstein's theorem}

\begin{thm}{Theorem}\label{thm:bernshtein}
Let $f\:\mathbb{R}^2\to\mathbb{R}$ be a smooth function.
Assume its graph $z=f(x,y)$ is a strictly saddle surface in $\mathbb{R}^3$.
Then $f$ is not bounded;
that is, there is no constant $C$ such that 
$|f(x,y)|\le C$ for any $(x,y)\in\mathbb{R}^2$.
\end{thm}

The theorem states that a saddle graph cannot lie between two horizontal planes;
combining this with \ref{ex:saddle-linear}, we get that saddle graphs cannot lie between two parallel planes.
The following exercise shows that the theorem does not hold for saddle surfaces that are not graphs.

\begin{thm}{Exercise}\label{ex:between-parallels}
Construct an open strictly saddle surface that lies between two parallel planes.
\end{thm}

Since $\exp(x-y^2)>0$,
the following exercise shows that there are strictly saddle graphs with functions bounded on one side;
that is, both (upper and lower) bounds are needed in the proof of the theorem.

\begin{thm}{Exercise}\label{ex:one-side-bernshtein}
Show that the graph
$z=\exp(x-y^2)$
is strictly saddle.
\end{thm}

The following exercise leads to a stronger version of Bernstein's theorem.

\begin{thm}{Advanced exercise}\label{ex:saddle-graph}
Let $\Sigma$ be an open smooth strictly saddle surface in $\mathbb{R}^3$.
Assume there is a compact subset $K\subset \Sigma$ such that the complement $\Sigma\setminus K$ is the graph $z=f(x,y)$ of a smooth function defined in an open domain of the $(x,y)$-plane.
Show that the surface $\Sigma$ is a graph.
\end{thm}

The following lemma gives an analogue of \ref{ex:length-of-bry}, and will be used in the proof of Theorem \ref{thm:bernshtein}.

\begin{thm}{Lemma}\label{lem:region}
There is no proper strictly saddle smooth surface that has its boundary line in a plane $\Pi$ and lies at a bounded distance from a line contained in $\Pi$.
\end{thm}

\parbf{Proof.}
By \ref{ex:saddle-linear}, the statement can be reformulated in the following way:
\textit{There is no proper strictly saddle smooth surface 
with its boundary line in the $(x,y)$-plane
and contained in a region of the form:}
\[R=\set{(x,y,z)\in\mathbb{R}^3}{0\le z\le r,0\le y\le r}.\]

Assume the contrary, let $\Sigma$ be such a surface.
Consider the projection $\hat \Sigma$ of $\Sigma$ to the $(x,z)$-plane.
It lies in the upper half-plane and below the line $z=r$.

Consider the open upper half-plane $H=\set{(x,z)\in \mathbb{R}^2}{z> 0}$.
Let $\Theta$ be the connected component of the complement $H\setminus \hat \Sigma$ that contains all the points above the line $z=r$.

The set $\Theta$ is convex.
If not, then there is a line segment $[p,q]$ for some $p,q\in \Theta$ that cuts from $\hat\Sigma$ a compact piece.
\begin{figure}[!ht]
\vskip-0mm
\centering
\includegraphics{mppics/pic-74}
\vskip0mm
\end{figure}
Consider the plane $\Pi$ thru $[p,q]$ that is perpendicular to the $(x,z)$-plane.
$\Pi$ cuts from $\Sigma$ a compact region $\Delta$.
By a general position argument (see \ref{lem:reg-section}),
we can assume that $\Delta$ is a compact surface with its boundary line in $\Pi$, and the remaining part of $\Delta$ lies strictly on one side from $\Pi$.
Since the plane $\Pi$ is convex, this statement contradicts \ref{lem:convex-saddle}.

Summarizing, $\Theta$ is an open convex set of $H$ that contains all points above $z=r$.
By convexity, together with any point $w$, the set $\Theta$ contains all points on the half-lines that start at $w$ and \textit{point up}; that is, in directions with positive $z$-coordinate.
In other words, with any point $w$,
the set $\Theta$ contains all points with larger $z$-coordinates.
\begin{figure}[!ht]
\vskip-0mm
\centering
\includegraphics{mppics/pic-75}
\vskip0mm
\end{figure}
Since $\Theta$ is open it can be described by an inequality $z>r_0$.
It follows that the plane $z=r_0$ supports $\Sigma$ at some point (in fact at many points).
By \ref{prop:surf-support}, the latter is impossible --- a contradiction.
\qeds

%Notice that the same proof goes through if one replaces the condition $0 \leq y \leq r$ in \ref{eq:R-bernstein} with $0 \leq y \leq h(x)$ for some positive continuous function $h$.

\parbf{Proof of \ref{thm:bernshtein}.}
Denote by $\Sigma$ the graph $z=f(x,y)$.
Assume the contrary; that is, $\Sigma$ lies between two planes $z=\pm C$.

The function $f$ cannot be constant.
It follows that the tangent plane $\T_p$ at some point $p\in\Sigma$ is not horizontal.

Denote by $\Sigma^+$ the part of $\Sigma$ that lies above $\T_p$.
Note that $\Sigma^+$ has at least two connected components which are approaching $p$ from both sides 
in the principal direction with positive principal curvature.
Otherwise, there is a curve that runs in $\Sigma^+$ and approaches $p$ from both sides.
This curve cuts a topological disc, say $\Delta$, from $\Sigma$ with the boundary line above or on $\T_p$ and some points strictly below $\T_p$;
the latter contradicts \ref{lem:convex-saddle}.
Summarizing, $\Sigma^+$ has at least two connected components.

\begin{figure}[!ht]
\vskip-0mm
\centering
\includegraphics{mppics/pic-76}
\caption*{The surface $\Sigma$ seen from above.}
\vskip0mm
\end{figure}

Let $z=h(x,y)=a\cdot x+b\cdot y+c$ be the equation of $\T_p$.
Notice that $\Sigma^{+} \z= \set{(x,y,f(x,y))\in \Sigma}{h(x,y) \leq f(x,y)}$, hence it contains the connected set
\[R_-=\set{(x,y,f(x,y))\in\Sigma}{h(x,y)< -C}\] 
and is disjoint from  
\[R_+=\set{(x,y,f(x,y))\in\Sigma}{h(x,y)> C}.\]
Whence one of the connected components, say $\Sigma^+_0$, lies in the strip
\[R_0=\set{(x,y,f(x,y))\in\Sigma}{|h(x,y)|\le  C}.\]
This set lies at a bounded distance from the line of intersection of $\T_p$ with the $(x,y)$-plane.

Let us move $\T_p$ slightly upward and cut from $\Sigma^+_0$ the piece above the obtained plane, say $\bar\Sigma^+_0$.
By the general position argument (\ref{lem:reg-section}),
we can assume that $\bar\Sigma^+_0$ is a surface with a smooth boundary line;
by construction the boundary line lies in the plane.
The obtained surface $\bar\Sigma^+_0$ still lies at a bounded distance to a line.
The latter is impossible by \ref{lem:region}.
\qeds

\section{Remarks}

Bernstein's theorem and the lemma in its proof do not hold for nonstrictly saddle surfaces;
counterexamples can be found among cylindrical surfaces over smooth curves; see the end of \ref{sec:shape}.
In fact, it can be shown that these are the only counterexamples;
a proof is based on the same idea, but it is more technical.

By \ref{prop:hat}, saddle surfaces can be defined as smooth surfaces without hats.
This definition can be used for arbitrary surfaces, not necessarily smooth.
Some results, for example, the characterization of saddle graphs, can be extended to these generalized saddle surfaces.
However, this class of surfaces is far from being understood; see \cite[Chapter 4]{alexander-kapovitch-petrunin2019} and the references therein.

%% file: part-geodesics.tex
\arxiv{\cleardoublepage
\phantomsection
\AddToShipoutPictureBG*{\includegraphics[width=\paperwidth]{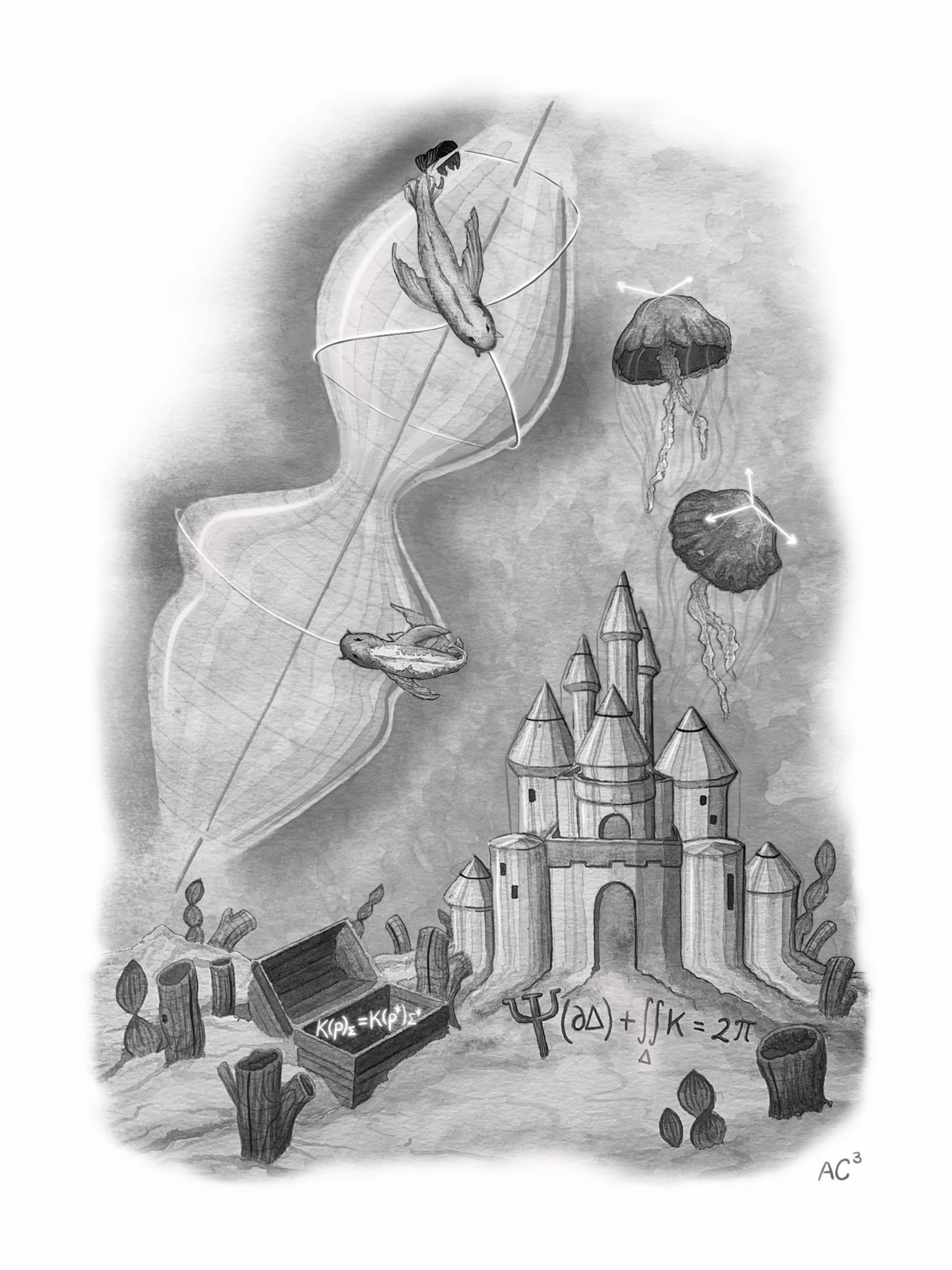}}
\cleardoublepage
\thispagestyle{empty}
\stepcounter{part}
\begin{center}
{\Huge\textbf{Part \Roman{part}:\qquad Geodesics}}\
\end{center}
\addcontentsline{toc}{part}{\Roman{part} Geodesics}
\clearpage
}
{\backgroundsetup{
scale=.95,
opacity=1,
angle=0,
vshift=10mm,
contents={%
  \includegraphics[width=\paperwidth%,height=\paperheight
  ]{pics/Geodesics}
  }%
}

\cleardoublepage
\phantomsection
\stepcounter{part}
\addcontentsline{toc}{part}{\Roman{part} Geodesics}
\thispagestyle{empty}
\begin{center}
{\Huge\textbf{Part \Roman{part}:\qquad Geodesics}}\
\end{center}
\BgThispage
}

%% file: shortest-path.tex
\chapter{Shortest paths}
\label{chap:shortest}

\section{Intrinsic geometry}

Geometry of surfaces is divided into intrinsic and extrinsic.
We say that something is \index{intrinsic}\emph{intrinsic} if it can be checked or measured within the surface;
for example, lengths of curves on the surface or angles between them.
Otherwise, if the definition of something essentially uses the ambient space, then this something is called \index{extrinsic}\emph{extrinsic};
examples include normal curvature.

The mean curvature, as well as the Gauss curvature, are defined via principal curvatures, which are extrinsic.
Later (\ref{thm:remarkable}) it will be shown that remarkably the Gauss curvature is intrinsic --- it can be calculated via measurements within the surface.
The mean curvature is not intrinsic; for example, the intrinsic geometry of the plane is not distinguishable from the intrinsic geometry of the graph $z\z=(x+y)^2$.
However, the mean curvature of the former vanishes at all points, while the mean curvature of the latter does not vanish, say at the origin.  

The following exercise should help you get in the right mood;
it might look like a tedious problem in calculus, but it actually is an amusing problem in geometry.

\begin{wrapfigure}[5]{r}{33 mm}
\vskip-8mm
\centering
\includegraphics{mppics/pic-77}
\vskip-0mm
\end{wrapfigure}

\begin{thm}{Exercise}\label{ex:lasso}
A cowboy stands at the bottom of a frictionless ice-mountain formed by a cone with a circular base with the angle of inclination~$\theta$.
He throws up his lasso which slips neatly over the tip of the cone, pulls it tight, and tries to climb.

What is the critical angle $\theta$ at which the cowboy cannot climb the ice-mountain?
\end{thm}

\section{Definition}

Let $p$ and $q$ be two points on a surface~$\Sigma$.
Recall that $\dist{p}{q}\Sigma$ denotes the induced length distance from $p$ to~$q$;
that is, the greatest lower bound of the lengths of paths in $\Sigma$ from $p$ to~$q$.

A path $\gamma$ from $p$ to $q$ in $\Sigma$ that minimizes the length is called a \index{shortest path}\emph{shortest path} from $p$ to~$q$.
The image of such a shortest path will be denoted by $[p,q]$ or $[p,q]_\Sigma$.\index{10aac@$[p,q]$, $[p,q]_\Sigma$ (shortest path)}
If we write $[p,q]_\Sigma$, we are assuming that such a shortest path exists, and we have chosen one of them.
In general, there might be no shortest path between two given points, and there might be many of them;
this is shown in the following two examples.

{

\begin{wrapfigure}[9]{r}{28 mm}
\vskip-6mm
\centering
\includegraphics{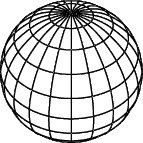}
\bigskip
\includegraphics{mppics/pic-79}
\end{wrapfigure}

\parbf{Nonuniqueness.} There are plenty of shortest paths between the poles of a round sphere --- each meridian is a shortest path.
The latter follows from \ref{obs:S2-length}.

\parbf{Nonexistence.} Let $\Sigma$ be the $(x,y)$-plane with the origin removed.
Consider two points $p=(1,0,0)$ and $q=(-1,0,0)$ in~$\Sigma$.
Let us show that \textit{there is no shortest path from $p$ to $q$ in~$\Sigma$.}

Note that $\dist{p}{q}\Sigma=2$. 
Indeed, given $\epsilon\z>0$, consider the point $s_\epsilon=(0,\epsilon,0)$.
The polygonal path $ps_\epsilon q$ lies in $\Sigma$, and its length $2\cdot\sqrt{1+\epsilon^2}$ approaches $2$ as $\epsilon\to0$.
It follows that $\dist{p}{q}\Sigma\le 2$.
Since $\dist{p}{q}\Sigma\z\ge \dist{p}{q}{\mathbb{R}^3}=2$, we get $\dist{p}{q}\Sigma=2$.

It follows that a shortest path from $p$ to $q$, if it exists, must have length~$2$.
By the triangle inequality, any curve of length $2$ from $p$ to $q$ must run along the line segment $[p,q]$;
in particular, it must pass thru the origin.
Since the origin does not lie in $\Sigma$, there is no shortest path from $p$ to $q$ in $\Sigma$.

\begin{thm}{Proposition}\label{prop:shortest-paths-exist}
Any two points in a proper smooth surface (possibly with boundary) can be joined by a shortest path. 
\end{thm}

\parbf{Proof.}
Choose two points $p$ and $q$ in a proper smooth surface~$\Sigma$.
Set $\ell=\dist{p}{q}\Sigma$.

Since $\Sigma$ is smooth, the points $p$ and $q$ can be joined by a piecewise smooth path in $\Sigma$.
Since any such path is rectifiable, $\ell$ has to be finite.

By the definition of induced length-metric (Section~\ref{sec:Length metric}),
there is a sequence of paths $\gamma_n$ from $p$ to $q$ in $\Sigma$ such that
\[\length\gamma_n\to \ell\quad\text{as}\quad n\to \infty.\]

Without loss of generality, we may assume that $\length\gamma_n<\ell+1$ for any $n$, and each $\gamma_n$ is parametrized proportional to its arc-length.
In particular, each path $\gamma_n\:[0,1]\to\Sigma$ is $(\ell+1)$-Lipschitz; 
that is,
\[|\gamma_n(t_0)-\gamma_n(t_1)|\le (\ell+1)\cdot|t_0-t_1|\]
for any $t_0,t_1\in[0,1]$ and $n$.

The image of $\gamma_n$ lies in the closed ball $\bar B[p,\ell+1]$ for any~$n$.
It follows that the coordinate functions of all $\gamma_n$ are uniformly equicontinuous and uniformly bounded.
By the Arzel\'{a}--Ascoli theorem (\ref{lem:equicontinuous}),
 there is a convergent subsequence of $\gamma_n$, and its limit, say $\gamma_\infty\:[0,1]\to\mathbb{R}^3$, is continuous;
that is, $\gamma_\infty$ is a path.
Evidently, $\gamma_\infty$ runs from $p$ to $q$;
in particular
\[\length\gamma_\infty\ge \ell.\]
Since $\Sigma$ is a closed set, $\gamma_\infty$ lies in~$\Sigma$.

Since length is semicontinuous (\ref{thm:length-semicont}), we get $\length\gamma_\infty\le \ell$.
Therefore, $\length\gamma_\infty\z= \ell$; equivalently, $\gamma_\infty$ is a shortest path from $p$ to~$q$.
\qeds

\section{Nearest-point projection}

The following statement belongs to convex geometry,
but we present its proof since it plays a significant role in the sequel.

\begin{thm}{Lemma}\label{lem:nearest-point-projection}
Let $R$ be a closed convex set in $\mathbb{R}^3$.
Then for every point $p\in\mathbb{R}^3$ there is a unique point $\bar p\in R$ that minimizes the distance to $R$;
that is, $|p-\bar p|\le |p-x|$ for any point $x\in R$.

Moreover, the map $p\mapsto \bar p$ is short;
that is,
\[|p-q|\ge|\bar p-\bar q| \eqlbl{eq:short-cpp}\]
for any pair of points $p,q\in \mathbb{R}^3$.
\end{thm}

The map $p\mapsto \bar p$ is called the \label{nearest-point projection}\index{nearest-point projection}\emph{nearest-point projection};
it maps the Euclidean space onto~$R$.
If $p\in R$, then $\bar p=p$.
In particular $\bar{\bar p}=\bar p$ for any point $p$.

\parbf{Proof.}
Fix a point $p$, and set 
\[\ell=\inf\set{|p-x|}{x\in R}.\]
Choose a sequence $x_n\in R$ such that $|p-x_n|\to \ell$ as $n\to\infty$.

Without loss of generality, we can assume that all the points $x_n$ lie in a ball of radius $\ell+1$ centered at~$p$.
Therefore, we can pass to a \index{partial limit}\emph{partial limit} $\bar p$ of $x_n$; that is, $\bar p$ is a limit of a subsequence of $x_n$.
Since $R$ is closed, $\bar p\in R$.
By construction 
\begin{align*}
|p-\bar p|&=\lim_{n\to\infty}|p-x_n|=\ell.
\end{align*}
Hence, the existence follows.

{

\begin{wrapfigure}{l}{22 mm}
\vskip-0mm
\centering
\includegraphics{mppics/pic-40}
\vskip-0mm
\end{wrapfigure}

Assume there are two distinct points $\bar p, \bar p'\in R$ that minimize the distance to~$p$.
Since $R$ is convex, their midpoint $m=\tfrac12\cdot (\bar p+\bar p')$ lies in~$R$.
Since $|p-\bar p|\z=|p-\bar p'|=\ell$, the triangle $[p\bar p\bar p']$ is isosceles; therefore the triangle $[p\bar p m]$ has a right angle at~$m$.
Since a leg of a right triangle is shorter than its hypotenuse, we have $|p-m|<\ell$, a contradiction. 

It remains to prove \ref{eq:short-cpp}.
We can assume that $\bar p\ne\bar q$; otherwise, there is nothing to prove.

}

If $\measuredangle \hinge{\bar p}{p}{\bar q}< \tfrac\pi2$, then $\dist{p}{x}{}\z<\dist{p}{\bar p}{}$ for a point $x\in [\bar p,\bar q]$ close to $\bar p$.
Since $[\bar p,\bar q]\subset R$,
the latter is impossible.
{

\begin{wrapfigure}{r}{37 mm}
\vskip-6mm
\centering
\includegraphics{mppics/pic-41}
\vskip-0mm
\end{wrapfigure}

Therefore, $\measuredangle \hinge{\bar p}{p}{\bar q}\ge \tfrac\pi2$ or $p=\bar p$.
In both cases, the orthogonal projection of $p$ to the line $\bar p\bar q$ lies behind $\bar p$, or coincides with $\bar p$.
The same way we show that the orthogonal projection of $q$ to the line $\bar p\bar q$ lies behind $\bar q$, or coincides with~$\bar q$.
This implies that the orthogonal projection of the line segment $[p,q]$ contains $[\bar p,\bar q]$.
Whence \ref{eq:short-cpp} follows.
\qeds

}

\section{Shortest paths on convex surfaces}

\begin{thm}{Theorem}\label{thm:shorts+convex}
Assume a surface $\Sigma$ bounds a closed convex region $R$, and $p,q\in \Sigma$.
Denote by $W$ the outer closed region of $\Sigma$;
in other words, $W$ is the union of $\Sigma$ and the complement of~$R$.
Then 
\[\length\gamma\ge \dist{p}{q}\Sigma\]
for any path $\gamma$ in $W$ from $p$ to~$q$.
Moreover, if  $\gamma$ does not lie in $\Sigma$, then the inequality is strict.
\end{thm}

\parbf{Proof.}
Let $\bar\gamma$ be the nearest-point projection of~$\gamma$ to $R$.
The projection $\bar\gamma$ lies in $\Sigma$ and connects $p$ to $q$; therefore, 
\[\length\bar\gamma\ge \dist{p}{q}\Sigma.\]

To prove the first statement, it remains to show that 
\[\length\gamma\ge\length\bar\gamma.\eqlbl{bar-gamma=<gamma}\]

Consider a polygonal line $\bar p_0\dots \bar p_n$ inscribed in $\bar\gamma$.
Let $p_0\dots p_n$ be the corresponding polygonal line inscribed in $\gamma$;
that is $p_i=\gamma(t_i)$ if $\bar p_i=\bar\gamma(t_i)$.
By \ref{lem:nearest-point-projection}, $|p_i-p_{i-1}|\z\ge|\bar p_i-\bar p_{i-1}|$ for any~$i$.
Therefore,
\[\length p_0\dots p_n\ge \length \bar p_0\dots \bar p_n.\]
Applying the definition of length, we get \ref{bar-gamma=<gamma}.

\begin{wrapfigure}{o}{37 mm}
\vskip-0mm
\centering
\includegraphics{mppics/pic-82}
\vskip-0mm
\end{wrapfigure}

It remains to prove the second statement.
Suppose there is a point $w\z=\gamma(t_1)\notin\Sigma$;
note that $w\notin R$.
Then there is a plane $\Pi$ that cuts $w$ from~$\Sigma$; see \ref{lem:separation}).

The curve $\gamma$ must intersect $\Pi$ at two points: one point before $t_1$ and one after.
Let $a\z=\gamma(t_0)$ and $b=\gamma(t_2)$ be these points.
The arc of $\gamma$ from $a$ to $b$ is strictly longer that $|a-b|$;
indeed its length is at least $|a-w|+|w-b|$, and $|a-w|+|w-b|>|a-b|$ since $w\notin[a,b]$.

Remove from $\gamma$ the arc from $a$ to $b$, and replace it with the line segment $[a,b]$.
Denote the obtained curve by $\gamma_1$, then
\[\length\gamma>\length \gamma_1.\]

Since $\gamma_1$ runs in~$W$, we get
\[\length \gamma_1\ge \dist{p}{q}\Sigma.\]
Hence the second statement follows.
\qeds

\begin{thm}{Exercise}\label{ex:length-dist-conv}
Suppose $\Sigma$ is an open or closed smooth surface with positive Gauss curvature, and $\Norm$ is the unit normal field on~$\Sigma$.
Show that 
\[\dist{p}{q}\Sigma\le 2\cdot \frac{|p-q|}{|\Norm(p)+\Norm(q)|}\]
for any two points $p,q\in \Sigma$.

\end{thm}

\begin{wrapfigure}{r}{27 mm}
\vskip-10mm
\centering
\includegraphics{mppics/pic-240}
\end{wrapfigure}

\begin{thm}{Exercise}\label{ex:hat-convex}
Suppose $\Sigma$ is a closed smooth surface that bounds a convex region $R$ 
in $\mathbb{R}^3$
and $\Pi$ is a plane that cuts a hat $\Delta$ from~$\Sigma$.
Assume the reflection of the interior of $\Delta$ across $\Pi$ lies in the interior of~$R$.
Show that $\Delta$ is \index{convex!set}\emph{convex} with respect to the intrinsic metric  of $\Sigma$;
that is, 
if both endpoints of a shortest path in $\Sigma$ 
lie in $\Delta$,
then the entire path lies in $\Delta$.
\end{thm}

Let us define the \index{intrinsic!diameter}\emph{intrinsic diameter} of a closed surface $\Sigma$ as the least upper bound of the lengths of shortest paths in the surface.

\begin{thm}{Exercise}\label{ex:intrinsic-diameter}
Assume a closed smooth surface $\Sigma$ with positive Gauss curvature lies in a unit ball~$B$.

\begin{subthm}{} Show that the intrinsic diameter of $\Sigma$ cannot exceed $\pi$.
\end{subthm}

\begin{subthm}{}
Show that the area of $\Sigma$ cannot exceed $4\cdot \pi$.
\end{subthm}

\end{thm}

%% file: geodesics-new.tex
\chapter{Geodesics}
\label{chap:geodesics}

\section{Definition}

A smooth curve $\gamma$ in a smooth surface is called a \index{geodesic}\emph{geodesic} if, for any~$t$, the acceleration $\gamma''(t)$ is perpendicular to the tangent plane $\T_{\gamma(t)}$.

\begin{thm}{Exercise}\label{ex:helix-geodesic}
Show that the helix $\gamma(t)\z\df(\cos t,\sin t, t)$ is a geodesic on the cylindrical surface $x^2+y^2=1$.
\end{thm}

Geodesics can be understood as the trajectories of particles that slide on $\Sigma$ without friction and external forces.
Indeed, since there is no friction, the force that keeps the particle on $\Sigma$ must be perpendicular to~$\Sigma$.
Therefore, by Newton's second law of motion,
we get that the acceleration $\gamma''$ is perpendicular to $\T_{\gamma(t)}$.

From a physics point of view, the following lemma is a corollary of the law of conservation of energy.

\begin{thm}{Lemma}\label{lem:constant-speed}
Let $\gamma$ be a geodesic in a smooth surface~$\Sigma$. 
Then $|\gamma'|$ is constant.
Moreover, for any $\lambda\in\mathbb{R}$, the curve 
$\gamma_{\lambda}(t)\df \gamma (\lambda\cdot t)$ is a geodesic as well. 
\end{thm}

In other words, any geodesic has a constant speed, and multiplying its parameter by a constant yields another geodesic.

\parbf{Proof.} 
Since $\gamma'(t)$ is a tangent vector at $\gamma(t)$,
we have that $\gamma''(t)\z\perp\gamma'(t)$, or equivalently $\langle\gamma'',\gamma'\rangle=0$ for any~$t$.
Hence 
$\langle\gamma',\gamma'\rangle'=2\cdot \langle\gamma'',\gamma'\rangle=0$.
That is, $|\gamma'|^2=\langle\gamma',\gamma'\rangle$ is constant.

The second part follows since 
$\gamma_{\lambda}''(t) =\lambda^2\cdot \gamma''(\lambda t)$.
\qeds

The statement in the following exercise is called \index{Clairaut's relation}\emph{Clairaut's relation};
it can be obtained from the conservation of angular momentum.

\begin{thm}{Exercise}\label{ex:clairaut}
Let $\gamma$ be a geodesic on a smooth surface of revolution.
Suppose $r(t)$ denotes the distance from $\gamma(t)$ to the axis of rotation
and $\theta(t)$ --- the angle between $\gamma'(t)$ and the latitudinal circle thru $\gamma(t)$. 

Show that the value $r(t)\cdot \cos\theta(t)$ is constant. 
\end{thm}

Recall that an {}\emph{asymptotic line} is a curve in the surface with vanishing normal curvature.

\begin{thm}{Exercise}\label{ex:asymptotic-geodesic}
Assume a curve $\gamma$ is a geodesic and, at the same time, an asymptotic line of a smooth surface.
Show that $\gamma$ is a straight-line segment.
\end{thm}

\begin{thm}{Exercise}\label{ex:reflection-geodesic}
Assume a smooth surface $\Sigma$ is mirror-symmetric with respect to a plane~$\Pi$.
Suppose $\Sigma$ and $\Pi$ intersect along a smooth curve~$\gamma$.
Show that $\gamma$, parametrized by arc-length, is a geodesic on~$\Sigma$.
\end{thm}

\section{Existence and uniqueness}

The following proposition says that the motion without friction and external forces depends smoothly on its initial data. 
The formulation uses the smoothness of the map $w$, which spits a point in $\mathbb{R}^3$ for a point $p$ on the surface $\Sigma$, a tangent vector $\vec{v}$ at $p$, and a real parameter $t$.
In a chart $s$ on $\Sigma$, the point $p$ is described by a pair of coordinates $(u,v)$, and the vector $\vec{v}$ can be expressed as the sum $a\cdot s_u+b\cdot s_v$.
Therefore, in local coordinates, $w$ is defined by the mapping $(u,v,a,b,t)\mapsto w(p,\vec v,t)$ from a subset of $\mathbb{R}^5$ to $\mathbb{R}^3$, to which the usual definition of smoothness applies.
Now, the map $(p,\vec v, t)\z\mapsto w(p,\vec v, t)$ is considered to be \index{smooth!map}\emph{smooth} if it is described by a smooth mapping $(u,v,a,b,t)\mapsto w(p,\vec v,t)$ in any local coordinates on $\Sigma$.

\begin{thm}{Proposition}\label{prop:geod-existence} 
Let $\Sigma$ be a smooth surface without boundary.
Given a tangent vector ${\vec v}$ to $\Sigma$ at a point $p$,
there is a unique geodesic $\gamma\:\mathbb{I}\to \Sigma$ defined on a maximal open interval $\mathbb{I}\ni 0$ that starts at $p$ with velocity vector ${\vec v}$;
that is, $\gamma(0)=p$ and $\gamma'(0)={\vec v}$.
Moreover, we have the following.
\begin{subthm}{prop:geod-existence:smooth}
The map $w\:(p,{\vec v},t)\mapsto \gamma(t)$ is smooth, and it has an open domain of definition%
\footnote{That is, if $w$ is defined for a triple $p\in \Sigma$, ${\vec v}\in \T_p$, and $t\in \mathbb{R}$,
then it is also defined for any triple $q\in \Sigma$, $\vec u\in \T_q$, and $s\in \mathbb{R}$, where $q$, $\vec u$, and $s$ are sufficiently close to $p$, ${\vec v}$, and $t$, respectively.}%
.
\end{subthm}

\begin{subthm}{prop:geod-existence:whole}
If $\Sigma$ is proper, then $\mathbb{I}=\mathbb{R}$; that is, $\gamma$ is defined on the entire real line.
\end{subthm}

\end{thm}

A surface that satisfies the conclusion of \ref{SHORT.prop:geod-existence:whole} for any tangent vector ${\vec v}$ is said to be \index{geodesically complete}\emph{geodesically complete}.
So part \ref{SHORT.prop:geod-existence:whole} says that \textit{any proper surface without boundary is geodesically complete}.
This is a part of the \index{Hopf--Rinow theorem}\emph{Hopf--Rinow theorem} \cite{hopf-rinow}.

In the proof, we will rewrite the definition of a geodesic using a differential equation and then apply Theorems \ref{thm:ODE-nth-order} and \ref{thm:ODE}.

\begin{thm}{Lemma}\label{lem:geodesic=2nd-order}
Let $f$ be a smooth function defined on an open domain in $\mathbb{R}^2$.
A smooth curve $t\mapsto \gamma(t)=(x(t),y(t),z(t))$ is a geodesic in the graph $z=f(x,y)$ if and only if $z(t)=f(x(t),y(t))$ for any $t$ and the functions $t\mapsto x(t)$ and $t\mapsto y(t)$
satisfy a system of differential equations
\[
\begin{cases}
x''=g(x,y,x',y'),
\\
y''=h(x,y,x',y'),
\end{cases}
\]
where $g$ and $h$ are smooth functions of four variables uniquely determined by~$f$.
\end{thm}

\parbf{Proof.} The first condition $z(t)=f(x(t),y(t))$ simply means that $\gamma(t)$ lies on the graph $z=f(x,y)$.

In the following calculations, we often omit the arguments;
in other words, we may use shortcuts $x=x(t)$, $f=f(x,y)=f(x(t),y(t))$, and so on.

First, let us express $z''$ in terms of $f$, $x$, and $y$.
\[
\begin{aligned}
z''&=f(x,y)''=
\\
&=\left(f_x\cdot x'+ f_y\cdot y'\right)'=
\\
&=
f_{xx}\cdot (x')^2
+
f_x\cdot x''
+2\cdot f_{xy}\cdot x'\cdot y'
+
f_{yy}\cdot (y')^2
+
f_y\cdot y''.
\end{aligned}
\eqlbl{eq:def-geod}
\]

The condition
\[\gamma''(t)\perp\T_{\gamma(t)}\] 
means that 
$\gamma''$ is perpendicular to the two basis vectors in $\T_{\gamma(t)}$; that is,
\[
\begin{cases}
\langle \gamma'',s_x\rangle=0,
\\
\langle\gamma'',s_y\rangle=0,
\end{cases}
\]
where $s(x,y)\df (x,y,f(x,y))$, $x=x(t)$, and $y=y(t)$.
Observe that 
$s_x=(1,0, f_x)$ 
and 
$s_y=(0,1, f_y)$.
Since $\gamma''\z=(x'',y'',z'')$, we can rewrite our system as
\[
\begin{cases}
x''+ f_x\cdot z''=0,
\\
y''+ f_y\cdot z''=0.
\end{cases}
\]
It remains to plug in \ref{eq:def-geod} for $z''$, combine the similar terms, and simplify.
\qeds

\parbf{Proof of \ref{prop:geod-existence}.}
Let $z=f(x,y)$ be a description of $\Sigma$ in tangent-normal coordinates at~$p$.
By Lemma \ref{lem:geodesic=2nd-order}, the condition $\gamma''(t)\perp\T_{\gamma(t)}$ can be written as a system of second-order differential equations.
From \ref{thm:ODE-nth-order} and \ref{thm:ODE} we get the existence and uniqueness of a geodesic $\gamma$ in an interval $(-\epsilon,\epsilon)$ for some $\epsilon>0$.

Let us extend the geodesic $\gamma$ to a maximal open interval $\mathbb{I}$.
Suppose there is another geodesic $\gamma_1$ with the same initial data that is defined on a maximal open interval $\mathbb{I}_1$.
Suppose $\gamma_1$ splits from $\gamma$ at some $t_0>0$;
that is, $\gamma_1$ and $\gamma$ coincides on the interval $[0,t_0)$, but they are different on any interval $[0,t_0+\delta)$ with $\delta>0$.
By continuity, $\gamma_1(t_0)=\gamma(t_0)$, and $\gamma_1'(t_0)=\gamma'(t_0)$.
Applying \ref{thm:ODE-nth-order} and \ref{thm:ODE} again, we get that $\gamma_1$ coincides with $\gamma$ in a small neighborhood of $t_0$ --- a contradiction.

The same argument shows that $\gamma_1$ cannot split from $\gamma$ at $t_0<0$.
It follows that $\gamma_1=\gamma$;
in particular, $\mathbb{I}_1=\mathbb{I}$.

If $\Sigma$ is the graph of a smooth function, then part \ref{SHORT.prop:geod-existence:smooth} follows from \ref{thm:ODE-nth-order}, \ref{thm:ODE}, and the lemma.
In this case, the mapping
\[\vec{w}(p,\vec v,t)\df\tfrac{\partial}{\partial t}w(p,\vec v,t)\]
is also smooth.
Note that $\vec{w}(p,\vec v,t_0) \in \T_{\gamma(t_0)}$ is the velocity vector of the geodesic $\gamma\:t\mapsto w(p,\vec v,t)$ at time $t_0$.

In the general case, suppose $w(p,{\vec v},b)$ is defined for $b\ge0$; that is, the geodesic $\gamma\:t\mapsto w(p,{\vec v},t)$ is defined on the interval $[0,b]$.
Then there exists a partition $0=t_0<t_1<\dots<t_n=b$ of the interval $[0,b]$ such that each geodesic segment $\gamma|_{[t_{i-1},t_i]}$ is covered by a chart defined by some tangent-normal coordinates.
Set $p_i=\gamma(t_i)$ and $\vec v_i=\gamma'(t_i)$, so $p_0=p$ and $\vec v_0=\vec v$.
Since $\gamma|_{[t_{i-1},t_i]}$ lies in a graph, by the previous reasoning, we conclude that for each $i$, the mappings
$w$ and $\vec w$ are defined and smooth in a neighborhood of the triples $(p_{i-1},\vec v_{i-1}, t_i-t_{i-1})$.
Note that $p_i=w(p_{i-1},\vec v_{i-1},t_i-t_{i-1})$ and $\vec v_i=\vec w(p_{i-1},\vec v_{i-1},t_i-t_{i-1})$ for each $i$.
Since the composition of smooth mappings is smooth, the mapping $w$ is smoothly defined in a neighborhood of the triple $(p,\vec v, b)$.

The case $b\le 0$ is similar; so, we have obtained the general case in \ref{SHORT.prop:geod-existence:smooth}.

Assume \ref{SHORT.prop:geod-existence:whole} does not hold;
that is, the maximal interval $\mathbb{I}$ is a proper subset of $\mathbb{R}$.
Without loss of generality, we may assume that $b=\sup\mathbb{I}\z<\infty$.
(If not, switch the direction of~$\gamma$.)

By \ref{lem:constant-speed}, $|\gamma'|$ is constant; in particular, $t\mapsto \gamma(t)$ is a uniformly continuous function.
Therefore, the limit point
$q=\lim_{t\to b}\gamma(t)$
is defined.
Since $\Sigma$ is a proper surface, $q\in \Sigma$. 

Applying the argument above in a tangent-normal coordinate chart at $q$, we conclude that $\gamma$ can be extended as a geodesic beyond~$q$.
Therefore, $\mathbb{I}$ is not maximal --- a contradiction.
\qeds

\begin{thm}{Exercise}\label{ex:round-torus}
Let $\Sigma$ be a smooth torus of revolution; that is,
a surface of revolution whose generatrix is a smooth closed curve.
Show that any closed geodesic on $\Sigma$ is noncontractible.

(In other words, if $s\:\mathbb{R}^2\to \Sigma$ is the natural bi-periodic parametrization of $\Sigma$, then
there is no closed curve $\gamma$ in $\mathbb{R}^2$ such that $s\circ\gamma$ is a geodesic.)
\end{thm}

\section{Exponential map}\label{sec:exp}

Let $\Sigma$ be a smooth surface and $p\in \Sigma$.
Given a tangent vector ${\vec v}\in \T_p$, consider a geodesic $\gamma_{\vec v}$ in $\Sigma$ that starts at $p$ with initial velocity~$\vec v$; 
that is, $\gamma(0)=p$ and $\gamma'(0)={\vec v}$.

Let us define the \index{exponential map}\emph{exponential map}%
\footnote{Explaining the reason for this term, would take us far away from the subject.}
at $p$ by
\[\exp_p\:\vec v\mapsto \gamma_{\vec v}(1).\]
By \ref{prop:geod-existence}, this map is smooth, and it is defined in a neighborhood of zero in the tangent plane $\T_p$;
moreover, if $\Sigma$ is proper, 
then $\exp_p$ is defined on the entire plane $\T_p$.

The exponential map
is defined on the tangent plane (or an open subset of it), which is a smooth surface,
and its target is another smooth surface.
We can identify the plane $\T_p$
with its tangent plane $\T_0\T_p$, so the differential $d_0(\exp_p)\:\vec v\mapsto D_{\vec v}\exp_p$ maps $\T_p$ to itself.
Furthermore, by Lemma \ref{lem:constant-speed}, this differential is the identity map; that is, $d_0\exp_p(\vec v)=
\vec v$ for any $\vec v\in \T_p$.

Summarizing, we get the following statement:

\begin{thm}{Observation}\label{obs:d(exp)=1}
Let $\Sigma$ be a smooth surface and $p\in \Sigma$.
Then:

\begin{subthm}{}
The exponential map $\exp_p$ is smooth, and its domain $\Dom(\exp_p)$ contains a neighborhood of the origin in $\T_p$.
Moreover, if $\Sigma$ is proper, then $\Dom(\exp_p)=\T_p$
\end{subthm}

\begin{subthm}{}
The differential $d_0(\exp_p)\:\T_p\to \T_p$ is the identity map.
\end{subthm}

\end{thm}

It is easy to check that $\Dom(\exp_p)$ is \index{star-shaped}\emph{star-shaped} in $\T_p$;
the latter means that if $\vec v\in \Dom(\exp_p)$, then $\lambda\cdot\vec v\in \Dom(\exp_p)$ for any $0\le \lambda\le 1$.

\section{Injectivity radius}

The \index{injectivity radius}\emph{injectivity radius} of $\Sigma$ at $p$ (briefly $\inj(p)$) is defined as least upper bound on radii $r_p\ge0$ such that the exponential map $\exp_p$ is defined on the open ball $B_p\df B(0,r_p)_{\T_p}$,
and the restriction $\exp_p|_{B_p}$ is a smooth regular parametrization of a neighborhood of $p$ in~$\Sigma$.

\begin{thm}{Proposition}\label{prop:exp}
Injectivity radius is positive at any point on a smooth surface $\Sigma$ (without boundary).
Moreover, it is {}\emph{locally bounded away from zero};
that is, for any $p\in\Sigma$ there is $\epsilon>0$ such that if $\dist{p}{q}\Sigma<\epsilon$ for some $q\in \Sigma$, then $\inj(q)\ge\epsilon$.
\end{thm}

The proof of the proposition uses \ref{obs:d(exp)=1} and the inverse function theorem (\ref{thm:inverse}).
In fact, it is true that \textit{the function $\inj\:\Sigma\to (0,\infty]$ is continuous} \cite[5.4]{gromoll-klingenberg-meyer}.

\parbf{Proof.}
Let $z=f(x,y)$ be a local graph representation of $\Sigma$ in tangent-normal coordinates at~$p$.
In this case, the $(x,y)$-plane coincides with the tangent plane $\T_p$.

Denote by $h$ the composition of $\exp_p$ with the projection $(x,y,z)\z\mapsto (x,y)$.
By \ref{obs:d(exp)=1}, the differential $d_0h$ is the identity;
in other words, the Jacobian matrix of $h$ at $0$ is the identity.
Applying the inverse function theorem (\ref{thm:inverse}) we get the first part of the proposition.

The proof of second part is similar, but more technical.

Denote by $h_q$ the composition of $\exp_q$ with the orthogonal projection $(x,y,z)\mapsto (x,y)$.
Consider the chart $s\:(u,v)\z\mapsto (u,v,f(u,v))$.
Let
\[m\:(u,v,a,b)\mapsto h_q(\vec v),\]
where $q=s(u,v)$ and $\vec v=a\cdot s_u+b\cdot s_v$.
By \ref{prop:geod-existence}, $m$ is a smooth map defined in a neighborhood of $0$.
Passing to a smaller neighborhood of $0$, we can assume that first and second partial derivatives of $m$ are bounded.
From above, the Jacobian matrix of $(a,b)\z\mapsto m(0,0,a,b)$ at $0$ is the identity.
It follows that for small fixed $u$ and $v$
the Jacobian matrix of $(a,b)\z\mapsto m(u,v,a,b)$ at $0$ is close to identity.
In particular, we can apply the second part of the inverse function theorem (\ref{thm:inverse}) to
the map $(a,b)\z\mapsto m(u,v,a,b)$  for small fixed $u$ and $v$, hence the result.
\qeds

The proof of the following statement will be indicated in \ref{ex:inj-rad}.

\begin{thm}{Proposition}\label{prop:inj-rad}
Let $p$ be a point on a smooth surface $\Sigma$ (without boundary).
If $\exp_p$ is injective in $B_p=B(0,r)_{\T_p}$, then the restriction $\exp_p|_{B_p}$ is a diffeomorphism between $B_p$ and its image in~$\Sigma$.

\end{thm}

In other words, $\inj(p)$ is the least upper bound on $r$ such that $\exp_p$ is injective in the ball $B(0,r)_{\T_p}$, which motivates the term \textit{injectivity radius}.

\section{Shortest paths are geodesics}

\begin{thm}{Proposition}\label{prop:gamma''}
Let $\Sigma$ be a smooth surface.
Then any shortest path $\gamma$ in $\Sigma$ parametrized proportionally  to its arc-length is a geodesic in~$\Sigma$.
In particular, $\gamma$ is a smooth curve.

A local converse also holds: any point $p$ in $\Sigma$ has a neighborhood $U$ such that any geodesic that lies entirely in $U$ is a shortest path.
\end{thm}

In particular, a sufficiently short arc of any geodesic is a shortest path.
A geodesic that is also a shortest path is called \index{minimizing geodesic}\emph{minimizing}.
As one can see from the following exercise, geodesics might fail to be minimizing.

\begin{thm}{Exercise}\label{ex:helix=geodesic}
Let $\Sigma$ be the cylindrical surface described by the equation $x^2\z+y^2=1$.
Show that the helix $\gamma\:[0,2\cdot\pi]\to \Sigma$ defined by $\gamma(t)\z\df(\cos t,\sin t, t)$
is a geodesic, but not a shortest path on~$\Sigma$.
\end{thm}

A formal proof of the proposition will be given in Section~\ref{sec:proof-of-gamma''}.
The following informal physical explanation might be sufficiently convincing.
In fact, if one assumes that $\gamma$ is smooth, then it is easy to convert this explanation into a rigorous proof.

\parbf{Informal explanation.}
Let us think about a shortest path $\gamma$ as a stable position of a stretched elastic thread that is forced to lie on a frictionless surface.
Since it is frictionless, the force density $\vec n=\vec n(t)$ that keeps $\gamma$ in the surface must be proportional to the normal vector to the surface at~$\gamma(t)$.

Denote the tension in the thread by $\tau$.
It has to be the same at all points;
otherwise, the thread would move back or forth, and it would not be stable.

We can assume that $\gamma$ has unit speed;
in this case, the net force from tension along the arc $\gamma_{[t_0,t_1]}$ is $\tau\cdot(\gamma'(t_1)-\gamma'(t_0))$.
Hence, the density of the net force from tension at $t_0$ is 
\begin{align*}
\vec f(t_0)&=\lim_{t_1\to t_0}\tau\cdot\frac{\gamma'(t_1)-\gamma'(t_0)}{t_1-t_0}=
\\
&=\tau\cdot\gamma''(t_0).
\end{align*}
By Newton's second law of motion,  
$\vec f+\vec n=0$; hence $\gamma''(t)\perp\T_{\gamma(t)}\Sigma$.
\qeds

\begin{thm}{Corollary}
Let $p$ be a point on a smooth surface $\Sigma$, and $r\z\le \inj(p)$.
Then the exponential map $\exp_p$ defines a diffeomorphism from $B(0,r)_{\T_p}$ to $B(p,r)_\Sigma$.
\end{thm}

\parbf{Proof.}
By \ref{prop:inj-rad}, the restriction of $\exp_p$ to $B_p={B(0,r)_{\T_p}}$ is a diffeomorphism to its image $\exp_p(B_p)\subset \Sigma$.

Evidently, $B(p,r)_\Sigma\supset\exp_p(B_p)$.
By \ref{prop:gamma''}, $B(p,r)_\Sigma\subset\exp_p(B_p)$, hence the result.
\qeds

By the corollary, the restriction $\exp_p|_{B(0,r)_{\T_p}}$ admits an inverse map called \index{logarithmic map}\emph{logarithmic};
it is denoted by \[\log_p\:B(p,r)_\Sigma\to B(0,r)_{\T_p}.\]

According to the proposition above, any shortest path parametrized by its arc-length is a smooth curve.
This observation helps to solve the following two exercises:

\begin{thm}{Exercise}\label{ex:two-min-geod}
Show that if two shortest paths have two distinct common points $p$ and $q$, then either $p$ and $q$ are the endpoints of both shortest paths, or both shortest paths share a common arc from $p$ to~$q$.

Show by an example that nonoverlapping geodesics can cross each other an arbitrary number of times.
\end{thm}

\begin{thm}{Exercise}\label{ex:min-geod+plane}
Assume a smooth surface $\Sigma$ is mirror-symmetric with respect to a plane $\Pi$.
Show that no shortest path in $\Sigma$ can {}\emph{cross} $\Pi$ more than once.

In other words, if you travel along a shortest path, then you change sides of $\Pi$ at most once.
\end{thm}

{

\begin{wrapfigure}{r}{40 mm}
\vskip-8mm
\centering
\includegraphics{mppics/pic-250}
\vskip-0mm
\end{wrapfigure}

\begin{thm}{Advanced exercise}\label{ex:milka}
Let $\Sigma$ be a smooth closed strictly convex surface 
in $\mathbb{R}^3$ 
and $\gamma\:[0,\ell]\z\to \Sigma$ be a unit-speed minimizing geodesic.
Set $p\z=\gamma(0)$, $q=\gamma(\ell)$, and 
\[p^s=\gamma(s)-s\cdot\gamma'(s).\]

Show that for any $s\in (0,\ell)$,
one cannot see $q$ from $p^s$;
that is, the line segment $[p^s,q]$ intersects $\Sigma$ at a point distinct from~$q$.

Show that the statement does not hold without assuming that $\gamma$ is minimizing.
\end{thm}

}

\begin{thm}{Exercise}\label{ex:round-sphere}
Let $\Sigma$ be a smooth closed surface.
Suppose that for any $p,q\z\in \Sigma$ the distance $\dist{p}{q}\Sigma$ depends only on the distance $\dist{p}{q}{\mathbb{R}^3}$.
Show that $\Sigma$ is a round sphere.
\end{thm}

\begin{thm}{Very advanced exercise}\label{ex:rad=2}
Let \(\Theta\) be a sphere of radius $2$ centered at $0\in\mathbb{R}^3$
and let $\Sigma$ be a smooth closed surface in the open ball bounded by \(\Theta\).
Assume that all normal curvatures of \(\Sigma\) do not exceed~$1$ by absolute value.

\begin{subthm}{ex:rad=2:a}
Show that there is a diffeomorphism \(\rho\:\Theta \to \Sigma \) (let us call it \emph{radial projection}) that sends a point \( p \in \Theta \) to a unique intersection point in \( \Sigma \cap [0,p]_{\mathbb{R}^3}\).
In particular, $\Sigma$ bounds a star-shaped region.
\end{subthm}

\begin{subthm}{ex:rad=2:b}
Show that the radial projection \(\rho\: \Theta \to \Sigma \) is length-nonincreasing.
\end{subthm}

\begin{subthm}{ex:rad=2:c}
Suppose that $x\in \Sigma$ maximize the distance to the origin.
Show that the ball with a diameter $[0, x]$ lies in the region bounded by $\Sigma$.
Conclude that the region bounded by $\Sigma$ contains a unit ball.
\end{subthm}

\end{thm}

\section{Liberman's lemma}

A version of the following lemma was used by Joseph Liberman \cite{liberman}.

\begin{thm}{Lemma}
\label{lem:liberman}
\index{Liberman's lemma}
Let $f$ be a smooth locally convex function defined on an open subset of the plane.
Suppose $t\mapsto \gamma(t)\z=(x(t),y(t),z(t))$ is a unit-speed geodesic on the graph $z=f(x,y)$.
Then $t\mapsto z(t)$ is a convex function; that is, $z''(t)\ge 0$ for any~$t$.
\end{thm}

\parbf{Proof.}
Choose the orientation on the graph so that the unit normal vector $\Norm$ always points up;
that is, $\Norm$ has a positive $z$-coordinate at each point.
Let us use the shortcut $\Norm(t)$ for $\Norm(\gamma(t))$.

Since $\gamma$ is a geodesic, we have $\gamma''(t)\perp\T_{\gamma(t)}$,
or equivalently the acceleration $\gamma''(t)$ is proportional to $\Norm(t)$ for any~$t$.
Furthermore,
\[\gamma''=k\cdot\Norm,\]
where $k=k(t)$ is the normal curvature of the graph at $\gamma(t)$ in the direction of $\gamma'(t)$.

Therefore,
\[z''=k\cdot\cos\theta,
\eqlbl{eq:z''}\]
where $\theta=\theta(t)$ denotes the angle between $\Norm(t)$ and the $z$-axis.

Since $\Norm$ points up, we have $\theta(t)<\tfrac\pi2$, or equivalently $\cos\theta>0$.

Since $f$ is convex, the tangent plane supports the graph from below at any point;
in particular, $k(t)\ge 0$ for any~$t$.
It follows that the right-hand side in \ref{eq:z''} is nonnegative;
whence the statement follows.
\qeds

\begin{thm}{Exercise}\label{ex:closed-liberman}
Let $\Sigma$ be the graph of a locally convex function defined on an open subset of the plane.
Show that $\Sigma$ has no closed geodesics.
\end{thm}

\begin{thm}{Exercise}\label{ex:rho''}
Assume $\gamma$ is a unit-speed geodesic on a smooth convex surface $\Sigma$, and a point $p$ lies in the interior of the convex set bounded by~$\Sigma$.
Set $\rho(t)=|p-\gamma(t)|^2$.
Show that $\rho''(t)\le 2$ for any~$t$.
\end{thm}

\section{Total curvature of geodesics}

Recall that $\tc\gamma$ denotes the total curvature of the curve~$\gamma$;
see~\ref{sec:Total curvature}.

\begin{thm}{Exercise}\label{ex:tc-spherical-image}
Let $\gamma$ be a geodesic on an oriented smooth surface $\Sigma$
with unit normal field $\Norm$.
Show that $\length(\Norm\circ\gamma)\ge \tc\gamma$.
\end{thm}

\begin{thm}{Theorem}\label{thm:usov}
Assume $\Sigma$ is the graph of a convex $\ell$-Lipschitz function $f$ defined on an open set in the $(x,y)$-plane.
Then the total curvature of any geodesic in $\Sigma$ is at most $2\cdot \ell$.
\end{thm}

This theorem was first proved by Vladimir Usov \cite{usov}.

\parbf{Proof.}
Let $t\mapsto\gamma(t)=(x(t),y(t),z(t))$ be a unit-speed geodesic on~$\Sigma$.
According to Liberman's lemma (\ref{lem:liberman}), the function
$t\mapsto z(t)$ is convex.

Since the slope of $f$ is at most $\ell$, we have
\[|z'(t)|\le \frac{\ell}{\sqrt{1+\ell^2}}\]
for any $t$.
We can assume that $\gamma$ is defined on the interval $[a,b]$.
Then
\[
\begin{aligned}
\int_a^b z''(t) dt&=z'(b)-z'(a)\le 
 2\cdot \frac{\ell}{\sqrt{1+\ell^2}}.
\end{aligned}
\eqlbl{eq:intz''}
\]

Also, notice that $z''$ is the projection of $\gamma''$ to the $z$-axis.
The slope of the tangent plane $\T_{\gamma (t)} \Sigma$ cannot be greater than $\ell$ for any~$t$.
Because $\gamma ''$ is perpendicular to that plane, we have that
\[|\gamma'' (t)| \le z''(t)\cdot\sqrt{1+ \ell ^2}.\]

By \ref{eq:intz''}, we get that
\begin{align*}
\tc\gamma&=\int_a^b|\gamma'' (t)|\cdot dt\le 
\\
&\le \sqrt{1+ \ell ^2}\cdot \int_a^b z''(t)\cdot dt\le 
\\
&\le 2\cdot \ell.
\end{align*}
\qedsf

By the following exercise, the estimate in the theorem is optimal.

\begin{thm}{Exercise}\label{ex:usov-exact}
Let $\Sigma$ be the graph $z=\ell\cdot\sqrt{x^2+y^2}$ with the origin removed.
Show that any both-side-infinite geodesic $\gamma$ in $\Sigma$ has total curvature exactly $2\cdot \ell$.
\end{thm}

\begin{thm}{Exercise}\label{ex:ruf-bound-mountain}
Assume $f$ is a smooth convex $\tfrac32$-Lipschitz function defined on the $(x,y)$-plane.
Show that any geodesic $\gamma$ on the graph $z\z=f(x,y)$ is simple;
that is, it has no self-intersections.

Construct a convex $2$-Lipschitz function defined on the plane
with a nonsimple geodesic $\gamma$ in its graph.
\end{thm}

\begin{thm}{Theorem}\label{thm:tc-of-mingeod}
Suppose a smooth surface $\Sigma$ bounds a convex set $K$.
Assume $B(0,\epsilon)\subset K\subset B(0,1)$.
Then the total curvatures of any shortest path in $\Sigma$ can be bounded in terms of~$\epsilon$. 
\end{thm}

\begin{wrapfigure}{r}{47 mm}
\vskip-0mm
\centering
\includegraphics{mppics/pic-83}
\vskip-0mm
\end{wrapfigure} 

\begin{thm}{Proof-guided exercise}\label{ex:bound-tc}
Let $\Sigma$ be as in the theorem, and $\gamma$ be a unit-speed shortest path in~$\Sigma$.
Denote by $\Norm(t)$ the unit normal vector at $\gamma(t)$ that points outside $\Sigma$;
denote by $\theta(t)$ the angle between $\Norm(t)$ and the direction from the origin to $\gamma(t)$.
Set $\rho(t)\z=|\gamma(t)|^2$; denote by $k(t)$ the curvature of $\gamma$ at~$t$.

\begin{subthm}{ex:bound-tc:a}
Show that $\cos(\theta(t))\ge \epsilon$ for any~$t$.
\end{subthm}

\begin{subthm}{ex:bound-tc:b}
 Show that $|\rho'(t)|\le 2$ for any~$t$.
\end{subthm}

\begin{subthm}{ex:bound-tc:c}
 Show that 
\[\rho''(t)=2-2\cdot k(t)\cdot \cos \theta(t)\cdot |\gamma(t)|\]
for any~$t$.
\end{subthm}

\begin{subthm}{ex:bound-tc:d}
 Use the nearest-point projection from the unit sphere to $\Sigma$ to show that 
\[\length \gamma\le \pi.\]
\end{subthm}

\begin{subthm}{ex:bound-tc:e}
Conclude that $\tc\gamma\le 100/\epsilon^2$.
\end{subthm}

\end{thm}

\parit{Remark.}
The obtained estimate tends to infinity as $\epsilon \to 0$,
but there is also an estimate that does not depend on $\epsilon$;
this is a result by Nina Lebedeva and the first author \cite{lebedeva-petrunin}.
Aleksei Pogorelov proposed the hypothesis that there exists a bound on the length of the spherical image of a shortest path \cite{pogorelov}.
According to \ref{ex:tc-spherical-image}, this hypothesis is stronger,
but counterexamples have been found for all of its possible variations \cite{zalgaller,milka,usov,pach}.

%% file: parallel-new.tex
\chapter{Parallel transport}
\label{chap:parallel-transport}

\section{Parallel tangent fields}

Let $\Sigma$ be a smooth surface and $\gamma\:[a,b]\z\to \Sigma$ be a smooth curve.
A smooth vector-valued function $t\mapsto {\vec v}(t) \in \mathbb{R}^3$ is called a \index{tangent!field}\emph{tangent field} along $\gamma$ if
the vector ${\vec v}(t)$ lies in the tangent plane $\T_{\gamma(t)}\Sigma$ for each~$t$.

A tangent field $\vec v$ on $\gamma$ is called \index{parallel!field}\emph{parallel} if ${\vec v}'(t)\perp\T_{\gamma(t)}$ for any~$t$.

In general, the family of tangent planes $\T_{\gamma(t)}\Sigma$ is not parallel.
Therefore, one cannot expect to have a \textit{truly} parallel field with ${\vec v}'(t)\equiv 0$.
The condition ${\vec v}'(t)\perp\T_{\gamma(t)}$ means that the family is \textit{as parallel as possible} --- it rotates together with the tangent plane but does not rotate inside the plane.

By the definition of geodesic, the velocity field ${\vec v}(t)\z=\gamma'(t)$ of any geodesic $\gamma$ is parallel on~$\gamma$.

\begin{thm}{Exercise}\label{ex:parallel}
Let $\Sigma$ be a smooth surface and let 
$\gamma\:[a,b]\to \Sigma$ be a smooth curve.
Suppose ${\vec v}(t)$, $\vec w(t)$ are parallel vector fields along~$\gamma$.

\begin{subthm}{ex:parallel:a} Show that $|{\vec v}(t)|$ is constant.
\end{subthm}

\begin{subthm}{ex:parallel:b} Show that the angle $\theta(t)=\measuredangle({\vec v}(t),\vec w(t))$ is constant.
\end{subthm}

\end{thm}

\section{Parallel transport}

\begin{thm}{Proposition-Definition}\label{prop:parallel}
Let $\gamma\:[a,b]\z\to \Sigma$ be a piecewise smooth curve on a smooth surface $\Sigma$.
Assume $p=\gamma(a)$ and $q=\gamma(b)$.

Given a tangent vector ${\vec w}\in\T_p$ there is a unique parallel field ${\vec w}(t)$ along $\gamma$ such that ${\vec w}(a)={\vec w}$.

The vector ${\vec w}(b)\in\T_q$ is called the \index{parallel!transport}\emph{parallel transport} of ${\vec w}(a)$ along~$\gamma$ in~$\Sigma$.
\end{thm}

The parallel transport along $\gamma$ will be denoted by $\iota_\gamma$;
so we can write $\vec w(b)=\iota_\gamma({\vec w}(a))$ or $\vec w(b)=\iota_\gamma({\vec w}(a))_\Sigma$ if we need to emphasize that $\gamma$ lies on the surface~$\Sigma$.
From Exercise~\ref{ex:parallel}, it follows that parallel transport $\iota_\gamma\:\T_p\z\to\T_q$ is an isometry.
In general, the parallel transport $\iota_\gamma\:\T_p\z\to\T_q$ depends on the choice of $\gamma$; that is, for another curve $\gamma_1$ connecting $p$ to $q$ in $\Sigma$, the parallel transports $\iota_{\gamma_1}$ and $\iota_{\gamma}$ might be different.

\parbf{Sketch of proof.}
Assume $\gamma$ is smooth and lies in a local chart $(u,v)\z\mapsto s(u,v)$ of $\Sigma$;
so $\gamma(t)\z=s(u(t),v(t))$ for smooth functions $t\mapsto u(t)$ and $t\mapsto v(t)$.
Set 
\[
\vec u(t)=s_u(u(t),v(t)),
\quad
\vec v(t)=s_v(u(t),v(t)),
\quad
\text{and}
\quad
\Norm(t)=\Norm(\gamma(t)).
\]

The conditions $\vec w(t)\in \T_{\gamma(t)}$ and $\vec w'(t)\perp \T_{\gamma(t)}$ can be written as a system of equations:
\[
\begin{cases}
\ \langle\vec w(t), \Norm(t)\rangle=0,
\\
\ \langle\vec w'(t), \vec u(t)\rangle=0,
\\
\ \langle\vec w'(t), \vec v(t)\rangle=0.
\end{cases}
\]
Rewriting this system in components of $\vec u$, $\vec v$, $\vec w$, and $\Norm$ produces a system of ordinary differential equations on the components of $\vec w$.
Applying \ref{thm:ODE}, we get a solution $\vec w(t)$.
By \ref{ex:parallel}, $|\vec w|$ is constant;
therefore \ref{thm:ODE} implies that $\vec w$ is defined on the whole interval $[a,b]$.

If $\gamma$ is only piecewise smooth and/or not covered by one chart we can subdivide it into smooth arcs $\gamma_1,\dots,\gamma_n$ such that each $\gamma_i$ lies in one chart.
Applying the statement consecutively to each $\gamma_i$ proves it for $\gamma$.
\qeds

Suppose $\gamma_1$ and $\gamma_2$ are two smooth curves in two smooth surfaces $\Sigma_1$ and $\Sigma_2$.
Denote by $\Norm_i\:\Sigma_i\to\mathbb{S}^2$ the Gauss maps of $\Sigma_1$ and $\Sigma_2$.
If $\Norm_1\circ\gamma_1(t)= \Norm_2\circ\gamma_2(t)$ for any $t$, then we say that the curves $\gamma_1$ and $\gamma_2$ have {}\emph{identical normals} in $\Sigma_1$ and $\Sigma_2$, respectively.

In this case, the tangent plane $\T_{\gamma_1(t)}\Sigma_1$ is parallel to $\T_{\gamma_2(t)}\Sigma_2$ for any~$t$.
Therefore we can identify $\T_{\gamma_1(t)}\Sigma_1$ with $\T_{\gamma_2(t)}\Sigma_2$.
In particular, if $\vec v(t)$ is a tangent vector field along $\gamma_1$,
then it is also a tangent vector field along $\gamma_2$.
Moreover, $\vec v'(t)\perp \T_{\gamma_1(t)}\Sigma_1$ is equivalent to $\vec v'(t)\perp \T_{\gamma_2(t)}\Sigma_2$; that is, if $\vec v(t)$ is parallel along $\gamma_1$,
then it is also parallel along $\gamma_2$.

The discussion above leads to the following observation that will play a key role in the sequel.

\begin{thm}{Observation}\label{obs:parallel=}
Let $\gamma_1$ and $\gamma_2$ be two smooth curves in two smooth surfaces $\Sigma_1$ and $\Sigma_2$.
Suppose $\gamma_1$ and $\gamma_2$ have identical normals (as the curves in $\Sigma_1$ and $\Sigma_2$, respectively).
Then the parallel transports $\iota_{\gamma_1}$ and $\iota_{\gamma_2}$ are identical. 
\end{thm}

\begin{thm}{Exercise}\label{ex:parallel-transport-support}
Let $\Sigma_1$ and $\Sigma_2$ be two surfaces with a common curve~$\gamma$.
Suppose $\Sigma_1$ bounds a region that contains $\Sigma_2$.
Show that the parallel transport along $\gamma$ in $\Sigma_1$
coincides with the parallel transport along $\gamma$ in $\Sigma_2$.
\end{thm}

\section{Bike wheel and projections}

Here we give two interpretations of parallel transport;
they might help to build the right intuition, but will not help in writing rigorous proofs.
The first is an experiment with a bike wheel suggested by Mark Levi \cite{levi}.
The second one uses orthogonal projections of tangent planes.

Suppose $\gamma\:[a,b]\to\Sigma$ is a smooth arc in a smooth surface $\Sigma$.
Think of walking along $\gamma$ and carrying a perfectly balanced bike wheel so that you keep the wheel's axis normal to $\Sigma$ and touch only its axis.
If the wheel is not spinning at the starting point $p=\gamma(a)$, then it will not be spinning after stopping at~$q=\gamma(b)$.
Indeed, by pushing the axis one cannot produce torque to spin the wheel.
Then the map that sends the initial position of the wheel to the final position is the parallel transport~$\iota_\gamma$.

\begin{figure}[ht!]
\vskip-0mm
\centering
\includegraphics[scale=.3]{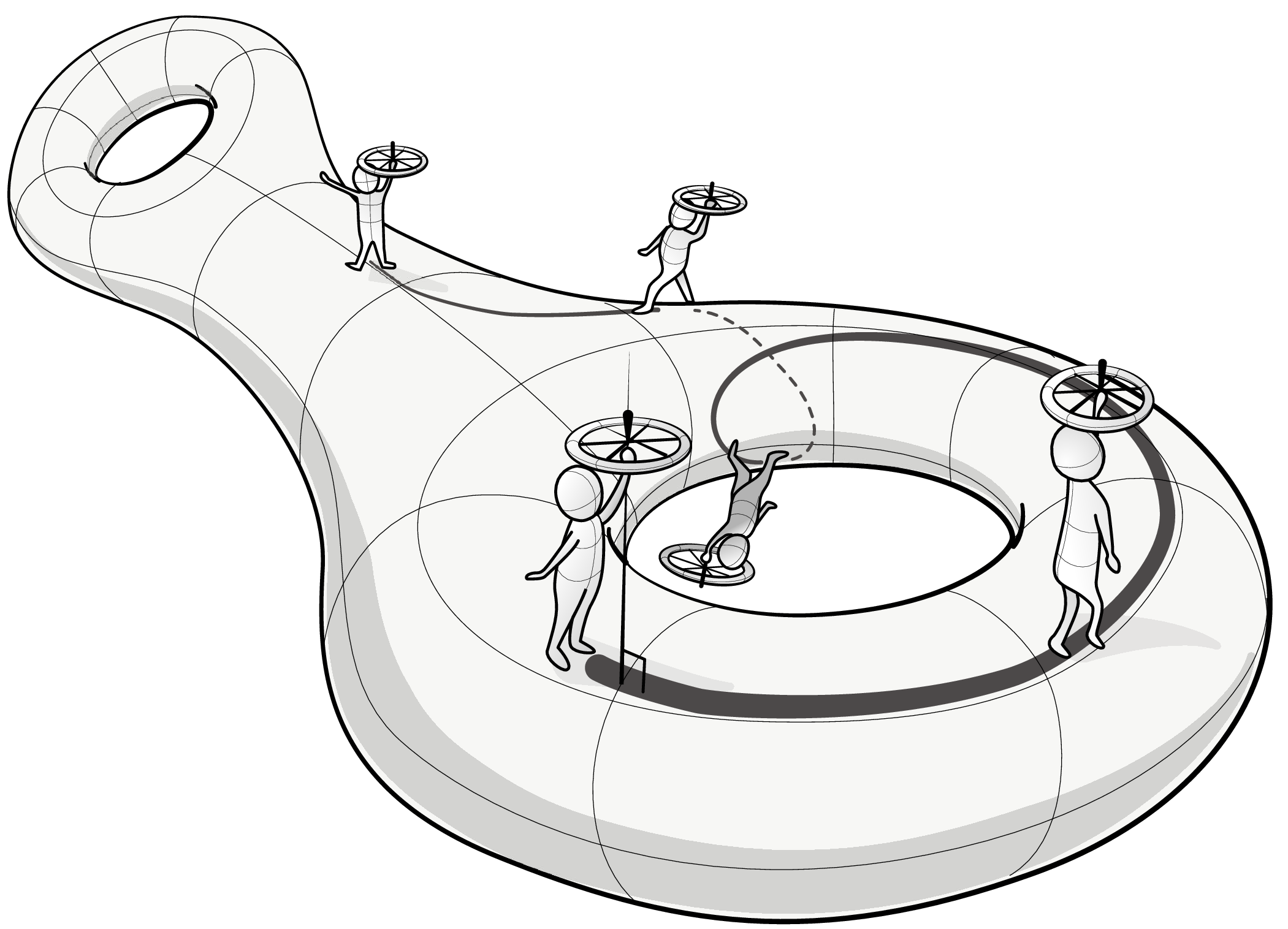}
\end{figure}

By the way, physical intuition should tell you that \textit{moving the axis of the wheel without changing its direction does not change the direction of the wheel's spokes};
this is essentially the statement of Observation~\ref{obs:parallel=}.

On a more formal level, one can choose a partition $a=t_0<\dots\z<t_n=b$ of $[a,b]$
and consider the sequence of orthogonal projections $\phi_i\:\T_{\gamma(t_{i-1})}\to\T_{\gamma(t_i)}$.
For a fine partition, the composition 
\[\phi_n\circ\dots\circ\phi_1\:\T_p\z\to\T_q\]
gives an approximation of $\iota_\gamma$.

Each $\phi_i$ does not increase the magnitude of a vector and neither the composition.
It is straightforward to see that if the partition is sufficiently fine, then the above composition is almost an isometry.
Therefore, the limit $\iota_\gamma$ is an isometry $\T_p\z\to\T_q$.

\begin{thm}{Exercise}\label{ex:holonomy=not0}
Construct a loop $\gamma$ in $\mathbb{S}^2$ such that the parallel transport $\iota_\gamma$ is not the identity map.
\end{thm}

\section{Total geodesic curvature}

Recall that geodesic curvature is defined in Section~\ref{sec:Darboux}.
It is denoted by $k_g$, and it measures how much a curve $\gamma$ on an oriented surface diverges from being a geodesic;
it is positive where $\gamma$ turns left and negative where $\gamma$ turns right.
In particular, geodesics have vanishing geodesic curvature.

Suppose $\Sigma$ is a smooth oriented surface, and $\gamma\:\mathbb{I}\to \Sigma$ is a smooth unit-speed curve.
The total geodesic curvature of $\gamma$ is defined by the integral 
\[\tgc\gamma
\df
\int_{\mathbb{I}} k_g(t)\cdot dt.\]

If $\Sigma$ is a plane, then the geodesic curvature of $\gamma$ coincides
with its signed curvature (see Section~\ref{sec:Total signed curvature}). 
In particular, its total geodesic curvature is equal to its total signed curvature.
For that reason, we use the same notation $\tgc\gamma$; if we need to emphasize that we consider $\gamma$ as a curve in $\Sigma$, we write $\tgc\gamma_\Sigma$.\index{10psi@$\tgc\gamma$, $\tgc\gamma_\Sigma$ (total geodesic curvature)}

If $\gamma$ is a piecewise smooth curve in $\Sigma$, then
its \index{total!geodesic curvature}\emph{total geodesic curvature} is defined as the sum of all total geodesic curvatures of its arcs plus the sum of the signed exterior angles of $\gamma$ at the joints.
More precisely, if $\gamma$ is a concatenation of smooth curves $\gamma_1,\dots,\gamma_n$, then
\[\tgc\gamma
\df
\tgc{\gamma_1}+\dots+\tgc{\gamma_n}+\theta_1+\dots+\theta_{n-1},\]
where $\theta_i$ is the signed external angle at the joint between $\gamma_i$ and $\gamma_{i+1}$;
it is positive if $\gamma$ turns left and negative if $\gamma$ turns right, it is undefined if $\gamma$ turns to the opposite direction.
If $\gamma$ is closed, then 
\[\tgc\gamma
\df
\tgc{\gamma_1}+\dots+\tgc{\gamma_n}+\theta_1+\dots+\theta_n,\]
where $\theta_n$ is the signed external angle at the joint between $\gamma_n$ and $\gamma_1$.

If each arc $\gamma_i$ in the concatenation is a minimizing geodesic, then $\gamma$ is called a \index{broken geodesic}\emph{broken geodesic}.
In this case, $\tgc{\gamma_i}=0$ for each~$i$.
Therefore the total geodesic curvature of $\gamma$ is the sum of its signed external angles.

\begin{thm}{Proposition}\label{prop:pt+tgc}
Assume $\gamma$ is a closed broken geodesic in a smooth oriented surface $\Sigma$ that starts and ends at the point~$p$.
Then the parallel transport $\iota_\gamma\:\T_p\to\T_p$ is a clockwise rotation of the plane $\T_p$ by the angle~$\tgc\gamma$.

Moreover, the same statement holds true for smooth closed curves and piecewise smooth curves.
\end{thm}

\parbf{Proof.}
Assume $\gamma$ is a cyclic concatenation of geodesics $\gamma_1,\dots,\gamma_n$.
Fix a tangent vector ${\vec v}$ at $p$, and extend it to a parallel vector field along~$\gamma$.
Since $\tan_i(t)=\gamma_i'(t)$ is parallel along $\gamma_i$, the oriented angle $\phi_i$ from $\tan_i$ to ${\vec v}$ stays constant along each~$\gamma_i$.

\begin{wrapfigure}{o}{22 mm}
\vskip-0mm
\centering
\includegraphics{mppics/pic-48}
\vskip4mm
\end{wrapfigure}

Suppose $\theta_i$ denotes the external angle at the joint between $\gamma_{i}$ and $\gamma_{i+1}$.
Then 
\[\phi_{i+1}=\phi_i-\theta_i \pmod{2\cdot\pi}.\]
Therefore, after going over all vertices we get that 
\[\phi_{n+1}-\phi_1=-\theta_1-\dots-\theta_n=-\tgc\gamma \pmod{2\cdot\pi}.\]
Hence, the first statement follows.

For a smooth unit-speed curve $\gamma\:[a,b]\to\Sigma$, the proof is analogous.
Denote by $\phi(t)$ the signed angle from $\tan (t)$ to ${\vec v} (t)$.
Let us show that 
\[\phi'(t)+k_g(t)\equiv0\eqlbl{eq:phi'+kg}\]

Recall that $\mu=\mu(t)$ denotes the counterclockwise rotation of $\tan \z=\tan(t)$ by angle $\tfrac\pi2$ in $\T_{\gamma(t)}$.
Denote by $\vec w=\vec w(t)$ the counterclockwise rotation of $\vec v=\vec v(t)$ by angle $\tfrac\pi2$ in $\T_{\gamma(t)}$.
Then
\begin{align*}
\tan&=\cos\phi\cdot \vec v-\sin\phi\cdot \vec w,
\\
\mu&=\sin\phi\cdot \vec v+\cos\phi\cdot \vec w.
\end{align*}

The vector fields $\vec v$ and $\vec w$ are parallel along $\gamma$; that is, $\vec v'(t),\vec w'(t)\z\perp\T_{\gamma(t)}$.
Therefore, $\langle\vec v',\mu\rangle\z=\langle\vec w',\mu\rangle\z=0$.
It follows that
\begin{align*}
k_g&=\langle\tan',\mu\rangle=
\\
&=-(\sin^2\phi+\cos^2\phi)\cdot \phi'
\\
&=-\phi'.
\end{align*}
It proves \ref{eq:phi'+kg}.

By \ref{eq:phi'+kg} we get that 
\begin{align*}
\phi(b)-\phi(a)&=\int_a^b \phi'(t)\cdot dt=
\\
&=-\int_a^b k_g\cdot dt=
\\
&=-\tgc\gamma
\end{align*}

In the case when $\gamma$ is a piecewise smooth curve, the result follows from a straightforward combination of the above two cases. 
\qeds

%% file: gauss-bonnet.tex
\chapter{Gauss--Bonnet formula}
\label{chap:gauss-bonnet}
\section{Formulation}

The following theorem was proved by Carl Friedrich Gauss \cite{gauss}
for geodesic triangles;
Pierre Bonnet and Jacques Binet independently 
generalized the statement for arbitrary curves.
\index{Gauss--Bonnet formula}

\begin{thm}{Theorem}\label{thm:gb}
Let $\Delta$ be a topological disc in a smooth oriented surface $\Sigma$ bounded by a piecewise smooth curve $\partial\Delta$.
Suppose $\partial \Delta$ is oriented so that $\Delta$ lies on its left.
Then 
\[\tgc{\partial\Delta}+\iint_\Delta K=2\cdot \pi,\eqlbl{eq:g-b}\]
where $K$ denotes the Gauss curvature.
\end{thm}

In this chapter,
we give an informal proof of this formula in a leading special case.
A formal computational proof will be given in Section~\ref{sec:gauss--bonnet:formal}.
Before going into the proofs, we suggest solving the following exercises using the Gauss--Bonnet formula.

\begin{thm}{Exercise}\label{ex:1=geodesic-curvature}
Assume $\gamma$ is a simple closed curve with constant geodesic curvature $1$ in a smooth closed surface $\Sigma$ with positive Gauss curvature.
Show that $\length\gamma<2\cdot\pi$;
that is, $\gamma$ is shorter than the unit circle.  
\end{thm}

\begin{thm}{Exercise}\label{ex:GB-hat}
Suppose that a disc $\Delta$ lies in a graph $z=f(x,y)$ of a smooth function,
the boundary of $\Delta$ lies in the $(x,y)$-plane 
and $\Delta$ meets this plane at a fixed angle $\alpha$.
Show that $\iint_\Delta K=2\cdot\pi\cdot(1-\cos\alpha)$.
\end{thm}

\begin{thm}{Exercise}\label{ex:geodesic-half}
Let $\gamma$ be a simple closed geodesic on a smooth closed surface $\Sigma$ with positive Gauss curvature.
Denote by $\Norm$ the spherical map on $\Sigma$.
Show that the curve $\alpha=\Norm\circ\gamma$ divides the sphere into two regions of equal area.

Conclude that $\length \alpha\ge 2\cdot\pi.$
\end{thm}

\begin{thm}{Exercise}\label{ex:closed-geodesic}
Let $\gamma$ be a closed geodesic possibly with self-intersections on a smooth closed surface $\Sigma$ with positive Gauss curvature.
Suppose $R$ is one of the regions that $\gamma$ cuts from~$\Sigma$.
Show that $\iint_R K\le 2\cdot\pi$.

Conclude that any two closed geodesics on $\Sigma$ have a common point.
\end{thm}

\begin{thm}{Exercise}\label{ex:self-intersections}
Let $\Sigma$ be a closed smooth surface with positive Gauss curvature. 
Show that in a single coordinate chart, a closed geodesic in $\Sigma$ cannot look like one of the curves on the following pictures.

\begin{figure}[h]
\begin{minipage}{.32\textwidth}
\centering
\includegraphics{mppics/pic-46}
\end{minipage}
\hfill
\begin{minipage}{.32\textwidth}
\centering
\includegraphics{mppics/pic-47}
\end{minipage}
\hfill
\begin{minipage}{.32\textwidth}
\centering
\includegraphics{mppics/pic-471}
\end{minipage}

\medskip

\begin{minipage}{.32\textwidth}
\centering
\caption*{\textit{(easy)}}
\end{minipage}
\hfill
\begin{minipage}{.32\textwidth}
\centering
\caption*{\textit{(tricky)}}
\end{minipage}
\hfill
\begin{minipage}{.32\textwidth}
\centering
\caption*{\textit{(hopeless)}}
\end{minipage}
\vskip-5mm
\end{figure}

\end{thm}

The following exercise optimizes \ref{ex:ruf-bound-mountain}.

\begin{thm}{Exercise}\label{ex:sqrt(3)}
Suppose $f\:\mathbb{R}^2\to\mathbb{R}$ is a $\sqrt{3}$-Lipschitz smooth convex function.
Show that any geodesic in the graph $z=f(x,y)$ has no self-intersections.
\end{thm}

A surface $\Sigma$ is called \index{simply-connected surface}\emph{simply-connected} if any simple closed curve in $\Sigma$ bounds a disc.
Equivalently, any closed curve in $\Sigma$ can be continuously 
deformed into a \index{trivial curve}\emph{trivial curve}; that is, a curve that stands at one point all the time.

Planes and spheres are examples of simply-connected surfaces;
tori and the cylinders are not simply-connected.

\begin{thm}{Exercise}\label{ex:unique-geod}
Suppose $\Sigma$ is a simply-connected open surface with nonpositive Gauss curvature.

\begin{subthm}{ex:unique-geod:unique}
Show that any two points in $\Sigma$ are connected by a unique geodesic.
Moreover, geodesics in $\Sigma$ do not have self-intersections.
\end{subthm}

\begin{subthm}{ex:unique-geod:diffeomorphism}
Conclude that for any point $p\in \Sigma$,
the exponential map $\exp_p$ is a diffeomorphism from the tangent plane $\T_p$ to~$\Sigma$.
In particular, $\Sigma$ is diffeomorphic to the plane.
\end{subthm}

\end{thm}

\section{Additivity}

Let $\Delta$ be a topological disc in a smooth oriented surface $\Sigma$ bounded by a simple piecewise smooth curve $\partial \Delta$.
As before we suppose $\partial \Delta$ is oriented in such a way that $\Delta$ lies on its left.
Set \index{10gb@$\GB$ (Gauss--Bonnet formula)}
\[\GB(\Delta)
\df
\tgc{\partial\Delta}+\iint_\Delta K-2\cdot \pi,
\eqlbl{eq:GB}\]
where $K$ denotes the Gauss curvature.
The Gauss--Bonnet formula can be written as
\[\GB(\Delta)=0,\]
and $\GB$ stands for Gauss--Bonnet.

{

\begin{wrapfigure}{r}{40 mm}
\vskip-16mm
\centering
\includegraphics{mppics/pic-1750}
\end{wrapfigure}

\begin{thm}{Lemma}\label{lem:GB-sum}
Suppose the disc $\Delta$ is subdivided into two discs $\Delta_1$ and $\Delta_2$ by a curve $\delta$.
Then
\[
\GB(\Delta)=\GB(\Delta_1)+\GB(\Delta_2).
\]
\end{thm}

\parbf{Proof.}
Let us subdivide $\partial \Delta$ into two curves $\gamma_1$ and $\gamma_2$ that share endpoints with $\delta$, so that $\Delta_i$ is bounded by the arc $\gamma_i$  and~$\delta$ for $i=1,2$.

}

Denote by $\phi_1$, $\phi_2$, $\psi_1$, and $\psi_2$ the angles between $\delta$ and $\gamma_i$ marked on the picture.
Suppose the arcs $\gamma_1$, $\gamma_2$, and $\delta$ are oriented as on the picture. 
Then
\begin{align*}
\tgc{\partial \Delta}&= \tgc{\gamma_1}-\tgc{\gamma_2}+(\pi-\phi_1-\phi_2)+(\pi-\psi_1-\psi_2),
\\
\tgc{\partial \Delta_1}&= \tgc{\gamma_1}-\tgc\delta+(\pi-\phi_1)+(\pi-\psi_1),
\\
\tgc{\partial \Delta_2}&= \tgc\delta-\tgc{\gamma_2}+(\pi-\phi_2)+(\pi-\psi_2),
\\
\iint_\Delta K&=\iint_{\Delta_1} K+\iint_{\Delta_2} K.
\end{align*}
It remains to plug in the results in the formulas for $\GB(\Delta)$, $\GB(\Delta_1)$, and $\GB(\Delta_2)$.
\qeds

\section{Spherical case}

If $\Sigma$ is a plane, then its Gauss curvature vanishes;
therefore the Gauss--Bonnet formula \ref{eq:g-b} can be written as 
\[\tgc{\partial\Delta}=2\cdot \pi,\]
and it follows from \ref{prop:total-signed-curvature}.
In other words, $\GB(\Delta)=0$ for any plane disc~$\Delta$ with piecewise smooth boundary.

If $\Sigma$ is the unit sphere, then $K\equiv1$;
in this case, Theorem~\ref{thm:gb} can be formulated in the following way:

\begin{thm}{Proposition}\label{prop:area-of-spher-polygon}
Let $P$ be a spherical polygon bounded by a simple closed broken geodesic $\partial P$.
Assume $\partial P$ is oriented such that $P$ lies on the left of $\partial P$.
Then 
\[\GB(P)=\tgc{\partial P}+\area P-2\cdot \pi=0.\]

Moreover, the same formula holds for any spherical region bounded by a piecewise smooth simple closed curve.
\end{thm}

This proposition will be used in the next section.

\parbf{Sketch of proof.}
Suppose a spherical triangle $\Delta$ has angles 
$\alpha$, $\beta$, and~$\gamma$.
According to \ref{lem:area-spher-triangle},
\[\area\Delta=\alpha+\beta+\gamma-\pi.\]

Recall that $\partial\Delta$ is oriented so that $\Delta$ lies on its left. 
Then its oriented external angles are  $\pi-\alpha$, $\pi-\beta$, and $\pi-\gamma$.
Therefore,
\[\tgc{\partial\Delta}=3\cdot\pi-\alpha-\beta-\gamma.\]
It follows that $\tgc{\partial\Delta}+\area \Delta=2\cdot\pi$ or, equivalently, $\GB(\Delta)=0$.
 
We can subdivide a given spherical polygon $P$ into triangles by cutting a polygon in two along a broken geodesic on each step.
By the additivity lemma (\ref{lem:GB-sum}), 
\[\GB(P)=0\]
for any spherical polygon~$P$.

The second statement can be proved by approximation.
One has to show that the total geodesic curvature of  
a piecewise smooth simple curve can be approximated by 
the total geodesic curvature of inscribed broken geodesics.
We omit the detailed proof,
but it is close to \ref{ex:total-curvature=}.
\qeds

\begin{thm}{Exercise}\label{ex:half-sphere-total-curvature}
Let $\gamma$ be a simple piecewise smooth loop on the unit sphere $\mathbb{S}^2$.
Assume $\gamma$ divides $\mathbb{S}^2$ in two regions of equal area.
Denote by $p$ the base point of~$\gamma$.
Show that the parallel transport $\iota_\gamma\:\T_p\mathbb{S}^2\to\T_p\mathbb{S}^2$ is the identity map.
\end{thm}

\section{Intuitive proof}\label{sec:gb-intuitive-proof}

In this section, we prove a special case of the Gauss--Bonnet formula.
This case is leading --- the general case can be proved similarly, using the signed area counted with multiplicity, but will do something else in \ref{sec:gauss--bonnet:formal}.

\parbf{Proof of \ref{thm:gb} for open and closed surfaces with positive Gauss curvature.}
Let $\Norm\:\Sigma\to\mathbb{S}^2$ be the spherical map.
By \ref{cor:intK}, we have
\[\GB(\Delta)=\tgc{\partial\Delta}+\area(\Norm(\Delta))-2\cdot \pi.
\eqlbl{eq:gb-area}\]

Choose a loop $\alpha$ that runs along $\partial\Delta$ so that $\Delta$ lies on the left from it; suppose that $p\in \partial\Delta$ is its base.
Consider the parallel transport $\iota_\alpha\:\T_p\to\T_p$ along $\alpha$.
According to \ref{prop:pt+tgc}, $\iota_\alpha$ is the clockwise rotation by angle $\tgc{\alpha}_\Sigma$.

Set $\beta=\Norm\circ\alpha$.
By \ref{obs:parallel=}, we have $\iota_\alpha=\iota_\beta$, where $\beta$ is considered as a curve in the unit sphere.

Further, $\iota_\beta$ is a clockwise rotation by angle $\tgc{\beta}_{\mathbb{S}^2}$.
By \ref{prop:area-of-spher-polygon} 
\[\GB(\Norm(\Delta))=\tgc{\beta}_{\mathbb{S}^2}+\area(\Norm(\Delta))-2\cdot \pi=0.\]
Therefore, 
$\iota_\beta$ is a counterclockwise rotation by $\area(\Norm(\Delta))$

Summarizing, the clockwise rotation by $\tgc{\alpha}_\Sigma$ is identical to a counterclockwise rotation by $\area(\Norm(\Delta))$.
The rotations are identical if the angles are equal modulo $2\cdot\pi$.
Therefore, 
\[
\begin{aligned}
\GB(\Delta)&=\tgc{\partial\Delta}_\Sigma+\area(\Norm(\Delta))-2\cdot \pi=
2\cdot n \cdot \pi
\end{aligned}
\eqlbl{eq:sum=2pin}\]
for an integer~$n$.

It remains to show that $n=0$.
By \ref{prop:total-signed-curvature}, this is so for a topological disc in a plane. 
One can think of a general disc $\Delta$ as the result of a continuous deformation of a plane disc. 
The integer $n$ cannot change during the process since the left hand side in \ref{eq:sum=2pin} is continuous along the deformation; whence $n=0$ for the result of the deformation.
\qeds

\section{Simple geodesic}

The following theorem provides an interesting application of the Gauss--Bonnet formula; it is proved by Stephan Cohn-Vossen \cite[Satz 9]{convossen}.

\begin{thm}{Theorem}\label{thm:cohn-vossen}
Any open smooth surface with positive Gauss curvature has a simple two-sided infinite geodesic.
\end{thm}

\parbf{Proof.}
Let $\Sigma$ be an open surface with positive Gauss curvature and $\gamma$ a two-sided infinite geodesic in~$\Sigma$.

If $\gamma$ has a self-intersection, then it contains a simple loop;
that is, for some closed interval $[a,b]$,
the restriction $\ell=\gamma|_{[a,b]}$ is a simple loop.

By \ref{ex:convex-proper-plane}, $\Sigma$ is parametrized by an open convex region $\Omega$ in the plane.
By Jordan's theorem (\ref{thm:jordan}), $\ell$ bounds a topological disc in $\Sigma$; denote it by~$\Delta$.
If $\phi$ is the internal angle at the base of the loop, then by the Gauss--Bonnet formula,
\[\iint_\Delta K=\pi+\phi.\] 

Recall that
\[\iint_\Sigma K\le 2\cdot\pi;
\eqlbl{intK=<2pi+}\]
see \ref{ex:intK:2pi}.
Therefore, $0<\phi<\pi$; that is, $\gamma$ has no concave simple loops.

Assume $\gamma$ has two simple loops, say $\ell_1$ and $\ell_2$ that bound discs $\Delta_1$ and $\Delta_2$, respectively.
Then the discs $\Delta_1$ and $\Delta_2$ have to overlap;
otherwise, the curvature of $\Sigma$ would exceed $2\cdot\pi$  contradicting \ref{intK=<2pi+}.

It follows that after leaving $\Delta_1$, the geodesic $\gamma$ has to enter it again before creating another simple loop.
\begin{figure}[h!]
\vskip-0mm
\centering
\includegraphics{mppics/pic-1550}
\end{figure}
Consider the moment when $\gamma$ enters $\Delta_1$ again;
two possible scenarios are shown in the picture.
On the left picture, we get two nonoverlapping discs which, as we know, is impossible.
The right picture is impossible as well --- in this case, we get a concave simple loop.

It follows that $\gamma$ contains only one simple loop.
This loop cuts a disc from $\Sigma$ 
and goes around it either clockwise or counterclockwise.
This way we divide all the self-intersecting geodesics on $\Sigma$
into two types which we call {}\emph{clockwise} and {}\emph{counterclockwise}.

The geodesic $t\mapsto \gamma(t)$ is clockwise 
if and only if the same geodesic traveled backwards
$t\mapsto \gamma(-t)$
is counterclockwise.
Let us shoot a unit-speed geodesic in each direction from a given point $p=\gamma(0)$.
This gives a one-parameter family of geodesics $\gamma_s$ for $s\in[0,\pi]$ connecting the geodesic $t\mapsto \gamma(t)$ with $t\mapsto \gamma(-t)$; that is, $\gamma_0(t)\z=\gamma(t)$, and $\gamma_\pi(t)=\gamma(-t)$.

The subset of values $s\in [0,\pi]$ such that $\gamma_s$ is clockwise (or counterclockwise) is open in $[0,\pi]$.
That is, if $\gamma_s$ is clockwise, then so is $\gamma_t$ for all $t$ sufficiently close to~$s$.%
\footnote{Informally speaking, this means that self-intersection cannot suddenly disappear.
Try to convince yourself of this.}

Since $[0,\pi]$ is connected, it cannot be subdivided into two open sets.
It follows that for some $s$, the geodesic $\gamma_s$ is neither  clockwise nor counterclockwise;
that is, $\gamma_s$ has no self-intersections.
\qeds

{

\begin{wrapfigure}{r}{17 mm}
\vskip-4mm
\centering
\includegraphics{mppics/pic-1575}
\end{wrapfigure}

\begin{thm}{Exercise}\label{ex:cohn-vossen}
Let $\Sigma$ be an open smooth surface with positive Gauss curvature.
Suppose $\alpha\:[0,1]\z\to \Sigma$ is a smooth loop such that $\alpha'(0)=-\alpha'(1)$.
Show that there is a simple two-sided infinite geodesic $\gamma$ that is tangent to $\alpha$.
\end{thm}

}

{

\begin{wrapfigure}[3]{r}{43 mm}
\vskip-0mm
\centering
\includegraphics{mppics/pic-1577}
\end{wrapfigure}

\begin{thm}{Exercise}\label{ex:3-curves}
Which of the curves in the figure can occur as geodesics in charts of open smooth surfaces with positive Gauss curvature, and why?
\end{thm}

}

\section{General domains}
\index{Gauss--Bonnet formula}

The following generalization of the Gauss--Bonnet formula was obtained by Walther von Dyck \cite{dyck}.

\begin{thm}{Theorem}\label{thm:GB-generalized}
Let $\Lambda$ be a compact domain on a smooth surface.
Assume $\Lambda$ is bounded by a finite (possibly empty) collection of closed piecewise smooth curves $\gamma_1,\dots,\gamma_n$.
Suppose that each $\gamma_i$ 
is oriented in such a way that $\Lambda$ lies on its left.
Then
\[\iint_\Lambda K=2\cdot \pi\cdot \chi-\tgc{\gamma_1}-\dots-\tgc{\gamma_n}\eqlbl{eq:g-b++}\]
for an integer $\chi=\chi(\Lambda)$.

Moreover, if a graph with $v$ vertices and $e$ edges embedded in $\Lambda$ subdivides it into $f$ discs and contains all $\gamma_i$, then $\chi=v-e+f$.
\end{thm}

The number $\chi=\chi(\Lambda)$ is called the \index{Euler characteristic}\emph{Euler characteristic} of $\Lambda$. 
The value $\chi$ does not depend on the choice of the subdivision since the remaining terms in \ref{eq:g-b++} do not.
Formula \ref{eq:g-b++} is deduced from the standard Gauss--Bonnet (\ref{thm:gb}).
The argument is similar to \ref{lem:GB-sum}, but the combinatorics gets trickier.

\begin{wrapfigure}{o}{38 mm}
\vskip-0mm
\centering
\includegraphics{mppics/pic-1580}
\end{wrapfigure}

Before reading the proof, try to modify the argument in \ref{lem:GB-sum} to derive the formula for the subdivision of the annulus $A$ on the picture;
its graph has $4$ vertices and $6$ edges (one of them is a loop) and it subdivides $A$ into two discs $\Delta_1$ and $\Delta_2$.
Therefore, $\chi(A)=4-6+2=0$ and  
\[\iint_A K=-\tgc{\gamma_1}-\tgc{\gamma_2}.\]

\parbf{Proof.}
Suppose that a graph with $v$ vertices and $e$ edges subdivides $\Lambda$ into $f$ discs $\Delta_1,\dots,\Delta_f$.
Let us apply the Gauss--Bonnet formula for each disc, and sum up the results:
\[
\begin{aligned}
\iint_\Lambda K&=\iint_{\Delta_1} K+\dots+\iint_{\Delta_f} K=
2\cdot f\cdot \pi-\tgc{\partial\Delta_1}-\dots-\tgc{\partial\Delta_f}.
\end{aligned}
\]
It remains to show that  
\[\tgc{\gamma_1}+\dots+\tgc{\gamma_n}-\tgc{\partial\Delta_1}-\dots-\tgc{\partial\Delta_f}
=
2\cdot\pi\cdot(v-e).
\eqlbl{eq:GB-sum}\]
To prove this identity, we compute the left-hand side adding up the contributions of each edge and each vertex separately.

Choose an edge $\sigma$ of the graph.
If $\sigma$ is not part of some $\gamma_i$,
then it appears twice in the boundary of the discs, say in  $\partial \Delta_i$ and $ \partial \Delta_j$. 
In the latter case, we may assume that $\Delta_i$ lies on the left from $\sigma$ and $\Delta_j$ is on its right, so 
 $\sigma$ contributes $\tgc\sigma$ to $\tgc{\partial\Delta_i}$ and $-\tgc\sigma$ to $\tgc{\partial\Delta_j}$; hence $\sigma$ contributes nothing to the left-hand side of \ref{eq:GB-sum}.
It might happen that $i=j$ (as for one of the edges on the picture above);
in this case, we have one disc on both sides of~$\sigma$, but still, it contributes nothing to \ref{eq:GB-sum}.
If $\sigma$ is part of some $\gamma_i$, then it also appears once in the boundary of some disc, say in $\partial \Delta_j$.
We may assume that both $\Lambda$ and $\Delta_j$ lie on the left from $\sigma$,
so it contributes $\tgc\sigma$ to both $\tgc{\gamma_i}$ and $\tgc{\partial\Delta_j}$.
These contributions cancel each other on the left-hand side of  \ref{eq:GB-sum}.

Summarizing, the contributions to the left-hand side of  \ref{eq:GB-sum} that come from the total geodesic curvatures of edges cancel out.

Let us now handle the external angles at the vertices.
Choose a vertex $p$.
Denote by $d$ the \index{degree of vertex}\emph{degree} of~$p$;
that is, the number of edges that end at~$p$.

\begin{wrapfigure}{r}{23 mm}
\vskip-6mm
\centering
\includegraphics{mppics/pic-1585}
\end{wrapfigure}

Suppose $p$ lies in the interior of $\Lambda$.
Denote by $\delta_1,\dots,\delta_d$ the internal angles at $p$ of the discs containing $p$ as a vertex. 
Let  $\phi_{i}=\pi-\delta_{i}$ be the corresponding external angles.
Then $p$ contributes  
$-\phi_1-\dots-\phi_d$ to the sum.
Note that $\delta_1+\dots+\delta_d=2\cdot\pi$;
so $p$ contributes
\[-\phi_1-\dots-\phi_d = (\delta_1+\dots+\delta_d) - d\cdot \pi=(2-d)\cdot \pi,\]
to the left-hand side of \ref{eq:GB-sum}.

\begin{wrapfigure}{r}{23 mm}
\vskip-0mm
\centering
\includegraphics{mppics/pic-1590}
\end{wrapfigure}

If $p$ lies on the boundary of $\Lambda$, then it is a vertex of $d-1$ internal angles
$\delta_1,\dots,\delta_{d-1}$, and $\phi_{i}\z=\pi-\delta_{i}$ are the corresponding external angles.
Note that
\[\delta_1+\dots+\delta_{d-1}\z=\pi-\theta,\]
where $\theta\in(-\pi,\pi)$ denotes the external angle of $\Lambda$ at $p$.
And again, $p$ contributes
\[\theta-\sum\phi_{i}=(2-d)\cdot \pi\]
to the left-hand side of \ref{eq:GB-sum}.

Summarizing, suppose $p_1,\dots,p_v$ are the vertices of the graph, and $d_1,\dots,d_v$ are their corresponding degrees.
Then the total contribution from external angles to the left-hand side of  \ref{eq:GB-sum} is 
\[2\cdot v\cdot \pi-(d_1+\dots+d_v)\cdot\pi.
\eqlbl{eq:GB-sum-d}\]
Since the edges contributed nothing, the left-hand side of  \ref{eq:GB-sum} equals \ref{eq:GB-sum-d}.

It remains to observe that $d_1+\dots+d_v=2\cdot e$.
Indeed, $d_1+\dots+d_v$ is the total number of ends of all edges in the graph, and each of $e$ edges has exactly two ends. 
This finishes the proof of \ref{eq:GB-sum} and \ref{eq:g-b++}.
\qeds

\begin{thm}{Exercise}\label{ex:g-b-chi}
Find the integral of the Gauss curvature on each of the following surfaces:

\setlength{\columnseprule}{0.4pt}
\begin{multicols}{2}

\begin{subthm}{ex:g-b-chi:torus}
A torus.
\end{subthm}

\begin{Figure}
\vskip-0mm
\centering
\includegraphics{mppics/pic-1595}
\end{Figure}

\begin{subthm}{ex:g-b-chi:moebius}
A Möbius band with geodesic boundary.
\end{subthm}

\begin{Figure}
\vskip-0mm
\centering
\includegraphics{mppics/pic-1605}
\end{Figure}

\begin{subthm}{ex:g-b-chi:pair-of-pants}
A pair of pants with geodesic boundary components.
\end{subthm}
\begin{Figure}
\vskip-0mm
\centering
\includegraphics{mppics/pic-1600}
\end{Figure}

\begin{subthm}{ex:g-b-chi:two-handles}
A sphere with two handles.
\end{subthm}

\begin{Figure}
\vskip-0mm
\centering
\includegraphics{mppics/pic-1610}
\end{Figure}

\end{multicols}

\begin{subthm}{ex:g-b-chi:cylinder}
A cylinder such that its boundary curves have flat neighborhoods.
\begin{Figure}
\vskip-0mm
\centering
\includegraphics{mppics/pic-1620}
\end{Figure}
\end{subthm}

\end{thm}

%% file: polar-chart.tex
\chapter{Semigeodesic charts}
\label{chap:semigeodesic}

This chapter contains computational proofs of several statements discussed above.
They include 
the alternative definition of injectivity radius (\ref{prop:inj-rad}),
that shortest paths are geodesics (\ref{prop:gamma''}),
and the Gauss--Bonnet formula (\ref{thm:gb}).
In addition, we discuss intrinsic isometries between surfaces and prove Gauss' remarkable theorem, stating that \textit{Gauss curvature is an intrinsic invariant}.

\section{Polar coordinates}\label{sec:Polar coordinates}

The property of the exponential map in \ref{prop:exp} can be used to define \index{polar coordinates}\emph{polar coordinates} in a smooth surface.

Namely, let $p$ be a point on a smooth surface $\Sigma$ and
let $(r,\theta) \in \mathbb{R}_{\ge0} \times \mathbb{S}^1$ be polar coordinates  on the tangent plane $\T_p$.
If $\vec v\in \T_p$ has coordinates $(r,\theta)$,
then we say that $s(r,\theta)\z=\exp_p\vec v$ is the point in $\Sigma$ with  polar coordinates $(r,\theta)$.

There might be many or no geodesics from $p$ to a given point $x$;
so, the point $x$ might be expressed in polar coordinates in multiple ways or it might have no polar coordinates.
However, from Proposition~\ref{prop:exp}, we get the following.

\begin{thm}{Observation}\label{obs:polar}
Let $(r,\theta)\mapsto s(r,\theta)$ describe polar coordinates with respect to a point $p$ in a  smooth surface~$\Sigma$.
Then there is $r_0>0$ such that $s$ is  regular at any pair $(r,\theta)$ with $0<r<r_0$.

Moreover, if $0\le r_1,r_2<r_0$, then $s(r_1,\theta_1) \z= s(r_2,\theta_2)$ if and only if
$r_1=r_2=0$ or $r_1=r_2$ and $\theta_1\z=\theta_2+2\cdot n\cdot\pi$ for an integer~$n$.
\end{thm}

Further we will always assume that polar coordinates on a surface are defined for $r<r_0$;
therefore, they behave in the usual way.

The following statement will play an important role in the formal proof that shortest paths are geodesics; see Section~\ref{sec:proof-of-gamma''}.

\begin{thm}{Gauss lemma}\label{lem:palar-perp}
Let $(r,\theta)\mapsto s(r,\theta)$ be polar coordinates with respect to a point $p$ in a smooth surface.
Then
$s_\theta\perp s_r$
for any $r$ and~$\theta$.
\end{thm}

\parbf{Proof.}
Fix $\theta \in \mathbb{S}^1$.
By the definition of the exponential map, the curve $\gamma(t)\z=s(t,\theta)$ is a unit-speed geodesic that starts at $p$.
\begin{enumerate}[(i)]
\item Since $\gamma$ has unit speed, we have $|s_r|=|\gamma'|=1$.
In particular,
 \[
 \tfrac{\partial}{\partial \theta}
 \langle s_r,s_r\rangle=0.\]
\item Since $\gamma$ is a geodesic, we have $s_{rr}(r,\theta)=\gamma''(r)\perp\T_{\gamma(r)}$.
Therefore 
\[
\langle s_\theta, s_{rr}\rangle=0.\]
\end{enumerate}
It follows that
\[
\begin{aligned}
\tfrac{\partial}{\partial r}
\langle s_\theta, s_r\rangle
&=
\langle s_{\theta r},s_r\rangle
+
\cancel{\langle s_\theta,s_{rr}\rangle}=
\\
&=
\tfrac12
\cdot 
\tfrac{\partial}{\partial \theta}
\langle s_r, s_r\rangle=
\\
&=0.
\end{aligned}
\eqlbl{eq:<s',s'>'=0}
\]

Since $s(0,\theta)=p$ for any $\theta$,
we have
$s_\theta(0,\theta)=0$;
in particular,
$\langle s_\theta, s_r\rangle=0$
if $r=0$.
By \ref{eq:<s',s'>'=0}, the value 
$\langle  s_\theta, s_r\rangle$ does not depend on~$r$ for fixed $\theta$.
Therefore
\[\langle s_\theta, s_r\rangle=0\]
for any $r$ and $\theta$.
\qeds

\section[Shortest paths are geodesics: a formal proof]{Shortest paths are geodesics:
\\
a formal proof}

\label{sec:proof-of-gamma''}

In this section, we use the construction of polar coordinates and the Gauss lemma (\ref{lem:palar-perp}) to prove Proposition~\ref{prop:gamma''}.

\parbf{Proof of \ref{prop:gamma''}.}
Let $\gamma\:[0,\ell]\to\Sigma$ be a shortest path parametrized by arc-length, and $p=\gamma(0)$.
Suppose $\ell=\length\gamma$ is sufficiently small, so $\gamma$ can be described in the polar coordinates at $p$;
say as $\gamma(t)=s(r(t),\theta(t))$ for functions $t\mapsto \theta(t)$ and $t\mapsto r(t)$ with $r(0)=0$.

By the chain rule, we have
\[\gamma'= s_\theta\cdot \theta'+ s_r\cdot r'
\eqlbl{eq:chain(gamma)}\]
whenever the left-hand side is defined and $t>0$.
By the Gauss lemma \ref{lem:palar-perp}, $s_\theta\perp s_r$, and by definition of polar coordinates, $|s_r|=1$.
Therefore, \ref{eq:chain(gamma)} implies
\[|\gamma'(t)|\ge r'(t).\eqlbl{eq:|gamma'|=r'}\]
for any $t>0$ where $\gamma'(t)$ is defined.

Since $\gamma$ is parametrized by arc-length, we have
\[|\gamma(t_2)-\gamma(t_1)|\le |t_2-t_1|.\]
In particular, $\gamma$ is Lipschitz.
Therefore, by Rademacher's theorem (\ref{thm:rademacher}) the derivative $\gamma'$ is defined almost everywhere.
By~\ref{adex:integral-length:a}, we have that
\begin{align*}
\length\gamma&=\int_0^\ell|\gamma'(t)|\cdot dt\ge
\\
&\ge\int_0^\ell r'(t)\cdot dt=
\\
&=r(\ell).
\end{align*}

By the construction of polar coordinates, there is a geodesic of length $r(\ell)$ from $p=\gamma(0)$ to $q=\gamma(\ell)$.
Since $\gamma$ is a shortest path, we get that $r(\ell)=\ell$, and, moreover, $r(t)=t$ for any~$t$.
This equality holds if and only if we have equality in \ref{eq:|gamma'|=r'} for almost all~$t$.
The latter implies that $\gamma$ is a geodesic.

It remains to prove the partial converse.

Fix a point $p\in\Sigma$.
Let $\epsilon>0$ be as in \ref{prop:exp}.
Assume a geodesic $\gamma$ of length less than $\epsilon$ from $p$ to $q$ does not minimize the length between its endpoints.
Then there is a shortest path from $p$ to $q$, which is distinct from $\gamma$,
and becomes a geodesic when it is parametrized by its arc-length.
That is, there are two geodesics from $p$ to $q$ of length smaller than $\epsilon$.
In other words, there are two vectors ${\vec v},\vec w\in\T_p$ such that $|{\vec v}|<\epsilon$, $|\vec w|<\epsilon$, and 
$q=\exp_p\vec v\z=\exp_p\vec w$.
But according to \ref{prop:exp}, the exponential map $\T_p \to \Sigma$ is injective in the $\epsilon$-neighborhood of zero --- a contradiction.\qeds

\section{Gauss curvature}\label{sec:jacobi-formula}

Let $s$ be a smooth map from a (possibly unbounded) coordinate rectangle in the $(u,v)$-plane to a smooth surface~$\Sigma$.
The map $s$ is called \index{semigeodesic}\emph{semigeodesic} if, for any fixed $v$, the map $u\mapsto s(u,v)$ is a unit-speed geodesic and $s_u\perp s_v$ for any $(u,v)$.

According to the Gauss lemma (\ref{lem:palar-perp}), the polar coordinates on $\Sigma$ are described by a semigeodesic map.

Let $\Norm=\Norm(u,v)$ be the unit vector normal to $\Sigma$ at $s(u,v)$.
For each pair $(u,v)$, consider an orthonormal frame $\Norm$, $\vec u\z=s_u$, and $\vec v\z=\Norm\times \vec u$.
Since the vector $s_v(u,v)$ is tangent to $\Sigma$ at $s(u,v)$, we get $s_v\perp \vec u$ and $s_v\perp \Norm$.
Therefore, we have that $s_v=b\cdot\vec v$ for a smooth function $(u,v)\z\mapsto b(u,v)$.%
\footnote{For a fixed value $v_0$, the vector field $s_v=b\cdot\vec v$ describes the difference between $\gamma_0$ and an \textit{infinitesimally close} geodesic $\gamma_1\:u\mapsto s(u,v_1)$.
The fields with this property are called \index{Jacobi field}\emph{Jacobi fields} along $\gamma_0$.}

\begin{thm}{Proposition}\label{prop:jaccobi}
Suppose $(u,v)\mapsto s(u,v)$ is a semigeodesic map to a smooth surface $\Sigma$ and $\Norm$, $\vec u$, $\vec v$, and $b$ are as above.
Then 
\[b\cdot K+b_{uu}=0,\]
where $K=K(u,v)$ is the Gauss curvature of $\Sigma$ at  $s(u,v)$.

Moreover, 
\[
\langle\vec u_u,\vec u\rangle=
\langle\vec u_u,\vec v\rangle=
\langle\vec u_v,\vec u\rangle=0,
\quad\text{and}\quad
\langle\vec u_v,\vec v\rangle=b_u.
\]

\end{thm}

The proof is lengthy but straightforward.

\parbf{Proof.}
Let $\ell=\ell(u,v)$, $m=m(u,v)$, and $n=n(u,v)$ be the components of the matrix describing the shape operator in the frame $\vec u, \vec v$;
that is,
\[
\begin{aligned}
\Shape\vec u&=\ell\cdot \vec u+ m\cdot \vec v,
&
\Shape\vec v&=m\cdot \vec u+ n\cdot \vec v.
\end{aligned}
\eqlbl{eq:Shape(u,v)}
\]
Recall that (see Section~\ref{sec:More curvatures})
\[K=\ell\cdot n-m^2.\]

The proposition follows from the identities
\[
\begin{aligned}
\vec u_u&=\ell\cdot \Norm,
&
\vec u_v&=\phantom{-}b_u\cdot \vec v+b\cdot m\cdot\Norm,
\\
\vec v_u&=m\cdot \Norm,
&
\vec v_v&=-b_u\cdot \vec u+b\cdot n\cdot\Norm.
\end{aligned}
\eqlbl{eq:uu-vv}
\]
Indeed 
\begin{align*}
b\cdot K&=b\cdot (\ell\cdot n-m^2)=
\\
&=\langle\vec u_u,\vec v_v\rangle-\langle\vec u_v,\vec v_u\rangle=
\tag{by \ref{eq:uu-vv}}
\\
&= 
\left(
\tfrac{\partial}{\partial v}
\langle\vec u_u,\vec v\rangle
-
\cancel{\langle\vec u_{uv},\vec v\rangle}
\right)-
\left(
\tfrac{\partial}{\partial u}
\langle\vec u_v,\vec v\rangle
-
\cancel{\langle\vec u_{uv},\vec v\rangle}
\right)=
\\
&=0-b_{uu}.
\tag{by \ref{eq:uu-vv}}
\end{align*}
It remains to prove the four identities in \ref{eq:uu-vv}.

\parit{Derivation of $\vec u_u=\ell\cdot \Norm$.}
Since the frame $\Norm$, $\vec u$, and $\vec v$ is orthonormal, this vector identity can be rewritten as the following three scalar identities:
\[
\begin{aligned}
\langle\vec u_u,\vec u\rangle&=0,
&
\langle\vec u_u,\vec v\rangle&=0,
&
\langle\vec u_u,\Norm\rangle&=\ell.
\end{aligned}
\]
Since $u\mapsto s(u,v)$ is a geodesic we have that $\vec u_u=s_{uu}(u,v)\perp\T_{s(u,v)}$.
Hence, the first two identities follow.
Now, 
\begin{align*}
\langle\vec u_u,\Norm\rangle
&=\langle s_{uu},\Norm\rangle=
\tag{by \ref{thm:shape-chart}}
\\
&=    \langle \Shape s_u,s_u\rangle=
\tag{since $\vec u=s_u$}
\\
&=\langle \Shape \vec u,\vec u\rangle=
\tag{by \ref{eq:Shape(u,v)}}
\\
&=\ell.
\end{align*}

\parit{Derivation of $\vec u_v= b_u\cdot \vec v+b\cdot m\cdot\Norm$.}
This vector identity can be rewritten as the following three scalar identities:
\[
\begin{aligned}
\langle\vec u_v,\vec u\rangle&=0,
&
\langle\vec u_v,\vec v\rangle&=b_u,
&
\langle\vec u_v,\Norm\rangle&=b\cdot m.
\end{aligned}
\eqlbl{eq:uu-vv:2}
\]

Since $\langle\vec u,\vec u\rangle=1$, we get 
$0=\tfrac{\partial}{\partial v}\langle\vec u,\vec u\rangle=2\cdot\langle\vec u_v,\vec u\rangle$; 
hence the first identity in \ref{eq:uu-vv:2} follows.
Further, 
\begin{align*}
\langle\vec u_v,\vec v\rangle&=\langle s_{vu},\vec v\rangle=
\tag{since $s_v=b\cdot \vec v$}
\\
&=\langle \tfrac{\partial}{\partial u} (b\cdot\vec v),\vec v\rangle =
\\
&=b_u\cdot \langle \vec v,\vec v\rangle+b\cdot \langle \vec v_u,\vec v\rangle=
\tag*{(since $\langle \vec v,\vec v\rangle=1$ and $\qquad$} 
\\
&=b_u; 
\tag*{$0=\tfrac{\partial}{\partial u}\langle \vec v,\vec v\rangle=2\cdot\langle \vec v_u,\vec v\rangle$)}
\end{align*}
hence the second identity in \ref{eq:uu-vv:2} follows.
Finally, 
\begin{align*}
\langle\vec u_v,\Norm\rangle
&=\langle s_{uv},\Norm\rangle=
\tag{by \ref{thm:shape-chart}}
\\
&=\langle \Shape s_u,s_v\rangle=
\tag{since $\vec u=s_u$ and $s_v=b\cdot \vec v$}
\\
&=\langle \Shape \vec u,b\cdot \vec v\rangle=
\tag{by \ref{eq:Shape(u,v)}}
\\
&=b\cdot m.
\end{align*}

\parit{Derivation of $\vec v_u=m\cdot \Norm$ and $\vec v_v=-b_u\cdot \vec u+b\cdot n\cdot\Norm$.}
Recall that $\vec v=\Norm\times \vec u$.
Therefore,
\[
\vec v_u=\Norm_u\times \vec u+\Norm\times \vec u_u,
\qquad\qquad
\vec v_v=\Norm_v\times \vec u+\Norm\times \vec u_v.
\eqlbl{eq:uu-vv:3+4}
\]
The expressions for $\vec u_u$ and $\vec u_v$ in \ref{eq:uu-vv} are proved already.
Further,
\begin{align*}
-\Norm_u&=\Shape s_u=
&
-\Norm_v&=\Shape s_v=
\\
&=\Shape \vec u=
&
&=b\cdot\Shape \vec v=
\\
&=\ell\cdot\vec u+m\cdot\vec v,
&
&=b\cdot(m\cdot \vec u+ n\cdot \vec v),
\end{align*}
It remains to plug into \ref{eq:uu-vv:3+4} the expressions for $\vec u_u$, $\vec u_v$, $\Norm_u$, and~$\Norm_v$. 
\qeds

A chart $(u,v)\mapsto s(u,v)$, as well as the corresponding local coordinates, will be called \index{semigeodesic}\emph{semigeodesic} if the map $(u,v)\z\mapsto s(u,v)$ is semigeodesic.
Note that the function $b=b(u,v)$ for a semigeodesic chart $s$ has a constant sign.
Therefore, by changing the sign of $\Norm$, we can (and always will) assume that $b>0$;
in other words, $b=|s_v|$.

\begin{thm}{Exercise}\label{ex:semigeodesc-chart}
Show that any point $p$ in a smooth surface $\Sigma$ can be covered by a semigeodesic chart.
\end{thm}

\begin{thm}{Exercise}\label{ex:inj-rad}
Let $p$ be a point on a smooth surface~$\Sigma$.
Assume $\exp_p$ is injective in the ball $B=B(0,r_0)_{\T_p}$.
Suppose the semigeodesic map $(r,\theta)\mapsto s(r,\theta)$ describes polar coordinates with respect to $p$, and the function $(r,\theta)\mapsto b(r,\theta)$ is as above.

Prove the following statements:

\begin{subthm}{ex:inj-rad:sign}
$b(r,\theta)$ does not change its sign for $0\z\le r\z<r_0$.
\end{subthm}

\begin{subthm}{ex:inj-rad:0}
$b(r,\theta)\ne0$ if $0< r<r_0$.
\end{subthm}

\begin{subthm}{ex:inj-rad:prop:inj-rad}
Apply \ref{SHORT.ex:inj-rad:sign} and \ref{SHORT.ex:inj-rad:0} to prove \ref{prop:inj-rad}.
\end{subthm}
 
\end{thm}

A chart $(u,v)\mapsto s(u,v)$ is called \index{orthogonal chart}\emph{orthogonal} if $s_u\perp s_v$ for any $(u,v)$.
For example, any semigeodesic chart is orthogonal.

A solution of the following exercise is similar to \ref{prop:jaccobi}.

\begin{thm}{Exercise}\label{lem:K(orthogonal)}
Let $(u,v)\mapsto s(u,v)$ be an orthogonal chart of a smooth surface~$\Sigma$.
Denote by $K=K(u,v)$ the Gauss curvature of $\Sigma$ at $s(u,v)$.
Set 
\begin{align*}
a=a(u,v)&\df|s_u|,&
b=b(u,v)&\df|s_v|,\\
\vec u=\vec u(u,v)&\df\tfrac{s_u}a,&
\vec v=\vec v(u,v)&\df\tfrac{s_v}b.
\end{align*}
Let $\Norm=\Norm(u,v)$ be the unit normal vector at $s(u,v)$.

\begin{subthm}{lem:K(orthogonal):uu-vv}
Show that 
\begin{align*}
\vec u_u
&=
-\tfrac1{b}\cdot a_v
\cdot
\vec v 
+
a\cdot \ell\cdot \Norm
,
&
\vec v_u
&=
\tfrac1{b}\cdot a_v
\cdot \vec u
+
a\cdot m\cdot \Norm
\\
\vec u_v
&=
\tfrac1{a}\cdot b_u\cdot\vec v
+
b\cdot m\cdot \Norm
,
&
\vec v_v
&=
-\tfrac1{a}\cdot b_u\cdot\vec u
+
b\cdot n\cdot \Norm,
\end{align*}
where $\ell=\ell(u,v)$, $m=m(u,v)$, and $n=n(u,v)$ are the components of the matrix describing the shape operator in the frame $\vec u, \vec v$.
\end{subthm}

\begin{subthm}{lem:K(orthogonal):K}
Show that
\[K=-\frac1{a\cdot b}\cdot
\left(
\frac{\partial}{\partial u}
\left(\frac{b_u}a \right)
+
\frac{\partial}{\partial v}
\left(\frac{a_v}b\right)
\right).\]
\end{subthm}
\end{thm}

\begin{thm}{Exercise}\label{ex:conformal}
Suppose $(u,v)\mapsto s(u,v)$ is a \index{conformal}\emph{conformal chart};
that is, there is a function $(u,v)\mapsto b(u,v)$ such that $b=|s_u|=|s_v|$ and $s_u\perp s_v$ for any $(u,v)$.
(The function $b$ is called a {}\emph{conformal factor} of the chart.)

Use \ref{lem:K(orthogonal)} to show that  
\[K=-\frac{\triangle (\ln b)}{b^2},\]
where $\triangle$ denotes the \index{Laplacian}\emph{Laplacian}; that is, $\triangle=\tfrac{\partial^2}{\partial u^2}+\tfrac{\partial^2}{\partial v^2}$, and 
 $K=K(u,v)$ is the Gauss curvature of $\Sigma$ at $s(u,v)$.
\end{thm}

It is useful to know that \textit{any point on a smooth surface can be covered by a conformal chart}; the corresponding coordinates are usually referred to as \index{isothermal}\emph{isothermal}.

\section{Rotation of a vector field}

Let $\gamma\:[0,1]\to\Sigma$ be a simple loop on a smooth oriented surface~$\Sigma$.
Suppose $\vec u$ is a field of unit tangent vectors on $\Sigma$ defined in a neighborhood of~$\gamma$.
Denote by $\vec v$ the field obtained from $\vec u$ by a counterclockwise rotation by $\tfrac{\pi}2$ of the tangent plane at each point; it could also be defined by $\vec v\df\Norm\times\vec u$, where $\Norm$ is the normal field on $\Sigma$.
Then the \index{rotation}\emph{rotation} of $\vec u$ around $\gamma$ is defined as the integral
\[\rot_\gamma\vec u
\df
\int_0^1\langle\vec u'(t),\vec v(t)\rangle\cdot dt.\]

\begin{thm}{Lemma}\label{lem:rotation-parallel}
Suppose $\gamma\:[0,1]\to\Sigma$ is a simple loop with the base at a point $p$ in a smooth oriented surface $\Sigma$, and $\vec u$ is a field of tangent unit vectors to $\Sigma$ defined in a neighborhood of~$\gamma$.
Then the parallel transport $\iota_\gamma\:\T_p\to\T_p$ is a clockwise rotation by the angle $\rot_\gamma\vec u$.

In particular, rotations of different vector fields around $\gamma$ may only differ by a multiple of $2\cdot\pi$.
\end{thm}

\parbf{Proof.}
As above, set $\vec v=\Norm\times\vec u$. 
Denote by $\vec u(t)$ and $\vec v(t)$ the vectors at $\gamma(t)$ of the fields $\vec u$ and $\vec v$, respectively.

Let $t\mapsto \vec x(t)\in \T_{\gamma(t)}$ be the parallel vector field along $\gamma$ with  $\vec x(0)\z=\vec u(0)$, and $\vec y\df\Norm\times\vec x$.

Note that there is a continuous function $t\mapsto \phi(t)$ such that 
$\vec u(t)$ is a counterclockwise rotation of $\vec x(0)$ by angle $\phi(t)$.
Since $\vec x(0)=\vec u(0)$, we can and will assume that $\phi(0)=0$.
Then
\begin{align*}
\vec u&=\cos\phi\cdot \vec x+\sin\phi\cdot \vec y
\\
\vec v&=-\sin\phi\cdot \vec x+\cos\phi\cdot \vec y
\end{align*}
It follows that 
\begin{align*}
\langle\vec u',\vec v\rangle
=\phi'\cdot\biggl(&(\sin \phi)^2\cdot \langle\vec x,\vec x\rangle+(\cos \phi)^2\cdot \langle\vec y,\vec y\rangle
\biggr)=
\phi'.
\end{align*}

Therefore,
\begin{align*}
\rot_\gamma\vec u&=\int_0^1\langle\vec u'(t),\vec v(t)\rangle\cdot dt=
\\
&=\int_0^1\phi'(t)\cdot dt=
\\
&=\phi(1).
\end{align*}

Observe that 
\begin{itemize}
\item $\iota_\gamma(\vec x(0))=\vec x(1)$,

\item  $\vec x (0) = \vec u (0) = \vec u (1),$ 

\item $\vec u(1)$ is a counterclockwise rotation of $\vec x(1)$ by the angle $\phi(1)\z=\rot_\gamma\vec u$,

\end{itemize}
It follows that $\vec x(1)$ is a \textit{clockwise} rotation of $\vec x(0)$ by angle $\rot_\gamma\vec u$, and the result follows.

The last statement follows from \ref{prop:pt+tgc}.
\qeds

The following lemma will play a key role in the proof of the Gauss--Bonnet formula.

\begin{thm}{Lemma}\label{lem:rotation-semigeoesic}
Let $(u,v)\mapsto s(u,v)$ be a semigeodesic chart on a smooth surface~$\Sigma$.
Suppose a simple loop $\gamma$ bounds a disc $\Delta$ that is covered completely by~$s$.
Then 
\[\rot_\gamma\vec u+\iint_\Delta K=0,\]
where $\vec u=s_u$, and $K$ denotes the Gauss curvature of~$\Sigma$.
\end{thm}

The calculations below use the so-called \index{Green formula}\emph{Green formula}, which can be formulated the following way.

\textit{Let $D$ be a compact region in the $(u,v)$-coordinate plane that is bounded by a piecewise smooth simple closed path $\alpha\:t\mapsto (u(t),v(t))$.
Suppose $\alpha$ is oriented in such a way that $D$ lies on its left.
Then}
\[\iint_D (Q_u- P_v)\cdot du\cdot dv=\int_\alpha (P\cdot du+Q\cdot dv)\df \int_0^1 (P\cdot u'+Q\cdot v')\cdot dt\]
\textit{for any two smooth functions $P$ and $Q$ defined on $D$.}

The Green and Gauss--Bonnet formulas are similar --- they relate the integral along a disc and its boundary curve.
So it shouldn't be  surprising that one helps to prove another.

\parbf{Proof.}
Let $\vec u$, $\vec v$, and $b$ be as in Section~\ref{sec:jacobi-formula}.
Let us write $\gamma$ in the $(u,v)$-coordinates: $\gamma(t)=s(u(t),v(t))$.
Then
\begin{align*}
\rot_\gamma \vec u&=\int_0^1\langle\vec u',\vec v\rangle\cdot dt=
\tag{by the chain rule}
\\
&=\int_0^1[\langle\vec u_u,\vec v\rangle\cdot u'+\langle\vec u_v,\vec v\rangle\cdot v']\cdot dt=\tag{by \ref{prop:jaccobi}}
\\
&=\int_0^1b_u\cdot v'\cdot dt=\int_{s^{-1}\circ\gamma}b_u\cdot dv=
\tag{by the Green formula}
\\
&=\iint_{s^{-1}(\Delta)}b_{uu}\cdot du\cdot dv=
\tag{$s_u\perp s_v \Longrightarrow \jac s=|s_u|\cdot|s_v|=b$}
\\
&=\iint_\Delta\frac{b_{uu}}{b}=
\tag{by \ref{prop:jaccobi}}
\\
&=-\iint_{\Delta}K.
\end{align*}
\qedsf

\section{Gauss--Bonnet formula: a formal proof}\label{sec:gauss--bonnet:formal}

Recall that the Gauss--Bonnet formula can be written as $\GB(\Delta)=0$,
where 
\[\GB(\Delta)
\df
\tgc{\partial\Delta}+\iint_\Delta K-2\cdot \pi,\]
$\Delta$ is a topological disc in a smooth oriented surface,
bounded by a piecewise smooth curve $\partial \Delta$ that is oriented in such a way that $\Delta$ lies on its left.

\parbf{Proof of the Gauss--Bonnet formula (\ref{thm:gb}).}
Assume $\Delta$ is covered by a semigeodesic chart.
From \ref{prop:pt+tgc},
\ref{lem:rotation-parallel},
and \ref{lem:rotation-semigeoesic} we get
\[\GB(\Delta)
=
2\cdot n\cdot \pi,
\eqlbl{eq:gb(n)}\]
where $n=n(\Delta)$ is an integer.

By \ref{ex:semigeodesc-chart}, any point can be covered by a semigeodesic chart.
Therefore, applying the additivity of $\GB$ (\ref{lem:GB-sum}) a finite number of times, we get
\ref{eq:gb(n)} for any disc $\Delta$ in~$\Sigma$.
More precisely, we can cut $\Delta$ by a smooth curve that runs from boundary to boundary
and repeat such subdivision recursively for the obtained discs;
see the picture.
\begin{figure}[!ht]
\vskip-0mm
\centering
\includegraphics{mppics/pic-1700}
\vskip-0mm
\end{figure}
After several such steps, every small disc is covered by a semigeodesic chart.
In particular, \ref{eq:gb(n)} holds for each small disc.
Then applying \ref{lem:GB-sum} several times we get \ref{eq:gb(n)} for the original disc.

It remains to show that $n=0$.
Assume $\Delta$ lies in a local graph realization $z = f(x,y)$ of $\Sigma$.
Consider the one-parameter family $\Sigma_t$ of graphs $z = t \cdot f(x,y)$;
denote by $\Delta_t$ the corresponding disc in $\Sigma_t$, so $\Delta_1 = \Delta$ and $\Delta_0$ is its projection onto the $(x,y)$-plane.
The function $h \: t \mapsto \GB (\Delta_t )$ is continuous.
By \ref{eq:gb(n)}, $h(t)$ is an integer multiple of  $2 \cdot \pi$ for any $t$.
It follows that $h$ is a constant function.
Therefore, 
\[\GB (\Delta)=\GB(\Delta_0)=0;\]
the last equality follows from \ref{prop:total-signed-curvature}.

We proved that 
\[\GB (\Delta)=0
\eqlbl{eq:GB=0}\]
if $\Delta$ lies in a graph $z = f(x,y)$ for some $(x,y,z)$-coordinate system. Any point of $\Sigma$ has a neighborhood that can be covered by such a graph. Therefore, applying Lemma~\ref{lem:GB-sum} the same way as above, we get that \ref{eq:GB=0} holds for any disc~$\Delta$ in~$\Sigma$.
\qeds

\section{Rauch comparison}

The following proposition is a special case of the so-called \index{Rauch comparison theorem}\emph{Rauch comparison theorem}.

\begin{thm}{Proposition}\label{prop:rauch}
Suppose $p$ is a point on a smooth surface $\Sigma$ and $r\le \inj(p)$.
Given a curve $\tilde\gamma$ in the $r$-neighborhood of $0$ in $\T_p$, set 
\[\gamma=\exp_p\circ\tilde\gamma
\quad
\text{or, equivalently}
\quad
\log_p\circ\gamma=\tilde\gamma
;\]
note that $\gamma$ is a curve in~$\Sigma$.

\begin{subthm}{prop:rauch:K=<0}
If $\Sigma$ has nonpositive Gauss curvature, then the logarithmic map $\log_p$ is length nonexpanding in the $r$-neighborhood of $p$ in $\Sigma$;
that is, 
\[\length \gamma\ge \length \tilde\gamma\]
for any curve $\gamma$ in the open ball $B(p,r)_{\Sigma}$.
\end{subthm}

\begin{subthm}{prop:rauch:K>=0}
If $\Sigma$ has nonnegative Gauss curvature, then the exponential map $\exp_p$ is length nonexpanding in the $r$-neighborhood of $0$ in $\T_p$;
that is, 
\[\length \gamma\le \length \tilde\gamma\]
for any curve $\tilde\gamma$ in the open ball $B(0,r)_{\T_p}$.
\end{subthm}

\end{thm}

\parbf{Proof.}
Suppose $(r(t),\theta(t))$ are the polar coordinates of $\tilde\gamma(t)$ and $s$ as in \ref{sec:Polar coordinates}.
Note that $\gamma(t)\z=s(r(t),\theta(t))$; that is, $(r(t),\theta(t))$ are the polar coordinates of $\gamma(t)$ with respect to $p$ on~$\Sigma$.

Let us show that $b(0,\theta)=0$, and $b_r(0,\theta)=1$ for any $\theta$.
We may assume that $\theta=0$.
Choose the standard basis $\vec v,\vec w$ in $\T_p$; so $\vec v$ points in the direction of curve $t\mapsto s(t,0)$.
Note that%
\footnote{Recall that $(D_{\vec w}\exp_p)(r\cdot \vec v)\df h'(0)$, where $h(t)=\exp_p(r\cdot \vec v+t\cdot \vec w)$; see Section~\ref{sec:dirder}.}
\[b(r,0)=r\cdot |(D_{\vec w}\exp_p)(r\cdot \vec v)|.\]
In particular, $b(0,0)=0$.
By \ref{obs:d(exp)=1}, $|(D_{\vec w}\exp_p)(0)|=1$.
In particular, $|(D_{\vec w}\exp_p)(r\cdot \vec v)|\z\ne 0$ for any small $r$.
Therefore, the function $r\mapsto|(D_{\vec w}\exp_p)(r\cdot \vec v)|$ is differentiable at $r=0$.
Taking the partial derivative of the expression for $b(r,0)$, we get $b_r(0,0)=1$.

Set $b(r,\theta)\df|s_\theta|$.
By \ref{prop:jaccobi}
\[b_{rr}=-K\cdot b.\]
If $K\ge 0$, then $r\mapsto b(r,\theta)$ is concave
and
if $K\le 0$, then $r\mapsto b(r,\theta)$ is convex for any fixed $\theta$.
Since $b(0,\theta)=0$, and $b_r(0,\theta)=1$,
\[
\begin{aligned}
b(r,\theta)\ge r\quad\text{if}\quad K&\le 0;
\\
b(r,\theta)\le r\quad\text{if}\quad K&\ge 0.
\end{aligned}
\eqlbl{eq:b-K}
\]

Without loss of generality, we may assume that $\tilde\gamma\:[a,b]\to \T_p$ is parametrized by arc-length;
in particular, it is a Lipschitz curve.
Note that
\begin{align*}
\length\tilde\gamma&=\int_a^b\sqrt{r'(t)^2+r(t)^2\cdot\theta'(t)^2}\cdot dt.
\shortintertext{Applying \ref{lem:palar-perp}, we get}
\length\gamma&=\int_a^b\sqrt{r'(t)^2+b(r(t),\theta(t))^2\cdot\theta'(t)^2}\cdot dt.
\end{align*}
Now, both statements \ref{SHORT.prop:rauch:K=<0} and \ref{SHORT.prop:rauch:K>=0} follow from \ref{eq:b-K}.
\qeds

\section{Intrinsic isometries}

Let $\Sigma$ and $\Sigma^{*}$ be two smooth surfaces.
A map $f\:\Sigma\to \Sigma^{*}$ is called \index{length-preserving}\emph{length-preserving} if, for any curve $\gamma$ in $\Sigma$, the curve $\gamma^{*}=f\circ\gamma$ in $\Sigma^{*}$ has the same length. 
If in addition $f$ is smooth and bijective, then it is called an  \index{intrinsic!isometry}\emph{intrinsic isometry}. 

An example of a length-preserving map can be obtained by wrapping a plane into a cylinder with the map $s\:\mathbb{R}^2\to\mathbb{R}^3$ given by 
$s(x,y)\z=(\cos x,\sin x,y)$.

\begin{thm}{Exercise}\label{ex:isom(geod)}
Show that an intrinsic isometry between smooth surfaces maps geodesics to geodesics.
\end{thm}

\begin{thm}{Exercise}\label{ex:K=0}
Suppose the Gauss curvature of a smooth surface $\Sigma$ vanishes.
Show that $\Sigma$ is \index{locally flat surface}\emph{locally flat};
that is, a neighborhood of any point in $\Sigma$ admits an intrinsic isometry to an open domain in the Euclidean plane.  
\end{thm}

\begin{thm}{Exercise}\label{ex:K=1}
Suppose a smooth surface $\Sigma$ has Gauss curvature  equal to 1 at every point.
Show that a neighborhood of any point in $\Sigma$ admits an intrinsic isometry to an open domain in the unit sphere.
\end{thm}

\begin{thm}{Exercise}\label{ex:deformation}
Given $a>0$, show that there is a smooth unit-speed curve 
$\gamma(t)=(x(t),y(t))$ with $y(t) = a\cdot \cos t$ and $y>0$.
Describe its interval of definition.

Let $\Sigma_a$ be the surface of revolution of $\gamma$ around the $x$-axis.
\begin{figure}[h!]
\vskip-0mm
\centering
\begin{lpic}[t(-0mm),b(6mm),r(0mm),l(0mm)]{asy/deformation(1.2)}
\lbl[t]{8,-.5;$a=2$}
\lbl[t]{24,3;$a=\sqrt{2}$}
\lbl[t]{41,4;$a=1$}
\lbl[t]{57,7;$a=\tfrac1{\sqrt{2}}$}
\lbl[t]{73,8;$a=\tfrac12$}
\end{lpic}
\vskip-0mm
\end{figure}
Show that the surface $\Sigma_a$ has a unit Gauss curvature at each point.

Use \ref{ex:K=1} to conclude that any small round disc $\Delta$ in $\mathbb{S}^2$ admits a smooth length-preserving deformation;
that is, there is a one-parameter family $\Delta_t$ of surfaces with boundary, such that $\Delta_0=\Delta$ and $\Delta_t$ is not congruent to $\Delta_0$ for $t\ne0$.%
\footnote{In fact, any disc in $\mathbb{S}^2$ admits a smooth length preserving deformation.
However, if the disc is larger than half-sphere, then the proof requires more;
it can be obtained as a corollary of two deep results of Alexandr Alexandrov: the gluing theorem and the theorem on the existence of a convex surface with an abstractly given metric \cite[p. 44]{pogorelov}.
}
\end{thm}

The following exercise illustrates the final step in the proof that \textit{any open surface with vanishing Gauss curvature a cylindrical surface}.
See the discussion after \ref{ex:flat-plane}.

\begin{thm}{Advanced exercise}\label{ex:line-cylinder} 
Suppose $(u,v)\mapsto f(u,v)$ describes an intrinsic isometry from the plane with standard $(u,v)$-coordinates to a surface $\Sigma$ in $\mathbb{R}^3$.
Assume that $f$ maps $v$-axis to the $z$-axis isometrically.
Show that $\Sigma$ is a cylindrical surface;
more precisely, it is a union of a family of lines parallel to the $z$-axis.
\end{thm}

\section{The remarkable theorem}

\begin{thm}{Theorem}\label{thm:remarkable}
Suppose $f\:\Sigma\to \Sigma^{*}$ is an intrinsic isometry between two smooth surfaces; $p\in \Sigma$ and $p^{*}\z=f(p)\in \Sigma^{*}$.
Then 
\[K(p)_{\Sigma}=K(p^{*})_{\Sigma^{*}};\]
that is, the Gauss curvature of $\Sigma$ at $p$ is the same as the Gauss curvature of $\Sigma^{*}$ at $p^{*}$.
\end{thm}

Recall that the Gauss curvature is defined as a product of principal curvatures which might be different at the points $p$ and $p^*$; however, according to the theorem, their products are the same.
In other words, the Gaussian curvature is an \textit{intrinsic invariant}.
This theorem was proved by Carl Friedrich Gauss \cite{gauss} who justifiably called it {}\emph{remarkable} ({}\emph{Theorema Egregium}).

In fact, the Gauss curvature can be found \textit{intrinsically},
by measuring the lengths of curves in the surface.
For example, the Gauss curvature $K(p)$ appears in the following formula for the circumference $c(r)$ of a geodesic circle centered at $p$ in a surface: 
\[c(r)=2\cdot\pi\cdot r-\tfrac\pi3\cdot K(p)\cdot r^3+o(r^3).\]

The theorem implies, for example, that there is no smooth length-preserving map that sends an open region in the unit sphere to the plane.%
\footnote{Smoothness is essential --- there are plenty of non-smooth length-preserving maps from the sphere to the plane; see \cite{petrunin-yashinski} and the references therein.}
This follows since the Gauss curvature of the plane is zero and the unit sphere has Gauss curvature~1. 
In particular, any geographic map has to have distortions.

\parbf{Proof.}
Choose a chart $(u,v)\mapsto s(u,v)$ on $\Sigma$, and set
$s^{*} =f\circ s$.
Note that $s^{*}$ is a chart of $\Sigma^{*}$, and 
\begin{align*}
\langle s_u,s_u\rangle
&=
\langle s_u^{*}, s_u^{*}\rangle,
&
\langle s_u, s_v\rangle
&=
\langle s_u^{*}, s_v^{*}\rangle,
&
\langle s_v, s_v\rangle
&=
\langle s_v^{*}, s_v^{*}\rangle
\end{align*}
at any $(u,v)$.
Indeed, since $f$ preserves the lengths of the coordinate lines $\gamma\:t\mapsto s(t,v)$ and  $\gamma\:t\z\mapsto s(u,t)$, we get the first and the third identities.
Now, since $f$ preserves the lengths of the curves $\gamma\:t\z\mapsto s(t,c-t)$ for any constant~$c$, the first and the third identities imply the second one.

By \ref{prop:gamma''}, if $s$ is a semigeodesic chart, then so is $s^{*}$.
It remains to apply \ref{prop:jaccobi} and \ref{ex:semigeodesc-chart}.
\qeds

%% file: comparison.tex
\chapter{Comparison theorems}
\label{chap:comparison}

The comparison theorems provide a powerful means of applying Euclidean intuition in differential geometry.

\section{Triangles and hinges}

Recall that a shortest path between points $x$ and $y$ in a surface $\Sigma$ is denoted by $[x,y]$ or $[x,y]_\Sigma$, and
$\dist{x}{y}\Sigma$ denotes the \index{intrinsic!distance}\emph{intrinsic distance} from $x$ to $y$ in~$\Sigma$.

A \index{geodesic!triangle}\emph{geodesic triangle} in a surface $\Sigma$ is defined as a triple of points $x,y,z\z\in \Sigma$ with a choice of minimizing geodesics $[x,y]_\Sigma$, $[y,z]_\Sigma$, and $[z,x]_\Sigma$.
The points $x,y,z$ are called the {}\emph{vertices} of the triangle,
the minimizing geodesics $[x,y]_\Sigma$, $[y,z]_\Sigma$, and $[z,x]_\Sigma$ are called its {}\emph{sides};
the triangle itself is denoted by $[xyz]$, or $[xyz]_\Sigma$;
the latter notation is used if we need to emphasize that the triangle lies on the surface~$\Sigma$.\index{10aad@$[xyz]$, $[xyz]_\Sigma$ (geodesic triangle)}

A triangle $[\tilde x\tilde y\tilde z]$ in the plane $\mathbb{R}^2$ is called a \index{model!angle and triangle}\emph{model triangle} of the triangle $[xyz]$
if its corresponding sides are equal;
that is,
\[\dist{\tilde x}{\tilde y}{\mathbb{R}^2}=\dist{x}{y}\Sigma,
\quad
\dist{\tilde y}{\tilde z}{\mathbb{R}^2}=\dist{y}{z}\Sigma,
\quad
\dist{\tilde z}{\tilde x}{\mathbb{R}^2}=\dist{z}{x}\Sigma.
\]
In this case, we write $[\tilde x\tilde y\tilde z]=\tilde\triangle xyz$.
\index{10aae@$\tilde\triangle$ (model triangle)}

A pair of minimizing geodesics $[x,y]$ and $[x,z]$ starting from one point $x$ is called a \index{hinge}\emph{hinge} and is denoted by $\hinge xyz$.\index{10aac@$\hinge yxz$ (hinge)}
The angle between these geodesics at $x$ is denoted by $\measuredangle\hinge xyz$.
The corresponding angle $\measuredangle\hinge {\tilde x}{\tilde y}{\tilde z}$ in a model triangle $[\tilde x\tilde y\tilde z]=\tilde\triangle xyz$ is denoted by $\modangle xyz$;
\index{10aab@$\modangle yxz$ (model angle)}
it is called {}\emph{model angle} of the triangle $[xyz]$ at~$x$.

By side-side-side congruence condition,
the model triangle $[\tilde x\tilde y\tilde z]$ is uniquely defined up to congruence.
Therefore, the model angle $\tilde\theta=\modangle xyz$ is uniquely defined as well.
By the cosine rule, we have 
\[\cos \tilde\theta=\frac{a^2+b^2-c^2}{2\cdot a \cdot b},\]
where $a=\dist{x}{y}{\Sigma}$, $b=\dist{x}{z}{\Sigma}$, and $c=\dist{y}{z}{\Sigma}$.

\begin{thm}{Exercise}\label{ex:wide-hinges}
Let $[x_ny_nz_n]$ be a sequence of triangles in a smooth surface~$\Sigma$.
Set $a_n=\dist{x_n}{y_n}{\Sigma}$,
$b_n=\dist{x_n}{z_n}{\Sigma}$,
and $c_n=\dist{y_n}{z_n}{\Sigma}$, and $\tilde\theta_n=\modangle {x_n}{y_n}{z_n}$.
Suppose the sequences $a_n$ and $b_n$ are bounded away from zero;
that is, $a_n>\epsilon$ and $b_n>\epsilon$ for a fixed $\epsilon>0$ and any~$n$.
Show that
\[(a_n+b_n-c_n)\to 0\qquad\iff\qquad \tilde\theta_n\to \pi\]
as $n\to\infty$
\end{thm}

\section{Formulations}

Part \ref{SHORT.thm:comp:cat} of the following theorem is called the {}\emph{Cartan--Hadamard theorem};
it was proved by 
Hans von Mangoldt \cite{mangoldt} and generalized by 
Elie Cartan \cite{cartan} and
Jacques Hadamard \cite{hadamard}.
Part \ref{SHORT.thm:comp:toponogov} is called the {}\emph{Toponogov comparison theorem} and sometimes the {}\emph{Alexandrov comparison theorem};
it was proved by Paolo Pizzetti \cite{pizzetti}, rediscovered by Alexandr Alexandrov \cite{aleksandrov}, and 
generalized by Victor Toponogov \cite{toponogov1957}.%, Mikhael Gromov, Yuri Burago and Grigory Perelman~\cite{burago-gromov-perelman}.

Recall that a surface $\Sigma$ is called \index{simply-connected surface}\emph{simply-connected} if any simple closed curve in $\Sigma$ bounds a disc.

\begin{thm}{Comparison theorems}
\label{thm:comp}
\index{comparison theorem}
\ 

\begin{subthm}{thm:comp:cat}
If $\Sigma$ is an open simply-connected  smooth surface with nonpositive Gauss curvature, then 
\[\measuredangle\hinge {x}{y}{z}\le\modangle xyz\]
for any geodesic triangle $[xyz]$ in $\Sigma$.
\end{subthm}

\begin{subthm}{thm:comp:toponogov}
If $\Sigma$ is a closed (or open) smooth surface with nonnegative Gauss curvature, then 
 \[\measuredangle\hinge {x}{y}{z}\ge\modangle xyz\]
for any geodesic triangle $[xyz]$.
\end{subthm}

\end{thm}

The proofs of parts \ref{SHORT.thm:comp:cat} and \ref{SHORT.thm:comp:toponogov} will be given in sections~\ref{sec:nonpos-comp} and~\ref{sec:nonneg-comp}, respectively.

\begin{thm}{Exercise}\label{ex:thm:comp:cat:nsc}
Show that the conclusion in \ref{thm:comp:cat} does not hold in the hyperboloid $\set{(x,y,z)\in\mathbb{R}^3}{x^2+y^2-z^2=1}$.
In particular, \ref{thm:comp:cat} does not hold without assuming that $\Sigma$ is simply-connected.
\end{thm}

Let us compare the Gauss--Bonnet formula with the comparison theorems.
Suppose a disc $\Delta$ is bounded by a geodesic triangle $[xyz]$ with internal angles $\alpha$, $\beta$, and~$\gamma$.
Then the Gauss--Bonnet implies that 
\[\alpha+\beta+\gamma-\pi=\iint_\Delta K;\]
in particular, both sides of the equation have the same sign.
It follows that
\begin{itemize}
\item if $K_\Sigma\ge 0$, then $\alpha+\beta+\gamma\ge\pi$, and
\item if $K_\Sigma\le 0$, then $\alpha+\beta+\gamma\le\pi$.
\end{itemize}

\begin{wrapfigure}{r}{35mm}
\centering
\vskip-10mm
\includegraphics{mppics/pic-2307}
\end{wrapfigure}

Now set 
$\hat\alpha\z=\measuredangle\hinge {x}{y}{z}$,
$\hat\beta\z=\measuredangle\hinge {y}{z}{x}$,
and $\hat\gamma\z=\measuredangle\hinge {z}{x}{y}$.
Note that $\hat\alpha,\hat\beta,\hat\gamma\in[0,\pi]$.
Since the angles of any plane triangle add up to $\pi$, from the comparison theorems, we get that
\begin{itemize}
\item if $K_\Sigma\ge 0$, then $\hat\alpha+\hat\beta+\hat\gamma\ge\pi$, and
\item if $K_\Sigma\le 0$, then $\hat\alpha+\hat\beta+\hat\gamma\le\pi$.
\end{itemize}

Now we see that despite the Gauss--Bonnet formula and the comparison theorems are closely related,
this relationship is not straightforward.

For example, suppose $K\ge 0$.
Then the Gauss--Bonnet formula does not forbid the internal angles $\alpha$, $\beta$, and $\gamma$ to be simultaneously close to $2\cdot\pi$.
But $\hat \alpha=\alpha$ if $\alpha\le \pi$, and otherwise $\hat \alpha=\pi-\alpha$;
that is,
\begin{align*}
\hat \alpha&=\min\{\,\alpha,2\cdot\pi-\alpha\,\},
&
\hat\beta &=\min\{\,\beta,2\cdot\pi-\beta\,\},
&
\hat\gamma&=\min\{\,\gamma,2\cdot\pi-\gamma\,\}.
\end{align*}
Therefore,
if $\alpha$, $\beta$, and $\gamma$ are each close to $2\cdot\pi$, then $\hat\alpha$, $\hat\beta$, and $\hat\gamma$ are close to $0$.
The latter is impossible by the comparison theorem.

\begin{thm}{Exercise}\label{ex:diam-angle}
Let $p$ and $q$ be points on a closed convex surface $\Sigma$ that lie at a maximal intrinsic distance from each other;
that is, for any $x,y\in \Sigma$ we have $\dist{p}{q}\Sigma\z\ge\dist{x}{y}\Sigma$.
Show that $\measuredangle\hinge xpq\ge \tfrac\pi3$
for any point $x\in \Sigma\setminus\{p,q\}$.
\end{thm}

\begin{thm}{Exercise}\label{ex:sum=<2pi}
Let $\Sigma$ be a closed (or open) surface with nonnegative Gauss curvature.
Show that 
\[\modangle pxy+\modangle pyz+\modangle pzx\le2\cdot \pi.\]
for any four distinct points $p,x,y,z$ in~$\Sigma$.
\end{thm}

\section{Local comparisons}\label{sec:loc-comp}

The following local theorem is the first step in the proof of the comparison theorems (\ref{thm:comp});
it will be deduced from the Rauch comparison (\ref{prop:rauch}).

\begin{thm}{Theorem}\label{thm:loc-comp}
The comparison theorem (\ref{thm:comp}) holds in a small neighborhood of any point.

Moreover, let $\Sigma$ be a smooth surface without boundary.
Then for any $p\in \Sigma$ there is $r>0$ such that if $\dist{p}{x}\Sigma<r$, then $\inj(x)_\Sigma>r$, and the following statements hold:

\begin{subthm}{thm:loc-comp:cba}
If $\Sigma$ has nonpositive Gauss curvature, then 
\[\measuredangle\hinge {x}{y}{z}\le\modangle xyz\]
for any geodesic triangle $[xyz]$ in $B(p,\tfrac r4)_\Sigma$.
\end{subthm}

\begin{subthm}{thm:loc-comp:cbb}
If $\Sigma$ has nonnegative Gauss curvature, then 
\[\measuredangle\hinge {x}{y}{z}\ge\modangle xyz\]
for any geodesic triangle $[xyz]$ in $B(p,\tfrac r4)_\Sigma$.
\end{subthm}

\end{thm}

\parbf{Proof.}
The existence of $r>0$ follows from \ref{prop:exp}.
Let $[xyz]$ be a geodesic triangle in $B(p,\tfrac{r}4)$.

Since $r<\inj(x)_\Sigma$, we can choose $\vec v,\vec w\in\T_x$ such that 
\begin{align*}
y&=\exp_x\vec v,
& 
z&=\exp_x\vec w,
\\
|\vec v|_{\T_x}&=\dist{x}{y}\Sigma,&
|\vec w|_{\T_x}&=\dist{x}{z}\Sigma,
&\measuredangle\hinge 0{\vec v}{\vec w}_{\T_x}&=\measuredangle\hinge xyz_\Sigma,
\end{align*}
where $0$ is the origin in the tangent plane $\T_x$.
Note that $|\vec v|, |\vec w|< \tfrac r2$.

\parit{\ref{SHORT.thm:loc-comp:cba}.}
Consider a minimizing geodesic $\gamma$ joining $y$ to $z$.
Since $\dist{x}{y}{\Sigma},\dist{x}{z}{\Sigma}\z<\tfrac r2$, the triangle inequality implies that $\gamma$ lies in the $r$-neighborhood of~$x$.
In particular, $\log_x\circ\gamma$ is defined, and the curve
$\tilde \gamma\df\log_x\circ\gamma$ lies in an $r$-neighborhood of zero in $\T_x$.
Note that $\tilde\gamma$ connects $\vec v$ to $\vec w$ in $\T_x$.

By the Rauch comparison (\ref{prop:rauch:K=<0}), we have
\[
\length \tilde \gamma \le \length \gamma.
\]
Since $\dist{\vec v}{\vec w}{\T_x} \le \length \tilde \gamma$ and $\length \gamma = \dist{y}{z}{\Sigma}$, it follows that
\[
\dist{\vec v}{\vec w}{\T_x} \le \dist{y}{z}{\Sigma}.
\]
By angle monotonicity (\ref{lem:angle-monotonicity}), we obtain
\[
\measuredangle \hinge{0}{\vec v}{\vec w}_{\T_x} \le \modangle xzy
\]
--- hence the result.

\parit{\ref{SHORT.thm:loc-comp:cbb}.}
Consider the line segment $\tilde \gamma$ joining $\vec v$ to $\vec w$ in the tangent plane $\T_x$, and set $\gamma\df\exp_x\circ\tilde \gamma$.
By Rauch comparison (\ref{prop:rauch:K>=0}), we have
\[\length \tilde \gamma\ge\length\gamma.\]
Since $\dist{\vec v}{\vec w}{\T_x}=\length\tilde \gamma$ and $\length\gamma\ge\dist{y}{z}\Sigma$, we get 
\[\dist{\vec v}{\vec w}{\T_x}\ge \dist{y}{z}\Sigma.\]
By the angle monotonicity (\ref{lem:angle-monotonicity}), we obtain
\[\measuredangle\hinge 0{\vec v}{\vec w}_{\T_x}\ge\modangle xzy\]
--- hence the result.
\qeds

\section{Nonpositive curvature}\label{sec:nonpos-comp}

\parbf{Proof of \ref{thm:comp:cat}.}
Since $\Sigma$ is simply-connected, \ref{ex:unique-geod} implies that
\[\inj(p)_\Sigma=\infty\]
for any $p\in\Sigma$.
Therefore, \ref{thm:loc-comp:cba} implies \ref{thm:comp:cat}.
\qeds

\section{Nonnegative curvature}\label{sec:nonneg-comp}

We will prove \ref{thm:comp:toponogov}, first assuming that $\Sigma$ is compact.
The general case requires only minor modifications; they are indicated in Exercise \ref{ex:open-comparison} at the end of the section.
The proof is taken from \cite{alexander-kapovitch-petrunin2027}, and it is close to the proof found independently by Urs Lang and Viktor Schroeder \cite{lang-schroeder}.

\parbf{Proof of \ref{thm:comp:toponogov} in the compact case.}\label{proof(thm:comp:toponogov)}
Assume $\Sigma$ is compact. 
From the local theorem (\ref{thm:loc-comp}), we get that there is $\epsilon>0$ such that the inequality 
\[\measuredangle\hinge {x}{p}{q}\ge\modangle xpq.\]
holds for any hinge $\hinge{x}{p}{q}$ with 
$\dist{x}{p}\Sigma+\dist{x}{q}\Sigma<\epsilon$.
The following lemma states that in this case, the same holds true for any hinge $\hinge{x}{p}{q}$ such that $\dist{x}{p}\Sigma+\dist{x}{q}\Sigma<\tfrac32\cdot\epsilon$.
Applying the key lemma (\ref{key-lem:globalization}) a few times we get that the comparison holds for an arbitrary hinge, which proves \mbox{\ref{thm:comp:toponogov}}.
\qeds

\begin{thm}{Key lemma}\label{key-lem:globalization} 
Let $\Sigma$ be an open or closed smooth surface.
Assume the comparison
\[\measuredangle\hinge x y z
\ge\modangle x y z\eqlbl{eq:key-lem:globalization}\]
holds for any hinge $\hinge x y z$ with 
$\dist{x}{y}\Sigma+\dist{x}{z}\Sigma
<
\frac{2}{3}\cdot\ell$.
Then the comparison \ref{eq:key-lem:globalization}
holds for any hinge $\hinge x y z$ with $\dist{x}{y}\Sigma+\dist{x}{z}\Sigma<\ell$.
\end{thm}

{

\begin{wrapfigure}{r}{35mm}
\centering
\vskip-0mm
\includegraphics{mppics/pic-2308}
\end{wrapfigure}

Given a hinge $\hinge x p q$ consider a triangle in the plane
with angle $\measuredangle\hinge x p q$ and two adjacent sides $\dist{x}{p}\Sigma$ and $\dist{x}{q}\Sigma$.
Let us denote by $\side \hinge x p q$ the third side of this triangle;
it will be called the \index{model!side}\emph{model side} of the hinge.

The next computational exercise plays a role in the following proof.

}

\begin{thm}{Exercise}\label{ex:s-r}
Suppose hinges $\hinge xpq$ and $\hinge xpy$ have a common side $[x,p]$ and $[x,y]\subset [x,q]$.
Show that 
\[\frac{\dist{x}{p}{}+\dist{x}{q}{}-\side\hinge xpq}{\dist{x}{q}{}}
\le
\frac{\dist{x}{p}{}+\dist{x}{y}{}-\side\hinge xpy}{\dist{x}{y}{}}.\]
\end{thm}

\parbf{Proof.} 
By the angle monotonicity (\ref{lem:angle-monotonicity}), we have
\[\measuredangle\hinge x p q\ge \modangle x p q\quad\iff\quad\side \hinge x p q
\ge\dist{p}{q}\Sigma.\]
Therefore, it is sufficient to prove that
\[\side \hinge x p q
\ge\dist{p}{q}\Sigma,
\eqlbl{eq:thm:=def-loc}\]
assuming that $\dist{x}{p}\Sigma+\dist{x}{q}\Sigma<\ell$.

Let us produce a new hinge $\hinge{x'}p q$ for a given hinge $\hinge x p q$ such that 
\[\tfrac{2}{3}\cdot\ell \le\dist{p}{x}\Sigma\z+\dist{x}{q}\Sigma< \ell.\]

\begin{wrapfigure}{r}{32mm}
\vskip0mm
\centering
\includegraphics{mppics/pic-2310}
\end{wrapfigure}

Assume $\dist{x}{q}\Sigma\ge\dist{x}{p}\Sigma$; otherwise, switch the roles of $p$ and $q$ in the construction.
Take $x'\in [x, q]$ such that 
\[\dist{p}{x}\Sigma+3\cdot\dist{x}{x'}\Sigma
=\tfrac{2}{3}\cdot\ell
\eqlbl{3|xx'|}\]
Choose a geodesic $[x', p]$, and consider the hinge $\hinge{x'}p q$ formed by $[x',p]$ and $[x',q]\subset [x,q]$.

By the triangle inequality, we have 
\[
\dist{p}{x}\Sigma+\dist{x}{q}\Sigma\ge\dist{p}{x'}\Sigma+\dist{x'}{q}\Sigma.
\eqlbl{eq:thm:=def-loc-fourstar}\]
Let us show that
\[\side \hinge x p q
\ge\side \hinge{x'}p q
\eqlbl{eq:thm:=def-loc-fivestar}\]

By \ref{3|xx'|}, we have that 
\[
\begin{aligned}
\dist{p}{x}{\Sigma}\z+\dist{x}{x'}{\Sigma}&<\tfrac{2}{3}\cdot\ell,
\\
\dist{p}{x'}{\Sigma}\z+\dist{x'}{x}{\Sigma}&<\tfrac{2}{3}\cdot\ell.
\end{aligned}
\]
Therefore, the assumption imply that 
\[\begin{aligned}
\measuredangle\hinge x p{x'}
\ge\modangle x p{x'}
\quad\text{and}\quad
\measuredangle\hinge {x'}p x
\ge\modangle {x'}p x.
  \end{aligned}
\eqlbl{eq:thm:=def-loc-threestar}
\]

Consider the model triangle
$[\tilde x\tilde x'\tilde p]\z=\modtrig xx'p$.
Take $\tilde q$ on the extension of $[\tilde x,\tilde x']$ beyond $x'$ such that $\dist{\tilde x}{\tilde q}\Sigma=\dist{x}{q}\Sigma$, and therefore $\dist{\tilde x'}{\tilde q}\Sigma\z=\dist{x'}{q}\Sigma$.

From \ref{eq:thm:=def-loc-threestar}, we get
\[\measuredangle\hinge x p q
=\measuredangle\hinge x p{x'}\ge\modangle x p{x'}.\]
Therefore,
\[
\side\hinge x q p
\ge
\dist{\tilde p}{\tilde q}{\mathbb{R}^2}.
\]
Further, since $\measuredangle\hinge{x'}p x+\measuredangle\hinge{x'}p q= \pi$,
the inequalities in \ref{eq:thm:=def-loc-threestar} imply
\[
\pi
-\modangle{x'}p x
\ge
\pi-\measuredangle\hinge{x'}p x
\ge
\measuredangle\hinge{x'}p q.
\]
Therefore,
\[\dist{\tilde p}{\tilde q}{\mathbb{R}^2}\ge\side \hinge{x'}q p \]
and \ref{eq:thm:=def-loc-fivestar} follows.

Set $x_0=x$; apply inductively the above construction to get a sequence of hinges $\hinge{x_n}p q$ with $x_{n+1}=x_n'$.
By \ref{eq:thm:=def-loc-fivestar} and \ref{eq:thm:=def-loc-fourstar}, both sequences
\[s_n=\side \hinge{x_n}pq\quad\text{and}\quad r_n=\dist{p}{x_n}\Sigma+\dist{x_n}{q}\Sigma\]
are nonincreasing.

The sequence might terminate at $x_n$ only if $r_n< \tfrac{2}{3}\cdot\ell $.
In this case, by the assumptions of the lemma, we have
\[s_n=\side \hinge{x_n}p q\ge\dist{p}{q}\Sigma.\]
Since the sequence $s_n$ is nonincreasing, we get
\[\side \hinge{x}p q=s_0\ge s_n\ge\dist{p}{q}\Sigma;\]
whence \ref{eq:thm:=def-loc} follows.

\begin{figure}[!ht]
\centering
\includegraphics{mppics/pic-2315}
\end{figure}

It remains to prove \ref{eq:thm:=def-loc} if the sequence $x_n$ does not terminate.
By \ref{3|xx'|}, we have 
\[
\dist{x_n}{x_{n-1}}\Sigma
\ge 
\tfrac1{100}\cdot \ell.
\eqlbl{eq:|x-x|><l}
\]
By \ref{eq:thm:=def-loc-fourstar}, $\dist{x_n}{p}{},\dist{x_n}{q}{}<
\ell$ for any~$n$.
In case $x_{n+1}\in [x_n,q]$, apply \ref{ex:s-r} for the hinges $\hinge{x_n}pq$ and $\hinge{x_n}p{x_{n+1}}$.
By \ref{eq:thm:=def-loc-threestar}, $\dist{p}{x_{n+1}}{}\le \side \hinge{x_n}{x_{n+1}}{p}$.
Therefore,
\[r_n-s_n\le 100\cdot (r_n-r_{n+1})\eqlbl{eq:r-s<100(r-r)}\]
In the case $x_{n+1}\in [x_n,p]$, inequality \ref{eq:r-s<100(r-r)} follows if one applies \ref{ex:s-r} for the hinges $\hinge{x_n}pq$ and $\hinge{x_n}{x_{n+1}}q$.

The sequences $r_n$ and $s_n$ are nonincreasing and nonnegative;
so, they have to converge.
In particular, $(r_n-r_{n+1})\to0$ as $n\to \infty$.
Therefore, \ref{eq:r-s<100(r-r)} implies that
\[\lim_{n\to\infty}s_n=\lim_{n\to\infty}r_n.\]
By the triangle inequality, $r_n\ge \dist{p}{q}\Sigma$ for any~$n$.
Since $s_n$ is nonincreasing, we get
\[\side \hinge{x}p q=s_0\ge \lim_{n\to\infty}s_n=\lim_{n\to\infty}r_n\ge \dist{p}{q}\Sigma\]
which finishes the proof of \ref{eq:thm:=def-loc}.
\qeds

\begin{thm}{Exercise}\label{ex:open-comparison}
Let $\Sigma$ be an open surface with nonnegative Gauss curvature.
Given $p\in\Sigma$, denote by $R_p$ 
(the {}\emph{comparison radius} at $p$) 
the maximal value (possibly $\infty$) such that the comparison 
\[\measuredangle\hinge x p y
\ge\modangle x p y\]
holds for any hinge $\hinge x p y$ with $\dist{p}{x}\Sigma+\dist{x}{y}\Sigma<R_p$.

\begin{subthm}{ex:open-comparison:positive}
Show that for any compact subset $K\subset \Sigma$, there is $\epsilon>0$ such that $R_p>\epsilon$ for any $p\in K$.
\end{subthm}

\begin{subthm}{ex:open-comparison:almost-min}
Use part \ref{SHORT.ex:open-comparison:positive} to show that 
there is a point $p\in\Sigma$ such that 
\[R_q>(1-\tfrac1{100})\cdot R_p,\]
for any $q\in B(p,100\cdot R_p)_\Sigma$.
\end{subthm}

\begin{subthm}{ex:open-comparison:proof}
Explain how to use \ref{SHORT.ex:open-comparison:almost-min} to extend the proof of \ref{thm:comp:toponogov} (page \pageref{proof(thm:comp:toponogov)}) to open surfaces. 
(That is, to show that $R_p=\infty$ for any $p\in\Sigma$.) 
\end{subthm}

\end{thm}

\section{Alexandrov's lemma}
\index{Alexandrov's lemma}

A reformulation of the following lemma (\ref{lem:alex-reformulation}) will be used in the next section to produce a few equivalent reformulations of the comparison theorems.

\begin{thm}{Lemma}
\label{lem:alex}
Consider two quadrangles in the Euclidean plane $pxyz$ and $p'x'y'z'$ with equal corresponding sides.
Assume the sides $[x',y']$ and $[y',z']$ extend each other; that is, $y'$ lies on the line segment $[x',z']$.
Then the following expressions have the same signs:
\begin{enumerate}[(i)]
 \item $|p-y|-|p'-y'|$;
 \item $\measuredangle\hinge xpy-\measuredangle\hinge {x'}{p'}{y'}$;
 \item $\pi-\measuredangle\hinge ypx-\measuredangle\hinge ypz$.
\end{enumerate}
\end{thm}

\parbf{Proof.} 
Take 
a point $\bar z$ on the extension of 
$[x,y]$ beyond $y$ so that $\dist{y}{\bar z}{}=\dist{y}{z}{}$ (and therefore $\dist{x}{\bar z}{}=\dist{x'}{z'}{}$). 
 
\begin{figure}[!ht]
\vskip-0mm
\centering
\includegraphics{mppics/pic-50}
\vskip-0mm
\end{figure}

By the angle monotonicity (\ref{lem:angle-monotonicity}), 
the following expressions have the same sign:
\begin{enumerate}[(i)]
\item $|p-y|-|p'-y'|$;
\item $\measuredangle\hinge{x}{y}{p}-\measuredangle\hinge{x'}{y'}{p'}=\measuredangle\hinge{x}{\bar z}{p}-\measuredangle\hinge{x'}{z'}{p'}$;
\item $|p-\bar z|-|p'-z'| = | p - \bar z | - | p-z | $;
\item $\measuredangle\hinge{y}{\bar z}{p}-\measuredangle\hinge{y}{z}{p}$.
\end{enumerate}
The statement follows since
\[\measuredangle\hinge{y'}{z'}{p'}+\measuredangle\hinge{y'}{x'}{p'}=\pi
\quad\text{and}\quad
\measuredangle\hinge{y}{\bar z}{p}+\measuredangle\hinge{y}{x}{p}=\pi.\]
\qedsf

\section{Reformulations}

For any triangle $[xyz]$ in a surface $\Sigma$, and its model triangle $[\tilde x \tilde y \tilde z]$, there is a natural map $p\mapsto \tilde p$ that isometrically sends the geodesics $[x,y]$, $[y,z]$, $[z,x]$ to the line segments $[\tilde x,\tilde y ]$, $[\tilde y , \tilde z]$, $[ \tilde z , \tilde x ]$, respectively.
The triangle $[xyz]$ is called \index{fat triangle}\emph{fat} (\index{thin triangle}\emph{thin})
if the inequality
\[\dist{p}{q}{\Sigma}\ge |\tilde p- \tilde q|_{\mathbb{R}^2}\qquad \text{(or, respectively,}\quad \dist{p}{q}{\Sigma}\le |\tilde p- \tilde q|_{\mathbb{R}^2})\]
holds for any two points $p$ and $q$ on the sides of $[xyz]$.

\begin{thm}{Proposition}\label{prop:comp-reformulations}
Let $\Sigma$ be an open or closed smooth surface.
Then the following three conditions are equivalent:

\begin{subthm}{mang>angk}
For any geodesic triangle $[xyz]$ in $\Sigma$ we have
 \[\measuredangle\hinge {x}{y}{z}\ge\modangle xyz.\]
\end{subthm}

\begin{subthm}{angk>angk} For any geodesic triangle $[pxz]$ in $\Sigma$ and $y$ on the side $[x,z]$ we have
 \[\modangle xpy \ge \modangle xpz.\]
 
\end{subthm}

\begin{subthm}{fat}
 Any geodesic triangle in $\Sigma$ is fat.
\end{subthm}

\medskip

Similarly, the following three conditions are equivalent:

\begin{subthmA}{mang<angk}
For any geodesic triangle $[xyz]$ in $\Sigma$ we have
 \[\measuredangle\hinge {x}{y}{z}\le\modangle xyz.\]
\end{subthmA}

\begin{subthmA}{angk<angk} For any geodesic triangle $[pxz]$ in $\Sigma$ and $y$ on the side $[x,z]$ we have
 \[\modangle xpy \le \modangle xpz.\]
\end{subthmA}

\begin{subthmA}{thin}
Any geodesic triangle in $\Sigma$ is thin.
\end{subthmA}

\end{thm}

Let us rewrite the Alexandrov lemma (\ref{lem:alex}) in the language of comparison triangles and angles.

\begin{thm}{Reformulation of Alexandrov lemma}\label{lem:alex-reformulation}
Let $[pxz]$ be a triangle in a surface $\Sigma$
and $y$ a point on the side $[x,z]$.
Consider its model triangle $[\tilde p\tilde x\tilde z]\z=\tilde\triangle pxz$, and let $\tilde y$ be the corresponding point on the side $[\tilde x,\tilde z]$.

\begin{wrapfigure}{r}{25mm}
\vskip-5mm
\centering
\includegraphics{mppics/pic-2305}
\end{wrapfigure}

Then the following expressions have the same signs:
\begin{enumerate}[(i)]
 \item $\dist{p}{y}\Sigma-\dist{\tilde p}{\tilde y}{\mathbb{R}^2}$;
 \item $\modangle xpy-\modangle {x}{p}{z}$;
 \item $\pi-\modangle ypx-\modangle ypz$;
\end{enumerate}
\end{thm}

\parbf{Proof of \ref{prop:comp-reformulations}.}
We will prove the implications \ref{SHORT.mang>angk}$\Rightarrow$\ref{SHORT.angk>angk}$\Rightarrow$\ref{SHORT.fat}$\Rightarrow$\ref{SHORT.mang>angk}.
The implications \ref{SHORT.mang<angk}$\Rightarrow$\ref{SHORT.angk<angk}$\Rightarrow$\ref{SHORT.thin}$\Rightarrow$\ref{SHORT.mang<angk} can be proved in the same way with all inequalities reversed.

\parit{\ref{SHORT.mang>angk}$\Rightarrow$\ref{SHORT.angk>angk}.}
Note that $\measuredangle\hinge ypx+\measuredangle\hinge ypz=\pi$.
By \ref{SHORT.mang>angk}, 
\[\modangle ypx+\modangle ypz\le \pi.\]
It remains to apply Alexandrov's lemma (\ref{lem:alex-reformulation}).

\parit{\ref{SHORT.angk>angk}$\Rightarrow$\ref{SHORT.fat}.}
Applying \ref{SHORT.angk>angk} twice, first for $y\in [x,z]$ and then for $w\in [p,x]$, we get that
\[\modangle xwy \ge \modangle xpy \ge \modangle xpz,\]
and therefore
\[\dist{w}{y}\Sigma\ge \dist{\tilde w}{\tilde y}{\mathbb{R}^2},\]
where $\tilde w$ and $\tilde y$ are the points corresponding to $w$ and $y$ points on the sides of the model triangle. 

\parit{\ref{SHORT.fat}$\Rightarrow$\ref{SHORT.mang>angk}.}
Since the triangle is fat, we have 
\[\modangle xwy \ge \modangle xpz\]
for any $w\in [x,p]\setminus \{x\}$ and $y\in [x,z]\setminus \{x\}$.
Note that $\modangle xwy\to \measuredangle\hinge xpz$ as $w,y\to x$.
Whence the implication follows.
\qeds

\begin{thm}{Exercise}\label{ex:geod-convexity}
Let $\Sigma$ be an open smooth surface, 
 $\gamma$ a unit-speed geodesic in $\Sigma$, and $p\in\Sigma$.

Consider the function
\[h(t)=\dist{p}{\gamma(t)}\Sigma^2-t^2.\]

\begin{subthm}{}
Show that if $\Sigma$ is simply-connected and the Gauss curvature of $\Sigma$ is nonpositive, then the function $h$ is convex.
\end{subthm}

\begin{subthm}{} Show that if the Gauss curvature of $\Sigma$ is nonnegative, then the function $h$ is concave.
\end{subthm}

\end{thm}

\begin{thm}{Exercise}\label{ex:midpoints}
Let $\bar x$ and $\bar y$ be the midpoints of minimizing geodesics $[p,x]$ and $[p,y]$ in an open smooth surface~$\Sigma$.

\begin{subthm}{}
Show that if $\Sigma$ is simply-connected and has nonpositive Gauss curvature, then 
\[2\cdot \dist{\bar x}{\bar y}\Sigma\le \dist{x}{y}\Sigma.\]
\end{subthm}
 
\begin{subthm}{} Show that if the Gauss curvature of $\Sigma$ is nonnegative, then 
 \[2\cdot \dist{\bar x}{\bar y}\Sigma\ge \dist{x}{y}\Sigma.\]
\end{subthm}

\end{thm}

\begin{thm}{Exercise}\label{ex:convex-dist}
Let $\gamma_1$ and $\gamma_2$ be two geodesics in an open smooth simply-connected surface $\Sigma$ with nonpositive Gauss curvature.
Show that the function $h(t)\df\dist{\gamma_1(t)}{\gamma_2(t)}\Sigma$
is convex.
\end{thm}

\begin{thm}{Exercise}\label{ex:disc+}
Let $\Sigma$ be an open or closed smooth surface with nonnegative Gauss curvature.
Show that area of any $R$-disc in the intrinsic metric of $\Sigma$ is at most $\pi\cdot R^2$.

\end{thm}

\begin{thm}{Exercise}\label{ex:disc-}
Let $\Delta$ be an $R$-disc in the intrinsic metric of an open simply-connected smooth surface $\Sigma$  with nonpositive Gauss curvature.

\begin{subthm}{ex:disc-:kg}
Show that the boundary of $\Delta$ is a smooth curve with geodesic curvature at least $\tfrac1R$.
\end{subthm}

\begin{subthm}{ex:disc-:area}
Show that the area of $\Delta$ is at least $\pi\cdot R^2$.
\end{subthm}

\end{thm}

The following exercise generalizes the Moon in a puddle theorem (\ref{thm:moon-orginal}).

\begin{thm}{Advanced exercise}\label{ex:moon-}
Let $\Delta$ be a disc in a smooth surface $\Sigma$ with nonpositive Gauss curvature.
Assume $\Delta$ is bounded by a smooth curve $\gamma$ with geodesic curvature at most 1 by absolute value.
Show that $\Delta$ contains a unit disc in the intrinsic metric of $\Sigma$.

Conclude that the area of $\Delta$ is at least $\pi$.
\end{thm}

%% file: prelim.tex
\bookmarksetupnext{startatroot,level=0}
\chapter[Appendix]{}
\chaptermark{}
%\addcontentsline{toc}{chapter}{Preliminaries}

\vskip0mm

{\footnotesize

This appendix should be used as a quick reference while reading the rest of the book;
it contains statements of necessary results and references to complete proofs.

\section{Metric spaces}\label{sec:metric-spcaes}

We assume familiarity with the notion of distance in the
Euclidean space.
In this section, we briefly discuss its generalization and fix notations that will be used further.
All these topics are discussed in detail in the introductory part of the book by Dmitri Burago, Yuri Burago, and Sergei Ivanov \cite{burago-burago-ivanov}.

\setlength{\multicolsep}{20pt}
\spell{\begin{multicols}{2}}{}

\subsection*{Definitions}

\emph{Metric} is a function that returns a real value $\Dist(x,y)$ for any pair of elements $x,y$ in a given set $\spc{X}$  and satisfies the following axioms for any triple $x,y,z\in \spc{X}$: \label{page:def:metric}
\begin{enumerate}[(a)]
\item\label{def:metric-space:a} Positiveness: 
$$\Dist(x,y)\ge 0.$$
\item\label{def:metric-space:b} $x=y$ if and only if 
$$\Dist(x,y)=0.$$
\item\label{def:metric-space:c} Symmetry: $$\Dist(x, y) = \Dist(y, x).$$
\item\label{def:metric-space:d} Triangle inequality: 
$$\Dist(x, z) \le \Dist(x, y) + \Dist(y, z).$$
\end{enumerate}

A set with a metric is called a \index{metric space}\emph{metric space}, and the elements of the set are called \index{point}\emph{points}.

\subsection*{Shortcut for distance}
Usually, we consider only one metric on a set, therefore we can denote the metric space and its underlying set by the same letter, say $\spc{X}$.
In this case, we also use the shortcut notations $\dist{x}{y}{}$ or $\dist{x}{y}{\spc{X}}$ for the {}\emph{distance} $\Dist(x,y)$ from $x$ to $y$ in $\spc{X}$.\index{10aaa@$\dist{x}{y}{}$, $\dist{x}{y}{\spc{X}}$ (distance)}
For example, the triangle inequality can be written as 
$$\dist{x}{z}{\spc{X}}\le \dist{x}{y}{\spc{X}}+\dist{y}{z}{\spc{X}}.$$

The Euclidean space and plane, as well as the real line, will be the most important examples of metric spaces.
In these examples, the introduced notation $\dist{x}{y}{}$ for the distance from $x$ to $y$ has perfect sense as the norm of the vector $x-y$.
However, let us emphasise that in a general metric space, the difference of points $x-y$ has no meaning.

\subsection*{More examples}

Usually, if we say {}\emph{plane} or {}\emph{space} we mean the {}\emph{Euclidean} plane or space.
However, the plane (as well as the space) admits many other metrics; for example, the so-called {}\emph{Manhattan metric} from the following exercise.

\begin{thm}{Exercise}\label{ex:ell-infty}
Consider the function
$$\Dist(p,q)=|x_1-x_2|+|y_1-y_2|,$$
where $p=(x_1,y_1)$ and $q=(x_2,y_2)$ are points in the coordinate plane $\mathbb{R}^2$.
Show that $\Dist$ is a metric on $\mathbb{R}^2$.
\end{thm}

Another example: the {}\emph{discrete space} --- an arbitrary nonempty set $\spc{X}$ with the metric defined as $\dist{x}{y}{\spc{X}}=0$ if $x=y$ and $\dist{x}{y}{\spc{X}}=1$ otherwise.

\subsection*{Subspaces}
Any subset of a metric space is also a metric space, by restricting the original metric to the subset;
the obtained metric space is called a {}\emph{subspace}.
In particular, all subsets of the Euclidean space are metric spaces.

\subsection*{Balls}
Given a point $p$ in a metric space $\spc{X}$ and a real number $R\ge 0$, the set of points $x$ at the distance less than $R$ (at most $R$) from $p$ is called the \index{open!ball}\emph{open} (respectively \index{closed!ball}\emph{closed}) {}\emph{ball} of radius $R$ with center~$p$.
The open ball is denoted as $B(p,R)$ or $B(p,R)_{\spc{X}}$;
the second notation is used if we need to emphasize that the ball is from $\spc{X}$.
Formally speaking
\begin{align*}
B(p,R)&=B(p,R)_{\spc{X}}=
\\
&=\set{x\in \spc{X}}{\dist{x}{p}{\spc{X}}< R}.
\end{align*}
\index{10b@$B(p,R)_{\spc{X}}$, $\bar B[p,R]_{\spc{X}}$ (ball)}
Analogously, the closed ball is denoted as $\bar B[p,R]$ or $\bar B[p,R]_{\spc{X}}$ and
\begin{align*}
\bar B[p,R]&=\bar B[p,R]_{\spc{X}}=
\\
&=\set{x\in \spc{X}}{\dist{x}{p}{\spc{X}}\le R}.
\end{align*}

\begin{thm}{Exercise}\label{ex:B2inB1}

\begin{subthm}{ex:B2inB1:a}
Let $p$ and $q$ be points in a metric space $\spc{X}$.
Show that if $\bar B[p,2]\z\subset \bar B[q,1]$, then $\bar B[p,2]\z=\bar B[q,1]$.
\end{subthm}

\begin{subthm}{ex:B2inB1:b}
Construct a metric space $\spc{X}$ with two points $p$ and $q$ such that the strict inclusion
$B(p,\tfrac32)\subsetneq B(q,1)$ holds.
\end{subthm}

\end{thm}

\subsection*{Isometries}

Let $\spc{X}$ and $\spc{Y}$ be two metric spaces. 
A~map $f\:\spc{X} \z\to \spc{Y}$
is
called \index{distance-preserving map}\emph{distance-preserving} if 
$$\dist{f(x)}{f(y)}{\spc{Y}}
 = \dist{x}{y}{\spc{X}}$$
for any $x,y\in {\spc{X}}$.

A bijective distance-preserving map is called an \index{isometry}\emph{isometry}. 
Two metric spaces are called {}\emph{isometric} if there exists an isometry between them.

\begin{thm}{Exercise}\label{ex:dist-preserv=>injective}
Show that any distance-preserving map $f\:\spc{X}\to\spc{Y}$ is \index{injective map}\emph{injective};
that is, $f(x)\ne f(y)$ for any pair of distinct points $x, y\in \spc{X}$.
\end{thm}

\subsection*{Continuity}

\begin{thm}{Definition}
 Let ${\spc{X}}$ be a metric space.
A sequence of points $x_1, x_2, \ldots$ in ${\spc{X}}$ {}\emph{converges}
if there is a point
$x_\infty\in {\spc{X}}$ such that $\dist{x_\infty}{x_n}{}\to 0$ as $n\to\infty$.  
That is, for every $\epsilon > 0$, there is a natural number $N$ such that for all $n \ge N$, we have
\[
\dist{x_\infty}{x_n}{\spc{X}}
<
\epsilon.
\]

In this case, we say that the sequence $x_n$ {}\emph{converges} to $x_\infty$, 
or $x_\infty$ is the {}\emph{limit} of the sequence $x_n$.
Notationally, we write $x_n\to x_\infty$ as $n\to\infty$
or $x_\infty=\lim_{n\to\infty} x_n$.
\end{thm}

\begin{thm}{Definition}\label{def:continuous}
Let $\spc{X}$ and $\spc{Y}$ be metric spaces.
A map $f\:\spc{X}\to \spc{Y}$ is called \index{continuous}\emph{continuous} if, for any convergent sequence $x_n\to x_\infty$ in ${\spc{X}}$,
we have $f(x_n) \to f(x_\infty)$ in $\spc{Y}$.

Equivalently, $f\:\spc{X}\to \spc{Y}$ is continuous if, for any $x\in {\spc{X}}$ and any $\epsilon>0$,
there is $\delta>0$ such that 
$$\dist{x}{y}{\spc{X}}<\delta\quad \Longrightarrow\quad \dist{f(x)}{f(y)}{\spc{Y}}<\epsilon.$$

\end{thm}

\begin{thm}{Exercise}\label{ex:shrt=>continuous}
Let $f\:\spc{X}\to \spc{Y}$ be a {}\emph{distance non-expanding} map between metric spaces; that is, 
\[\dist{f(x)}{f(y)}{\spc{Y}}\le \dist{x}{y}{\spc{X}}\]
for any $x,y\in \spc{X}$.
Show that $f$ is continuous.
\end{thm}
\spell{\end{multicols}}{}

\section{Topology}\label{sec:topology}

The following material is covered in any introductory text to topology; 
one of our favorites is a textbook by Czes Kosniowski \cite{kosniowski}.

\spell{\begin{multicols}{2}}{}

\subsection*{Closed and open sets}

\begin{thm}{Definition}
A subset $C$ of a metric space $\spc{X}$ is called \index{closed!set}\emph{closed} if whenever a sequence $x_n$ of points from $C$ converges in $\spc{X}$, we have that $\lim_{n\to\infty} x_n \in C$.

A set $\Omega \subset \spc{X}$ is called \index{open!set}\emph{open} if, for any $z\in \Omega$, 
there is $\epsilon>0$ such that $B(z,\epsilon)\subset\Omega$.
\end{thm}

\begin{thm}{Exercise}\label{ex:close-open}
Let $Q$ be a subset of a metric space $\spc{X}$.
Show that $Q$ is closed if and only if its complement $\Omega=\spc{X}\setminus Q$ is open.
\end{thm}

An open set $\Omega$ that contains a given point $p$ is called a \index{neighborhood}\emph{neighborhood of~$p$}.
A closed subset $C$ that contains $p$ together with its neighborhood is called a {}\emph{closed neighborhood of~$p$}.

A point $p$ lies on the \index{boundary}\emph{boundary} of a set $Q$ (briefly, $p\in\partial Q$) if any neighborhood of $p$ contains points in $Q$ and in the complement of~$Q$.

\subsection*{Compact sets}

A subset $K$ of a metric space is called \index{compact}\emph{compact} if any sequence of points in $K$ has a subsequence that converges to a point in~$K$.

The following properties follow directly from the definition:

\begin{itemize}
\item A closed subset of a compact set is compact.
\item An image of a compact set under a continuous map is compact.
\end{itemize}

\begin{thm}{Heine--Borel lemma}\label{thm:Heine--Borel}
A subset of Euclidean space is compact if and only if it is closed and bounded.
\end{thm}

\subsection*{Homeomorphisms\\ and embeddings}

A bijection $f\:\spc{X}\to\spc{Y}$ between metric spaces is called a \index{homeomorphism}\emph{homeomorphism} if $f$ and its inverse $f^{-1}$ are continuous.
A homeomorphism to its image is called an \index{embedding}\emph{embedding}.

If there exists a homeomorphism $f\:\spc{X}\z\to \spc{Y}$,
we say that ${\spc{X}}$ is {}\emph{homeomorphic} to $\spc{Y}$,
or $\spc{X}$ and $\spc{Y}$ are {}\emph{homeomorphic}.

If a metric space $\spc{X}$ is homeomorphic to a known space, for example, plane, sphere, disc, circle, and so on,
then we may also say that $\spc{X}$ is a \index{topological surface}\emph{topological} plane, sphere, disc, circle, and so on.

The following theorem characterizes homeomorphisms between compact spaces:

\begin{thm}{Theorem}\label{thm:Hausdorff-compact}
A continuous bijection $f$ between compact metric spaces has a continuous inverse.
In particular, we have the following.

\begin{subthm}{}
Any continuous bijection between compact metric spaces
is a homeomorphism.
\end{subthm}

\begin{subthm}{}
Any continuous injection from a compact metric spaces to another metric space
is an embedding.
\end{subthm}

\end{thm}

\subsection*{Connected sets}

Recall that a continuous map $\alpha$ from the unit interval $[0,1]$ to a Euclidean space is called a \index{path}\emph{path}.
If $p=\alpha (0)$ and $q = \alpha (1)$, then we say that $\alpha$ connects $p$ to~$q$.

A nonempty set $X$ in the Euclidean space is called \index{path-connected set}\emph{path-connected} if any two points $x,y\in X$ can be connected by a path lying in~$X$.

A nonempty set $X$ in the Euclidean space is called \index{connected}\emph{connected} if one cannot cover $X$ with two disjoint open sets $V$ and $W$ such that both intersections $X\cap V$ and $X\cap W$ are nonempty.

Notice that path-connected and connected sets are nonempty by the definition.

\begin{thm}{Proposition}
Any path-connected set is connected.

Moreover, any open connected set in the Euclidean space or plane is path-connected.
\end{thm}

Given a point $x\in X$, the maximal connected subset of $X$ containing $x$ is called the {}\emph{connected component} of $x$ in~$X$.

\subsection*{Jordan's theorem}
\index{Jordan's theorem}

The first part of the following theorem was proved by Camille Jordan, the second part is due to Arthur Schoenflies:

\begin{thm}{Theorem}\label{thm:jordan}
The complement of any simple closed curve $\gamma$ in $\mathbb{R}^2$ has exactly two connected components. 

Moreover, there is a homeomorphism $h\:\mathbb{R}^2\to \mathbb{R}^2$ that maps the unit circle to~$\gamma$.
In particular, $\gamma$ bounds a topological disc.
\end{thm}

This theorem is known for its simple formulation and quite hard proof.
For the first statement, a very short proof based on a somewhat developed technique is given by Patrick Doyle \cite{doyle},
among elementary proofs, one of our favorites is the proof given by Aleksei Filippov \cite{filippov}.

We use mostly the smooth case of this theorem which is simpler.
An amusing proof of this case was given by Gregory Chambers and Yevgeny Liokumovich \cite{chambers-liokumovich}.

\spell{\end{multicols}}{}

\section{Elementary geometry}

%%%???orientation is not defined???
\spell{\begin{multicols}{2}}{}
\subsection*{Internal angles}

Polygon is defined as a compact set bounded by a closed polygonal line. 
Recall that the internal angle of a polygon $P$ at a vertex $v$
is defined as the angular measure of the intersection of $P$ with a small circle centered at~$v$.

\begin{thm}{Theorem}\label{thm:sum=(n-2)pi}
The sum of all the internal angles of an $n$-gon is $(n\z-2)\cdot\pi$. 
\end{thm}

A clean proof of this theorem can be found, for example, in \cite{meisters}.
It uses induction on $n$ and is based on the following nontrivial statement:

\begin{thm}{Claim}
Suppose $P$ is an $n$-gon with $n\ge 4$.
Then a diagonal of $P$ lies completely in~$P$.
\end{thm}

\subsection*{Angle monotonicity}

The {}\emph{measure} of angle with sides $[p,x]$ and $[p,y]$ will be denoted by $\measuredangle\hinge pxy$\index{10aab@$\measuredangle\hinge yxz$ (angle measure)};
it takes a value in the interval $[0,\pi]$.

The following corollary of the cosine rule is simple and useful.
It says that the angle of a triangle monotonically depends on the opposite side, assuming we keep the other two sides fixed.

\begin{thm}{Monotonicity lemma}\label{lem:angle-monotonicity}
Let $x$, $y$, $z$, $x^{*}$, $y^{*}$, and $z^{*}$ be 6 points such that $\dist{x}{y}{}\z=|x^{*}-y^{*}|>0$ and $|y-z|=|y^{*}-z^{*}|>0$.
Then 
\[\measuredangle\hinge yxz
\ge
\measuredangle\hinge {y^{*}}{x^{*}}{z^{*}}
\ \Longleftrightarrow\
|x-z|\ge |x^{*}-z^{*}|.\]
\end{thm}

\subsection*{Spherical triangle inequality}

The following theorem says that the triangle inequality holds for angles between half-lines from a fixed point.
In particular, it implies that a sphere with the angle metric is a metric space.

\begin{thm}{Theorem}\label{thm:spherical-triangle-inq}
The following inequality holds for any three line segments $[o,a]$, $[o,b]$, and $[o,c]$ in the Euclidean space:
\[\measuredangle\hinge oab
+
\measuredangle\hinge obc
\ge
\measuredangle\hinge oac\]

\end{thm}

Most of the authors use this theorem without mentioning, but the proof is not that simple.
A short elementary proof can be found in the classical textbook in Euclidean geometry by Andrey Kiselyov \cite[§ 47]{kiselyov}.

\subsection*{Area of spherical triangle}

\begin{thm}{Lemma}\label{lem:area-spher-triangle}
Let $\Delta$ be a spherical triangle;
that is, an intersection of three closed half-spheres in the unit sphere $\mathbb{S}^2$.
Then 
\[\area\Delta=\alpha+\beta+\gamma-\pi,\eqlbl{eq:area(Delta)}\]
where $\alpha$, $\beta$, and $\gamma$ are the angles of $\Delta$.
\end{thm}

The value $\alpha+\beta+\gamma-\pi$ is called the \index{excess of triangle}\emph{excess} of the triangle $\Delta$,
so the lemma says that the area of a spherical triangle equals its excess.

This lemma will be important for the intuitive understanding of the Gauss--Bonnet formula.
By that reason, we present its proof.

\begin{wrapfigure}{r}{17 mm}
\vskip-0mm
\centering
\includegraphics{mppics/pic-43}
\vskip2mm
\end{wrapfigure}

\parbf{Proof.}
Recall that 
\[\area\mathbb{S}^2=4\cdot\pi.\eqlbl{eq:area(S2)}\]

Note that the area of a spherical slice $S_\alpha$ between two meridians meeting at angle $\alpha$ is proportional to $\alpha$.
Since $S_\pi$ is a half-sphere, from \ref{eq:area(S2)}, we get $\area S_\pi\z=2\cdot\pi$.
Therefore, the coefficient is $2$; that is,
\[\area S_\alpha=2\cdot \alpha
\eqlbl{eq:area(Sa)}\]
for any $\alpha$.

Extending the sides of $\Delta$ we get 6 slices: two $S_\alpha$, two $S_\beta$, and two $S_\gamma$.
They cover most of the sphere once,
but the triangle $\Delta$ and its centrally symmetric copy $\Delta^{*}$ are covered 3 times.
It follows that
\begin{align*}
2\cdot \area S_\alpha &+2\cdot \area S_\beta+2\cdot \area S_\gamma=
\\
&=\area\mathbb{S}^2+4\cdot\area\Delta.
\end{align*}
It remains to apply \ref{eq:area(S2)} and \ref{eq:area(Sa)}.
\qeds

\spell{\end{multicols}}{}

\section{Convex geometry}

\spell{\begin{multicols}{2}}{}
A set $X$ in the Euclidean space is called \index{convex!set}\emph{convex} if, for any two points $x,y\in X$, any point $z$ between $x$ and $y$ lies in~$X$.
It is called  {}\emph{strictly convex} if, for any two points $x,y\in X$, any point $z$ between $x$ and $y$ lies in the interior of~$X$.

From the definition, it is easy to see that the intersection of an arbitrary family of convex sets is convex. 
The intersection of all convex sets containing $X$ is called the \index{convex!hull}\emph{convex hull} of $X$;
it is the minimal convex set containing the set~$X$.

These definitions and the following statements should appear on the first pages of any introductory text in convex geometry;
see for example the book by Roger Webster \cite{webster}.

\subsection*{Separating and supporting planes}

We will use the following corollary of the so-called \index{hyperplane separation theorem}\emph{hyperplane separation theorem}.

\begin{thm}{Lemma}\label{lem:separation}
Let $K\subset \mathbb{R}^3$ be a closed convex set.
Then for any point $p\notin K$ there is a plane $\Pi$ that separates $K$ from $p$;
that is, $K$ and $p$ lie on opposite open half-spaces separated by $\Pi$.

Moreover, for any boundary point $p\in\partial K$, there is a plane $\Pi$ \index{supporting!plane}\emph{supporting} $K$ at $p$;
that is, $\Pi\ni p$ and $K$ lies in a closed half-space bounded by $\Pi$
\end{thm}

\spell{\end{multicols}}{}

\begin{figure*}[h!]
\begin{minipage}{.48\textwidth}
\centering
\includegraphics{mppics/pic-3540}
\end{minipage}\hfill
\begin{minipage}{.48\textwidth}
\centering
\includegraphics{mppics/pic-3542}
\end{minipage}
\end{figure*}

\section{Linear algebra}

\spell{\begin{multicols}{2}}{}

The following theorem can be found in any textbook in linear algebra;
the book by Sergei Treil \cite{treil} will do.

\begin{thm}{Spectral theorem}\label{thm:spectral}
Any symmetric matrix is diagonalizable  by an orthogonal matrix.
\end{thm}

We will use this theorem only for $2{\times}2$ matrices.
In this case, it can be restated as follows:
Consider a function 
\begin{align*}
f(x,y)&=
\begin{pmatrix}
x&y
\end{pmatrix}
\cdot
\begin{pmatrix}
\ell&m
\\
m&n
\end{pmatrix}
\cdot
\begin{pmatrix}
x\\y
\end{pmatrix}=
\\
&=\ell\cdot x^2+2\cdot m\cdot x\cdot y+n\cdot y^2,
\end{align*}
that is defined on a $(x,y)$-coordinate plane.
Then after proper rotation of the coordinates, 
the expression for $f$ in the new $(x,y)$-coordinates will be
\begin{align*}
\bar f(x,y)&=
\begin{pmatrix}
x&y
\end{pmatrix}
\cdot
\begin{pmatrix}
k_1&0
\\
0&k_2
\end{pmatrix}
\cdot
\begin{pmatrix}
x\\y
\end{pmatrix}=
\\
&=k_1\cdot x^2+k_2\cdot y^2.
\end{align*}

\spell{\end{multicols}}{}

\section{Analysis}\label{sec:analysis}

The following material is discussed in any course of real analysis, the classical book by Walter Rudin \cite{rudin} is one of our favorites.

\spell{\begin{multicols}{2}}{}

\subsection*{Measurable functions}

A function is called \index{measurable function}\emph{measurable} if the inverse image of any Borel set is a Borel set.
Virtually all functions that naturally appear in geometry are measurable.

The following theorem makes it possible to extend many statements about continuous functions to measurable functions.

\begin{thm}{Lusin's theorem}\label{thm:lusin}
Let $\phi\:[a,b]\to \mathbb{R}$ be a measurable function.
Then for any $\epsilon>0$, there is a continuous function $\psi_\epsilon\:[a,b]\to \mathbb{R}$ that coincides with $\phi$ outside a set of measure at most $\epsilon$.
Moreover, if $\phi$ is bounded above and/or below, then we may assume that so is~$\psi_\epsilon$.  
\end{thm}

\subsection*{Lipschitz condition}

Recall that a function $f$ between metric spaces is called \index{Lipschitz function}\emph{Lipschitz} if there is a constant $L$ such that 
\[\dist{f(x)}{f(y)}{}\le L\cdot\dist{x}{y}{}\]
for all values $x$ and $y$ in the domain of definition of~$f$.

The following theorem makes it possible to extend many results about smooth functions to Lipschitz functions.
Recall that {}\emph{almost all} means all values, with exceptions in a set of zero Lebesgue measure.

\begin{thm}{Rademacher's theorem}\label{thm:rademacher}
Let $f\:[a,b]\to\mathbb{R}$ be a Lipschitz function.
Then its derivative $f'$ is a bounded measurable function defined almost everywhere in $[a,b]$.
Moreover, it satisfies the fundamental theorem of calculus; that is, the following identity 
\[f(b)-f(a)=\int_a^b f'(x)\cdot dx,\]
holds if the integral is understood in the sense of Lebesgue.
\end{thm}

\begin{figure*}[t!]
\begin{minipage}{.48\textwidth}
\centering
\includegraphics{mppics/pic-320}
\end{minipage}\hfill
\begin{minipage}{.48\textwidth}
\centering
\includegraphics{mppics/pic-321}
\end{minipage}
\end{figure*}

\begin{figure*}[t!]
\begin{minipage}{.48\textwidth}
\centering
\includegraphics{mppics/pic-325}
\end{minipage}\hfill
\begin{minipage}{.48\textwidth}
\centering
\includegraphics{mppics/pic-326}
\end{minipage}
\end{figure*}

\subsection*{Uniform continuity and convergence}

Let $f\:{\spc{X}}\to \spc{Y}$ be a map between metric spaces.
If  for any $\epsilon>0$ there is $\delta>0$ such that 
\[\dist{x_1}{x_2}{\spc{X}}<\delta\ \Longrightarrow\ \dist{f(x_1)}{f(x_2)}{\spc{Y}}<\epsilon,\]
then $f$ is called \index{uniformly continuous}\emph{uniformly continuous}.

Evidently, every uniformly continuous function is continuous;
the converse does not hold.
For example, the function $f(x)=x^2$ is continuous, but not uniformly continuous.

\begin{thm}{Heine--Cantor theorem}
Any continuous function defined on a compact metric space is uniformly continuous.
\end{thm}

If the condition above holds for any function $f_n$ in a sequence, and $\delta$ depends solely on $\epsilon$,
then the sequence $f_n$ is called \index{uniformly continuous}\emph{uniformly equicontinuous}.
More precisely, 
a sequence of functions $f_n:{\spc{X}}\to \spc{Y}$ is called {}\emph{uniformly equicontinuous} if 
for any $\epsilon>0$ there is $\delta>0$ such that 
\[\dist{x_1}{x_2}{\spc{X}}<\delta\ \Longrightarrow\ \dist{f_n(x_1)}{f_n(x_2)}{\spc{Y}}<\epsilon\]
for any~$n$.

We say that a sequence of functions $f_i\: {\spc{X}} \to \spc{Y}$ \index{uniform convergence}\emph{converges uniformly} to a function $f_{\infty}\: {\spc{X}} \to \spc{Y}$ if, for any 
$\epsilon>0$, there is a natural number $N$ such that for all $n \ge N$, we have $\dist{f_{\infty}(x)}{f_n (x)}{}<\epsilon$
for all $x  \in {\spc{X}}$.

\begin{thm}{Arzel\'{a}--Ascoli theorem}\label{lem:equicontinuous}
Suppose $\spc{X}$ and $\spc{Y}$ are compact metric spaces. 
Then any uniformly equicontinuous sequence of functions $f_n\:\spc{X}\z\to \spc{Y}$ has a subsequence that converges uniformly to a continuous function $f_\infty\:\spc{X}\z\to \spc{Y}$.
\end{thm}

\subsection*{Cutoffs and mollifiers}

Here we explain how to construct smooth functions that mimic the behavior of certain model functions.
These functions are used to smooth model geometric objects keeping their shape nearly unchanged.

As an example, consider the following functions
\begin{align*}
h(t)&\df
\begin{cases}
0&\text{if}\ t\le 0,
\\
t&\text{if}\ t> 0;
\end{cases}
\\
f(t)&\df
\begin{cases}
0&\text{if}\ t\le 0,
\\
\frac{t}{e^{1\!/\!t}}&\text{if}\ t> 0.
\end{cases}
\end{align*}
The functions $h$ and $f$ behave alike ---
both vanish at $t\le 0$ and grow to infinity for positive~$t$.
The function $h$ is not smooth --- its derivative at $0$ is undefined.
Unlike $h$, the function $f$ is smooth.
Indeed, the existence of all derivatives $f^{(n)}(x)$ at $x\ne 0$ is evident, and direct calculations show that $f^{(n)}(0)=0$ for all~$n$.

Other useful examples of that type are the so-called \index{bell function}\emph{bell function} --- a smooth function that is positive in an $\epsilon$-neighborhood of zero and vanishing outside this neighborhood.
Such functions can be obtained from the function $f$ constructed above, say 
\[b_\epsilon(t)\df c\cdot f(\epsilon^2-t^2);\]
typically one chooses the constant $c$ so that $\int b_\epsilon=1$.

Another useful example is a \index{sigmoid}\emph{sigmoid} --- nondecreasing function that vanishes for $t\le -\epsilon$ and takes value $1$ for any $t\ge \epsilon$.
It mimics the step function and can be defined the following way: \label{page:sigma-function}
\[\sigma_\epsilon(t)
\df 
\int_{-\infty}^t b_\epsilon(x)\cdot dx.\]

\spell{\end{multicols}}{}

\section{Multivariable  calculus}\label{sec:Multivariable calculus}

The following material is discussed in any course of multivariable calculus;
for example, in the already mentioned book by Walter Rudin \cite{rudin}.

\spell{\begin{multicols}{2}}{}

\subsection*{Regular values}

Let $\Omega\subset \mathbb{R}^m$ be an open subset.
Any map $\bm{f}\:\Omega\to\mathbb{R}^n$ can be thaut of as an array of its coordinate functions
\[f_1,\dots,f_n\:\Omega\to \mathbb{R}.\]
The map $\bm{f}$ is called \index{smooth!map}\emph{smooth} if each function $f_i$ is smooth;
that is, all partial derivatives of $f_i$ are defined in $\Omega$.

The \emph{Jacobian matrix} of $\bm{f}$ at $\bm{x}\in\mathbb{R}^m$ is defined as \index{10j@$\Jac$ (Jacobian matrix)}
\[\Jac_{\bm{x}}\bm{f}=
\begin{pmatrix}
\dfrac{\partial f_1}{\partial x_1} & \cdots & \dfrac{\partial f_1}{\partial x_m}\\
\vdots & \ddots & \vdots\\
\dfrac{\partial f_n}{\partial x_1} & \cdots & \dfrac{\partial f_n}{\partial x_m} \end{pmatrix};\]
we assume the right-hand side is evaluated at $\bm{x}=(x_1,\dots,x_m)$.

If the Jacobian matrix defines a surjective linear map $\mathbb{R}^m\to\mathbb{R}^n$ (that is, if $\rank(\Jac_{\bm{x}}\bm{f})=n$), then we say that 
$\bm{x}$ is a \index{regular!point}\emph{regular point} of~$\bm{f}$.

If $\bm{x}$ is regular anytime $\bm{f}(\bm{x})=\bm{y}$,
then we say that $\bm{y}$ is a \index{regular!value}\emph{regular value} of $\bm{f}$.
The following lemma states that \textit{most} values of a smooth map are regular.

\begin{thm}{Sard's lemma}\label{lem:sard}
Given a smooth map $\bm{f}\colon \Omega \z\to \mathbb{R}^n$ defined on an open set $\Omega\subset \mathbb{R}^m$, almost all values in $\mathbb{R}^n$ are regular.
\end{thm}

The words \index{almost all}\emph{almost all} mean all values, with the possible exceptions belonging to a set with vanishing Lebesgue measure.
In particular, if one chooses a random value equidistributed in an arbitrarily small ball $B\z\subset \mathbb{R}^n$, then it is a regular value of $\bm{f}$ with probability~$1$.

Note that if $m<n$, then any point $\bm{y}=\bm{f}(\bm{x})$ is \emph{not} a regular value of $\bm{f}$.
Therefore, the only regular values of $\bm{f}$ are the points in the complement of the image $\Im \bm{f}$.
In this case, the theorem states that almost all points in $\mathbb{R}^n$ do \textit{not} belong to $\Im \bm{f}$.

\subsection*{Inverse function theorem}

The \index{inverse function theorem}\emph{inverse function theorem} gives an if condition for a map $\bm{f}$ to be invertible in a neighborhood of a given point $\bm{x}$.
The condition is formulated in terms of $\Jac_{\bm{x}}\bm{f}$ --- the Jacobian matrix of $\bm{f}$ at $\bm{x}$.

The \index{implicit function theorem}\emph{implicit function theorem} is its close relative and a corollary.
It is used when we need to pass from parametric to implicit descriptions of curves and surfaces.

Both theorems reduce the existence of a map satisfying a certain equation to a question in linear algebra.
We use these two theorems only for $n\le 3$.

\begin{thm}{Inverse function theorem}\label{thm:inverse}
Let $\bm{f}\z=(f_1,\dots,f_n)\:\Omega\to\mathbb{R}^n$ be a smooth map
defined on an open set $\Omega\subset \mathbb{R}^n$.
Assume that $\Jac_{\bm{x}}\bm{f}$
is invertible at $\bm{x}\in \Omega$.
Then there is a smooth map $\bm{h}\:\Phi\to\mathbb{R}^n$ defined in an open neighborhood $\Phi$ of ${\bm{y}}\z=\bm{f}(\bm{x})$ that is a {}\emph{local inverse of $\bm{f}$ at $\bm{x}$};
that is, there is a neighborhood $\Psi\ni \bm{x}$ such that
$\bm{f}$ defines a homeomorphism $\Psi\leftrightarrow \Phi$, and
$\bm{h} \circ \bm{f}$ is an identity map on $\Psi$.

Moreover, if $|\det[\Jac_{\bm{x}}\bm{f}]|>\epsilon>0$, the domain $\Omega$ contains an $\epsilon$-neighborhood of $\bm{x}$, 
and the first and second partial derivatives $\tfrac{\partial f_i}{\partial x_j}$, $\tfrac{\partial^2 f_i}{\partial x_j\partial x_k}$ are bounded by a constant $C$ for all $i$, $j$, and $k$, then we can assume that $\Phi$ is a $\delta$-neighborhood of $\bm{y}$, for some $\delta>0$ that depends only on $\epsilon$ and~$C$. 
\end{thm}

\begin{thm}{Implicit function theorem}\label{thm:imlicit}
Let $\bm{f}\z=(f_1,\dots,f_n)\:\Omega\to\mathbb{R}^n$ be a smooth map defined on an open subset $\Omega\subset\mathbb{R}^{n+m}$, where
$m,n\z\ge 1$.
Regard $\mathbb{R}^{n+m}$ as a product space $\mathbb{R}^n\times \mathbb{R}^m$ with coordinates
$(\bm{x},\bm{y})\z=(x_1,\dots,x_n,y_1,\dots,y_m)$.
Consider the following matrix 
\[
M=\begin{pmatrix}
\dfrac{\partial f_1}{\partial x_1} & \cdots & \dfrac{\partial f_1}{\partial x_n}\\
\vdots & \ddots & \vdots\\
\dfrac{\partial f_n}{\partial x_1} & \cdots & \dfrac{\partial f_n}{\partial x_n} \end{pmatrix}\]
formed by the first $n$ columns of $\Jac_{(\bm{x},\bm{y})}\bm{f}$.
Assume $M$ is invertible at $(\bm{x},\bm{y})\z\in \Omega$ and $\bm{f}(\bm{x},\bm{y})\z=0$.
Then there is a smooth function $\bm{h}\:\Phi\z\to\mathbb{R}^n$ defined on a neighborhood $\Phi\ni \bm{y}$ in $\mathbb{R}^{m}$
and a neighborhood $\Psi\ni (\bm{x},\bm{y})$ in $\mathbb{R}^n\times \mathbb{R}^m$
such that
for any $(x_1,\dots,x_n,y_1,\dots, y_m)\z\in \Psi$, the equality
\[\bm{f}(x_1,\dots,x_n,y_1,\dots, y_m)=\bm{0}\]
holds if and only if 
\[(x_1,\dots,x_n)=\bm{h}(y_1,\dots, y_m).\]

\end{thm}

\subsection*{Multiple integral}

Let $\bm{f}\:\mathbb{R}^n\to\mathbb{R}^n$ be a smooth map (maybe partially defined).
Set
\[\jac_{\bm{x}}\bm{f}\df|\det[\Jac_{\bm{x}}\bm{f}]|;
\index{10j@$\jac$ (Jacobian determinant)}\]
that is, $\jac_{\bm{x}}\bm{f}$ is the absolute value of the determinant of the Jacobian matrix of $\bm{f}$ at $\bm{x}$.

The following theorem plays the role of a substitution rule for multiple variables.

\index{Borel subsets}\emph{Borel subsets} are defined as the class of subsets that are generated from open sets by applying the following operations recursively: countable union, countable intersection, and complement.
Since the complement of a closed set is open and the other way around, these sets can be also generated from all closed sets.
This class of sets includes virtually all sets that naturally appear in geometry but does not include pathological examples that create problems with integration.

\begin{thm}{Theorem}\label{thm:mult-substitution} 
Let $h\:K\to\mathbb{R}$ be a continuous function on a Borel subset $K\subset \mathbb{R}^n$.
Assume $\bm{f}\:\Omega\to \mathbb{R}^n$ is an injective smooth map that is defined on an open set $\Omega\supset K$.
Then 
\[\idotsint_{\bm{x}\in K} h(\bm{x})\cdot \jac_{\bm{x}}\bm{f}
=
\idotsint_{\bm{y}\in \bm{f}(K)} h\circ \bm{f}^{-1}(\bm{y}).\]

\end{thm}

\subsection*{Convex functions}

The following statements will be used only in dimensions $\le 3$.

A function $f\:\mathbb{R}^n\to \mathbb{R}$ is called \index{convex!function}\emph{(strictly) convex} (respectively, \emph{concave}) if
its epigraph $z\ge f(\bm{x})$ (respectively, its subgraph $z\le f(\bm{x})$) is a (strictly) convex set in $\mathbb{R}^n\times \mathbb{R}$.

Let $f\:\mathbb{R}^n\to \mathbb{R}$ be a smooth function (maybe partially defined).
Choose a vector $\vec w\in \mathbb{R}^n$.
Given a point $p\in\mathbb{R}^n$ consider the function $\phi(t)=f(p+t\cdot \vec w)$.
Then the \index{directional derivative}\emph{directional derivative} $(D_{\vec w}f)(p)$ of $f$ at $p$ with respect to vector $\vec w$ is defined by
\[(D_{\vec w}f)(p)=\phi'(0).\]

\begin{thm}{Theorem}\label{thm:Jensen}
A smooth function $f\:K\z\to \mathbb{R}$ defined on a convex subset $K\subset\mathbb{R}^n$ is convex if and only if one of the following equivalent conditions holds:

\begin{subthm}{}
The second directional derivative of $f$ at any point in the direction of any vector is nonnegative; that is,
\[(D_{\vec w}^2f)(p)\ge 0\]
for any $p\in K$ and $\vec w\in\mathbb{R}^n$.
\end{subthm}

\begin{subthm}{}
The so-called \index{Jensen's inequality}\emph{Jensen's inequality}
\begin{align*}
f ((1-t)\cdot x_0 &+ t\cdot x_1 ) \le
\\
&\le (1-t)\cdot f(x_0)+ t\cdot f(x_1)
\end{align*}
holds for any $x_0,x_1\in K$ and $t\in[0,1]$.

\end{subthm}

\begin{subthm}{}
For any $x_0,x_1\in K$, we have 
\[f \left (\frac{x_0 + x_1}2 \right ) \le \frac{f(x_0) + f(x_1)}2.\]
\end{subthm}

\end{thm}

\spell{\end{multicols}}{}

\section{Ordinary differential equations}

The following material is discussed at the beginning of any course of ordinary differential equations; the classical book by Vladimir Arnold \cite{arnold} is one of our favorites.

\spell{\begin{multicols}{2}}{}

\subsection*{First-order equations}

The following theorem guarantees existence and uniqueness of solutions of an initial value problem
for a system of ordinary first-order differential equations
\[
\begin{cases}
x_1'&=f_1(x_1,\dots,x_n,t),
\\
&\,\,\vdots
\\
x_n'&=f_n(x_1,\dots,x_n,t),
\end{cases}
\]
where each $t\mapsto x_i=x_i(t)$ is a real-valued function defined on a real interval $\mathbb{J}$
and each $f_i$ is a smooth function defined on an open subset $\Omega\subset \mathbb{R}^n\times \mathbb{R}$.

The array of functions $(f_1,\dots,f_n)$ can be packed into one vector-valued function 
$\bm{f}\:\Omega\to \mathbb{R}^n$;
the same way the function array $(x_1,\dots,x_n)$ can be packed into a vector-valued function  $\bm{x}\:\mathbb{J}\to\mathbb{R}^n$.
Therefore, the system can be rewritten as one vector equation 
\[\bm{x}'=\bm{f}(\bm{x}, t).\] 

\begin{thm}{Theorem}\label{thm:ODE}
Suppose $\bm{f}\:\Omega\to \mathbb{R}^n$ is a smooth function defined on an open subset $\Omega\z\subset \mathbb{R}^n\times \mathbb{R}$.
Then for any initial data $\bm{x}(t_0)=\bm{u}$ such that $(\bm{u},t_0)\in\Omega$ the differential equation
\[\bm{x}'=\bm{f}(\bm{x},t)\]
has a unique solution $t\mapsto \bm{x}(t)$ defined at a maximal interval $\mathbb{J}$ that contains $t_0$.
Moreover,
\begin{enumerate}[(a)]
\item  if $\mathbb{J}\ne \mathbb{R}$ (that is, if an end $b$ of $\mathbb{J}$ is finite), then $\bm{x}(t)$ does not have a limit point in $\Omega$ as $t\to b$;%
\footnote{In other words, if $\bm{x}(t_n)$ converges for a sequence $t_n\to b$, then its limit does not lie in $\Omega$.}
\item  the function $w\:(\bm{u},t_0,t)\mapsto \bm{x}(t)$ has an open domain of definition in $\Omega\times \mathbb{R}$ that contains all points $(\bm{u},t_0,t_0)$ for $(\bm{u},t_0)\z\in \Omega$, and $w$ is smooth in this domain.
\end{enumerate}

\end{thm}

\subsection*{Higher order}

Consider an ordinary differential equation of order $k$
\[\bm{x}^{(k)}=\bm{f}(\bm{x},\bm{x}',\dots,\bm{x}^{(k-1)},t),\eqlbl{eq:nth-order}\]
where $t\mapsto\bm{x}=\bm{x}(t)$ is a function from a real interval to $\mathbb{R}^n$.

This equation can be rewritten as $k$ first-order equations  with $k-1$ new vector variables 
$\bm{y}_1=\bm{x}'$,
$\bm{y}_2=\bm{x}'',\dots,\bm{y}_{k-1}=\bm{x}^{(k-1)}$:
\[
\begin{cases}
\bm{x}'(t)&=\bm{y}_1(t),
\\
\bm{y}_1'(t)&=\bm{y}_2(t),
\\
&\,\,\vdots
\\
\bm{y}_{k-2}'(t)&=\bm{y}_{k-1}(t),
\\
\bm{y}_{k-1}'(t)&=\bm{f}(\bm{x},\bm{y}_{1},\dots,\bm{y}_{k-1},t).
\end{cases}
\eqlbl{eq:nth-order-new}
\]

Thus, we have derived the following.

\begin{thm}{Theorem}\label{thm:ODE-nth-order}
The $k$-th order equation \ref{eq:nth-order} is equivalent to the system \ref{eq:nth-order-new} consisting of $k$ first-order equations.
\end{thm}

This trick reduces a higher-order ordinary differential equation to a first-order equation. 
In particular, we get local existence and uniqueness for solutions of higher-order equations as in Theorem \ref{thm:ODE};
one only has to assume that
$\Omega\z\subset \mathbb{R}^{n\cdot k}\times \mathbb{R}$, and the initial data consists of $\bm{x}(t_0)$, $\bm{x}'(t_0),\dots,\bm{x}^{(k-1)}(t_0)$.

\spell{\end{multicols}}{}
}

%% file: sol-new.tex
\stepcounter{chapter}
\setcounter{eqtn}{0}

%\raggedcolumns\setlength{\multicolsep}{-3mm}
\spell{\begin{multicols}{2}}{}

\chapter[Semisolutions]{\vspace*{-23mm} Semisolutions\vspace*{-12mm}}

%\end{multicols}
%\par\noindent\rule{\textwidth}{0.4pt}
%\begin{multicols}{2}

\stepcounter{chapter}
\setcounter{eqtn}{0}

\parbf{\ref{ex:9}}; \ref{SHORT.ex:9:compact}. Use that a continuous injection defined on a compact domain is an embedding (\ref{thm:Hausdorff-compact}). 

\parit{\ref{SHORT.ex:9:9}.} The image of $\gamma$ might have the shape of the digit $8$ or $9$.

\parbf{\ref{ex:mono}.}
Push neighborhoods of the ends to the ends.

\parbf{\ref{aex:simple-curve}.}
Let $\alpha$ be a path, connecting $p$ to~$q$.

Passing to a subarc of $\alpha$,
we can assume that $\alpha(t)\ne p,q$ for $t\ne0,1$.

An open set $\Omega$ in $(0,1)$ will be called {}\emph{suitable}
if, for any connected component $(a,b)$ of $\Omega$, we have $\alpha(a)=\alpha(b)$.
Show that the union of nested suitable sets is suitable.
Therefore, we can find a maximal suitable set $\hat \Omega$.

Define $\beta(t)=\alpha(a)$ for any $t$ in a connected component $(a,b)\subset\hat \Omega$, and $\beta (t) = \alpha (t) $ for $t\notin\hat{\Omega}$.
Note that for any $x\in [0,1]$ the set $\beta^{-1}\{\beta(x)\}$ is connected.

It remains to show that a reparametrization of $\beta$ is a simple path.
In order to do that, we need to construct a non-decreasing surjective function $\tau\:[0,1]\z\to[0,1]$ such that 
$\tau(t_1)\z=\tau(t_2)$ if and only if there is a connected component $(a,b)\subset\hat \Omega$ such that $t_1,t_2\z\in [a,b]$.

The required function $\tau$ can be constructed similarly to the so-called {}\emph{devil's staircase} --- learn this construction and modify it.

\parbf{\ref{ex:L-shape}.}
Denote the union of two half-axes by~$L$.

Observe that $f(t)\to\infty$ as $t\to \infty$.
Since $f(0)=0$, the intermediate value theorem implies that $f(t)$ takes all nonnegative values for $t\ge 0$.
Use it to show that $L$ is the range of~$\alpha$.

Further, show that the function $f$ is strictly increasing for $t> 0$.
Use this to show that the map $t\mapsto \alpha(t)$ is injective.

Summarizing, we get that $\alpha$ is a smooth parametrization of~$L$. 

Now suppose $\beta\:t\z\mapsto (x(t),y(t))$ is a smooth parametrization of~$L$.
Without loss of generality, we may assume that $x(0)\z=y(0)=0$.
Note that $x(t)\ge 0$ for any $t$ therefore $x'(0)=0$.
The same way we get that $y'(0)\z=0$.
That is, $\beta'(0)=0$;
so $L$ does not admit a smooth \textit{regular} parametrization.

\parbf{\ref{ex:cycloid}.}
Apply the definitions.
For \ref{SHORT.ex:cycloid:regular} you need to check that $\gamma'_\ell\ne 0$.
For \ref{SHORT.ex:cycloid:simple} you need to check that $\gamma_\ell(t_0)\z=\gamma_\ell(t_1)$ only if $t_0=t_1$.

\parbf{\ref{ex:nonregular}.}
Note that the parametrization $t\mapsto (t,t^3)$ is smooth and regular.
Modify it so it has zero speed at a point.

\begin{wrapfigure}{r}{20 mm}
\vskip-6mm
\centering
\includegraphics{mppics/pic-270}
\vskip-3mm
\end{wrapfigure}

\parbf{\ref{ex:y^2=x^3}.}
This is the so-called \index{semicubical parabola}\emph{semicubical parabola}; it is shown on the picture.
Try to argue similarly to \ref{ex:L-shape}.

\parbf{\ref{ex:viviani}.}
For $\ell=0$ the system describes a pair of points $(0,0,\pm1)$, so we can assume that $\ell\ne 0$.

The first equation describes the unit sphere centered at the origin, and the second equation describes a cylinder over the circle in the $(x,y)$-plane with center at $(-\tfrac\ell2,0)$ and radius~$|\tfrac\ell2|$.

Find the gradients $\nabla f$ and $\nabla h$ for the functions
\begin{align*}
 f(x,y,z)&=x^2+y^2+z^2-1,
 \\
 h(x,y,z)&=x^2+\ell\cdot x+y^2.
\end{align*}
Show that for $\ell\ne 0$,
the gradients are linearly dependent only on the $x$-axis.
Conclude that for $\ell\ne\pm 1$ each connected component of the set of solutions is a smooth curve.

Show that 
\begin{itemize}
\item if $|\ell|<1$, then the set has two connected components with $z>0$ and $z<0$.
\item if $|\ell|\ge1$, then the set is connected.
\end{itemize}

The linear independence of the gradients provides only a sufficient condition.
Therefore, the case $\ell=\pm1$ has to be checked by hand.
In this case, a neighborhood of $(\pm1,0,0)$ does not admit a smooth regular parametrization --- try to prove it. 
The case $\ell=1$ is shown on the picture.

\begin{Figure}
\centering
\vskip-0mm
\begin{lpic}[t(2mm),b(0mm),r(0mm),l(0mm)]{asy/viviani(1)}
\lbl[r]{-.5,18;$x$}
\lbl[l]{41,22;$y$}
\lbl[r]{18,54;$z$}
\end{lpic}
\end{Figure}

\parit{Remark.}
In the case $\ell=\pm1$, the curve is called {}\emph{Viviani's curve}.
It admits the following smooth regular parametrization with a self-intersection at $(\pm1,0,0)$
\[t\mapsto(\pm(\cos t)^2,\cos t\cdot\sin t,\sin t).\]

\parbf{\ref{ex:open-curve}.}
Assume that $|\gamma(t)|\z\to\infty$ as $t\to\pm\infty$.
Choose a compact set $K\subset \mathbb{R}^3$.
Show that $\gamma^{-1}(K)$ is a bounded closed set in $\mathbb{R}$,
and apply the Heine--Borel lemma (\ref{thm:Heine--Borel}).
Conclude that $\gamma$ is proper.

Now assume that $\gamma(t_n)$ converges for some sequence $t_n\to \pm \infty$; let $p$ be its limit and $K$ be a closed ball centered at~$p$.
Show that the inverse image $\gamma^{-1}(K)$ is not compact.
Conclude that $\gamma$ is not proper.

\parbf{\ref{ex:proper-closed}.}
Show and use that a set $C\subset \mathbb{R}^3$ is closed if and only if the intersection $K\cap C$ is compact for any compact $K\subset \mathbb{R}^3$.

\parbf{\ref{ex:proper-curve}.}
Without loss of generality, we may assume that the origin does not lie on the curve.

Show that inversion of the plane $(x,y)\z\mapsto (\tfrac{x}{x^2+y^2},\tfrac{y}{x^2+y^2})$ maps our curve to a closed curve with the origin removed.
Apply Jordan's theorem for the obtained curve, and use the inversion again.

%\end{multicols}
%\par\noindent\rule{\textwidth}{0.4pt}
%\begin{multicols}{2}

\stepcounter{chapter}
\setcounter{eqtn}{0}

\parbf{\ref{ex:integral-length-0}.}
Show that if we take the least upper bound in \ref{def:length} for all sequences
$a=t_0\le t_1\le\z\dots\le t_k=b$, then the result is the same.

Suppose $\gamma_2$ is reparametrization of $\gamma_1$ by $\tau\:[a_1,b_1]\to [a_2,b_2]$;
without loss of generality, we may assume that $\tau$ is nondecreasing.
Set $\theta_i=\tau(t_i)$.
Observe that $a_2=\theta_0\z\le\theta_1\le \z\dots\le\theta_k=b_2$ if and only if 
$a_1=t_0\le t_1\le\z\dots\le t_k=b_1$.
Make a conclusion.

\parbf{\ref{ex:length-chain}.}
Show that for any inscribed polygonal line $\beta$ and any $\epsilon>0$ we have
\[\length\beta_n>\length\beta-\epsilon.\]
for all large $n$.
Conclude that
\[\liminf_{n\to\infty}\length\beta_n\ge \length \gamma.\]
Use the definition of length to show that 
\[\limsup_{n\to\infty}\length\beta_n\le \length \gamma.\]
Observe that the two obtained inequalities imply the required statement.

\parbf{\ref{ex:length-image}.}
Choose a partition $0=t_0<\dots <t_n=1$ of $[0,1]$.
Set $\tau_0=0$ and 
\[\tau_i=\max\set{\tau \in[0,1]}{\beta(\tau_i)=\gamma(t_i)}\]
for $i>0$.
Show that $(\tau_i)$ is a partition of $[0,1]$;
that is, $0=\tau_0<\tau_1<\dots<\tau_n=1$.

By construction 
\begin{align*}
|\gamma(t_0)&-\gamma(t_1)|+|\gamma(t_1)-\gamma(t_2)|+\dots
\\
&\qquad\dots+|\gamma(t_{n-1})-\gamma(t_n)|=
\\
&=
|\beta(\tau_0)-\beta(\tau_1)|+|\beta(\tau_1)-\beta(\tau_2)|+\dots
\\
&\qquad\dots+|\beta(\tau_{n-1})-\beta(\tau_n)|.
\end{align*}
Since the partition $(t_i)$ is arbitrary, we get 
\[\length \beta\ge \length \gamma.\]

\parit{Remarks.}
It is instructive to compare this exercise with \ref{obs:S2-length}.

The inequality might be strict.
It happens if $\beta$ runs back and forth along $\gamma$.
In this case the partition $(\tau_i)$ above is \textit{not} arbitrary.

If one assumes only $\beta([0,1])\supset\gamma([0,1])$, then the problem becomes harder.
It follows from the following inequality \cite[2.6.1+2.6.2]{burago-burago-ivanov}:
\[h\le \length\gamma.\]
where $h$ denotes the 1-dimensional Hausdorff measure of the image $\gamma([0,1])$.
Moreover the equality holds if and only if $\gamma$ is simple.

\parbf{\ref{ex:integral-length}}; \ref{SHORT.ex:integral-length>}.
Apply the fundamental theorem of calculus for each segment in a given partition.

\parit{\ref{SHORT.ex:integral-length<}.} Consider a partition such that the velocity vector $\alpha'(t)$ is nearly constant on each of its segments.

\parbf{\ref{adex:integral-length}.}
Use the theorems of Rademacher and Lusin (\ref{thm:rademacher} and \ref{thm:lusin}).

\parbf{\ref{ex:nonrectifiable-curve}}; \ref{SHORT.ex:nonrectifiable-curve:a}.
Look at the picture, and guess the parametrization of an arc of the snowflake by $[0,1]$.
Extend it to the whole snowflake.
Show that it indeed describes an embedding of the circle into the plane.

\begin{Figure}
\vskip-0mm
\centering
\includegraphics{mppics/pic-226}
\vskip0mm
\end{Figure}

\parit{\ref{SHORT.ex:nonrectifiable-curve:b}.}
Suppose $\gamma\:[0,1]\to\mathbb{R}^2$ is a rectifiable curve and let $\gamma_k$ be a scaled copy of $\gamma$ with factor $k>0$;
that is, $\gamma_k(t)=k\cdot\gamma(t)$ for any~$t$.
Show that 
\[\length\gamma_k=k\cdot\length \gamma.\]

Now suppose the arc $\gamma$ of the Koch snowflake shown on the picture is rectifiable; denote its length by $\ell$.
Observe that $\gamma$ can be divided into 4 arcs such that each one is a scaled copy of $\gamma$ with factor $\tfrac13$.
It follows that $\ell=\tfrac43\cdot\ell$.
Evidently, $\ell>0$ --- a contradiction.

\parbf{\ref{ex:cont-length}.}
Use \ref{thm:length-semicont} to show that $s$ is lower semicontinuous;
that is, if $t_n\to t_\infty$ as $n\to\infty$, then 
\[\liminf_{n\to\infty} s(t_n)\ge s(t_\infty).\]
Observe that
\[s(t)=s(b)-\length(\gamma|_{[t,b]}).\]
Apply \ref{thm:length-semicont} to show that $s$ is upper semicontinuous.
Conclude that $s$ is continuous.
Finally, show and use that $s$ is nondecreasing.

\parbf{\ref{ex:arc-length-helix}.} 
We have to assume $a\ne 0$ or $b\ne0$;
otherwise, we get a constant curve.

Show that $|\gamma'(t)|\equiv \sqrt{a^2+b^2}$;
in particular, the velocity is constant.
Therefore, $s\z=t/\sqrt{a^2+b^2}$ is an arc-length parameter.
It remains to substitute $s\cdot \sqrt{a^2+b^2}$ for~$t$.

\parbf{\ref{ex:convex-hull}.}
Choose a closed polygonal line $p_1\dots p_n$ inscribed in $\beta$.
By \ref{cor:convex=>rectifiable}, we can assume that its length is arbitrarily close to the length of $\beta$;
that is, given $\epsilon>0$ 
\[\length (p_1\dots p_n)>\length\beta-\epsilon.\]

Show that we may assume in addition that each point $p_i$ lies on $\alpha$.

Since $\alpha$ is simple, the points $p_1,\dots,p_n$ appear on $\alpha$ in the same cyclic order;
that is, the polygonal line $p_1\dots p_n$ is also inscribed in $\alpha$.
In particular,
\[\length\alpha\ge \length (p_1\dots p_n).\]
It follows that 
\[\length\alpha>\length\beta-\epsilon.\]
for any $\epsilon>0$.
Whence 
\[\length\alpha\ge\length\beta.\]

\begin{wrapfigure}{r}{25 mm}
\vskip-6mm
\centering
\includegraphics{mppics/pic-275}
\vskip0mm
\end{wrapfigure}

If $\alpha$ has self-intersections, then the points $p_1,\dots, p_n$ might appear on $\alpha$ in a different cyclic order, say $p_{i_1},\dots,p_{i_n}$.
Apply the triangle inequality to show that 
\[\length(p_{i_1}\dots p_{i_n})\ge \length (p_1\dots p_n)\]
and use it to modify the given proof.

\parbf{\ref{ex:convex-croftons}.} 
Denote by $\ell_{\vec u}$ the line segment 
obtained by orthogonal projection of $\gamma$ to the line in the direction ${\vec u}$.
Since $\gamma_{\vec u}$ runs along $\ell_{\vec u}$ back and forth, we get 
\[\length\gamma_{\vec u}\ge 2\cdot\length\ell_{\vec u}.\]
By the Crofton formula, 
\[\length\gamma\ge \pi\cdot \overline{\length\ell_{\vec u}}.\]

In the case of equality, the curve $\gamma_{\vec u}$ runs exactly back and forth along $\ell_{\vec u}$ without additional zigzags for almost all (and therefore for all) ${\vec u}$.

Let $K$ be a closed set bounded by~$\gamma$.
The last statement implies that every line may intersect $K$ only along a closed segment or a single point.
It follows that $K$ is convex.

\parbf{\ref{adex:more-croftons}.}
The proof is identical to the proof of the standard Crofton formula.
To find the coefficients it is sufficient to check it on a unit interval.
The latter can be done by integration:
\begin{align*}
\frac1{k_1}&=\frac{1}{\area \mathbb{S}^2}\cdot\iint_{\mathbb{S}^2} |x|;
\\
\frac1{k_2}&=\frac{1}{\area \mathbb{S}^2}\cdot\iint_{\mathbb{S}^2} \sqrt{1-x^2}.
\end{align*}
The answers are $k_1=2$ and $k_2=\tfrac4\pi$.

\parbf{\ref{ex:induced-is-length}.}
Let $d$ be the original metric on $\spc{X}$.
Suppose that $\gamma\:[a,b]\to (\spc{X},\ell)$ is continuous.
Show that so is $\gamma\:[a,b]\to (\spc{X},d)$.

Choose a curve $\gamma$ in $(\spc{X},\ell)$.
Denote by $r$ and $s$ its $d$-length and $\ell$-length, respectively.
Show that $s\ge r\ge s-\epsilon$ for any $\epsilon>0$.
Deduce that any curve in $(\spc{X},\ell)$ has equal $d$-length and $\ell$-length.
Make a conclusion.

\parit{Remark.} The converse to the first statement does not hold in general.
Consider for example, the following metric on the upper half-plane 
\[d(x,y)=\min\{|x|+|y|,\bigl||x|-|y|\bigr|+\sqrt{\measuredangle(x,y)}\}.\]  

\parbf{\ref{ex:intrinsic-convex}.}
The only-if part is trivial.
To show the if part, assume $A$ is not convex;
that is, there are points $x,y\in A$ and a point $z\notin A$ that lies between $x$ and~$y$.

Since $A$ is closed, its complement is open.
That is, the ball $B(z,\epsilon)$ does not intersect $A$ for some $\epsilon>0$.

Show that there is $\delta>0$ such that any path from $x$ to $y$ of length at most $\dist{x}{y}{\mathbb{R}^3}+\delta$ passes thru $B(z,\epsilon)$.
It follows that $\dist{x}{y}A\ge \dist{x}{y}{\mathbb{R}^3}+\delta$; 
in particular, $\dist{x}{y}A\ne \dist{x}{y}{\mathbb{R}^3}$.

\begin{wrapfigure}{r}{23 mm}
\vskip-0mm
\centering
\includegraphics{mppics/pic-280}
\vskip0mm
\end{wrapfigure}

\parbf{\ref{ex:antipodal}.}
The spherical curve shown on the picture does not have antipodal pairs of points.
However, there are three points $x,y,z$ on an equator that lie on one of its sides and their antipodal points $-x,-y,-z$ on the other.
(We assume that the points $x, -y, z, -x, y, -z$ lie in the same order on the equator.)

Show that such a curve cannot lie in a hemisphere.

\parbf{\ref{ex:bisection-of-S2}.}
Assume the contrary, then $\gamma$ lies in an open hemisphere (\ref{lem:hemisphere}).
In particular, it cannot divide $\mathbb{S}^2$ into two regions of equal area --- a contradiction.

\parbf{\ref{ex:flaw}.}
The first sentence is wrong --- it is \textit{not} sufficient to show that the diameter is at most~2.
For example, if an equilateral triangle has circumradius slightly above $1$,
then its diameter (which is defined as the maximal distance between its points) is slightly bigger than $\sqrt3$, so it is smaller than $2$.

On the other hand, it is easy to modify the proof of the hemisphere lemma (\ref{lem:hemisphere}) to get a correct solution.
That is, (1) choose two points $p$ and $q$ on $\gamma$ that divide it into two arcs of the same length;
(2) set $z$ to be a midpoint of $p$ and $q$,
and (3) show that $\gamma$ lies in the unit disc centered at~$z$.

\parbf{\ref{adex:crofton}}; \ref{SHORT.adex:crofton:crofton}.
Modify the proof of the original Crofton formula
(Section~\ref{sec:crofton}).

\parit{\ref{SHORT.adex:crofton:hemisphere}.}
Assume $\length \gamma<2\cdot\pi$.
By \ref{SHORT.adex:crofton:crofton},
\[\overline{\length \gamma^*_{\vec u}}<2\cdot\pi.\]
Therefore, we can choose ${\vec u}$ so that 
\[\length \gamma^*_{\vec u}<2\cdot\pi.\]

Observe that $\gamma^*_{\vec u}$ runs in a semicircle~$h$.
Therefore $\gamma$ lies in a hemisphere with $h$ as a diameter.

%\end{multicols}
%\par\noindent\rule{\textwidth}{0.4pt}
%\begin{multicols}{2}

\stepcounter{chapter}
\setcounter{eqtn}{0}

\parbf{\ref{ex:zero-curvature-curve}.}
Let $\gamma$ be a unit-speed curve with vanishing curvature.
Then $\gamma''\equiv 0$; therefore, the velocity $\vec v=\gamma'$ is constant.
Hence, $\gamma(t)=p+(t-t_0)\cdot \vec v$, where $p=\gamma(t_0)$.

\parbf{\ref{ex:scaled-curvature}.} 
Observe that $\alpha(t)\df\gamma_{\lambda}(t/\lambda)$ is a unit-speed parametrization of the curve $ \gamma_{\lambda}$.
Apply the chain rule twice.

\parbf{\ref{ex:curvature-of-spherical-curve}.} Differentiate the identity $\langle\gamma(s),\gamma(s)\rangle=1$ a couple of times.

\parbf{\ref{ex:curvature-formulas}.} 
Let $\tan=\tfrac{\gamma'}{|\gamma'|}$.
Prove and use the following identities: 
\begin{align*}
\gamma''-(\gamma'')^\perp&=\tan\cdot\langle\gamma'',\tan\rangle,
&
|\gamma'|&=\sqrt{\langle \gamma',\gamma'\rangle}.
\end{align*}

\parbf{\ref{ex:curvature-graph}.} 
Apply \ref{ex:curvature-formulas:a} for the parametrization $t\z\mapsto (t,f(t))$.

\parbf{\ref{ex:approximation-const-curvature}.}
Without loss of generality, we may assume that $\gamma$ has a unit-speed parametrization.

Consider the tangent indicatrix $\tan(s)\z=\gamma'(s)$.
Note that $\tan$ is a spherical curve and $|\tan'|\le 1$.
Use it to construct a sequence of unit-speed spherical curves $\tan_n\:\mathbb{I}\to\mathbb{S}^2$ such that $\tan_n(s)\to \tan(s)$ as $n\to\infty$ for any~$s$.
Show that the following sequence of curves solves the problem:
\[\gamma_n(s)=\gamma(a)+\int_a^s\tan_n(t)\cdot dt.\]

\parbf{\ref{ex:no-parallel-tangents}.}
Reuse the construction in \ref{ex:antipodal} to get a closed smooth simple spherical curve $\tan\:\mathbb{S}^1\z\to\mathbb{S}^2$ that contains no pair of antipodal points and has zero average.
Then construct a curve with tangent indicatrix $\tan$.

\parit{Source:} This curve was constructed by Beniamino Segre \cite{segre}.

\parbf{\ref{ex:helix-curvature}.}
Show that $\gamma_{a,b}''\perp \gamma'_{a,b}$, and apply \ref{ex:curvature-formulas:a}.

\parbf{\ref{ex:length>=2pi}.} Apply Fenchel's theorem.

\parbf{\ref{ex:gamma/|gamma|}.}
We can assume that $\gamma$ is unit-speed.
Set $\theta(s)=\measuredangle(\gamma(s),\gamma'(s))$.
Since $\langle\tan,\tan\rangle=1$, we have $\tan'\perp \tan$.
Observe that $\langle \tan,\sigma\rangle=\cos\theta$;
therefore
\begin{align*}
\kur\cdot \sin\theta
&=|\tan'|\cdot \sin\theta\ge
-\langle \tan',\sigma\rangle=
\\
&=
\langle \tan,\sigma'\rangle-\langle \tan,\sigma\rangle'=
\\
&=
(|\sigma'|+\theta')\cdot \sin\theta.
\end{align*}
Whence $\kur\ge |\sigma'|+\theta'$
if $\theta\ne0,\pi$.
It remains to integrate this inequality and show that the set with $\theta\ne0$ or $\pi$ does not create a problem.

\parit{An alternative solution} can be built on \ref{prop:inscribed-total-curvature}.

\parbf{\ref{ex:DNA}}; \ref{SHORT.ex:DNA:c''c>=k}.
Since $|\gamma(s)|\le 1$, we have
\begin{align*}
\langle\gamma''(s),\gamma(s)\rangle&\ge -|\gamma''(s)|\cdot|\gamma(s)|\ge-\kur(s)
\end{align*}
for all~$s$.

\parit{\ref{SHORT.ex:DNA:int>=length-tc}.}
Since $\gamma$ is unit-speed, we have $\ell=\length\gamma$ and $\langle\gamma',\gamma'\rangle\equiv1$.
Therefore,
\begingroup
\allowdisplaybreaks
\begin{align*}
\int_0^\ell\langle\gamma(s),\gamma'(s)\rangle'\cdot ds&=\\
=\int_0^\ell\langle\gamma'(s),\gamma'(s)\rangle\cdot ds&+\int_0^\ell\langle\gamma(s),\gamma''(s)\rangle\cdot ds\ge
\\
&\ge \length\gamma-\tc\gamma.
\end{align*}
\endgroup

\parit{\ref{SHORT.ex:DNA:end}.}
By the fundamental theorem of calculus, we have
\begin{align*}
\int_0^\ell\langle\gamma(s),\gamma'(s)\rangle'\cdot ds
&=\langle\gamma(s),\gamma'(s)\rangle\bigg|_0^\ell.
\end{align*}
Since $\gamma(0)=\gamma(\ell)$ and $\gamma'(0)=\gamma'(\ell)$, the right-hand side vanishes.

We can assume the curve in \ref{thm:DNA} is described by a loop $\gamma\:[0,\ell]\to\mathbb{R}^3$ parametrized by arc-length.
We can also assume that the origin is the center of the ball; that is $|\gamma|\le 1$.
Since $\gamma$ is a smooth closed curve, we have 
$\gamma'(0)=\gamma'(\ell)$ and $\gamma(0)=\gamma(\ell)$.
Therefore, \ref{SHORT.ex:DNA:int>=length-tc} and \ref{SHORT.ex:DNA:end} imply \ref{thm:DNA}.

\parbf{\ref{ex:tangent-support}.}
Show that no line distinct from $\ell$ can support $F$ at $p$. 
Apply \ref{lem:separation} to show that there is a line supporting $F$ at $p$.
Make a conclusion.

{

\begin{wrapfigure}{r}{22 mm}
\vskip-6mm
\centering
\includegraphics{mppics/pic-255}
\vskip0mm
\end{wrapfigure}

\parbf{\ref{ex:anti-bow}.}
Start with the curve $\gamma_1$ shown on the picture.
To obtain $\gamma_2$, slightly unbend (that is, decrease the curvature of) the dashed arc.

}

\parbf{\ref{ex:bow'}}; \ref{SHORT.ex:bow'+}.
We can assume that $p=\gamma_1(a)\z=\gamma_2(a)$ is the origin, 
$\vec u=\tan_1(a)\z=\tan_2(a)$ points in the direction of $x$-axis,
and the points $q_1=\gamma_1(b)$ and $q_2=\gamma_2(b)$ lie in the upper half-plane of the $(x,y)$-plane.

Suppose $\alpha_1<\alpha_2$.
Find a point $x$ such that 
$\dist{x}{q_1}{}<\dist{x}{q_2}{}$, and
$\gamma_1$ with the line segments $[p,x]$ and $[x,q_1]$ bound a convex region in the $(x,y)$-plane.
(Here you have to use that $\beta_1\le\tfrac\pi2$.)

Modify the proof of the bow lemma (\ref{lem:bow}) so that it works for concatenations of $[x,p]$ with $\gamma_1$ and $\gamma_2$, and arrive at a contradiction.

\begin{Figure}
\vskip-0mm
\centering
\includegraphics{mppics/pic-257}
\vskip0mm
\end{Figure}

\parit{\ref{SHORT.ex:bow'-}.} Look at the picture and think.

\parit{Remark.} It is instructive to compare the exercise and the chord lemma (\ref{lem:chord}).

\parbf{\ref{ex:length-dist}}; \ref{SHORT.ex:length-dist:>}.
Choose a value $s_0\in[a,b]$ that splits the total curvature of $\gamma$ in half.
Observe that $\measuredangle(\gamma'(s_0),\gamma'(s))\le \theta$ for any~$s$.
Use this inequality as in the bow lemma.

\begin{wrapfigure}{r}{23 mm}
\vskip-0mm
\centering
\includegraphics{mppics/pic-290}
\vskip-0mm
\end{wrapfigure}

\parit{\ref{SHORT.ex:length-dist:self-intersection:>pi}.} Apply \ref{SHORT.ex:length-dist:>}, and see the picture.

\parit{\ref{SHORT.ex:length-dist:=}.} Start with a polygonal line made from two equal segments with the external angle $2\cdot\theta$ and smooth its joint, using cutoffs and mollifiers (\ref{sec:analysis}).

\parbf{\ref{ex:schwartz}.}
Let $\ell=\length\gamma$.
Suppose $\ell_1<\ell<\ell_2$.
Let $\gamma_1$ be an arc of the unit circle with length~$\ell$.

Show that the distance between the endpoints of $\gamma_1$ is larger than $|p-q|$, and apply the bow lemma (\ref{lem:bow}).

\parit{Source:} The exercise is generally attributed to Hermann Schwarz \cite{shur}.

\parbf{\ref{ex:loop}.}
Suppose $\length\gamma<2\cdot\pi$. 
Apply the bow lemma (\ref{lem:bow}) to $\gamma$ and an arc of the unit circle of the same length.

\parbf{\ref{ex:bow-upper}.}
Choose a smooth spherical curve $\alpha$ that runs near one point;
for example, a small spherical circle centered at this point.
Consider a curve $\tan$ that runs along $\alpha$ with speed $\kappa(s)$ at any $s\in [0,\ell]$.
Show that a curve with tangent indicatrix $\tan$ solves the problem.

\parbf{\ref{ex:gromov-twist}.}
We can assume that $\gamma$ is unit-speed.
Suppose the inequality fails at $t=0$.
We may assume that $\alpha(0)\le\tfrac\pi2$;
otherwise, revert the parametrization.

In the plane spanned by $\gamma(0)$ and~$\gamma'(0)$, choose a unit-speed circle arc (or line segment) $\sigma$ from $0$ to $\gamma(0)$ that comes to $\gamma(0)$ in the direction opposite to $\gamma'(0)$.
Consider a unit-speed semicircle $\tilde\gamma$ with curvature $2$
that starts at $\gamma(0)$ in the direction $\gamma'(0)$ so that the concatenation $\sigma*\tilde\gamma$ is an arc of a convex plane curve; see the figure.

\begin{Figure}
\vskip-1mm
\centering
\includegraphics{mppics/pic-282}
\vskip-1mm
\end{Figure}

Show that if $|\gamma(0)|> \sin(\alpha(0))$, then $\tilde\gamma$ leaves the unit ball; that is, $|\tilde\gamma(t_0)|>1$ for some~$t_0$.

The concatenations $\sigma*\gamma$ and $\sigma*\tilde\gamma$ fail to be smooth at the joint, but it is differentiable at this point.
Check that the proof of the bow lemma works in this more general case.

Apply this bow lemma to $\sigma*\gamma$ and $\sigma*\tilde\gamma$, to get $|\gamma(t_0)|\ge |\tilde\gamma(t_0)|$, and arrive at a contradiction.

\parit{Source:} This statement was used by the first author \cite{petrunin2023}.

%\end{multicols}
%\par\noindent\rule{\textwidth}{0.4pt}
%\begin{multicols}{2}

\parbf{\ref{ex:chord-lemma-optimal}.} 
Set $\alpha=\measuredangle(\vec u,\vec w)$ 
and $\beta=\measuredangle(\vec v,\vec w)$.
Try to guess the example from the picture.

The shown curve is divided into three arcs: I, II, and III. 
Arc I turns from $\vec u$ to $\vec w$;
it has total curvature $\alpha$.
Similarly, arc III turns from $\vec w$ to $\vec v$ and has total curvature $\beta$. 
Arc II goes very close and almost parallel to the chord $pq$, and its total curvature can be made arbitrarily small.

\begin{Figure}
\vskip-1mm
\centering
\includegraphics{mppics/pic-285}
\vskip-1mm
\end{Figure}

\parbf{\ref{ex:monotonic-tc}.}
Use that an exterior angle of a triangle equals the sum of the two remote interior angles.
For the second part apply induction on the number of vertices.

\parbf{\ref{ex:sef-intersection}.} Assume $x$ is a point of self-intersection.
Show that we may choose two points $y$ and $z$ between the self-intersections on $\gamma$ so that the triangle $xyz$ is nondegenerate.
In particular, 
$\measuredangle\hinge yxz
\z+
\measuredangle\hinge zyx
\z<\pi$, or, equivalently, for the \textit{nonclosed} inscribed polygonal line $xyzx$ we have $\tc{xyzx}\z>\pi$.
It remains to apply \ref{prop:inscribed-total-curvature}.

\parbf{\ref{ex:quadrisecant}.}
Consider the closed polygonal line $acbd$.
Observe that $\tc{acbd}=4\cdot\pi$.
It remains to apply \ref{prop:inscribed-total-curvature}.

\parit{Remarks.}
Lines crossing a curve as in the exercise are called \index{alternating quadrisecants}\emph{alternating quadrisecants}.
It is known that any {}\emph{nontrivial knot} admits an alternating quadrisecant \cite{denne};
according to the exercise, the latter implies the so-called {}\emph{F\'ary--Milnor theorem} --- \textit{the total curvature of any knot exceeds~$4\cdot \pi$}; see \cite{petrunin-stadler} and the references therein.

\parbf{\ref{ex:total-curvature=}.}
By \ref{prop:inscribed-total-curvature}, $\tc\gamma\ge \tc\beta$;
it remains to show that
$\tc\gamma\le\sup\{\tc\beta\}$.
In other words, 
given $\epsilon>0$ and a polygonal line $\sigma=\vec u_0\dots \vec u_k$ inscribed in the tangent indicatrix $\tan$ of $\gamma$, 
we need to construct a polygonal line $\beta$ inscribed in $\gamma$ such that
\[\length\sigma<\tc\beta+\epsilon.
\eqlbl{eq:tc=<tc}\]

Suppose $\vec u_i=\tan(s_i)$.
Choose an inscribed polygonal line $\beta=p_0\dots p_{2\cdot k+1}$ such that $p_{2\cdot i}$ and $p_{2\cdot i+1}$ lie sufficiently close to $\gamma(s_i)$; so we can assume that the direction of $p_{2\cdot i+1}-p_{2\cdot i}$ is sufficiently close to $\vec u_i$ for each~$i$.
Show that \ref{eq:tc=<tc} holds true for the constructed polygonal line~$\beta$.

\parbf{\ref{ex:tc-rectifiable}.}
Show that for any polygonal line $\beta$ in a ball of radius $R$, we have
\[\tc{\beta}+2\cdot \pi\ge\frac{\length\beta}R.\]
Observe that $\gamma$ lies in a ball and apply this inequality.

An example for the second question can be found among logarithmic spirals.

%\end{multicols}
%\par\noindent\rule{\textwidth}{0.4pt}
%\begin{multicols}{2}

\stepcounter{chapter}
\setcounter{eqtn}{0}

\parbf{\ref{ex:helix-torsion}.} 
The arc-length parameter $s$ is already found in   \ref{ex:arc-length-helix}.
It remains to calculate the Frenet frame, curvature, and torsion.

{

\begin{wrapfigure}{r}{25 mm}
\vskip-3mm
\centering
\begin{lpic}[t(-0mm),b(0mm),r(0mm),l(0mm)]{asy/helix(1)}
\lbl[br]{8,24;$\norm$}
\lbl[b]{2,26;$\bi$}
\lbl[wl]{15,25;$\tan$}
\end{lpic}
\vskip-0mm
\end{wrapfigure}

One should get 
\begin{align*}
\tan(t)&=\tfrac{(-a\cdot\sin t, a\cdot\cos t,b)}{\sqrt{a^2+b^2}},
\\
\norm(t)&=(-\cos t,-\sin t,0),
\\
\bi(t)&=\tfrac{(b\cdot \sin t,-b\cdot \cos t, a)}{\sqrt{a^2+b^2}},
\\
\kur &\equiv \tfrac{a}{a^2+b^2},
\\
\tor &\equiv \tfrac{b}{a^2+b^2}.
\end{align*}

It remains to show that the map $(a,b) \z\mapsto (\frac{a}{a^2+b^2}, \frac{b}{a^2+b^2})$ sends bijectively the half plane $a>0$ onto itself.

}

\parbf{\ref{ex:beta-from-tau+nu}.} By the product rule, we get
\begin{align*}
\bi'&=(\tan\times \norm)'=
\tan'\times \norm+\tan\times\norm'.
\end{align*}
It remains to substitute the values from \ref{eq:frenet-tau} and \ref{eq:frenet-nu} and simplify.

\parbf{\ref{ex:torsion=0}.}
This is a consequence of the equation $\bi' = - \tor\cdot \norm $.

\parbf{\ref{ex:+B}.}
We can assume that $\gamma_0$ is unit-speed.
Show that
$|\gamma_1'|=|\tan-\tor\cdot \norm|\ge 1$ and use it.

\parbf{\ref{ex:frenet}.}
Observe that $\tfrac{\gamma'\times\gamma''}{|\gamma'\times\gamma''|}$ is a unit vector perpendicular to the plane spanned by $\gamma'$ and $\gamma''$, so, up to sign, it has to be equal to $\bi$.
It remains to check that the sign is right.

\parbf{\ref{ex:moment-curve}.} Apply \ref{ex:curvature-formulas:b} and \ref{ex:frenet}.

\parbf{\ref{ex:bow-converse}.}
We can assume that $t_0=0$ and $\gamma_i(0)\z=0$.
Consider the functions
\begin{align*}
\rho_i(t)&\df|\gamma_i(t)|^2=\langle \gamma_i(t),\gamma_i(t)\rangle.
\end{align*}
Observe that $\rho_1\ge \rho_2$ and $\rho_i(0)=0$.
Show that 
\begin{align*}
\rho_i'(0)&=0,
&
\rho_i''(0)&=2,
\\
\rho_i'''(0)&=0,
&
\rho_i''''(0)&=-2\cdot\kur(0)^2_{\gamma_i}.
\end{align*}
Make a conclusion. (Compare to \ref{ex:const-dist}.)

\parbf{\ref{ex:torsion-indicatrix}.} Assume that the tangent indicatrix has no self-intersection.
Show that it lies in an open hemisphere and argue as in Fenchel's theorem.

\parbf{\ref{ex:lancret}}; \ref{SHORT.ex:lancret:a}.
Observe that 
$\langle \vec w,\tan\rangle'=0$.
Show that it implies that $\langle \vec w, \norm\rangle =0$.
Further, show and use that $\langle \vec w, \tan\rangle^2+\langle \vec w, \bi\rangle^2=\langle \vec w, \vec w\rangle$.
To prove the last identity, apply the Frenet formula for~$\norm'$.

\parit{\ref{SHORT.ex:lancret:b}.}
Show that $\vec w'=0$;
it implies that $\langle \vec w,\tan\rangle\z=\tfrac\tor\kur$.
In particular, the velocity vector makes a constant angle with $\vec w$; that is, $\gamma$ has a constant slope.

\parbf{\ref{ex:evolvent-constant-slope}.}
Suppose $\alpha$ is an evolvent of $\gamma$, and $\vec w$ is a fixed vector.
Show and use that $\langle \vec w,\alpha\rangle$ is constant if $\gamma$ makes a constant angle with $\vec w$.

\parbf{\ref{ex:spherical-frenet}}.
Part \ref{SHORT.ex:spherical-frenet:tau} follows since $(\tan,\norm,\bi)$ is an orthonormal basis.
For \ref{SHORT.ex:spherical-frenet:nu}, take the first and second derivatives of the identity $\langle\gamma,\gamma\rangle=1$ and simplify them using the Frenet formulas.
Part \ref{SHORT.ex:spherical-frenet:beta} follows from \ref{SHORT.ex:spherical-frenet:nu} and the Frenet formulas.
By \ref{SHORT.ex:spherical-frenet:beta}, $\int\tfrac\tor\kur=0$, hence \ref{SHORT.ex:spherical-frenet:beta+} follows.
Part \ref{SHORT.ex:spherical-frenet:kur-tor} is proved by algebraic manipulations.
For \ref{SHORT.ex:spherical-frenet:f},
use the Frenet formulas to show that $(\gamma+\tfrac1\kur\cdot \norm+\tfrac{\kur'}{\kur^2\cdot\tor}\cdot\bi)'=0$.

\parit{Remark.}
Part \ref{SHORT.ex:spherical-frenet:beta} implies that \textit{$\int\tfrac\tor\kur=0$ for any closed smooth spherical curve}.
It is known that this property characterizes spheres and planes; this was proved by Beniamino Segre \cite{segre}.

\parbf{\ref{ex:cur+tor=helix}.} Use the second statement in \ref{ex:helix-torsion}.

\parbf{\ref{ex:const-dist}.} The function
\begin{align*}
\rho(\ell)&=|\gamma(t+\ell)-\gamma(t)|^2=
\\
&=\langle \gamma(t+\ell)-\gamma(t),\gamma(t+\ell)-\gamma(t)\rangle
\end{align*}
is smooth and does not depend on~$t$.
Express speed, curvature, and torsion of $\gamma$ in terms of the derivatives $\rho^{(n)}(0)$.
Be patient, you will need two derivatives for the speed,
four for the curvature,
and six for the torsion.
Once it is done, apply \ref{ex:cur+tor=helix}.

%\end{multicols}
%\par\noindent\rule{\textwidth}{0.4pt}
%\begin{multicols}{2}

\stepcounter{chapter}
\setcounter{eqtn}{0}

\parbf{\ref{ex:bike}.}
Without loss of generality, we may assume that $\gamma_0$ is parametrized by its arc-length.
Then
\begin{align*}
|\gamma_1'|&=|\gamma_0'+\tan'|=|\tan+\kur\cdot\norm|.
\end{align*}
It follows that
\[|\gamma_1'(t)|\ge|\gamma_0'(t)|
\quad\text{and}\quad
|\gamma_1'(t)|\ge|\kur(t)_{\gamma_0}|
\]
for any $t\in[a,b]$.
Integrate these inequalities and apply 
\ref{ex:integral-length}.

\parbf{\ref{ex:trochoids}.}
Observe that 
\[\gamma'_a(t)=(1+a\cdot \cos t, -a\cdot \sin t);\]
that is, $\gamma'_a$ runs clockwise along a circle with center at $(1,0)$ and radius $\vert a \vert$.

\parit{Case $|a|>1$.}
Notice that $\tan_a(t)=\gamma'_a/|\gamma'_a|$ runs clockwise and makes a full turn in time $2\cdot\pi$.
Therefore, $\tgc{\gamma_a}=-2\cdot\pi$ and $\tc{\gamma_a}\z=|\tgc{\gamma_a}|=2\cdot\pi$.

\parit{Case $|a|<1$.}
Set $\theta_a=\arcsin a$.
Show that $\tan_a(t)=\gamma'_a/|\gamma'_a|$ starts with the horizontal direction $\tan_a(0)=(1,0)$, turns monotonically to angle $\theta_a$, then monotonically to $-\theta_a$ and then monotonically back to $\tan_a(2\cdot\pi)=(1,0)$.
It follows that 
$\tgc{\gamma_a}=0$ and $\tc{\gamma_a}=4\cdot\theta_a$.

\parit{Case $a=-1$.}
The velocity $\gamma'_{-1}(t)$ vanishes at $t=0$ and $2\cdot\pi$.
Nevertheless, the curve admits a smooth regular parametrization --- find it.
You should get $\tc{\gamma_{-1}}=-\tgc{\gamma_{-1}}=\pi$.

\parit{Case $a=1$.}
The velocity $\gamma'_1(t)$ vanishes at $t=\pi$.
At $t=\pi$ the curve has a cusp;
that is, $\gamma_1$ turns exactly back at the time $\pi$.
So $\gamma_1(t)$ has undefined total signed curvature.
The curve is a joint of two smooth arcs with external angle $\pi$, and
the total curvature of each arc is $\tfrac\pi2$, so 
$\tc{\gamma_{1}}=\tfrac\pi2+\pi+\tfrac\pi2=2\cdot\pi$.

\parbf{\ref{ex:zero-tsc}.}
The answers are shown.
In the last picture, we assume that the two marked points have parallel tangent lines.

\begin{Figure}
\begin{minipage}{.27\textwidth}
\centering
\includegraphics{mppics/pic-260}
\end{minipage}\hfill
\begin{minipage}{.42\textwidth}
\centering
\includegraphics{mppics/pic-261}
\end{minipage}
\hfill
\begin{minipage}{.27\textwidth}
\centering
\includegraphics{mppics/pic-262}
\end{minipage}
\end{Figure}

\parbf{\ref{ex:length'}}; \ref{SHORT.ex:length':reg}.
Show that
\[
\gamma_\ell'(t)=(1-\ell\cdot\skur(t))\cdot \gamma'(t).
\]
By regularity of $\gamma$, $\gamma'\ne0$.
So if $\gamma_\ell'(t)=0$, then $\ell\cdot \skur(t)=1$.

\parit{\ref{SHORT.ex:length':formula}.} We can assume that $\gamma$ is parametrized by its arc-length, so $\gamma'(t)=\tan(t)$.
Suppose $|\ell|<\frac{1}{\kur(t)}=\frac{1}{|\skur(t)|}$ for any~$t$.
Then 
\[
|\gamma_\ell'(t)|=(1-\ell\cdot\skur(t)).
\]
Therefore,
\begin{align*}
L(\ell)
&=
\int_a^b(1-\ell\cdot\skur(t))\cdot dt=
\\&=
\int_a^b1\cdot dt-\ell\cdot \int_a^b\skur(t)\cdot dt=
\\
&=
L(0)-\ell\cdot\tgc\gamma.
\end{align*}

\parit{\ref{SHORT.ex:length':antiformula}.}
Consider the unit circle $\gamma(t)=(\cos t,\sin t)$ for $t\in[0,2\cdot\pi]$ and $\gamma_\ell$ for $\ell=2$.

\parbf{\ref{ex:inverse}.}
Use the definition of osculating circle via order of contact, and that inversions send circles to circles or lines. 

\parbf{\ref{ex:evolute}.}
We assume that $\gamma$ is unit-speed.
Show that $\omega'=\tfrac{\skur'}{\skur^2}\cdot \norm$. 
Conclude that $\omega$ has Frenet frame either $(\norm,-\tan)$ or $(-\norm,\tan)$, and its curvature is $|\tfrac{\kur^3}{\kur'}|$.

\begin{wrapfigure}{r}{19 mm}
\vskip-4mm
\centering
\includegraphics{mppics/pic-296}
\vskip-0mm
\end{wrapfigure}

\parbf{\ref{ex:3D-spiral}.}
Start with a plane spiral curve as shown on the picture.
Applying \ref{thm:fund-curves}, increase the torsion of the dashed arc without changing the curvature until a self-intersection appears.

You may think that the dashed arc is very short, so the tangent vector $\tan$ is nearly constant on this arc. 
Increasing the torsion can rotate normal vector $\norm$ arbitrary around $\tan$.
An intersection appears if $\norm$ is rotated by certain angle near $\pi$.
(Compare to \ref{ex:approximation-const-curvature}.)

\parbf{\ref{ex:double-tangent}.} 
Observe that if a line or circle is tangent to $\gamma$,
then it is tangent to the osculating circle at the same point.
Then apply the spiral lemma (\ref{lem:spiral}).

\parbf{\ref{ex:spherical-spiral}.}
Note that osculating circles of a spherical curve lie in the sphere.
Prove an analog of \ref{lem:spiral} for these circles.
(Compare to \ref{ex:spherical-frenet:beta+}.)

%\end{multicols}
%\par\noindent\rule{\textwidth}{0.4pt}
%\begin{multicols}{2}

\stepcounter{chapter}
\setcounter{eqtn}{0}

\parbf{\ref{ex:vertex-support}.}
Apply the spiral lemma (\ref{lem:spiral}).

\parit{Computational solution.} 
We will assume that the curvature does not vanish at $p$, the remaining case is simpler.
We may assume that $\gamma$ is unit-speed, $p=\gamma(0)$,
and the center of curvature of $\gamma$ at $p$ is the origin;
in other words, $\kur(0)\cdot\gamma(0)+ \norm(0)=0$.

Consider the function $f\:t\mapsto \langle\gamma(t),\gamma(t)\rangle$.
Direct calculations show that 
\begin{align*}
f'
&=\langle\gamma,\gamma\rangle'
=2\cdot\langle\tan,\gamma\rangle,
\\
f''
&=2\cdot\langle\tan,\gamma\rangle'
=2\cdot\kur\cdot\langle\norm,\gamma\rangle+2,
\\
f'''
&=2\cdot\kur'\cdot\langle\norm,\gamma\rangle
+2\cdot\kur\cdot\langle\norm',\gamma\rangle+\cancel{2\cdot\kur\cdot\langle\norm,\tan\rangle}.
\end{align*}

Observe that $\norm'(0)\perp\gamma(0)$.
Therefore, $f'(0)=0$, $f''(0)=0$, and $f'''(0)\z=-2\cdot\kur'/\kur$.

If the osculating circle supports $\gamma$ at $p$,
then the function $f$ has a local maximum or minimum at $0$.
Therefore $f'''(0)=0$ and $\kur'=0$.

\parbf{\ref{ex:support}.}
Consider the coordinate system with the origin at $p$ and the common tangent line to $\gamma_1$ and $\gamma_2$ as the $x$-axis.
We may assume that $\gamma_1$ and $\gamma_2$ are defined in $(-\epsilon,\epsilon)$ for small $\epsilon>0$,
and they run almost horizontally.

Given $t\in[0,1]$ consider the curve $\gamma_t$ that is tangent and cooriented to the $x$-axis at $\gamma_t(0)\z=p$ and has signed curvature defined by $\skur_t(s)\z\df(1\z-t)\cdot\skur_0(s)+t\cdot\skur_1(s)$.
It exists by \ref{thm:fund-curves-2D}.

Choose $s\approx 0$.
Consider the curve $\alpha_s\:t\z\mapsto \gamma_t(s)$.
Show that $\alpha_s$ moves almost vertically up.
Use that $\gamma_t$ moves almost horizontally to show that in a small neighborhood of $p$, the curve $\gamma_1$ lies above $\gamma_0$,
whence the statement follows.

\parbf{\ref{ex:in-circle}.}
Move the unit circle straight until its center gets to the base of the loop,
and then reduce its radius to zero.
Show that at the moment when the circle first touches the loop, it supports the loop at a point distinct from the base.
Apply~\ref{prop:supporting-circline}.

\parbf{\ref{ex:between-parallels-1} $\bm{+}$ \ref{ex:in-triangle} $\bm{+}$ \ref{ex:lens}.}
Observe that one of the arcs of curvature 1 in the families shown on the picture supports $\gamma$, and apply \ref{prop:supporting-circline}.

\begin{Figure}
\begin{minipage}{.35\textwidth}
\centering
\includegraphics{mppics/pic-265}
\end{minipage}
\hfill
\begin{minipage}{.3\textwidth}
\centering
\includegraphics{mppics/pic-266}
\end{minipage}
\hfill
\begin{minipage}{.25\textwidth}
\centering
\includegraphics{mppics/pic-267}
\end{minipage}
\end{Figure}

The second part in \ref{ex:between-parallels-1} can be reduced to \ref{ex:in-circle} using the shown family and another family of arcs curved in the opposite direction.

\parit{Remark.} Compare to \ref{ex:moon-rad}.

\begin{wrapfigure}{r}{34 mm}
\vskip-6mm
\centering
\includegraphics{mppics/pic-268}
\vskip0mm
\end{wrapfigure}

\parbf{\ref{ex:convex small}.}
Note that we can assume that $\gamma$ bounds a convex figure~$F$.
Otherwise, by \ref{prop:convex} its curvature changes sign.
Therefore, $\gamma$ has zero curvature at some point.

Choose two points $x$ and $y$ surrounded by $\gamma$ such that $|x\z-y|\z>2$.
Look at the maximal lens bounded by two arcs with a common chord $xy$ that lies in~$F$.
Apply the supporting test (\ref{prop:supporting-circline}).

\parbf{\ref{ex:convex-lens}.}
Apply the lens lemma (\ref{lem:lens}) to show that $\gamma$ lies on one side from the line $pq$.
Conclude that the union of the arc $\gamma$ and its chord $[p,q]$ is a simple closed curve;
by Jordan's theorem, it bounds a figure, say~$F$.

Assume that $F$ is not convex and arrive at a contradiction as in \ref{prop:convex}. 

\parbf{\ref{ex:diameter-of-simple-curve}.}
Choose two points $p,q\in\gamma$ on the maximal distance.
Apply \ref{ex:convex-lens} to the arc of $\gamma$ from $p$ to $q$.
After that use the bow lemma (\ref{lem:bow}), or argue as in \ref{ex:convex small}.

For the last part of the problem, consider the smallest circle $\sigma$ that surrounds $\gamma$ and argue as before for one of the arcs of $\gamma$ between common points with $\sigma$.

\parit{Remark.}
If self-intersections are allowed, then diameter of the curve can be arbitrarily large.
An example can be guessed from the following picture.
In fact the diameter can be made arbitrarily close to the length of the curve.
\begin{Figure}
\vskip-0mm
\centering
\includegraphics{mppics/pic-269}
\vskip-0mm
\end{Figure}

\parbf{\ref{ex:moon-rad}.}
Show that $\gamma$ contains a simple loop $\gamma_1$.
Apply \ref{thm:moon-orginal} to $\gamma_1$.

\parbf{\ref{ex:2-squares}.}
We may assume that the smaller square is on the left of $\gamma$.
Use \ref{thm:4-vert-supporting} to show that there are two points on $\gamma$ with signed curvature $\ge 1$ and $\le 1$, respectively.
Apply the intermediate value theorem for signed curvature to show that there is a point with curvature $1$.

\parit{More challenging question:} What is the minimal number of such points?

\parbf{\ref{ex:moon-area}.} 
Applying rescaling, we may assume that $a=\pi$; that is, our loop bounds area of the unit disc.
Then apply \ref{thm:moon} to show that the curvature is at least 1 at some point.
Argue as at the end of the proof of \ref{thm:4-vert-supporting}, showing that the curvature is at most 1 at some point.

Apply the intermediate value theorem to show that curvature is $1$ at some point.

\parbf{\ref{ex:curve-crosses-circle}.}
Repeat the proof of \ref{thm:moon} for each cyclic concatenation of an arc of $\gamma$ from $p_i$ to $p_{i+3}$ and the arc of the circle with points $p_{i+4},\z\dots,p_{i-1}$ on it.

\begin{wrapfigure}{r}{20 mm}
\vskip-2mm
\centering
\includegraphics{mppics/pic-305}
\vskip-2mm
\end{wrapfigure}

An example for the second part can be guessed from the picture.

%\parbf{\ref{ex:order}.}
% given two points $p,q\in\gamma$ write (p,q)\in in or (p,q)\in out
% We can assume that if (p,q)\in in, then (p,q)\notin out and other way around.
% Otherwise $\gamma$ is a circle and the statement is trivial.

% Let $A, B\subset\gamma$ be the set of all points at which osculating circle supports $\gamma$ from inside and outside respectively.
%We can assume that $A\cap B=\emptyset$, otherwise $\gamma$ is a circle and the statement is trivial.

%Notice that $A$ and $B$ are closed sets.
%If the statement does not hold $B$ lies in one open arc, say $\gamma_1$, of $\gamma$ in the complement of $A$.
%Show that there is 

\parbf{\ref{ex:berk}.}
The if part is trivial; let us show the only-if part.
Choose two points $p,q\in\gamma$ such that $\dist{p}{q}{}\ge 4$.
Denote by $\sigma_p$, $\sigma_q$, $r_p$, $r_q$ the incircles and inradii at $p$ and $q$ respectively.

Suppose $r_p\ge 1$ and $r_q\ge1$.
Show that two unit discs tangent to $\gamma$ from inside at $p$ and $q$ solve the problem.

Suppose $r_p<1$ and $r_q<1$.
Observe that $\sigma_p$ and $\sigma_q$ do not intersect.
Note that we can choose two nonoverlapping arcs $\gamma_p$ and $\gamma_q$ of $\gamma$ with ends on  $\sigma_p$ and $\sigma_q$ respectively.
Start with these two arcs and argue as in \ref{thm:moon}.
Show that the obtained supporting osculating circles do not intersect;
make a conclusion.

In the remaining case $r_p\ge 1$ and $r_q<1$, combine both arguments.

\parit{Source:} The problem is due to Berk Ceylan \cite{ceylan}.
It is unknown whether \textit{any simple closed plane curve with curvature at most 1 and length at least $4\cdot\pi$ surrounds two disjoint unit discs.}

\parbf{\ref{ex:4x0-torsion}.}
Observe that the tangent indicatrix $\tan$ of $\gamma$ is a smooth simple closed spherical curve.
Show that each point of $\gamma$ with vanishing torsion corresponds to an {}\emph{inflection point} of $\tan$;
that is, a point with vanishing geodesic curvature (see Section~\ref{sec:Darboux}).

Try to prove that if $\tan$ has less than $4$ inflection points, then it lies in an open hemisphere,
or extract this statement from the proof of the tennis ball theorem \cite[§ 20]{arnold1994}.
After that argue as in \ref{thm:fenchel}.

%\end{multicols}
%\par\noindent\rule{\textwidth}{0.4pt}
%\begin{multicols}{2}

\stepcounter{chapter}
\setcounter{eqtn}{0}

\parbf{\ref{ex:hyperboloids}.} 
Denote the level set by $\Sigma_\ell$.

Show that $\nabla_p f=0$ if and only if $p\z=(0,0,0)$.
Use \ref{prop:implicit-surface} to conclude that if $\ell\ne 0$, then $\Sigma_\ell$ is a union of disjoint smooth surfaces.

Show that $\Sigma_\ell$ is connected if and only if $\ell\ge 0$.
It follows that if $\ell>0$, then $\Sigma_\ell$ is a smooth surface, and if $\ell<0$ --- it is not.

The case $\ell=0$ has to be done by hand --- it does not satisfy the sufficient condition in \ref{prop:implicit-surface}, but that does not by itself imply that $\Sigma_0$ is not a smooth surface.

Show that any neighborhood of the origin in $\Sigma_0$ cannot be described by a graph in any coordinate system;
so by the definition (Section~\ref{sec:def-smooth-surface}) $\Sigma_0$ is not a smooth surface.

\parbf{\ref{ex:9-surf}.} Modify the solution of \ref{ex:9:9}.

\parbf{\ref{ex:smooth-fun(surf)}.}
To prove the if part, use that the composition of smooth functions is smooth.

Now, choose a graph representation $z\z=f(x,y)$ of a neighborhood of $p$ in $\Sigma$.
Note that $g$ is smooth in the neighborhood if and only if the function $\hat g\:(x,y)\mapsto g(x,y,f(x,y))$ is smooth in its domain of definition.
Define $h(x,y,z)\z\df \hat g(x,y)$; it solves the only-if part.

For the last part, consider function $g\:(x,y)\mapsto\tfrac1{x^2+y^2}$ defined on the $(x,y)$-plane with removed origin.
Observe that $g$ cannot be extended to a continuous function on $\mathbb{R}^3$.

\parit{Remark.} Suppose $\Sigma$ is a smooth proper surface.
Then \textit{any smooth function $g\:\Sigma\to\mathbb{R}$ can be extended to a smooth function on $\mathbb{R}^3$}.
A proof uses the so-called \emph{partition of unity};
read about it and try to prove this statement.

\parbf{\ref{ex:inversion-chart}.} 
Check that the image of $s$ lies in the unit sphere centered at $(0,0,1)$;
that is, show that 
$\left(\tfrac{2\cdot u}{1+u^2+v^2}\right)^2
\z+
\left(\tfrac{2\cdot v}{1+u^2+v^2}\right)^2
\z+
\left(\tfrac{2}{1+u^2+v^2}-1\right)^2=1.
$
for any $u$ and~$v$.

Show that the map 
$(x,y,z)\z\mapsto \tfrac{2\cdot(x,y)}{x^2+y^2+z^2}$
describes the inverse of $s$, which is continuous away from the origin.
In particular, $s$ is an embedding that covers the entire sphere except for the origin.

Show that $s$ is regular; that is, $s_u$ and $s_v$ are linearly independent at all points of the $(u,v)$-plane.

\parit{Remark.}
The map $s$ in the exercise can be visualized as
\[(u,v,1)\mapsto (\tfrac{2\cdot u}{1+u^2+v^2},\tfrac{2\cdot v}{1+u^2+v^2},\tfrac{2}{1+u^2+v^2})\]
which is called \index{stereographic projection}\emph{stereographic projection} from the plane $z=1$ to the unit sphere with center at $(0,0,1)$.
Note that the point $(u,v,1)$ and its image lie on the same half-line emerging from the origin. 

\begin{Figure}
\vskip-0mm
\centering
\includegraphics{mppics/pic-750}
\vskip0mm
\end{Figure}

\parbf{\ref{ex:revolution}.}
Set
$s\:(t,\theta)\mapsto (x(t), y(t)\cdot\cos\theta,y(t)\cdot\sin\theta)$.
Show that $s$ is regular; that is, $s_t$ and $s_\theta$ are linearly independent.
(It might help to observe that $s_t\perp s_\theta$).

Show that $s$ is an embedding;
that is, any $(t_0,\theta_0)$ admits a neighborhood $U$ in the $(t,\theta)$-plane such that the restriction $s|_U$ has a continuous inverse.
It remains to apply \ref{cor:reg-parmeterization}.

\parbf{\ref{ex:inv-diffeomorphism}.} Apply the inverse function theorem (\ref{thm:inverse}) in charts for both surfaces. 

\parbf{\ref{ex:star-shaped-disc}.} 
The solutions of these exercises are built on the following general construction known as the \index{Moser trick}\emph{Moser trick}.

Suppose $\vec u_t$ is a smooth time-dependent vector field on a plane.
Consider the ordinary differential equation $x'(t)=\vec u_t(x(t))$.
Consider the map $\iota\:x(0)\mapsto x(1)$ where $t\mapsto x(t)$ is a solution of the equation.
The map $\iota$ is called the \index{flow}\emph{flow} of the vector field $\vec u_t$ for the time interval $[0,1]$.
By \ref{thm:ODE}, the map $\iota$ is smooth in its domain of definition.
Moreover, the same holds for its inverse $\iota^{-1}$;
indeed, $\iota^{-1}$ is the flow of the vector field $-\vec u_{1-t}$.
That is, $\iota$ is a diffeomorphism from its domain of definition to its image.

Therefore, to construct a diffeomorphism from one open subset of the plane to another, it is sufficient to construct a smooth vector field such that its flow maps one set to the other;
such a map is automatically a diffeomorphism.

\parit{\ref{SHORT.ex:plane-n}.}
Consider two sets $\mathbb{R}^2\setminus\{p_1,\dots,p_n\}$ and $\mathbb{R}^2\setminus\{q_1,\dots,q_n\}$.
Choose smooth paths $\gamma_i\:[0,1]\z\to \mathbb{R}^2$ such that $\gamma_i(0)=p_i$,
$\gamma_i(1)=q_i$, and $\gamma_i(t)\ne \gamma_j(t)$ if $i\ne j$.

Choose a smooth vector field $\vec v_t$ such that $\vec v_t(\gamma_i(t))=\gamma'_i(t)$ for any $i$ and~$t$.
We can assume in addition that $\vec v_t$ vanishes outside a sufficiently large disc; this can be arranged by multiplying the vector field by an appropriate bell function;
for example, $\sigma_1(R-|x|)$, where $R$ is large and $\sigma_1$ is defined on page \pageref{page:sigma-function}.

It remains to apply the Moser trick to the constructed vector field. 

\parit{\ref{SHORT.ex:star-shaped-disc:smooth}--\ref{SHORT.ex:star-shaped-disc:star-shaped}.}
In \ref{SHORT.ex:star-shaped-disc:smooth}--\ref{SHORT.ex:star-shaped-disc:nonsmooth}, we can assume that the origin belongs to both sets, and in \ref{SHORT.ex:star-shaped-disc:star-shaped} that the sets are star-shaped with respect to the origin. 

In each case show that there is a vector field $\vec v$ defined on $\mathbb{R}^2$ that flows one surface to an other.
In fact, one can choose radial fields of that type,
but be careful with the cases \ref{SHORT.ex:star-shaped-disc:nonsmooth} and \ref{SHORT.ex:star-shaped-disc:star-shaped} --- they are not as easy as one might think.

%\end{multicols}
%\par\noindent\rule{\textwidth}{0.4pt}
%\begin{multicols}{2}

\stepcounter{chapter}
\setcounter{eqtn}{0}

\parbf{\ref{ex:tangent-normal}.}
Let $\gamma$ be a smooth curve in~$\Sigma$.
Observe that $f\circ\gamma(t)\equiv 0$.
Differentiate this identity, and apply the definition of tangent vector (\ref{def:tangent-vector}).

\parbf{\ref{ex:vertical-tangent}.}
Assume a neighborhood of $p$ in $\Sigma$ is a graph $z=f(x,y)$.
In this case, $s\:(u,v)\z\mapsto (u,v,f(u,v))$ is a smooth chart at~$p$.
Show that the plane spanned by $s_u$ and $s_v$ is not vertical;
together with \ref{def:tangent-plane}, this proves the if part.

For the only-if part, fix a chart 
$s\:(u,v)\z\mapsto(x(u,v),y(u,v),z(u,v))$,
and apply the inverse function theorem to the map $(u,v)\z\mapsto(x(u,v),y(u,v))$.

\parbf{\ref{ex:tangent-single-point}.}
Choose $(x,y,z)$-coordinates so that $\Pi$ is the $(x,y)$-plane and $p$ is the origin.
Let $(u,v)\mapsto s(u,v)$ be a chart of $\Sigma$ such that $p=s(0,0)$.
Denote by $\vec k$ the unit vector in the direction of the $z$-axis.

Show that we can assume that $\langle s(u,v),\vec k\rangle>0$ in a punctured neighborhood of $0$.
Conclude that $s_u\perp \vec k$ and $s_v\perp \vec k$ at
$0$, hence the result.

\parbf{\ref{ex:lin-ind-chart}.}
By \ref{thm:ODE}, there is a curve $\alpha$ in $\Sigma$ such that $\alpha(0)=p$ and $\alpha'(x)=\vec x_{\alpha(x)}$.
The same way construct a curve $\beta_x$ in $\Sigma$ such that $\beta_x(0)=\alpha(x)$ and $\beta_x'(y)=\vec y_{\beta_x(y)}$.
Apply \ref{thm:inverse} to show that $(x,y)\mapsto \beta_x(y)$ describes a chart of a small neighborhood $W$ of $p$.

Show that $u\:\beta_x(y)= x$ meets the required conditions in $W$.
The same way one can construct $v$.

It remains to apply \ref{thm:inverse} again.

\parbf{\ref{ex:const-normal}.}
Show and use that for any smooth curve $\gamma$ in $\Sigma$ the function $t\mapsto \langle\nu_0,\gamma(t)\rangle$ is constant.

\parbf{\ref{ex:implicit-orientable}.}
By \ref{ex:tangent-normal}, $\Norm=\tfrac{\nabla h}{|\nabla h|}$ defines a unit normal field on~$\Sigma$.

\parbf{\ref{ex:plane-section}.}
Use cutoffs and mollifiers (Appendix~\ref{sec:analysis}) to construct a smooth nonnegative function $f$ on the $(x,y)$-plane such that $f(x,y)=0$ if and only if $(x,y)\in A$.
The graph $z=f(x,y)$ describes the required surface.

%\end{multicols}
%\par\noindent\rule{\textwidth}{0.4pt}
%\begin{multicols}{2}

\stepcounter{chapter}
\setcounter{eqtn}{0}

\parbf{\ref{ex:line-of-curvature}.}
Fix a point $p$ on~$\gamma$.
Since $\Sigma$ is mirror-symmetric with respect to $\Pi$,
we have $\T_p\perp \Pi$.

Choose $(x,y)$-coordinates on $\T_p$ so that the $x$-axis is the intersection $\Pi\cap \T_p$.
Suppose the osculating paraboloid is described by the graph 
$z=\tfrac12\cdot(\ell\cdot x^2+2\cdot m\cdot x\cdot y+n\cdot y^2)$.
Since $\Sigma$ is mirror-symmetric, so is the paraboloid;
that is, changing $y$ to $(-y)$ does not change the value 
$\ell\cdot x^2+2\cdot m\cdot x\cdot y+n\cdot y^2$.
In other words $m=0$, or equivalently, the $x$-axis points in a direction of curvature.

\parbf{\ref{ex:gauss+orientable}.}
Notice that the principal curvatures have the same sign at each point.
Therefore, we can choose a unit normal $\Norm$ at each point so that both principal curvatures are positive.
Show that it defines a global field on the surface.

\parbf{\ref{ex:re-scale-surface-curvature}.} Compute the Taylor series of the function $g(x,y)= \lambda \cdot f( x/ \lambda , y/\lambda)$.

\parbf{\ref{ex:self-adjoint}.}
Apply \ref{thm:shape-chart} to a map $s$ such that $s_u(0,0)=\vec u$ and $s_v(0,0)=\vec v$.

\parbf{\ref{ex:normal-curvature=const}}; \ref{SHORT.ex:normal-curvature=const:a}.
Observe that $\Sigma$ has unit Hessian matrix at each point, and apply the definition of the shape operator.

\parit{\ref{SHORT.ex:normal-curvature=const:b}.}
Choose a chart $s$ in~$\Sigma$.
Show that
\[\tfrac{\partial }{\partial u}(s+\Norm)
=
\tfrac{\partial }{\partial v}(s+\Norm)
=
0.\]
Make a conclusion.

\parbf{\ref{ex:normal-curvature=0}.} 
The problem might seem trivial until one realizes that a connected set can be quite peculiar.
For example, it may contain no nontrivial curves.

\medskip

Apply Sard's lemma (\ref{lem:sard}) to show that $\Norm$ is constant on $Z_0$, denote its value $\Norm_0$.
Apply Sard's lemma again to show that $x\mapsto \langle \Norm_0,x\rangle$ is constant on $Z_0$.
Make a conclusion.

\parit{Remark.} This result was used by Richard Sacksteder \cite[Lemma 6]{sacksteder}.
Try to solve the following slightly more challenging problem:
\textit{Let $\Sigma$ be a smooth surface with orientation defined by a unit normal field $\Norm$,
and let $Z_1\subset \Sigma$ be a connected set with unit shape operator.
Show that $Z_1$ lies in a unit sphere.}

\parbf{\ref{ex:shape-curvature-line}.} 
We can assume that $\gamma$ is parametrized by arc-length.
Denote by $\Norm_1(s)$ and $\Norm_2(s)$ the unit normal vectors to $\Sigma_1$ and $\Sigma_2$ at $\gamma(s)$.

By the assumption, $\langle \Norm_1,\Norm_2\rangle= c$ for some constant $c$.
We know that
$\Norm_1'$ is proportional to $\gamma'$ and need to show that $\Norm_2'$ is proportional to $\gamma'$.
If $c=\pm1$, then $\Norm_1\equiv\pm \Norm_2$, and the statement trivially holds.

Note that $\langle\Norm_1',\Norm_2\rangle=0$;
use that $c'=0$ to show that $\langle\Norm_1,\Norm_2'\rangle=0$.
Show that $\langle\Norm_2,\Norm_2'\rangle=0$.
Conclude that $\Norm_2'$ is proportional to $\gamma'$.

\parit{Source:} This is a result by Ferdinand Joachimsthal \cite{joachimsthal} generalized by Ossian Bonnet \cite{bonnet}.

\parbf{\ref{ex:equidistant}};
\ref{SHORT.ex:equidistant:smooth}.
Fix~$t$.
Set $f_t\:p\mapsto p+t\cdot\Norm(p)$; it maps $\Sigma$ to $\Sigma_t$.

Apply the definition of the shape operator to show that $d_pf_t(\vec v)=\vec v-t\cdot\Shape_p(\vec v)$.
Since $\Sigma$ is closed, the norm of $\Shape$ is bounded.
Whence $f_t$ is regular if $t\approx 0$.

Further show that $f_t$ is injective if $t\approx0$;
conclude that $f_t$ is a smooth embedding.

\parit{\ref{SHORT.ex:equidistant:area}.} By the area formula (\ref{prop:surface-integral}), we have
\[\area\Sigma_t=\int_{p\in \Sigma}\jac_pf_t.\]
Show and use that $\jac_pf_t=1-t\cdot H+t^2\cdot K$.

\parbf{\ref{ex:flat-plane};} \ref{SHORT.ex:flat-plane:orthonormal}.
Use that principal directions are orthogonal; see \ref{sec:Principal curvatures}.

\parit{\ref{SHORT.ex:flat-plane:depend}.} 
Apply \ref{cor:Shape(ij)} to show that $\Norm_u=0$; use it to show that $\Norm$ depends only on~$v$.

Observe that $\Norm_{uv}=\Norm_{vu}=0$.
Use \ref{cor:Shape(ij)} to show that $\Norm_v=-k \cdot s_v$.
Conclude that $\vec v_u=0$ and therefore $\vec v$ depends only on~$v$.

Apply \ref{SHORT.ex:flat-plane:orthonormal} to show that $\vec u$ depends only on~$v$.
Make a conclusion.

\parit{\ref{SHORT.ex:flat-plane:depend-u}.}
The first part follows from the conclusion of \ref{SHORT.ex:flat-plane:depend}.

By \ref{SHORT.ex:flat-plane:orthonormal}, $\langle s_u,s_v\rangle=0$.
Since $s_{uu}$ is proportional to $s_u$, we also get $\langle s_{v},s_{uu}\rangle=0$.
Therefore,
\begin{align*}
\tfrac{\partial}{\partial v}\langle s_u,s_u\rangle&=2\cdot \langle s_{vu},s_u\rangle=
\\
&=2\cdot \tfrac{\partial}{\partial u}\langle s_v,s_u\rangle-2\cdot \langle s_{v},s_{uu}\rangle=0.
\end{align*}
Make a conclusion.

\parit{\ref{SHORT.ex:flat-plane:linear}.} 
Since $\langle s_u,s_u\rangle=1$,
\[0=\tfrac{\partial}{\partial u}\langle s_u,s_u\rangle=2\cdot\langle s_{uu},s_u\rangle.\]
By \ref{SHORT.ex:flat-plane:depend-u}, $s_{uu}\parallel s_u$; therefore, $s_{uu}=0$.

Conclude that $s_v$ depends linearly on $u$.
Use this together with the equations
\[\Norm_v=-k\cdot s_v,\quad \Norm_{uv}=0.\]

\parit{Remarks.}
This exercise is based on the proof given by Sergei Ivanov \cite[3$^\text{d}$ Sem. Lect. 13]{ivanov}.
Another proof can be built on the Peterson--Codazzi formulas.

\parbf{\ref{ex:mean-curvature}.}
Apply \ref{obs:k1-k2} and the definition of mean curvature.

\parbf{\ref{ex:average}.}
Recall that $H=k_1+k_2$ and $K=k_1\cdot k_2$, where $k_1$ and $k_2$ are the principal curvatures.
By Euler's identity, we need to find the average value of 
\[(k_1\cdot (\cos\theta)^2+k_2\cdot (\sin\theta)^2)^2.\]
Show and use that $\tfrac38$, $\tfrac38$, and $\tfrac18$ are the average values of $(\cos\theta)^4$, $(\sin\theta)^4$, and $(\cos\theta\cdot \sin\theta)^2$, respectively.

\parbf{\ref{ex:meusnier}.}
Use Meusnier's theorem (\ref{thm:meusnier}), to find the center and radius of curvature of $\gamma$ in terms of its normal curvature at $p$.
Make a conclusion.

\parit{Source:}
This statement, as well as \ref{thm:meusnier} were proved by Jean Baptiste Meusnier \cite{meusnier}.

\parbf{\ref{ex:principal-revolution}.}
Use \ref{ex:line-of-curvature} and Meusnier's theorem (\ref{thm:meusnier}).

\parbf{\ref{ex:catenoid-is-minimal}.}
Use \ref{ex:principal-revolution:a},  \ref{ex:curvature-graph}, and \ref{thm:meusnier}.

\parbf{\ref{ex:helicoid-is-minimal}.} Apply Meusnier's theorem (\ref{thm:meusnier}) to show that the coordinate curves $\alpha_v\:u\mapsto s(u,v)$ and $\beta_u\:v\mapsto s(u,v)$ are asymptotic; that is, they have vanishing normal curvature.

Observe that these two families are orthogonal to each other.
Apply \ref{ex:mean-curvature}.

\parbf{\ref{ex:rev(sin)}.}
By \ref{ex:principal-revolution:a}, parallels and meridians are lines of curvature.
One principal curvature is the curvature of generatrix.
Use \ref{ex:curvature-formulas:a} to show that it is at most $a\cdot \sin x$.
Note that the curvature of parallel is $\tfrac1{a\cdot \sin x}$;
apply Meusnier's theorem, to show that the principal curvature in its direction cannot be larger.
Make a conclusion.

\parbf{\ref{ex:rev(lin)}.} Apply \ref{ex:principal-revolution:formula}.

\parbf{\ref{ex:moon-revolution}.}
Use \ref{ex:line-of-curvature} and \ref{thm:moon-orginal}.

\parbf{\ref{ex:lagunov-genus4}.}
Drill an extra hole or combine two examples together.
For the second part, combine $V_2$ with $V_3$.

\parit{Comment.}
It is possible to create such an example bounded by a torus, which is a very advanced exercise.
A body bounded by a sphere must include a ball of radius a bit larger than $r_3=\sqrt{3/2}-1>r_2$; see \cite{lagunov-1960, lagunov-1961, lagunov-fet-1963, lagunov-fet-1965}.

%\end{multicols}
%\par\noindent\rule{\textwidth}{0.4pt}
%\begin{multicols}{2}

\stepcounter{chapter}
\setcounter{eqtn}{0}

\parbf{\ref{ex:supp>tan}.}
We may assume that $\Sigma_1$ is defined implicitly by a smooth function $f$ in a small neighborhood $U$ of $p$, so $f |_{\Sigma_1 \cap U} \equiv 0$.
Since $\Sigma_1$ locally supports $\Sigma_2$ at $p$, we may assume that $f(x)\ge 0$ for all $x\in \Sigma_2\cap U$.

Show that $(f\circ\gamma_2)'(0)=0$ for any smooth curve $\gamma_2\:(-\varepsilon,\varepsilon)\to\Sigma_2$ such that $\gamma_2(0)=p$.
Deduce from this that the velocity vector $\gamma_2'(0)$ is tangent to $\Sigma_1$ at $p$ and make the last step.

\parbf{\ref{ex:surf-support}.}
Choose curvatures such that 
\[k_2(p)_{\Sigma_1}\z>k_2(p)_{\Sigma_2}> k_1(p)_{\Sigma_1}> k_1(p)_{\Sigma_2},\]
and suppose the first principal direction of $\Sigma_1$ coincides with the second principal direction of $\Sigma_2$ and vice-versa.

\parbf{\ref{ex:positive-gauss-0} $\bm{+}$ \ref{ex:positive-gauss}.}
Apply the same reasoning as in \ref{ex:between-parallels-1}--\ref{ex:lens}, but use families of spheres instead.

\parbf{\ref{ex:convex-surf}.}
Show that any tangent plane $\T_p$ supports $\Sigma$ at~$p$.
Apply \ref{prop:surf-support}.

\parbf{\ref{ex:convex-lagunov}.}
Assume a maximal ball in $R$ has radius $r$, and it touches the boundary of $R$ at the points $p$ and~$q$.

\columnbreak

\begin{wrapfigure}{r}{18 mm}
\vskip-0mm
\centering
\includegraphics{mppics/pic-1050}
\vskip-0mm
\end{wrapfigure}

Consider the projection $S$ of $R$ to a plane thru $p$, $q$ and the center of the ball.
Show that $S$ is bounded by a smooth closed convex curve with curvature at most $1$.
Argue as in \ref{ex:between-parallels-1} to
show that $r\ge 1$.

\parbf{\ref{ex:section-of-convex}.}
Choose a plane $\Pi$.
Suppose a point $p$ lies in the intersection $\Pi\cap\Sigma$.

Show that if $\Pi$ is tangent to surface at $p$,
then $p$ is an isolated point of the intersection $\Pi\cap\Sigma$.

It follows that if $\gamma$ is a connected component of the intersection $\Pi\cap\Sigma$ that is not an isolated point,
then $\Pi$ intersects $\Sigma$ \index{transversality}\emph{transversally} along $\gamma$;
that is, $\T_p\Sigma\ne\Pi$ for any point $p \in \gamma$.
Apply the implicit function theorem to show that $\gamma$ is a smooth curve.

Finally, the curvature of $\gamma$ cannot be smaller than the normal curvature of $\Sigma$ in the same direction.
Hence $\gamma$ has no points with vanishing curvature.
Therefore its signed curvature has a constant sign.

\parbf{\ref{ex:surrounds-disc}.}
Assume the contrary, then by \ref{thm:convex-embedded} the surface is convex.
Show that the surface has a supporting spherical dome with the unit circle as the boundary.
Arrive at a contradiction with~\ref{cor:surf-support}.

\parbf{\ref{ex:small-gauss}.}
We can assume that the Gauss curvature of the surface is positive; otherwise, the statement is evident.
By \ref{thm:convex-embedded}, the surface bounds a convex region that contains a line segment of length~$\pi$.

\begin{Figure}
\vskip-0mm
\centering
\includegraphics{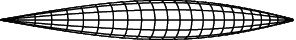}
\vskip-0mm
\end{Figure}

By \ref{ex:rev(sin)}, the Gauss curvature of the surface of revolution of the graph $y=a\cdot \sin x$ for $x\in(0,\pi)$ cannot exceed $1$.
Try to support the surface $\Sigma$ from the inside by a surface of revolution of the described type.
(Compare to \ref{ex:convex small}.)

\parit{Remark.}
The exercise can be deduced from the following deeper result: \textit{if the Gauss curvature of $\Sigma$ is at least $1$,
then
the intrinsic diameter of $\Sigma$ cannot exceed $\pi$}.
The latter means that any two points in $\Sigma$ can be connected by a path in $\Sigma$ that has length at most~$\pi$.
This theorem was proved by Heinz Hopf and Willi Rinow \cite{hopf-rinow} and 
named after Sumner Myers who generalized it \cite{myers}.

\parbf{\ref{ex:convex-proper-sphere}}; \ref{SHORT.ex:convex-proper-sphere:single}.
Use the convexity of~$R$.

\parit{\ref{SHORT.ex:convex-proper-sphere:smooth}.}
Observe that the map $\Sigma\to\mathbb{S}^2$ is smooth and regular, then applying the inverse function theorem, show that its inverse is smooth as well.

\begin{wrapfigure}{r}{29 mm}
\vskip-3mm
\centering
\includegraphics{mppics/pic-1182}
\vskip-3mm
\end{wrapfigure}

\parbf{\ref{ex:convex-proper-plane}}; \ref{SHORT.ex:convex-proper-plane:a}
(The argument is similar to \ref{prop:convex-monotone:open}.)
We can assume that the origin lies on~$\Sigma$.
Consider a sequence of points $x_n\in \Sigma$ such that $|x_n|\z\to \infty$ as $n\to \infty$.
Denote by $\vec u_n$ the unit vector in the direction of $x_n$; that is $\vec u_n\z=\tfrac{x_n}{|x_n|}$.

Since the unit sphere is compact, we can pass to a subsequence of $x_n$ such that $\vec u_n$ converges to a unit vector, say $\vec u$.
Show that the half-line from the origin in the direction of $\vec u$ can be taken as $\ell$.

\parit{\ref{SHORT.ex:convex-proper-plane:b} $+$ \ref{SHORT.ex:convex-proper-plane:c} $+$ \ref{SHORT.ex:convex-proper-plane:d}.}
Since $R$ is convex, so is its projection $\Omega$.

Note and use that for any $q\in \Sigma$, the directions $\vec v_n=\tfrac{x_n-q}{|x_n-q|}$ converge to $\vec u$ as well.

Show that $\Sigma$ has no vertical tangent planes.
Conclude that the projection map from $\Sigma$ to the $(x,y)$-plane is regular.
Use the inverse function theorem to show that $\Omega$ is open.

\parit{\ref{SHORT.ex:convex-proper-plane:e}.} 
Arguing by contradiction, suppose for some sequence $(x_n,y_n)\z\to(x_\infty,y_\infty)\in \partial\Omega$ the sequence $f(x_n,y_n)$ stays bounded above.
We can pass to a subsequence such that either $f(x_n,y_n)$ converges to a finite value, say $z_\infty$, or it diverges to $-\infty$.

In the first case, show that the point $(x_\infty, y_\infty,z_\infty)$ does not lie on $\Sigma$, but it has arbitrarily close points on~$\Sigma$.
That is $\Sigma$ is not proper --- a contradiction.

If $f(x_n,y_n)\to -\infty$, use the convexity of $f$ to show that $f(\tfrac{x_n}2 ,\tfrac{y_n}2)\z\to -\infty$.
Note that the origin belongs to $\Omega$;
use it to show that $(\tfrac{x_\infty}2, \tfrac{y_\infty}2)\in\Omega$.
Arrive at a contradiction.

\parbf{\ref{ex:open+convex=plane}.}
By \ref{ex:convex-proper-plane:d}, $\Sigma$ is parametrized by an open convex plane domain $\Omega$.
It remains to show that $\Omega$ can parametrize the whole plane.

We may assume that the origin of the plane lies in $\Omega$.
Show that in this case, the boundary of $\Omega$ can be written in polar coordinates as $(\theta,f(\theta))$ where $f\:\mathbb{S}^1\to\mathbb{R}$ is a positive continuous function.
Then a homeomorphism from $\Omega$ to the plane can be described in polar coordinates by changing only the radial coordinate;
for example, as 
$(\theta,r)\z\mapsto (\theta,
\tfrac{r}{1-r/f(\theta)})$.

In the second part, one may apply \ref{ex:star-shaped-disc:nonsmooth}.

\parbf{\ref{ex:circular-cone}.}
Choose a coordinate system so that the $(x,y)$-plane supports $\Sigma$ at the origin, so $\Sigma$ lies in the upper half-space. 

Show that there is $\epsilon>0$ such that any line starting from the origin with slope at most $\epsilon$ may intersect $\Sigma$ only in the unit ball centered at the origin;
we may assume that $\epsilon$ is small, say $\epsilon<1$.
Consider the cone formed by the half-lines from the origin with slope $\epsilon$ shifted down by $1$.
Observe that the entire surface lies inside this cone.

\parbf{\ref{ex:intK}}.
Choose distinct points $p,q\in\Sigma$.
Apply \ref{thm:convex-embedded} to show that the angle 
$\measuredangle(\Norm(p),p-q)$ is acute and $\measuredangle(\Norm(q),p-q)$ is obtuse.
Conclude that $\Norm(p)\ne\Norm(q)$;
that is, $\Norm\:\Sigma\to\mathbb{S}^2$ is injective.

\parit{\ref{SHORT.ex:intK:4pi}.}
Given a unit vector $\vec u$, consider a point $p\in \Sigma$ that maximizes the scalar product $\langle p,\vec u\rangle$.
Show that $\Norm(p)=\vec u$.
Conclude that the spherical map $\Norm\:\Sigma\to\mathbb{S}^2$ is onto, and therefore it is a bijection.

Applying \ref{cor:intK}, we get that the integral is $4\cdot\pi=\area\mathbb{S}^2$.

\parit{\ref{SHORT.ex:intK:2pi}.} Choose an $(x,y,z)$-coordinate system provided by \ref{ex:convex-proper-plane:d}.
Observe that for any $p$ the normal vector $\Norm(p)$ forms an obtuse angle with the $z$-axis.
It follows that the image $\Norm(\Sigma)$ lies in the southern hemisphere.

Applying \ref{cor:intK}, we get that the integral does not exceed $2\cdot\pi=\tfrac12\cdot\area\mathbb{S}^2$.

\parit{Remark.} This exercise reminds Corollary~\ref{cor:fenchel=convex}.

\parbf{\ref{ex:convex-revolution}.}
Apply \ref{ex:principal-revolution}.

\parbf{\ref{ex:ruled=>saddle}.}
Prove and use that each point on the surface has a direction with vanishing normal curvature.

\parbf{\ref{ex:saddle-convex}.} Suppose $p\in \Sigma$ is a point of local maximum of~$f$.
Show that $\Sigma$ is supported by its tangent plane at~$p$.
Arrive at a contradiction.

\parbf{\ref{ex:panov}.}
Denote by $\Pi_t$ the affine tangent plane to $\Sigma$ at $\gamma(t)$ and by $\ell_t$ the tangent line to $\gamma$ at time~$t$.

Note that $\Pi_t$ is the graph of a linear function, say $h_t$, defined on the $(x, y)$-plane.
Denote by $\bar\ell_t$ the projection of $\ell_t$ to the $(x, y)$-plane.
Show that the derivative $\tfrac{d}{dt}h_t(w)$ vanishes at the point $w$ if and only if $w\in \bar\ell_t$ 
and the derivative changes sign if $w$ goes from one side of $\bar\ell_t$ to the other.

If $\bar\gamma$ is star-shaped with respect to a point $w$, then $w$ cannot cross $\bar\ell_t$.
Therefore, the function $t\mapsto h_t(w)$ is monotone on $\mathbb{S}^1$.
Observe that this function cannot be constant, and arrive at a contradiction.

\parit{Source:}
The statement follows from the theorem of Galina Kovaleva \cite{kovaleva}, which was rediscovered by Dmitri Panov \cite{panov-curves}.

\parbf{\ref{ex:crosss}.}
This problem follows easily from the so-called \index{Morse lemma}\emph{Morse lemma}.
The following sketch is a slightly stripped version of it.
A more conceptual proof \cite{abraham-marsden-ratiu}
can be built on the Moser trick; see the solution of \ref{ex:star-shaped-disc}.

\medskip

Choose tangent-normal coordinates at $p$ so that the coordinate axes point in the principal directions;
let $z=f(x,y)$ be the local graph representation of~$\Sigma$.
We need to show that the solution of $f(x,y)=0$ is a union of two smooth curves that intersect transversely at~$p$.

Prove the following claim:
\textit{Suppose $x\z\mapsto h(x)$ is a smooth function defined in an open interval $\mathbb{I}\ni0$ such that $h(0)=h'(0)=0$ and $h''(0)>0$.
Then, for a smaller interval $\mathbb{J}\ni0$ there is a unique smooth function $a\:\mathbb{J}\to\mathbb{R}$ such that $h=a^2$, $a(0)=0$ and $a'(0)> 0$.}

Note that passing to a small square domain $|x|,|y|<\epsilon$, we can assume that $f_{xx}\z>\epsilon$ and $f_{yy}<-\epsilon$. 
Show that if $\epsilon$ is small, then for every $x$ there is unique $y(x)$ such that $f_x(x,y(x))=0$; 
moreover, the function $x\z\mapsto y(x)$ is smooth.

Set $h(x)=f(x,y(x))$.
Note that $h(0)\z=h'(0)=0$ and $h''>0$.
Applying the claim, we get a function $a$ such that $h=a^2$, where $a(0)=0$ and $a'(0)>0$.

Observe that $g(x,y)=h(x)-f(x,y)\ge 0$, $g_y(x,y(x))\z=g(x,y(x))\z=0$, and $g_{yy}>0$.
Applying the claim to each function $y\mapsto g(x,y)$ with fixed $x$, we get that $g(x,y)=b(x,y)^2$ for a smooth function $b$ such that 
$b(x,y(x))=0$ and $b_y(x,y(x))>0$.

It follows that 
\begin{align*}
f(x,y)&=a(x)^2-b(x,y)^2=
\\
&=
(a(x)-b(x,y))\cdot (a(x)+b(x,y)).
\end{align*}
That is $f(x,y)=0$ if $a(x)\pm b(x,y) =0$.

It remains to observe that the two functions $g_\pm(x,y)=a(x)\pm b(x,y)$ have distinct non-zero gradients at $0$.
Therefore, each equation $a(x)\pm b(x,y) =0$ defines a smooth curve in a neighborhood of $p$;
see Section~\ref{sec:implicit-curves}.

\parbf{\ref{ex:proper-saddle}.} Look at the pseudosphere described in \ref{ex:principal-revolution:pseudosphere}.

\parbf{\ref{ex:length-of-bry}.}
Use \ref{lem:convex-saddle} and the hemisphere lemma (\ref{lem:hemisphere}).

For the second part, consider a thin tube bounded by two closed spherical curves.

\parbf{\ref{ex:circular-cone-saddle}.}
Assume $\Sigma$ is an open saddle surface that lies in a cone~$K$.
Show that there is a plane $\Pi$ that cuts $\Sigma$ and cuts from $K$ a compact region.
Conclude that $\Pi$ cuts from $\Sigma$ a compact region as well. 

By \ref{lem:reg-section} one can move the plane $\Pi$ slightly so that it cuts from $\Sigma$ a compact surface with boundary.
Apply \ref{lem:convex-saddle}.

\parbf{\ref{ex:disc-hat}.} 
Consider the radial projection of $F_\epsilon$ to the sphere $\Sigma$ with center at $p=(0,0,\epsilon)$;
that is, a point $q\in F_\epsilon$ is mapped to a point $s(q)$ on the sphere that lies on the half-line $pq$.

Show that $s$ is a diffeomorphism from $F_\epsilon$ to the southern hemisphere of~$\Sigma$.
It remains to observe that the unit disc is diffeomorphic to the hemisphere.

\parbf{\ref{ex:saddle-linear}.} 
By the fundamental theorem of affine geometry, any affine transformation is smooth.
It remains to apply \ref{prop:hat}.

\parbf{\ref{ex:between-parallels}.}
Find an example among the surfaces of revolution;
use \ref{ex:principal-revolution}.

\parbf{\ref{ex:one-side-bernshtein}.} Look at the sections of the graph by planes parallel to the $(x,y)$-plane and to the $(x,z)$-plane, then apply Meusnier’s theorem (\ref{thm:meusnier}).

\parbf{\ref{ex:saddle-graph}.}
Suppose the orthogonal projection of $\Sigma$ to the $(x,y)$-plane is not injective.
Show that there is a point $p\in\Sigma$ with a vertical tangent plane;
that is, $\T_p\Sigma$ is perpendicular to the $(x,y)$-plane.

Let $\Gamma_p$ be the connected component of $p$ in the intersection of $\Sigma$ and the affine tangent plane $\T_p\Sigma$.
Use \ref{ex:crosss} to show that $\Gamma_p$ is a union of smooth curves that can cross each other transversely.
Moreover, two of these curves pass thru $p$, and $\Gamma_p$ does not bound a compact region on~$\Sigma$.

Show that $\Gamma_p$ must have at least 4 ways to escape to infinity.
On the other hand, since $\Sigma$ is a graph outside the compact set $K$, we have that $\Gamma_p\setminus K$ is a graph of a real-to-real function that has only two ways to escape to infinity --- a contradiction.

%\end{multicols}
%\par\noindent\rule{\textwidth}{0.4pt}
%\begin{multicols}{2}

\stepcounter{chapter}
\setcounter{eqtn}{0}

\parbf{\ref{ex:lasso}.}
Cut the lateral surface of the ice-mountain by a line from the cowboy to the tip.
Unfold it on the plane (see the picture) and try to figure out what is the image of the strained lasso.

\begin{wrapfigure}{r}{24 mm}
\vskip-6mm
\centering
\includegraphics{mppics/pic-1250}
\vskip-0mm
\end{wrapfigure}

\parit{Source:}
We learned this problem from Joel Fine, who attributed it to Frederic Bourgeois \cite{fine}; see also \cite[Problem 12]{khesin-tabachnikov}

\parbf{\ref{ex:length-dist-conv}.} 
Note that by \ref{thm:convex-embedded}, $\Sigma$ bounds a strictly convex region.
Therefore, we can assume that $\Norm(p)\ne\Norm(q)$; otherwise, $p=q$, and the inequality is evident.

Further, we can assume that $\Norm(p)+\Norm(q)\z\ne 0$; otherwise, the right-hand side is undefined.

In the remaining case, the tangent planes $\T_p$ and $\T_q$ intersect along a line, say $\ell$.
Set $\alpha=\tfrac12\cdot\measuredangle(\Norm(p),\Norm(q))$.
Show that $2\cdot\cos\alpha\z= |\Norm(p)+\Norm(q)|$.
Let $x\in \ell$ be the point that minimizes the sum $|p-x|\z+|x-q|$.
Show that $\measuredangle\hinge xpq\ge \pi-2\cdot\alpha$.
Conclude that 
\[|p-x|+|x-q|\le \tfrac{|p-q|}{\cos\alpha}.\]
Apply \ref{thm:shorts+convex} to show that
$\dist{p}{q}\Sigma\le |p-x|\z+|x-q|$.

\parbf{\ref{ex:hat-convex}.}
Suppose there is a minimizing curve $\gamma\not\subset\Delta$ with endpoints $p$ and $q$ in $\Delta$.

Without loss of generality, we may assume that only one arc of $\gamma$ lies outside $\Delta$.
Reflection of this arc with respect to $\Pi$ together with the remaining part of $\gamma$ forms another curve $\hat\gamma$ from $p$ to $q$;
it runs partly along $\Sigma$ 
and partly outside $\Sigma$,
but does not get inside.
Note that
\[\length\hat\gamma=\length\gamma.\]

Denote by $\bar\gamma$ the nearest-point projection of $\hat\gamma$ on~$\Sigma$.
The curve $\bar\gamma$ lies in $\Sigma$, 
it has the same endpoints as $\gamma$,
and by \ref{thm:shorts+convex}
\[\length\bar\gamma<\length\gamma.\]
This means that $\gamma$ is not length minimizing --- 
a contradiction.

\parbf{\ref{ex:intrinsic-diameter}.} 
Show that the nearest-point projection $\mathbb{S}^2\to\Sigma$ is surjective.
Apply \ref{lem:nearest-point-projection} and \ref{thm:area-axioms:monotonicity}.

%\end{multicols}
%\par\noindent\rule{\textwidth}{0.4pt}
%\begin{multicols}{2}

\stepcounter{chapter}
\setcounter{eqtn}{0}

\parbf{\ref{ex:helix-geodesic}.} Show that $\gamma''(t)$ is proportional to $\nabla_{\gamma(t)} f$, where $f=x^2+y^2$. 
Apply \ref{ex:tangent-normal}.

\parbf{\ref{ex:clairaut}.} 
We can assume that the origin lies on the axis of revolution, and $\vec i$ points in the direction of such axis.
Use \ref{lem:constant-speed} to show that it is sufficient to prove that 
$\langle\gamma'\times \gamma,\vec i\rangle$
is constant.

Since $\gamma''(t)\perp\T_{\gamma(t)}$, the three vectors $\vec i$, $\gamma$, and $\gamma''$ lie in the same plane.
In particular, $\langle\gamma''\times \gamma,\vec i\rangle=0$.
Therefore,
\[
\langle\gamma'\times \gamma,\vec i\rangle'
=
\langle\gamma'\times \gamma',\vec i\rangle+\langle\gamma''\times \gamma,\vec i\rangle =0
.\]

\parbf{\ref{ex:asymptotic-geodesic}.} By Lemma \ref{lem:constant-speed},
we can assume that $\gamma$ is parametrized by arc-length.
By the definition of geodesic, we have that $\gamma''(s)\perp\T_{\gamma(s)}$. 
Therefore, 
\[\gamma''(s)=k_n(s)\cdot\Norm(\gamma(s)),\]
where $k_n(s)$ is the normal curvature of $\gamma$ at time~$s$.
Since $\gamma$ is asymptotic, $k_n(s)\equiv 0$;
that is, $\gamma''(s)\equiv 0$, therefore $\gamma'$ is constant and hence $\gamma$ runs along a line segment; see \ref{ex:zero-curvature-curve}.

\parbf{\ref{ex:reflection-geodesic}.}
Denote by $\mu$ a unit vector perpendicular to $\Pi$.
Since $\gamma$ lies in $\Pi$, we have that $\gamma''$ is parallel to $\Pi$, or equivalently $\gamma''\perp \mu$.
Since $\gamma$ is parametrized by arc-length, \ref{prop:a'-pertp-a''} implies that $\gamma''\perp\gamma'$.

Since $\Sigma$ is mirror-symmetric with respect to the plane $\Pi$,
the tangent plane $\T_{\gamma(t)}\Sigma$ is also mirror-symmetric with respect to $\Pi$.
It follows that $\T_{\gamma(t)}\Sigma$ is spanned by $\mu$ and $\gamma'(t)$.
Hence, $\gamma''\perp \mu$ and $\gamma''\perp\gamma'$ imply $\gamma''\perp\T_{\gamma(t)}\Sigma$;
that is, $\gamma$ is a geodesic.

\parbf{\ref{ex:round-torus}.}
By \ref{ex:reflection-geodesic}, any meridian of $\Sigma$ is a closed geodesic.
Consider an arbitrarily closed geodesic~$\gamma$.

If $\gamma$ is tangent to a meridian at some point, then by the uniqueness part of Proposition \ref{prop:geod-existence}, $\gamma$ runs along that meridian;
in particular, it is non-contractible.

In the remaining case, $\gamma$ can intersect meridians only transversely.
Therefore, the longitude of $\gamma$ is monotone.
Whence again $\gamma$ is non-contractible.

\parbf{\ref{ex:helix=geodesic}.}
Check that $\Norm(\gamma(t))=(\cos t,\sin t,0)$.
Calculate $\gamma''(t)$, and show that it is proportional to $\Norm(\gamma(t))$.
The latter is equivalent to $\gamma''(t)\perp\T_{\gamma(t)}$. 

Note that the line segment from $\gamma (0) $ to $\gamma (2{\cdot}\pi) $ is contained in~$\Sigma$.

\parbf{\ref{ex:two-min-geod}.}
Assume two shortest paths $\alpha$ and $\beta$ have two common points $p$ and~$q$.
Denote by $\alpha_1$ and $\beta_1$ the arcs of $\alpha$ and $\beta$ from $p$ to~$q$.
Suppose $\alpha_1$ is distinct from $\beta_1$.

Note that $\alpha_1$ and $\beta_1$ are shortest paths with the same endpoints;
in particular, they have the same length.
Exchanging $\alpha_1$ in $\alpha$ to $\beta_1$ produces a shortest path, say $\hat\alpha$, that is distinct from $\alpha$.
By \ref{prop:gamma''}, $\hat\alpha$ is a geodesic.

\begin{Figure}
\vskip-0mm
\centering
\includegraphics{mppics/pic-308}
\vskip0mm
\end{Figure}

Suppose $\alpha_1$ is a proper subarc of $\alpha$;
that is, $\alpha_1\ne\alpha$, or equivalently, $p$ or $q$ is not an endpoint of $\alpha$.
Then $\alpha$ and $\hat\alpha$ share one point and velocity vector at this point.
By \ref{prop:geod-existence}, $\alpha$ coincides with $\hat\alpha$ --- a contradiction.

It follows that $p$ and $q$ are the endpoints of~$\alpha$.
Analogously, $p$ and $q$ are the endpoints of~$\beta$.

For the second part, one could consider two distinct geodesics of the form 
\[ \gamma_b(t) = ( \cos t , \sin t , b\cdot t ) , t \in \mathbb{R} \]
in the cylinder $x^2 + y^2 =1$.

\parbf{\ref{ex:min-geod+plane}.}
Assume a shortest path $\alpha$ crosses $\Pi$ at least twice.
In this case, there is a subarc $\alpha_1$ of $\alpha$ that lies strictly on one side of $\Pi$;
only its endpoints are on $\Pi$, and these endpoints are distinct from the endpoints of $\alpha$.
 
\begin{Figure}
\vskip-1mm
\centering
\includegraphics{mppics/pic-310}
\vskip-1mm
\end{Figure}

Let us remove the arc $\alpha_1$ from $\alpha$ and replace it with the reflection of $\alpha_1$ across $\Pi$.
Note that the obtained curve, say $\beta$, lies on the surface;
it has the same length as $\alpha$, and it connects the same pair of points, say $p$ and~$q$.
Therefore, $\beta$ is a shortest path from $p$ to $q$ in $\Sigma$ that is distinct from~$\alpha$.

We may assume that both $\alpha$ and $\beta$ have arc-length parametrization.
By \ref{prop:gamma''}, $\alpha$ and $\beta$ are geodesics.
Since $\alpha$ and $\beta$ have a common subarc, they share one point and velocity vector at this point;
by \ref{prop:geod-existence} $\alpha$ coincides with $\beta$ --- a contradiction.

\parbf{\ref{ex:milka}.}
Let $W$ be the closed region outside~$\Sigma$.
Show that the distance $\dist{p^s}{q}W$ is constant with respect to~$s$.
This implies that the concatenation of the line segment $[p^s,\gamma(s)]$ and the arc $\gamma|_{[s,\ell]}$ is a shortest path from $p^s$ to $q$ in~$W$.

Since $\Sigma$ is strictly convex, 
\[ \dist{\gamma (s)}{q}W > \dist{\gamma(s)}{q}{\mathbb{R}^3} \]
for all $s < \ell$.
Hence,
\begin{align*}
\dist{p^s}{q}W&=\dist{p^s}{\gamma(s)}W+\dist{\gamma (s)}{q}W 
> 
\\
&>\dist{p^s}{\gamma(s)}{\mathbb{R}^3} + \dist{ \gamma (s)}{q}{\mathbb{R}^3} 
\ge
\\
&\ge\dist{p^s}{q}{\mathbb{R}^3}. 
\end{align*}
If the segment $[p^s,q]$ was entirely contained in $W$ for some $s<\ell$, then $\dist{p^s}{q}W = \dist{p^s}{q}{\mathbb{R}^3} $, which would be a contradiction.

\parbf{\ref{ex:round-sphere}.}
Argue as in \ref{ex:const-dist} to show that all geodesics in $\Sigma$ have the constant curvature.
Conclude that $\Sigma$ has constant normal curvature and argue as in \ref{ex:normal-curvature=const}.

\parit{Source:} The problem is due to Ali Taghavi \cite{taghavi}.

\parbf{\ref{ex:rad=2}}; \ref{SHORT.ex:rad=2:a}.
Consider a geodesic $\gamma$ on $\Sigma$.
Observe that its curvature does not exceed $1$.
Argue as in \ref{ex:gromov-twist} to show that there is $\epsilon>0$ such that
\[\epsilon<|\gamma(t)|<2\cdot\sin(\measuredangle(\gamma'(t),\gamma(t)))\]
for any $t$.
Conclude that we can choose a unit normal field $\Norm$ on $\Sigma$ such that
\[\epsilon<|x|<2\cdot\cos(\measuredangle(\Norm(x),x)) \eqlbl{eq:|gamma(t)|}\]
for any point $x\in\Sigma$.

Use \ref{eq:|gamma(t)|} to show that the map $\Sigma\to \Theta$ defined by $x\mapsto 2\cdot\tfrac x{|x|}$ is locally invertible.
Conclude that it is invertible.
Its inverse is the radial projection $\rho\:\Theta\to\Sigma$.

\parit{\ref{SHORT.ex:rad=2:b}.}
Suppose $x=\rho(p)$ and $\vec v\in \T_p\Theta$.
Show that
\[2\cdot\cos(\measuredangle(\Norm(x),x))\cdot|d_p\rho(\vec v)|
\le
|x|\cdot|\vec v|.\]
Apply \ref{eq:|gamma(t)|} to show that $|d_p\rho(\vec v)|
\le|\vec v|$; conclude that $\rho$ is length-nonincreasing.

\parit{\ref{SHORT.ex:rad=2:c}.}
Assume $y\in \Sigma$ minimizes distance to~$\tfrac x2$.
Let $x=\rho(p)$ and $y=\rho(q)$.

Arguing by contradiction, suppose that
\[\dist{\tfrac x2}{y}{}\z<|\tfrac x2|.\]
Note that $\Norm\:\Sigma\to \mathbb{S}^2$ is length-nonincreasing.
Use all this and \ref{SHORT.ex:rad=2:b} to deduce the following
\begin{align*}
\measuredangle(\Norm(x),\Norm(y))&= \measuredangle(x,y-\tfrac x2)>
\\
&>2\cdot\measuredangle\hinge 0pq=
\\
&=\dist{p}{q}{\Theta}\ge
\\
&\ge\dist{x}{y}{\Sigma}\ge
\\
&\ge\measuredangle(\Norm(x),\Norm(y));
\end{align*}
which is a contradiction.

Finally, we can move $\Sigma$ so that $x$ gets arbitrarily close to $\Theta$.
Applying the above argument, we get that $\Sigma$ surrounds a unit ball that touches $\Sigma$ at $x$.

\parit{Source:} Problem taken from the paper by Hongda Qiu \cite{qiu2025}.

\parbf{\ref{ex:closed-liberman}.} Assume $t\mapsto \gamma(t)=(x(t),y(t),z(t)$ is a closed unit-speed geodesic in the graph $z=f(x,y)$.
By Liberman's lemma, $t\mapsto z'(t)$ is monotonic.
Conclude that $z$ is constant on $\gamma$, so $\gamma$ runs in a horizontal plane.
Show that $\gamma$ is straight and arrive at a contradiction.

\parbf{\ref{ex:rho''}.}
Equip $\Sigma$ with unit normal field $\Norm$ that points inside.
Denote by $k(t)$ the normal curvature of $\Sigma$ at $\gamma(t)$ in the direction of $\gamma'(t)$.
Since $\Sigma$ is convex, $k(t)\ge 0$ for any~$t$.

Since $\gamma$ is a unit-speed geodesic, we have $\gamma''(t)=k(t)\cdot\Norm(\gamma(t))$ and $\langle\gamma'(t),\gamma'(t)\rangle=1$ for any~$t$.
Without loss of generality, we can assume that $p$ is the origin of~$\mathbb{R}^3$.
Since $p$ is inside $\Sigma$, we have that $\langle\gamma(t),\Norm(\gamma(t))\rangle\le 0$ for any~$t$.
It follows that 
\[\langle\gamma''(t),\gamma(t)\rangle=k(t)\cdot \langle\gamma(t),\Norm(\gamma(t))\rangle\le 0\]
for any~$t$.
Conclude that
\begin{align*}
&\rho''(t)
=\langle\gamma(t),\gamma(t)\rangle''=
\\
&=2\cdot\langle\gamma''(t),\gamma(t)\rangle+2\cdot\langle\gamma'(t),\gamma'(t)\rangle\le 2.
\end{align*}

\parbf{\ref{ex:tc-spherical-image}.}
We can assume $\gamma$ is unit-speed.
Set $\Norm (t) = \Norm ( \gamma (t)) $. 
Then $\langle \gamma'(t) , \Norm (t) \rangle \equiv 0$. Also, since $\gamma$ is a geodesic, we have that $\gamma''(t) \parallel \Norm (t)$.
Therefore,
\[
\begin{aligned}
|\gamma''(t)|
&=|\langle\gamma''(t),\Norm(t)\rangle|=
\\
&=|\langle\gamma'(t),\Norm(t)\rangle'-\langle\gamma'(t),\Norm'(t)\rangle|=
\\
&=|\langle\gamma'(t),\Norm'(t)\rangle|\le
\\
&\le|\gamma'(t)|\cdot|\Norm'(t)|=|\Norm'(t)|.
\end{aligned}
\]

Integrating with respect to $t$, we get 
\[\tc\gamma\le\length(\Norm\circ\gamma).\]

\parbf{\ref{ex:usov-exact}.} 
Suppose $\gamma(t)=(x(t),y(t),z(t))$. 
Show that
\[|\gamma'' (t)| = z''(t)\cdot\sqrt{1+ \ell ^2}\eqlbl{eq:gamma''=z''}\]
for any~$t$.

Observe that $z'(t)\to\pm \tfrac\ell{\sqrt{1+ \ell ^2}}$ as $t\to\pm\infty$.
Conclude that 
\[\int_{-\infty}^{+\infty}z''(t)\cdot dt
=
\frac{2\cdot\ell}{\sqrt{1+ \ell ^2}}.\eqlbl{eq:int z''}\]
By \ref{eq:gamma''=z''} and \ref{eq:int z''}, we have
\begin{align*}
\tc\gamma&=\int_{-\infty}^{+\infty}|\gamma''(t)|\cdot dt=
\\
&=
\sqrt{1+ \ell ^2}\cdot \int_{-\infty}^{+\infty}z''(t)\cdot dt=
\\
&=2\cdot \ell.
\end{align*}

\parbf{\ref{ex:ruf-bound-mountain}.}
Use \ref{thm:usov} and \ref{ex:sef-intersection}.
For the second part, consider a geodesic on a cone with slope $2$, and mollify its tip.

\parit{Remark.}
The statement still holds for $\sqrt{3}$-Lipschitz functions, and $\sqrt{3}=\tg\tfrac\pi3$ is the optimal bound; see \ref{ex:sqrt(3)}.
It is the same slope as in the problem about the cowboy and the lasso (\ref{ex:lasso}).

\parbf{\ref{ex:bound-tc}}; \ref{SHORT.ex:bound-tc:a}.
Show and use that $|\gamma|\le 1$, and
\[\langle\Norm,\epsilon\cdot \Norm-\gamma\rangle\le0.\]

\parit{\ref{SHORT.ex:bound-tc:b}.}
Show and use that $|\gamma|\le1$, $|\gamma'|=1$, and
\[\rho'=2\cdot \langle\gamma,\gamma'\rangle.\]

\parit{\ref{SHORT.ex:bound-tc:c}.}
Show and use that $\langle\gamma,\Norm\rangle=|\gamma|\cdot\cos\theta$, $|\gamma'|=1$, $\gamma''=-k\cdot \Norm$, and
\[\rho''=2\cdot \langle\gamma',\gamma'\rangle+2\cdot \langle\gamma,\gamma''\rangle.\]

\parit{\ref{SHORT.ex:bound-tc:d}.}
Suppose $x,y$ are the points in the unit sphere that project to the ends of $\gamma$. 
Connect $x$ to $y$ by an arc of length $\le \pi$ in the sphere, project it to $\Sigma$.
Apply \ref{lem:nearest-point-projection} and use that $\gamma$ is a shortest path.

\parit{\ref{SHORT.ex:bound-tc:e}.}
Suppose $\gamma\:[0,\ell]\to\Sigma$ is a unit-speed parametrization.
Use \ref{SHORT.ex:bound-tc:d} to show that $\ell\le \pi$.
Use \ref{SHORT.ex:bound-tc:a}, \ref{SHORT.ex:bound-tc:b}, \ref{SHORT.ex:bound-tc:c}, and $|\gamma|\ge \epsilon$, to show that 
\[2\cdot \ell-2\cdot \epsilon^2\cdot \tc{\gamma}
\ge
\int_0^\ell\rho''(t)\cdot dt.\]
Use \ref{SHORT.ex:bound-tc:b} to show that $\int\rho''\ge -4$.
Make a conclusion.

%\end{multicols}
%\par\noindent\rule{\textwidth}{0.4pt}
%\begin{multicols}{2}
 
\stepcounter{chapter}
\setcounter{eqtn}{0}

\parbf{\ref{ex:parallel}}; \ref{SHORT.ex:parallel:a}.
Show and use that $\langle\vec v(t),\vec v'(t)\rangle=0$.

\parit{\ref{SHORT.ex:parallel:b}}
Show that $|\vec v(t)|$, $|\vec w(t)|$, and
$\langle\vec v(t),\vec w(t)\rangle$,
are constants; it can be done the same way as \ref{SHORT.ex:parallel:a}.
Then use that 
$\langle\vec v(t),\vec w(t)\rangle\z=|\vec v(t)|\cdot|\vec w(t)|\cdot\cos\theta$.

\stepcounter{chapter}
\setcounter{eqtn}{0}

\parbf{\ref{ex:parallel-transport-support}.}
Observe that $\Sigma_1$ supports $\Sigma_2$ at any point of~$\gamma$.
Conclude that $\gamma$ has identical normals as a curve in $\Sigma_1$ and in $\Sigma_2$.
Apply \ref{obs:parallel=}.

\parbf{\ref{ex:holonomy=not0}.}
Nearly any loop solves the problem.
For example, consider the right-angled spherical triangle that an octant of $\mathbb{R}^3$ cuts from the sphere.
Argue that parallel transport around it rotates the tangent plane by the angle $\tfrac\pi 2$. 

\stepcounter{chapter}
\setcounter{eqtn}{0}

\parbf{\ref{ex:1=geodesic-curvature}.}
By \ref{ex:convex-proper-sphere}, $\Sigma$ is a smooth sphere.
By Jordan's theorem (\ref{thm:jordan}) the curve $\gamma$ divides $\Sigma$ into two discs.
Let us denote by $\Delta$ the disc that lies on the left from~$\gamma$.

Observe that $\tgc\gamma=\length\gamma$, and apply the Gauss--Bonnet formula (\ref{thm:gb}) for $\Delta$.

\parbf{\ref{ex:GB-hat}.}
Show that $\skur=\cos\alpha\cdot\skur_0$,
where $\skur$ is the geodesic curvature of $\partial \Delta$ in $\Delta$
and $\skur_0$ is the geodesic curvature of $\partial \Delta$ in the $(x,y)$-plane.
Apply \ref{prop:total-signed-curvature} and the Gauss--Bonnet formula (\ref{thm:gb}).

\parit{Remark.}
By \ref{ex:shape-curvature-line}, $\partial\Delta$ is a line of curvature on the graph.

\parbf{\ref{ex:geodesic-half}.}
Apply \ref{cor:intK}, \ref{ex:intK:4pi}, and the Gauss--Bonnet formula (\ref{thm:gb}).

For the last part, apply \ref{ex:bisection-of-S2}.

\parbf{\ref{ex:closed-geodesic}.} For the first part apply the Gauss--Bonnet formula (\ref{thm:gb}).

For the second part, arguing by contradiction, assume two closed geodesics $\gamma_1$ and $\gamma_2$ do not intersect. 
This would imply that $\gamma_2$ lies in one of the regions, say $R_1$, that $\gamma_1$ cuts from~$\Sigma$.
Similarly, $\gamma_1$ lies in one of the regions, say $R_2$, that $\gamma_2$ cuts from~$\Sigma$.

Observe that $R_1$ and $R_2$ cover $\Sigma$ with overlap.
Therefore, the first part implies that the integral of the Gauss curvature over $\Sigma$ is less than $4\cdot\pi$.
The latter contradicts \ref{ex:intK:4pi}.

\parbf{\ref{ex:self-intersections}}; \textit{(easy)}.
Consider the 4 regions bounded by loops.
Apply the Gauss--Bonnet formula (\ref{thm:gb}) to show that the integral of the Gauss curvature on each of these regions exceeds $\pi$.
The latter contradicts \ref{ex:intK:4pi}.

\parit{(tricky)}.
Denote by $\alpha$, $\beta$, and $\gamma$ the angles of the triangle.
Apply the Gauss--Bonnet formula (\ref{thm:gb}) to show that the loops surround regions over which the integral of the Gauss curvature is $\pi+\alpha$, $\pi+\beta$, and $\pi+\gamma$, respectively.

Apply the Gauss--Bonnet formula to the triangular region to show that $\alpha+\beta+\gamma>\pi$.
Use \ref{ex:intK:4pi} to get at a contradiction.

\parit{(hopeless)}.
It is difficult.
A solution based on \ref{thm:comp:toponogov} is given in \cite{petrunin2021}.

\parbf{\ref{ex:sqrt(3)}.}
It is sufficient to show that the surface has no geodesic loops.
Assume that there is a loop, estimate the integrals of the Gauss curvature over the entire surface and a disc surrounded by a geodesic loop.

\parbf{\ref{ex:unique-geod}}; \ref{SHORT.ex:unique-geod:unique}.
By \ref{prop:shortest-paths-exist} and \ref{prop:gamma''}, any two points in $\Sigma$ can be connected by a geodesic.
Suppose points $p$ and $q$ can be connected by two distinct geodesics $\gamma_1$ and $\gamma_2$.
Passing to their subarcs, we may assume that $\gamma_1$ and $\gamma_2$ share only their endpoints.
Since the surface is simply-connected, $\gamma_1$ and $\gamma_2$ together bound a disc, say $\Delta\subset\Sigma$.
It remains to apply the Gauss--Bonnet formula to $\Delta$ and make a conclusion.
 
\parit{\ref{SHORT.ex:unique-geod:diffeomorphism}.} 
Use part \ref{SHORT.ex:unique-geod:unique} and \ref{prop:inj-rad}.

\parbf{\ref{ex:half-sphere-total-curvature}.}
Apply \ref{prop:pt+tgc} and \ref{prop:area-of-spher-polygon}.

\parbf{\ref{ex:cohn-vossen}.}
Repeat the end of the proof of \ref{thm:cohn-vossen} for a one-parameter family of geodesics $\gamma_\tau$ defined by $\gamma_\tau(0)=\alpha(\tau)$ and $\gamma'_\tau(0)=\alpha'(\tau)$.

\parbf{\ref{ex:3-curves}.}
The first curve corresponds to a geodesic on a smoothed cone.
Check that the second and third are ruled out by the argument in the proof of the theorem.

\parbf{\ref{ex:g-b-chi}.}
Subdivide each surface into discs,
count number of edges, vertices, and discs in the subdivision.
Calculate the Euler characteristic and apply \ref{thm:GB-generalized}.
For the Möbius band and pair of pants, use that the boundary has vanishing geodesic curvature.
For the cylinder, use that each component of the boundary is a simple plane curve and apply \ref{prop:total-signed-curvature}. 

%\end{multicols}
%\par\noindent\rule{\textwidth}{0.4pt}
%\begin{multicols}{2}

\stepcounter{chapter}
\setcounter{eqtn}{0}

\parbf{\ref{ex:semigeodesc-chart}.}
By the Gauss lemma (\ref{lem:palar-perp}), polar coordinates with respect to $q$ produce a semigeodesic chart at any nearby point.
Therefore, it is sufficient to find a point $q\ne p$ such that polar coordinates on $\Sigma$ with respect to $q$ cover~$p$.
By \ref{prop:exp}, any $q$ sufficiently close to $p$ does the job.

\parbf{\ref{ex:inj-rad}}; \ref{SHORT.ex:inj-rad:sign}.
Show that we can choose an orientation on $\Sigma$ so that $b_r(0,\theta)\z=1$ for any $\theta$.
Conclude that we can assume that $b(r,\theta)>0$ for all small $r>0$.

Suppose $b(r_1,\theta_1)<0$ at a pair $(r_1,\theta_1)$ such that $0<r_1<r_0$.
Observe that if $\theta_2$ is sufficiently close to $\theta_1$, then the radial curves $r\mapsto b(r,\theta_1)$ and $r\mapsto b(r,\theta_2)$ defined on the interval $(0,r_1)$ intersect.
Therefore, $\exp_p$ is not injective in $B$ --- a contradiction.

\parit{\ref{SHORT.ex:inj-rad:0}.}
Suppose $b(r_1,\theta_1)=0$.
Apply \ref{SHORT.ex:inj-rad:sign} to show that $b_r(r_1,\theta_1)=0$.

Apply \ref{prop:jaccobi} to conclude $b(r,\theta_1)=0$ for any~$r$.
The latter contradicts that $b_r(0,\theta_1)=1$.

\parit{\ref{SHORT.ex:inj-rad:prop:inj-rad}.}
We need to show that $\exp_p$ is regular in~$B$.
Suppose vector $\vec v\in B$ has polar coordinates $(r,\theta)$ for some $r>0$.
Show that $\exp_p$ is regular at $\vec v$ if $b(r,\theta)\ne 0$.
Conclude that $\exp_p$ is regular in $B\setminus \{0\}$.

By \ref{obs:d(exp)=1}, $\exp_p$ is regular at $0$.
Whence $\exp_p|_B$ is a regular injective smooth map.
Use the inverse function theorem (\ref{thm:inverse}) to show that the restriction $\exp_p|_B$ is a diffeomorphism to its image.

\parbf{\ref{lem:K(orthogonal)}}; \ref{SHORT.lem:K(orthogonal):uu-vv}.
Since the frame $\vec u, \vec v,\Norm$ is orthonormal,
the first two vector identities are equivalent to the following six real identities:
\[
\begin{aligned}
\langle\vec u_u,\vec u\rangle
&=0,
&
\langle\vec v_u,\vec v\rangle
&=0,
\\
\langle\vec u_u,\vec v\rangle
&=-\tfrac1{b}\cdot a_v,
&
\langle\vec v_u,\vec u\rangle
&=
\tfrac1{b}\cdot a_v,
\\
\langle\vec u_u,\Norm\rangle
&=a\cdot \ell,
&
\langle\vec v_u,\Norm\rangle
&=
a\cdot m.
\end{aligned}
\eqlbl{eq:key-orthogonal/2}
\]

Taking the partial derivatives of the identities
$\langle\vec u,\vec u\rangle=1$ and
$\langle\vec v,\vec v\rangle=1$ with respect to $u$,
we get the first two identities in \ref{eq:key-orthogonal/2}.

Further, observe that
\[
\begin{aligned}
\vec v_u
&=
\tfrac{\partial}{\partial v}
(\tfrac1b\cdot s_v)=\tfrac1b\cdot s_{uv}
+
\tfrac{\partial}{\partial u}(\tfrac1b)
\cdot
 s_v.
\end{aligned}
\eqlbl{eq:dv/du}
\]
Since $s_u\perp s_v$, it follows that
\begin{align*}
\langle\vec v_u,\vec u\rangle
&=
\tfrac1{a\cdot b}\cdot \langle s_{vu}, s_u\rangle=
\\
&=
\tfrac1{2\cdot a\cdot b}\cdot \tfrac{\partial}{\partial v}\langle s_u, s_u\rangle=
\\
&=
\tfrac1{2\cdot a\cdot b}\cdot \tfrac{\partial a^2}{\partial v}=
\\
&=
\tfrac1{b}\cdot a_v.
\end{align*}
Taking the partial derivative of $\langle\vec u,\vec v\rangle=0$ with respect to $u$
we get
\begin{align*}
\langle\vec v_u,\vec u\rangle+
\langle\vec v,\vec u_u\rangle
&=0.
\end{align*}
Hence, we get two more identities in \ref{eq:key-orthogonal/2}.

Since $\vec u$ and $\vec v$ are orthonormal, \ref{thm:shape-chart} implies
\[
\begin{aligned}
\tfrac1{a^2}
\cdot
\langle s_{uu},\Norm\rangle
&=\ell,
&
\tfrac1{a\cdot b}
\cdot
\langle s_{uv},\Norm\rangle
&=m,
\\
\tfrac1{a\cdot b}
\cdot
\langle s_{vu},\Norm\rangle
&=m,
&
\tfrac1{b^2}
\cdot
\langle s_{vv},\Norm\rangle
&=n.
\end{aligned}
\eqlbl{eq:shape-lmn}
\]

Applying \ref{eq:dv/du}, \ref{eq:shape-lmn}, and $s_v\perp\Norm$ we get
\begin{align*}
\langle\vec u_u,\Norm\rangle
&=
\tfrac1{a}\cdot \langle s_{uu},\Norm\rangle
=a\cdot \ell,
\\
\langle\vec v_u,\Norm\rangle
&=
\tfrac1{a}\cdot \langle s_{uv},\Norm\rangle
=a\cdot m.
\end{align*}
It implies the last two equalities in \ref{eq:key-orthogonal/2}.

Therefore, the first two identities in \ref{SHORT.lem:K(orthogonal):uu-vv} are proved;
the remaining two identities can be proved along the same lines.

\parit{\ref{SHORT.lem:K(orthogonal):K}.}
Recall that the Gauss curvature equals the determinant of the matrix $
(\begin{smallmatrix}
\ell&m
\\
m&n
\end{smallmatrix}
)
$;
that is, $K=\ell\cdot n-m^2$.
Therefore, 
\begin{align*}
a\cdot b\cdot K
&=
a\cdot b\cdot(\ell\cdot n-m^2)
=
\\
&=
\langle\vec u_u,\vec v_v\rangle 
-
\langle\vec u_v,\vec v_u\rangle= 
\\
&= 
\left(
\tfrac{\partial}{\partial v}
\langle\vec u_u,\vec v\rangle
-
\cancel{\langle\vec u_{uv},\vec v\rangle}
\right)-
\\
&-
\left(
\tfrac{\partial}{\partial u}
\langle\vec u_v,\vec v\rangle
-
\cancel{\langle\vec u_{uv},\vec v\rangle}
\right)=
\\
&=\tfrac{\partial}{\partial v}(-\tfrac1{b}\cdot a_v)
-
\tfrac{\partial}{\partial u}(\tfrac1{a}\cdot b_u).
\end{align*}

\parbf{\ref{ex:conformal}.}
Apply \ref{lem:K(orthogonal)} assuming that $a=b$, and simplify.

\parbf{\ref{ex:isom(geod)}.}
Verify that an intrinsic isometry maps shortest paths to shortest paths and apply \ref{prop:gamma''}.

\parbf{\ref{ex:K=0}.}
Note and use that by \ref{prop:rauch}, $\exp_p$ is length-preserving.

\parbf{\ref{ex:K=1}.} Modify the proof of \ref{prop:rauch} to show that $K\equiv 1$ implies that $b(\theta,r)\z=\sin r$ for all small $r\ge 0$.

\parbf{\ref{ex:deformation}.}
To find $x(t)$, one may solve the differential equation
$x'(t)=\sqrt{1-|y'(t)|^2}$.
Its solution is defined while $y(t)>0$ and $|y'(t)|\z=|a\cdot \sin t|<1$.

Apply \ref{ex:principal-revolution:formula} to calculate Gauss curvature of~$\Sigma_a$.

By \ref{ex:K=1}, small discs $\Delta_a$ in $\Sigma_a$ admit intrinsic isometry to a disc in a unit sphere.
To show that $\Delta_a$ is not congruent to a spherical disc $\Delta_1$ for $a\ne 1$, show that its principal curvature is not 1 at some point.

\parbf{\ref{ex:line-cylinder}.}
Observe that for any two points $p$ and $q$ in the $(u,v)$ plane we have $|f(p)-f(q)|\le |p-q|$.

Let $f(u,v)=(x(u,v),y(u,v),z(u,v))$;
we can assume that $z(0,v)=v$ for any $v$.

Apply the observation and the triangle inequality to points $(0,-a)$, $(u,v_1)$, $(u,v_2)$, and $(0,a)$ for $a\to \infty$ to show that 
\begin{align*}
z(u,v_1)&=v_1,
&
z(u,v_2)&=v_2,
\\
x(u,v_1)&=x(u,v_2),
&
y(u,v_1)&=y(u,v_2).
\end{align*}
Make a conclusion.

%\end{multicols}
%\par\noindent\rule{\textwidth}{0.4pt}
%\begin{multicols}{2}

\stepcounter{chapter}
\setcounter{eqtn}{0}

\parbf{\ref{ex:wide-hinges}.}
Apply the triangle inequality $c_n\le a_n+b_n$ and the bound $a_n,b_n>\epsilon$ to show that the sequence 
\[s_n=\frac{a_n+b_n+c_n}{2\cdot a_n\cdot b_n}\]
is bounded.
Conclude that 
\[\frac{(a_n+b_n)^2-c_n^2}{2\cdot a_n\cdot b_n}=s_n\cdot (a_n+b_n-c_n)\to 0\]
as $n\to\infty$.
Use the last statement together with the cosine rule
\[a_n^2+b_n^2-2\cdot a_n\cdot b_n\cdot\cos\tilde\theta_n -c_n^2=0\]
to show that $\cos\tilde\theta_n\to -1$ as $n\to\infty$;
conclude that $\tilde\theta_n\to \pi$.

\parbf{\ref{ex:thm:comp:cat:nsc}.}
Consider the triangle on the hyperboloid with vertices $(1,0,0)$, $(-\tfrac{1}2, \tfrac{\sqrt{3}}2,0)$, and $(-\tfrac{1}2,-\tfrac{\sqrt{3}}2,0)$.
Note that all its angles are $\pi$, but all model angles are $\tfrac{\pi}3$.

\parbf{\ref{ex:diam-angle}.}
Show that $\dist{p}{x}\Sigma\le \dist{p}{q}\Sigma$ and $\dist{q}{x}\Sigma\le \dist{p}{q}\Sigma$.
Conclude that
\[\modangle xpq\ge \tfrac\pi3.\]
Apply \ref{thm:comp:toponogov}.

\parbf{\ref{ex:sum=<2pi}.}
Show that 
$\measuredangle\hinge pxy+\measuredangle\hinge pyz+\measuredangle\hinge pzx\le2\cdot \pi$.
Apply \ref{thm:comp:toponogov}.

\parbf{\ref{ex:s-r}.}
Choose vectors $\vec u$, and $\vec w$ 
such that $|\vec w|\z=\dist{p}{x}{}$, $|\vec u|=1$, and $\measuredangle(\vec u,\vec w)=\measuredangle\hinge xpq$.
Consider the function
$f\:t\mapsto t+|\vec w|-|\vec w-t\cdot \vec u|$.
Observe that $f(0)\z=0$,
\begin{align*}
f(\dist{x}{y}{})&=\dist{x}{p}{}+\dist{x}{y}{}-\side\hinge xpy,
\\
f(\dist{x}{q}{})&=\dist{x}{p}{}+\dist{x}{q}{}-\side\hinge xpq.
\end{align*}
Show and use that $f$ is concave.

\parbf{\ref{ex:open-comparison}}; \ref{SHORT.ex:open-comparison:positive}.
Apply the Rauch comparison (\ref{prop:rauch}) and the properties of the exponential map in \ref{prop:exp}.

\parit{\ref{SHORT.ex:open-comparison:almost-min}.} Argue by contradiction.
If the statement does not hold, then for any $p\in\Sigma$ there is a point $q=q(p)\in \Sigma$ such that 
$\dist{q}{p}\Sigma<100\cdot R_p$
and
$R_q<(1-\tfrac1{100})\cdot R_p$.

Start with any point $p_0$, and consider a sequence $p_n$ defined by $p_{n+1}\z\df q(p_n)$.
Show that $p_n$ converges, and $R_{p_n}\to 0$ as $n\to\infty$.
Arrive at a contradiction using \ref{SHORT.ex:open-comparison:positive}.

\parit{\ref{SHORT.ex:open-comparison:proof}.} Repeat the proof assuming that $p$ is provided by \ref{SHORT.ex:open-comparison:almost-min}.

\parbf{\ref{ex:geod-convexity}.}
Apply \ref{angk>angk} and \ref{angk<angk}.

\parbf{\ref{ex:midpoints}.} Use \ref{angk>angk} or \ref{angk<angk} twice:
first --- for the triangle $[pxy]$ and $\bar x\in [p,x]$;
second --- for the triangle $[p\bar xy]$ and $\bar y\in [p,y]$.
Then apply the angle monotonicity (\ref{lem:angle-monotonicity}).

\parbf{\ref{ex:convex-dist}.}
It is sufficient to prove the Jensen inequality;
that is, 
\begin{align*}
\dist{\gamma_1(\tfrac12)}{\gamma_2(\tfrac12)}{}
\le
\tfrac12\cdot&\dist{\gamma_1(0)}{\gamma_2(0)}{}
+
\\
+
\tfrac12\cdot&\dist{\gamma_1(1)}{\gamma_2(1)}{}.
\end{align*}

Let $\delta$ be the geodesic path from $\gamma_2(0)$ to $\gamma_1(1)$.
From \ref{ex:midpoints}, we have
\begin{align*}
\dist{\gamma_1(\tfrac12)}{\delta(\tfrac12)}{}
&\le
\tfrac12\cdot\dist{\gamma_1(0)}{\delta(0)}{},
\\
\dist{\delta(\tfrac12)}{\gamma_2(\tfrac12)}{}
&\le
\tfrac12\cdot\dist{\delta(1)}{\gamma_2(1)}{}.
\end{align*}
Sum it up, and apply the triangle inequality.

\begin{Figure}
\vskip-0mm
\centering
\includegraphics{mppics/pic-1650}
\vskip1mm
\end{Figure}

\parit{Remark.} Modulo the comparison theorem, 
the case of the Euclidean plane is just as hard.

\parbf{\ref{ex:disc+}.} Consider a ball $B=\bar B[p,R]_\Sigma$.
Given a point $x\in B$ choose a geodesic $[p,x]$ and denote by $\tilde x$ the point in $\T_p$ that lies at the distance $|p-x|_\Sigma$ from $p$ in the direction of $[p,x]$.
The map $x\mapsto \tilde x$ sends $B$ into the $R$-ball in $\T_p$;
by \mbox{\ref{mang>angk}}, this map is distance-noncontracting.
Whence the statement follows.

\parbf{\ref{ex:disc-}}.
Let $p\in\Sigma$ be the center of $\Delta$. 

\parit{\ref{SHORT.ex:disc-:kg}.}
Consider geodesic triangle $[pxy]_\Sigma$ with $|p-x|_\Sigma=|p-y|_\Sigma=R$;
set $\alpha=\measuredangle\hinge xpy$, 
$\beta=\measuredangle\hinge ypx$ and $\ell=|x-y|_\Sigma$.
Use \mbox{\ref{mang<angk}} to show that 
\[\pi-\alpha-\beta>\tfrac\ell R+o(\ell).\]

Denote by $\sigma$ the arc of $\partial \Delta$ from $x$ to $y$.
Note that $[x,y]_\Sigma$ lies in the region bounded by $[p,x]$, $[p,y]$, and $\sigma$.
Apply the Gauss--Bonnet formula to show that $\tgc\sigma>\pi-\alpha-\beta$.
Pass to a limit as $\ell\to0$, and make a conclusion.

\parit{\ref{SHORT.ex:disc-:area}.}
By \ref{mang<angk}, the map $\exp_p\:\T_p\to\Sigma$ is distance-noncontracting.
By \ref{prop:gamma''}, it maps the $R$-ball in $\T_p$ centered at the origin onto $\Delta$.
Whence the statement follows.

\parbf{\ref{ex:moon-}.}
Mimic the proof of \ref{thm:moon-orginal} using \ref{ex:disc-:kg}.
For the last statement, apply \ref{ex:disc-:area}.

\parit{Source:}
Suggested by Dmitri Burago.

%\end{multicols}
%\par\noindent\rule{\textwidth}{0.4pt}
%\begin{multicols}{2}

\stepcounter{chapter}
\setcounter{eqtn}{0}

\parbf{\ref{ex:ell-infty}.}
Check all the conditions in the definition of metric, page \pageref{page:def:metric}.

\parbf{\ref{ex:B2inB1}}; \ref{SHORT.ex:B2inB1:a}.
Observe that $\dist{p}{q}{\spc{X}}\le 1$.
Apply the triangle inequality to show that $\dist{p}{x}{\spc{X}}\le 2$ for any $x\in B[q,1]$.
Make a conclusion.

\parit{\ref{SHORT.ex:B2inB1:b}.}
Take $\spc{X}$ to be the half-line $[0,\infty)$ with the standard metric; $p=0$ and $q=\tfrac45$.

\parbf{\ref{ex:dist-preserv=>injective}.}
If $x\ne y$, then $\dist{x}{y}{\spc{X}}>0$.
Since $f$ is distance-preserving,
$$\dist{f(x)}{f(y)}{\spc{Y}}=\dist{x}{y}{\spc{X}}.$$
Conclude that $\dist{f(x)}{f(y)}{\spc{Y}}>0$, and hence $f(x)\z\ne f(y)$.

\parbf{\ref{ex:shrt=>continuous}.}
Show that the conditions in \ref{def:continuous} hold for $\delta=\epsilon$.

\parbf{\ref{ex:close-open}.}
Suppose the complement $\Omega=\spc{X}\setminus Q$ is open.
Then for each point $p\in \Omega$ there is $\epsilon>0$ such that $\dist{p}{q}{\spc{X}}>\epsilon$ for any $q\in Q$.
It follows that $p$ is \textit{not} a limit point of any sequence $q_n\in Q$.
That is, any limit of a sequence in $Q$ lies in $Q$;
by the definition, $Q$ is closed.

Now suppose $\Omega=\spc{X}\setminus Q$ is not open.
Show that there is a point $p\in \Omega$ and a sequence $q_n\in Q$ such that $\dist{p}{q_n}{\spc{X}}\z<\tfrac 1n$ for any~$n$.
Conclude that $q_n\to p$ as $n\to \infty$;
therefore $Q$ is not closed.

\spell{\end{multicols}}{}

%% file: afterword.tex
\chapter{Afterword}

After our book, if you continue with Riemannian geometry,
then half of the material in this subject will look familiar.
But first, you need to read a text on tensor calculus;
the book by Richard Bishop and Samuel Goldberg \cite{bishop-goldberg} is one of our favorites.

Let us list a few introductory texts we know and love that range from very detailed to very condensed and challenging:
\begin{itemize}
\item ``Riemannian manifolds''  \cite{lee2006riemannian} by John Lee.
\item ``Riemannian geometry'' \cite{carmo1992riemannian} by Manfredo do Carmo.
\item ``Riemannsche Geometrie im Großen'' \cite{gromoll-klingenberg-meyer} by 
Detlef Gromoll,
Wilhelm Klingenberg, 
and  Werner Meyer;
it is in German, plus there is a Russian translation.
In Russian, there is a more elementary textbook by Yurii Burago and Viktor Zalgaller \cite{burago-zalgaller}.
\item ``Comparison geometry'' \cite{cheeger-ebin} by Jeff Cheeger and David Ebin. 
\item ``Sign and geometric meaning of curvature'' \cite{gromov-1991} by Mikhael Gromov.
\end{itemize}
Good luck.

\begin{flushright}
Anton Petrunin and Sergio Zamora Barrera,\\
State College, Pennsylvania, May 11, 2021.
\end{flushright}